\documentclass[envcountsame,envcountsect,a4paper]{svmonociro}

\usepackage{makeidx,multicol,graphics,graphicx,amssymb,amsfonts,amsthm,amsmath,mathrsfs,epsfig}
\usepackage[sort]{cite}

\usepackage[all]{xy}
\usepackage[english]{babel}
\usepackage[bottom]{footmisc}

\newcommand{\ax}{\operatorname{Ax}}
\newcommand{\op}{\operatorname{op}}
\newcommand{\edm}{\operatorname{End}_\QQ}
\newcommand{\obj}{\operatorname{Obj}}

\newcommand{\dom}{\operatorname{dom}}
\newcommand{\cod}{\operatorname{cod}}
\newcommand{\id}{\operatorname{id}}

\newcommand{\luk}{\operatorname{\text{\L}}}

\newcommand{\frsl}{\operatorname{Free}_\SL}

\newcommand{\fr}{\operatorname{Free}}
\newcommand{\Quot}{\operatorname{Qt}}
\newcommand{\Cl}{\operatorname{Cl}}
\newcommand{\Mor}{\operatorname{Mor}}
\newcommand{\Tp}{\operatorname{Tp}}

\newcommand{\lang}{\mathcal{L}}
\newcommand{\fml}{\mathit{Fm}_\mathcal{L}}
\newcommand{\Fml}{\mathbf{Fm}_\mathcal{L}}
\newcommand{\fmll}{\mathit{Fm}_\mathcal{L'}}
\newcommand{\Fmll}{\mathbf{Fm}_\mathcal{L'}}

\newcommand{\sfm}{\Sigma_{\cat L}}
\newcommand{\Sfm}{\mathbf{\Sigma_{\cat L}}}
\newcommand{\sfmm}{\Sigma_{\cat L'}}
\newcommand{\Sfmm}{\mathbf{\Sigma_{\cat L'}}}

\newcommand{\SL}{\mathcal{SL}}
\newcommand{\Q}{\mathcal{Q}}
\newcommand{\Op}{\mathit{Op}}

\newcommand{\RL}{\mathcal{RL}}
\newcommand{\Mod}{\textrm{-}\mathcal{M}\!\!\:\mathit{od}}
\newcommand{\Modd}{\mathcal{M}\!\!\:\mathit{od}\textrm{-}}
\newcommand{\cat}{\mathcal}
\newcommand{\Edm}{\mathbf{End}_\QQ}

\renewcommand{\I}{\mathrm{I}}
\newcommand{\MM}{\mathbf{M}}
\renewcommand{\AA}{\mathbf{A}}
\newcommand{\NN}{\mathbf{N}}
\newcommand{\LL}{\mathbf{L}}
\newcommand{\PP}{\mathbf{P}}

\newcommand{\QQ}{\mathbf{Q}}
\newcommand{\RR}{\mathbf{R}}
\newcommand{\EQ}{\mathbf{Eq}}
\newcommand{\Eq}{\mathit{Eq}}
\renewcommand{\SS}{\mathbf{S}}
\renewcommand{\hom}{\operatorname{Hom}}
\newcommand{\Hom}{\mathbf{Hom}}
\newcommand{\R}{\mathbb{R}}
\newcommand{\N}{\mathbb{N}}
\newcommand{\Z}{\mathbb{Z}}
\newcommand{\U}{\mathbb{U}}
\renewcommand{\wp}{\mathscr{P}}
\renewcommand{\th}{\mathit{Th}}
\newcommand{\Th}{\mathbf{Th}}
\newcommand{\V}{\mathit{Var}}
\newcommand{\TT}{\mathbf{T}}
\newcommand{\VV}{\mathbf{V}}

\renewcommand{\phi}{\varphi}
\newcommand{\e}{\varepsilon}
\newcommand{\g}{\gamma}
\renewcommand{\d}{\delta}

\newcommand{\restr}{\upharpoonright}
\newcommand{\la}{\langle}
\newcommand{\ra}{\rangle}
\newcommand{\liff}{\Longleftrightarrow}
\newcommand{\upm}{_{\uparrow m}}

\newcommand{\under}{\backslash}
\newcommand{\ost}{{}_\star/}
\newcommand{\lto}{\longrightarrow}
\newcommand{\To}{\Rightarrow}
\newcommand{\lTo}{\Longrightarrow}
\newcommand{\lmapsto}{\longmapsto}
\newcommand{\ust}{\under_\star}
\newcommand{\ov}{\overline}
\newcommand{\ul}{\underline}

\newcommand{\tensor}{\otimes}

\newcommand{\eq}{{\approx}}
\newcommand{\da}{\downarrow}
\newcommand{\ua}{\uparrow}
\newcommand{\Dashv}{= \!\! \!\! \: |}

\newcommand{\Sse}{\Dashv  \! \! \models}
\renewcommand{\vDash}{\models}

\makeindex

\begin{document}

\author{Ciro Russo}
\title{Quantale Modules, with Applications to Logic and Image Processing}
\subtitle{Doctoral dissertation \\ Ph.D. in Mathematics}
\date{November 14, 2007}
\maketitle


\frontmatter

\thispagestyle{empty}
\vspace*{3.5cm}
\begin{flushright}

{\em To the memory of my father, Francesco. \\
To my newborn nephew and godson, \\
wishing he will share with his grandad \\
much more than the bare name.}

\end{flushright}

\ack

Most of the people who read a doctoral dissertation are academics. Then, as it often understandably happens, they may underestimate its importance for the author. A Ph.D. thesis represents a sort of finishing line, of a run begun more than twenty years before. So there is no reason for being sparing of thanks and gratitude.

A few years ago, while writing my degree thesis, I read a booklet by Umberto Eco, entitled ``Come si fa una tesi di laurea'' (that is ``How to make a degree thesis''). One of the first hints he gives in that book is that it is {\bf inelegant} to thank your advisor in the acknowledgements of your thesis, because he's simply doing his job, nothing more, and in most cases your acknowledgement would be nothing but an act of obedience. It is probably true in many situations and, on the other hand, I have to admit that there were many useful suggestions in that book. Nevertheless, I believe (and I challenge everybody to find a disproof) that the way we can do our job is never unique, and this fact implies that not all the advisors and teachers are the same.

That's why I would rather ignore an advice from a famous and highly esteemed man of letters and be ``inelegant''. My advisor Antonio Di Nola deserves to be thanked not for his role, but for how he plays it. First of all, he bet on me in a moment of my life when even I wouldn't have been so brave, and during this three years~--- for some mysterious reasons~--- he kept on doing it, by giving me many opportunities. His constant, kind or harsh, encouragement and spur have been indispensable for me.

For similar reasons I would like to thank my co-advisor Constantine Tsinakis, I feel lucky and honoured of working with him, as well. His precision and his careful analysis of my results left a deep mark on this thesis, and prevented me from making mistakes and slips. Last, I'm grateful to him for the opportunity of teaching a course at the Vanderbilt University that he (again, for mysterious reasons) offered me, and for his kind hospitality during my first stay in Nashville. In this connection, I also want to thank Costas' wife, Annell, but she deserves a special thank~--- too~--- for the care she spent in helping me to find the Gibson Les Paul of my dreams.

There are also four further people who interpreted and/or interprete their institutional roles with uncommon devotion, always helpful in solving any kind of issue: the Director of the Ph.D. programme in Mathematics at the University of Salerno, Patrizia Longobardi, her predecessor Giangiacomo Gerla (who also provided me with relevant and rare scientific works), as well as the current and the former Directors of the C.A.D. in Mathematics, Mariella Transirico and Mercede Maj.

Another person I never had the opportunity to thank is the one who, before anybody else, taught me the beauty of this discipline: my Mathematics teacher at the Secondary School, Speranza Fedele. After my advisors and her, I want to give special thanks to Nick Galatos and Enrico Marchioni because, despite of their young age,\footnote{At least I hope their age can be considered ``young'', since they are more or less as old as me\ldots} they gave me extremely important suggestions during the preparation of my thesis.

Then I cannot fail to mention all the ``senior'' mathematicians I had the pleasure to meet in the last years and stroke me for their intellectual and human stature. Namely: Vito Di Ges\`u, Anatolij Dvurec\v{e}nskij, Francesc Esteva, Llu\'is Godo, Revaz Grigolia, Petr H\'ajek, Charles W. Holland, Bjarni J\'onsson, Ada Lettieri, Ralph McKenzie, Franco Montagna, Daniele Mundici, Hiroakira Ono, Francesco Paoli, Enric Trillas, Esko Turunen and Yde Venema.

Then there are, of course, all the Ph.D. students, Post-Doc fellows and researchers who act within the orbit of Office 6/12 of our department. Even if they're so many (and that's why our office is also called the ``Hilbert infinite office''), I want to mention each of them: Pina Albano, Marilia Amendola and Maria Ferro,\footnote{Because they always criticize me, whatever happens.} Michael B\"achtold, Serena Boccia, Ivana Bochicchio, Diego Catalano Ferraioli, Dajana Conte, Cristina Coppola,\footnote{For her willingness and her uncorrectable fuzziness, but also for her patience with somebody we know\ldots} Raffaele D'Ambrosio, Sergio Di Martino, Christian Di Pietro,\footnote{The best horoscope reader.} Daniel Donato, Daniele Genito, Diana Imperatore, Annamaria Lucibello,\footnote{It's simply fantastic to eat and drink wine with her.} Giovanni Moreno,\footnote{For his absurd mathematical quips and photomontages.} Tiziana Pacelli,\footnote{For all the times she helped me in standing Cristina.} Rino Paolillo, Rossella Piscopo, Pasquale Ronca, Carmela Sica,\footnote{Besides all her merits listed at the end of these acknowledgements, she is the best confectioner within a radius of hundreds light-years; e.g. a Sacher Torte like hers can hardly be found even in Vienna!} Luca Spada,\footnote{For all the battles he fights against the obtuse bureaucracy and the fruitful discussions we had during the preparation of this thesis.} Antonio Tortora, Maria Tota,\footnote{!} Alessandro Vignes and Luca Vitagliano.

I wish to thank also all the young researchers with whom I had a good time either discussing Mathematics or just relaxing: Stefano Aguzzoli, Agata Ciabattoni, Pietro Codara, Tommaso Flaminio, Brunella Gerla, Peppe Greco, Dan Guralnik, Antonio Ledda, Vincenzo Marra, George Metcalfe (with whom I'm particularly indebted, for several reasons), Norbert Preining and Annika Wille. And I don't forget special people, like Nia Stephens and Akram Aldroubi, that I had the pleasure to meet in Nashville.

Of course, nothing would have been possible without the constant support and encouragement of all my large and enlarged family: my father Francesco who, even from some other world, still keeps on guiding and protecting all of us, my mother (and colleague) Virginia who gave me the ``gene of Mathematics'', my brother Salvatore and his wife Andreina who gave us my adorable newborn nephew Francesco, the best present of the last twenty years, my elder sister Carmen and her husband Davide, my little sister (and future colleague) Anna and her eternal boyfriend Francesco, and all the Sica family.

I have mentioned and thanked so many people that it would be unfair to forget my best friends: Nicola Allocca, Lello Aviello, Antonio Bianco, Geremia De Stefano, Enzo Esposito, Domenico Lieto, Nando Gelo, Giampiero Guarino, Domenico Iavarone, Tino Nigro, Stefano Sgambati and Ester Sposato. Fabio Corrente deserves a special thank for his redemption. After several years studying Physics, he finally found ``the straight and narrow path'': a Ph.D. programme in Applied (well, it's enough, I don't want to expect the impossible...) Mathematics; moreover he made me notice that the acknowledgements of a thesis can be completely free and informal (this is a meta-acknowledgement!).

I wish {\bf not} to thank Italian governments and politicians who, year by year, are increasingly dismantling our education and research systems. Contrarywise, there is an Italian foundation whose financial support, at the beginning of my Ph.D., has been determinant: the Fondazione ONAOSI, a precious backing for children of retired or dead physicians and doctors. 

Last, just like the sweet at the end of a good meal, I want to express my gratitude to my beloved fianc\'ee Carmela. A valuable colleague, a precious, funny, pleasant and patient friend, the sweetest, kindest and most beautiful woman I've ever met. I'm sure: a wonderful lifemate. Thank you for standing always by me. 

\vspace{\baselineskip}
\begin{flushright}\noindent
November 14, 2007 \hfill {\it Ciro Russo} \\
\end{flushright}

\tableofcontents



\mainmatter


\introduction
\setcounter{footnote}{0}

The formal concept of deductive system can be probably traced back to Friedrich Ludwig Gottlob Frege and David Hilbert; afterwards, several extensions and refinements of those definitions have been proposed. For example the existence of algebraic semantics for certain logics, together with the birth of Model Theory, can be ascribed to Alfred Tarski's work (although the first discoveries on the relationship between Algebra and Logic had been made by George Boole in the middle of nineteenth century), and the classes of algebras that constituted such semantics were endowed on their turn with consequence relations. In 1934, Gerhard Karl Erich Gentzen proposed --- in \cite{gentzen} --- a further definition of deductive system that included both the aforementioned ones.

Since then, many authors have approached logical systems essentially by studying their deducibility relations. The fact that such relations can be treated from three different point of view --- namely as binary relations, closure operators or unary operations --- is well known, and those different perspectives are equally common in literature; the choice of one of them rather than the others usually depends either on how plainly they allow the authors to present the results or on the genesis of the results themselves or, simply, on the style of the single authors. However the wide literature on consequence relations pass through the twentieth century, and a systematic (and not too old) collection of the results in this area can be found in the book by Ryszard W\'ojcicki, \cite{woj}. The connection of Hilbert and Gentzen systems with equational systems, in the wake of Tarski's approach, has been deeply investigated in the last two decades by Willem Johannes Blok, Janusz Czelakowski, Bjarni J\'onsson and Don Pigozzi (besides many other authors), and important results have been obtained on this subject; see \cite{blokjonsson2,blokpigozzi,blokpigozzi2}. On the contrary, the relations between Hilbert and Gentzen systems cannot boast the same interest from logicians. This is probably due to the fact that, of course, the possibility of reducing the study of a formal logical system to that of a class of algebraic structures is more useful than moving from a formal system to another with different rules of calculus. By the way, a detailed work, \cite{raftery} by James G. Raftery, on the connection between Hilbert and Gentzen systems was published in 2006.

The approach to deductive systems by means of complete residuated lattices and complete posets has been proposed by Nikolaos Galatos and Constantine Tsinakis in~\cite{galatostsinakis}, and we will revisit it in terms of quantales and quantale modules, adding a few new contributions. This theory is very recent\footnote{As the reader will see from the bibliographical references, the cited paper has not been published yet, at the time this thesis is being written.}; nonetheless the results shown in~\cite{galatostsinakis}, that will be recalled in Chapters~\ref{logicchap} and~\ref{consrelchap}, together with our contributions presented in Chapter~\ref{consrelchap}, are unquestionably promising and~--- in our opinion~--- open a new and fruitful perspective on Mathematical Logic, besides proving, once more, its strong relationship with Algebra.

Briefly, in~\cite{galatostsinakis} the authors prove that equivalences and similarities between deductive systems, over a propositional language, can be treated with categorical and algebraic tools, by representing each deductive system $\cat S = \la \lang, \vdash \ra$ as a pair $\la \mathbf{P}, \g_\vdash \ra$, where $\mathbf{P}$ is a complete poset over which an action from a complete residuated lattice is defined, and $\g_\vdash$ is a \emph{structural closure operator}, i.e. a closure operator that is invariant, in a precise sense, under the action of the residuated lattice.

More explicitly, the complete poset $\mathbf{P} = \la P, \leq \ra$ is the powerset of the set of formulas, equations or sequents over the language $\lang$, ordered by set inclusion, and the complete residuated lattice is the powerset $\wp(\sfm)$ of the set of substitutions over $\lang$, again with the set inclusion as order relation. The structural closure operator is defined as the operator that sends each set $\Phi$ of formulas, equations or sequents to the set of all the formulas, equations or sequents that are deducible, with respect to $\vdash$, from $\Phi$. Starting from this representation, $\wp(\Sfm) = \la \wp(\sfm), \cup, \circ, \varnothing, \{\id\}\ra$ turns out to be a (unital) quantale and $\wp(\Fml) = \la \wp(\fml), \cup, \varnothing \ra$, $\wp(\EQ_\lang) = \la \wp(\Eq_\lang), \cup, \varnothing \ra$ and $\wp(\mathbf{Seq}_\lang) = \la \wp(\mathit{Seq}_\lang), \cup, \varnothing \ra$ quantale modules over $\wp(\Sfm)$. Then the images of such modules under the operator $\g_\vdash$~--- that are easily seen to be the lattices of theories, $\th_\vdash$, of $\vdash$~--- have a structure of $\wp(\Sfm)$-module as well: $\Th_\vdash = \la \th_\vdash, \g_\vdash \circ \cup, \g(\varnothing)\ra$, where the action of $\wp(\Sfm)$ is defined as the composition of $\g_\vdash$ with the action of $\wp(\Sfm)$ on $\wp(\Fml)$, $\wp(\EQ_\lang)$ or $\wp(\mathbf{Seq}_\lang)$. In this setting, $\Th_\vdash$ is homomorphic image of the original module in the category of left $\wp(\Sfm)$-modules.

The second motivation that stimulated our investigation on quantale modules comes from the area of Image Processing. Indeed, in the literature of Image Processing, several suitable representations of digital images as $[0,1]$-valued maps are proposed. Such representations are the starting point for defining both compression and reconstruction algorithms based on fuzzy set theory, also called \emph{fuzzy algorithms}, and mathematical morphological operators, used for shape analysis in digital images. Most of the fuzzy algorithms use a suitable pair of operators, one for compressing the image and the other one for approximating the original image starting from the compressed one; see, for instance,~\cite{dinolarusso,loiasessa,perfilieva}. The idea is similar to that of the so-called ``integral transforms'' in Mathematical Analysis: every map can be discretized by means of the \emph{direct transform} and then approximated through the application of a suitable \emph{inverse transform}. Moreover, as well as a direct integral transform is defined as an integral-product composition, a fuzzy compression operator is defined as a join-product composition (where the product is actually a left-continuous triangular norm), and its inverse operator has the form of a meet-division composition, where the division is the residual operation of the same triangular norm. In mathematical morphology, the operators of \emph{dilation} and \emph{erosion}~--- whose action on an image can easily be guessed by their names~--- that are translation invariant can be expressed, again, as compositions join-product and meet-residuum respectively. All these methods can be placed under a common umbrella by essentially abstracting their common properties. Indeed they are all examples of \emph{$\Q$-module transforms}, that we will define in Chapter~\ref{mqchap} and that turn out to be precisely the homomorphisms between free $\Q$-modules.

All these considerations show that a deep categorical and algebraic study of quantale modules, will probably push progresses in the fields of Logic and Image Processing and the aim of this thesis is right to study quantale modules keeping their possible applications as a constant inspiration and a further intention. Our hope is also to give impulse to an extensive study of quantale modules, by showing a glimpse of their great potential.

\smallskip
\smallskip

\begin{center}
{\sc Structure of the Work}
\end{center}
The thesis is organized in three distinct parts and seven chapters.

\paragraph{{\sc Part~\ref{prelpart}} \textnormal{contains most of the preliminary notions and results, and is divided in three chapters.}}
\noindent {\bf Chapter~\ref{logicchap}} \quad We recall some definitions and results regarding algebraizable logics and deductive systems, useful for motivating the study of quantale modules and necessary for the comprehension of the results of Logic that follow as applications of the properties of quantale modules.

\noindent {\bf Chapter~\ref{catschap}} \quad This chapter is dedicated to a brief overview of the categorical notions that will be involved in the theory developed.

\noindent {\bf Chapter~\ref{orderchap}} \quad We show some results on residuated maps, sup-lattices and quantales. A quantale module can be thought of as an object similar to a ring module, where we have a quantale instead of a ring and a sup-lattice instead of an Abelian group. Then, in order to make the thesis as self-contained as possible, it is necessary to include some preliminary notions and results, also considering that such structures may not be as familiar to the reader as rings and groups.

\paragraph{{\sc Part~\ref{modupart}} \textnormal{can be considered as the main (theoretical) part, and is divided in two chapters.}}

\noindent {\bf Chapter~\ref{mqchap}} \quad This chapter contains all the main results and constructions on quantale modules. In Sections~\ref{mqbasnotsec} and~\ref{freemqsec} we present the categories of quantale modules and start establishing the first results. So we define objects, morphisms, subobjects, free objects, hom-sets, and we show that the product and the coproduct of a family of quantale modules are both the Cartesian product equipped with coordinatewise defined operations and with, respectively, canonical projections and their left adjoints as the associated families of morphisms. Sections~\ref{nucleisec} and~\ref{mqtransec} are devoted to the study of two classes of operators, that we already mentioned, on quantale modules: the \emph{structural closure operators} (also called \emph{nuclei}) and the \emph{transforms}. The importance of such operators for applications has been already underlined, but it will be clear soon also their centrality for the theory of quantale modules. Indeed, among other things, we show that each transform is a $\Q$-module homomorphism of free modules and vice versa. On the other hand, $\Q$-module nuclei are strongly connected with homomorphisms as well; we can say that, somehow, $\Q$-module morphisms, transforms and nuclei are three different points of view of the same concept.

Projective (and injective) quantale modules are investigated in Section~\ref{projinjsec}. Apart from free modules, that are obviously projective, we show a characterization of projective cyclic modules and we prove that the product of projective (respectively, injective) objects is projective (resp., injective); several results of this section are due to N. Galatos and C. Tsinakis. In Section~\ref{amalsec}, we prove that the categories of quantale modules have the strong amalgamation property while, in Section~\ref{mqtensorsec}, we show the existence of tensor products of quantale modules. Their construction and properties are similar to the analogues for ring modules, and this analogy include also the use of tensor products for extending the set of scalars of a module. Later on, we show that any module obtained by extending the quantale of scalars of a coproduct of cyclic projective modules, is the coproduct of cyclic projectives itself, thus projective.

\noindent {\bf Chapter~\ref{consrelchap}} \quad The exposition of how quantale modules are connected to Logic is the content of this chapter. At the beginning, we abstract the definitions of consequence relations and deductive systems in the algebraic frameworks of sup-lattices and quantale modules. This approach follows, even if from a slightly different point of view, the aforementioned work by N. Galatos and C. Tsinakis. The main novelty, here, consists of the algebraic treatment of the concepts of translation and interpretation between logics over different languages.

\smallskip

\paragraph{{\sc Part~\ref{applpart}} \textnormal{contains the applications of quantale modules~--- and especially of $\Q$-module transforms~--- to Image Processing; it is composed of two chapters.}}

\noindent {\bf Chapter~\ref{imagechap}} \quad After a brief overview on the literature on fuzzy image compression and mathematical morphology, we show how parts of these areas fall within the formal theory we have established in Part~\ref{modupart}.

\noindent {\bf Chapter~\ref{ltbchap}} \quad We present an example of $\Q$-module transform together with a concrete application of it. The results of the application have been compared with those obtained by using JPEG, the most famous algorithm for image compression. The operator shown, called \L ukasiewicz transform, is defined between free modules on the quantale reduct of the MV-algebra $[0,1]$; the algorithm based on it is called, not surprisingly, \L TB~--- ``\L ukasiewicz Transform Based'' (see~\cite{dinolarusso,dinolarusso2}).


\part{Preliminaries}
\label{prelpart}

\chapter{On Propositional Deductive Systems}
\label{logicchap}

This chapter shall be considered as an overview of the main notions of Mathematical Logic we will deal with in this thesis. At the present time, as we anticipated in the Introduction, the main applications of quantale modules to Mathematical Logic are limited to the case of logics on a propositional language. Therefore we will not go beyond this level, since it suffices to introduce all the notions we need.

On the other hand we will see, both in this chapter and in Chapter~\ref{consrelchap}, that many of the results presented hold for any kind of deductive system defined on a propositional language, regardless of whether it is a propositional logic, an equational system or a sequent-based system.

In Section~\ref{proplogicsec}, once recalled some basic notions, we will start by defining the most simple kind of deductive system~--- the propositional one~--- showing some of its best known examples. In Section~\ref{eqlogicsec}, we will define an equational deductive system, and we will recall the notion and the characterization of algebraizable logics, according to W. J. Blok and D. Pigozzi~\cite{blokpigozzi}. Moreover we will give a reformulation, due to N. Galatos and C. Tsinakis, of algebraizability, that is amenable to a generalization to ``abstract'' deductive systems. A further step~--- namely, the introduction of Gentzen-style deductive systems~--- is the subject of Section~\ref{seqlogicsec}.

Last, in Section~\ref{conspwsetsec}, we will walk a first step toward the algebraic approach to consequence relations proposed in~\cite{galatostsinakis}. In particular, we will first extend the notions of asymmetric and symmetric consequence relation to arbitrary sets, and then we will show that symmetric and asymmetric consequence relations are essentially the same concept, thus the study of consequence relations can be limited to one of them. We will see, in Chapter~\ref{consrelchap}, that an abstract definition of consequence relation is possible essentially thanks to this result.

\section{Propositional deductive systems}
\label{proplogicsec}

In this first section we will recall the very basic definitions regarding \emph{propositional deductive systems}\index{Propositional!-- deductive system} (or \emph{propositional logics}\index{Propositional!-- logic}, for short). In order to define a propositional logic, we need several preliminary concepts; so, starting from the notion of ``language'', we will now introduce all the essential consituents of a propositional logic.

A \emph{propositional}~--- or \emph{algebraic}~--- \emph{language}\index{Propositional!--language}\index{Algebraic language}\index{Language!Propositional --}\index{Language!Algebraic --} is a pair  $\lang = \la L, \nu \ra$ of a set $L$ and a map $\nu: L \lto \N_0$. The elements of $L$ are called (\emph{primitive}) \emph{connectives}\index{Connective}\index{Connective!Primitive --}, or \emph{operation symbols}\index{Operation symbol}, and the image of a connective under $\nu$ is called the \emph{arity}\index{Arity of a connective}\index{Connective!Arity of a --} of the connective; nullary connectives, i.e. connectives whose arity is zero, are also called \emph{constant symbols}\index{Constant symbol}\index{Connective!Nullary --}.

Given a propositional language $\lang$ and a denumerable set of variables $V = \{x_n \mid n \in \N\}$, the $\lang$-formulas are strings of connectives and variables that respect certain constraints; more precisely, the $\lang$-formulas are defined recursively by means of the following conditions:
\begin{enumerate}
\item[(F1)] every propositional variable is an $\lang$-formula,
\item[(F2)] every constant symbol is a formula,
\item[(F3)] if $f$ is a connective of arity $\nu(f) > 0$ and $\phi_1, \ldots, \phi_{\nu(f)}$ are $\lang$-formulas, then $f(\phi_1, \ldots, \phi_{\nu(f)})$ is an $\lang$-formula,
\item[(F4)] all $\lang$-formulas are built by iterative applications of (F1), (F2) and (F3).
\end{enumerate}

In what follows we will use the word ``connective'' only for those operation symbols whose arity is greater than zero, while we will always refer to nullary connectives as ``constant symbols'', or ``constants'' for short.

For $\phi \in \fml$ we write $\phi = \phi(x_{i_1}, \ldots, x_{i_n})$ to indicate that the variables of $\phi$ are all included in the set $\{x_{i_1}, \ldots, x_{i_n}\}$. We denote the set of all $\lang$-formulas by $\fml$. If $\sigma: V \lto \fml$ is a map that assigns an $\lang$-formula to each variable, then $\sigma$ can be naturally extended to a map from $\fml$ into itself~--- also denoted by $\sigma$~--- by setting
$$\sigma(\phi(x_{i_1}, \ldots, x_{i_n})) = \phi(x_{i_1}/\sigma x_{i_1}, \ldots, x_{i_n}/\sigma x_{i_n}),$$
where $\phi(x_{i_1}/\sigma x_{i_1}, \ldots, x_{i_n}/\sigma x_{i_n})$ is the formula obtained from $\phi(x_{i_1}, \ldots, x_{i_1})$ by substituting each variable $x_{i_k}$ with its image $\sigma x_{i_k}$ under $\sigma$. Since $\phi$ is a formula, $\sigma(\phi)$ is a formula as well, by virtue of (F1--F4) above.

A map like $\sigma$ is called a \emph{substitution}\index{Substitution}. Of course, $\id_{\fml}$ is a substitution and the composition of two substitutions is again a substitution; thus, once denoted by $\sfm$ the set of all the substitutions over $\fml$, the structure $\Sfm = \la \sfm, \circ, \id_{\fml} \ra$ is a monoid, called the \emph{substitution monoid over $\fml$}\index{Substitution!-- monoid}.

A \emph{(finitary) inference rule}\index{Inference rule} over $\lang$ is a pair $\la \Phi,\psi \ra$ where $\Phi$ is a finite set of formulas and $\psi$ is a single formula. We may think of an inference rule as a law telling us that from a set of formulas that are similar, in a precise sense, to those in $\Phi$, we can infer a formula that is similar, in the same sense, to $\psi$. Indeed, a formula $\phi$ is \emph{directly derivable} from a set $\Psi$ of formulas by the rule $\la \Phi,\psi \ra$ if there is a substitution $\sigma$ such that $\sigma \psi = \phi$ and $\sigma[\Phi] \subseteq \Psi$. An inference rule $\la \Phi,\psi \ra$ is usually denoted by $\frac{\Phi}{\psi}$.

An \emph{axiom}\index{Axiom} in the language $\lang$ is simply a formula of $\fml$. Given a deductive system $\cat S$ (see the definition below), we will denote by $\ax_\cat S$ the set of axioms of $\cat S$.

\begin{definition}\label{consrel}
A \emph{propositional deductive system}\index{Propositional!-- deductive system}\index{Deductive system!Propositional --}, or a \emph{propositional logic}\index{Propositional!-- logic} for short, $\cat S$ over a given language $\lang$, is defined by means of a (possible infinite) set of inference rules and axioms. It consists of the pair $\cat S = \la \lang, \vdash \ra$, where $\vdash$ is a subset of $\wp(\fml) \times \fml$ defined by the following condition: $\Phi \vdash \psi$ iff $\psi$ is contained in the smallest set of formulas that includes $\Phi$ together with all substitution instances of the axioms of $\cat S$, and is closed under direct derivability by the inference rules of $\cat S$. The relation $\vdash$ is called the (\emph{asymmetric}) \emph{consequence relation}\index{Consequence relation}\index{Consequence relation!Asymmetric --} of $\cat S$.
\end{definition}
It can be proved that $\vdash$ satisfies the following conditions for all $\Phi, \Psi \subseteq \fml$ and $\phi, \psi \in \fml$
\begin{eqnarray}
&&\textrm{if $\psi \in \Phi$ then $\Phi \vdash \psi$,} \label{entailin}\\
&&\textrm{if $\Phi \vdash \psi$ and $\Psi \vdash \phi$ for all $\phi \in \Phi$, then $\Psi \vdash \psi$,} \label{entailmp} \\
&&\textrm{if $\Phi \vdash \psi$ then $\Phi_0 \vdash \psi$ for some finite $\Phi_0 \subseteq \Phi$,} \label{entailfinitary} \\
&&\textrm{if $\Phi \vdash \psi$ then $\sigma[\Phi] \vdash \sigma \psi$ for every substitution $\sigma \in \sfm$,} \label{entailstructural}  
\end{eqnarray}
(\ref{entailfinitary}) holding because inference rules are assumed to be finitary.
Reciprocally, it has been proved by \L o\'s and Suszko in~\cite{lossuszko} that, given a language $\lang$, any subset of $\wp(\fml) \times \fml$ satisfying conditions (\ref{entailin}--\ref{entailstructural}) is the consequence relation for some deductive system $\cat S$ over $\lang$.

Thanks to this result, we can give an equivalent definition of a deductive system on a propositional language $\lang$.
\begin{definition}\label{aconsrel}
A subset of $\wp(\fml) \times \fml$ that satisfies (\ref{entailin},\ref{entailmp}) is called an \emph{asymmetric consequence relation} over $\lang$. An asymmetric consequence relation is said to be \emph{finitary}\index{Consequence relation!Finitary --} if it satisfies (\ref{entailfinitary}) and \emph{structural}\index{Consequence relation!Structural --}\index{Structural!-- consequence relation}, or \emph{substitution invariant}\index{Consequence relation!Substitution invariant --}, if it satisfies (\ref{entailstructural}). Then a deductive system over $\lang$ can be defined as a pair $\la \lang, \vdash \ra$, where $\vdash$ is a finitary and substitution invariant consequence relation over $\lang$.
\end{definition}
A \emph{theory}\index{Theory} of a consequence relation $\vdash$ over $\fml$ is a subset $T$ of $\fml$ closed under $\vdash$; i.e. $T$ is a theory of $\vdash$ iff, for all $\phi \in \fml$, $T \vdash \phi$ implies $\phi \in T$. The set of theories of $\vdash$ forms a lattice usually denoted by $\Th_\vdash$. Next, we define \emph{symmetric} consequence relations.
\begin{definition}\label{sconsrel}
A \emph{symmetric consequence relation}\index{Consequence relation}\index{Consequence relation!Symmetric --} $\vdash$ over $\lang$ is a binary relation over $\wp(\fml)$, i.e. a subset of $\wp(\fml) \times \wp(\fml)$, such that, for all $\Phi, \Psi, \Xi \in \wp(\fml)$,
\begin{eqnarray}
&&\textrm{if $\Psi \subseteq \Phi$, then $\Phi \vdash \Psi$,} \label{sentailin} \\
&&\textrm{if $\Phi \vdash \Psi$ and $\Psi \vdash \Xi$, then $\Phi \vdash \Xi$,} \label{sentailimp} \\
&&\Phi \vdash \bigcup_{\Phi \vdash \Psi} \Psi. \label{sentailjoin}
\end{eqnarray}
A symmetric consequence relation $\vdash$ over $\lang$ is called \emph{finitary}\index{Consequence relation!Finitary --}, if for all subsets $\Phi, \Psi$ of $\fml$, with $\Psi$ finite,
\begin{eqnarray}
&&\textrm{if $\Phi \vdash \Psi$, there exists a finite $\Phi_0 \subseteq \Phi$ such that $\Phi_0 \vdash \Psi$;} \label{sentailfinitary}
\end{eqnarray}
it is called \emph{substitution invariant}\index{Consequence relation!Substitution invariant --} or \emph{structural}\index{Structural!-- consequence relation}\index{Consequence relation!Structural --}, if for every substitution $\sigma \in \sfm$ and for all $\Phi, \Psi \subseteq \fml$,
\begin{eqnarray}
&&\textrm{$\Phi \vdash \Psi$ \quad implies \quad $\sigma[\Phi] \vdash \sigma[\Psi]$.} \label{sentailstructural}
\end{eqnarray}
\end{definition}
It will be shown in Section~\ref{conspwsetsec} that symmetric and asymmetric consequence relations, over the same language, are interdefinable, i.e. they are in one-one correspondence. Indeed an asymmetric consequence relation $\vdash^a$ naturally induces a binary relation, $\vdash^s$ on $\wp(\fml)$. It is defined, for all $\Phi, \Psi \in \fml$, by $\Phi \vdash^s \Psi$ iff $\Phi \vdash^a \psi$ for all $\psi \in \Psi$, and it is easy to verify that $\vdash^s$ is really a symmetric consequence relation over $\lang$. When such a correspondence will be shown, we will see that it also preserves both finitarity and substitution invariance. In other words, a symmetric consequence relation $\vdash^s$ is finitary (respectively, structural) if and only if its corresponding asymmetric relation $\vdash^a$ is finitary (respectively, structural). This result allows another equivalent definition of propositional logic: a \emph{deductive system}\index{Deductive system} $\cat S$ is a pair $\la \lang, \vdash \ra$, where $\lang$ is a propositional language and $\vdash$ is a binary relation over $\wp(\fml)$ satisfying (\ref{sentailin}--\ref{sentailstructural}).

We conclude this section giving three well-known examples of propositional logics: the Classical Propositional Logic ($\mathit{CPL}$), the \L ukasiewicz Propositional Logic ($\mathit{\textrm{\emph{\L}}PL}$) and the Modal Logic $S3$.
\begin{example}\label{cpl}
The language of Classical Propositional Logic is $\lang_{\mathit{CPL}} = \la \{\bot, \to\}, \nu_{\mathit{CPL}} \ra$, with $\nu_{\mathit{CPL}}(\bot) = 0$ and $\nu_{\mathit{CPL}}(\to) = 2$. The set $\ax_{\mathit{CPL}}$ of axioms is composed of the three following schemes of formulas
\begin{enumerate}
\item[] $\phi \to (\psi \to \phi)$,
\item[] $(\phi \to (\psi \to \chi)) \to ((\phi \to \psi) \to (\phi \to \chi))$,
\item[] $((\phi \to \bot) \to (\psi \to \bot)) \to (\psi \to \phi)$.
\end{enumerate}
The only rule of inference for \emph{CPL} is \emph{Modus Ponens}\index{Modus Ponens}
$$\mathit{MP} \qquad \frac{\phi \qquad \phi \to \psi}{\psi}.$$
If we define the \emph{derived} unary connective $\lnot$ by setting $\lnot \phi = \phi \to \bot$ for all $\phi \in \mathit{Fm_{CPL}}$, the axioms above take the well known form
\begin{enumerate}
\item[$\mathit{CPL}1$] \quad $\phi \to (\psi \to \phi)$,
\item[$\mathit{CPL}2$] \quad $(\phi \to (\psi \to \chi)) \to ((\phi \to \psi) \to (\phi \to \chi))$,
\item[$\mathit{CPL}3$] \quad $(\lnot \phi \to \lnot \psi) \to (\psi \to \phi)$.
\end{enumerate}
\end{example}

\begin{example}\label{lpl}
The language of \L ukasiewicz Propositional Logic is $\lang_{\luk} = \la \{\lnot, \to\}, \nu_{\luk} \ra$, with $\nu_{\luk}(\lnot) = 1$ and $\nu_{\luk}(\to) = 2$, and the set $\ax_{\luk}$ of axioms is composed of the four following schemes of formulas
\begin{enumerate}
\item[$\luk1$] \quad $\phi \to (\psi \to \phi)$,
\item[$\luk2$] \quad $(\phi \to \psi) \to ((\psi \to \chi) \to (\phi \to \chi))$,
\item[$\luk3$] \quad $((\phi \to \psi) \to \psi) \to ((\psi \to \phi) \to \phi)$,
\item[$\luk1$] \quad $(\lnot \phi \to  \lnot \psi) \to (\psi \to \phi)$,
\end{enumerate}
with $\mathit{MP}$, again, as the only inference rule.
\end{example}

\begin{example}\label{s3}
The language of $S3$ Modal Logic is $\lang_{S3} = \la \{\bot, \Box, \to\}, \nu_{S3} \ra$, with $\nu_{S3}(\bot) = 0$, $\nu_{S3}(\Box) = 1$ and $\nu_{S3}(\to) = 2$. The set $\ax_{S3}$ of axioms is composed of the following formulas
\begin{enumerate}
\item[$S3_1$] \quad $\Box(\phi \to (\psi \to \phi))$,
\item[$S3_2$] \quad $\Box((\phi \to (\psi \to \chi)) \to ((\phi \to \psi) \to (\phi \to \chi)))$,
\item[$S3_3$] \quad $\Box(((\phi \to \bot) \to (\psi \to \bot)) \to (\psi \to \phi))$,
\item[$S3_4$] \quad $\Box(\Box \phi \to \phi)$,
\item[$S3_5$] \quad $\Box(((\Box(\phi \to \psi) \to (\Box(\psi \to \chi) \to \bot)) \to \bot) \to \Box(\phi \to \chi))$,
\item[$S3_6$] \quad $\Box(\Box \phi \to \Box((\phi \to \bot) \to \psi))$,
\item[$S3_7$] \quad $\Box(\Box(\phi \to \psi) \to \Box(\Box \phi \to \Box \psi))$.
\end{enumerate}
The Modus Ponens is an inference rule also for $S3$, but two further rules are added:
$$\frac{\Box \phi}{\phi} \qquad \textrm{ and } \qquad \frac{\Box(\phi \to \psi) \qquad \Box(\psi \to \phi)}{\Box(\Box \phi \to \Box \psi)}.$$
If we define the derived binary connectives $\vee$ and $\&$ by setting
$$\phi \vee \psi = (\phi \to \bot) \to \psi \qquad \textrm{ and } \qquad \phi \& \psi = (\phi \to (\psi \to \bot)) \to \bot,$$
for all $\phi, \psi \in \mathit{Fm_{S3}}$, the axioms $S3_5$ and $S3_6$ above take the simpler and better known form
\begin{enumerate}
\item[$S3_5$] \quad $\Box((\Box(\phi \to \psi) \& \Box(\psi \to \chi)) \to \Box(\phi \to \chi))$,
\item[$S3_6$] \quad $\Box(\Box \phi \to \Box(\phi \vee \psi))$.
\end{enumerate}
\end{example}

\begin{example}\label{bck}
Let us consider the language $\la \{\to\}, 2\ra$, and denote simply by $\{\to\}$ such language and by $\vdash_{\mathit{BCK}}$ the substitution invariant consequence relation on $\mathbf{Fm}_{\{\to\}}$ having $\mathit{MP}$ as the only rule of inference and the following axioms:  
\begin{enumerate}  
\item[$B$] \quad $(\phi \to \psi) \to ((\psi \to \chi) \to (\phi \to \chi))$,
\item[$C$] \quad $(\phi \to (\psi \to \chi)) \to (\psi \to (\phi \to \chi))$,
\item[$K$] \quad $\phi \to (\psi \to \psi)$.
\end{enumerate}
The deductive system $\la \{\to\}, \vdash_{\mathit{BCK}} \ra$ is called \emph{BCK-logic}.
\end{example}

In the next section, once given the definition of \emph{equivalent algebraic semantics} for a deductive system and the one of \emph{algebraizable logics}, we will see that Examples~\ref{cpl},~\ref{lpl} and~\ref{bck} are also examples of algebraizable logics whose respective equivalent algebraic semantics are some famous classes of algebras. On the contrary, the Modal Logic $S3$ in Example~\ref{s3} is an example of a propositional logic that is not algebraizable.

\section{Equational deductive systems and algebraizability}
\label{eqlogicsec}

In the present section we recall the definition of \emph{equational deductive system} over an algebraic language $\lang$, and the notion of \emph{algebraizable logics}, given by W.~J. Blok and D. Pigozzi in~\cite{blokpigozzi}, together with its characterization.

We will see that the concept of equational deductive system allows new definitions of equivalent algebraic semantics and algebraizable logics, in terms of equivalence between a propositional logic and an equational deductive system over the same language. These definitions and results, mainly due to W.~J. Blok, B. J\'onsson and D. Pigozzi (see~\cite{blokjonsson1,blokjonsson2,blokpigozzi,blokpigozzi2}) and, more recently, to N. Galatos and C. Tsinakis (see~\cite{galatostsinakis}), find a further natural generalization to the case of \emph{deductive systems on sets of sequents} (also called \emph{Gentzen-style systems}) that we will present in Section~\ref{seqlogicsec}.

We start with the definition of $\lang$-algebra.
\begin{definition}\label{lalgebra}
Let $\lang = \la L, \nu \ra$ be an algebraic language. An \emph{$\lang$-algebra}\index{L-algebra@$\lang$-algebra} is a pair $\AA = \la A, \Op[L] \ra$, where $A$ is a set, $\Op$ is a map that assigns an operation $\Op(f) = f^\AA$ on $A$ of arity $\nu(f)$ to each connective $f$ of $L$, and an element of $A$ to each constant symbol of $L$. If $L$ is finite, the elements of $\Op[L]$ in the expression $\la A, \Op[L] \ra$ are usually listed individually. The endomorphism monoid of an $\lang$-algebra $\AA$ is denoted by $\Sfm = \la \sfm, \circ, \id_A \ra$.
\end{definition}
Before we continue, we need to make two remarks. First of all, we observe that an $\lang$-algebra is, in a certain sense, a formal object; in other words, the operations on the algebra need not satisfy, a priori, any property. The second comment concerns the notation introduced for the endomorphism monoids. In Definition~\ref{lalgebra}, we denote by $\Sfm$ the endomorphism monoid of an $\lang$-algebra, exactly as we denoted the substitution monoid of $\fml$ in the previous section. This notation does not generate confusion; indeed the pair $\Fml = \la \fml, L \ra$ is an example of $\lang$-algebra~--- the absolutely free $\lang$-algebra~--- and its endomorphism monoid is precisely the substitution monoid $\Sfm = \la \sfm, \circ, \id_{\fml} \ra$.

Now we can start treating equational deductive systems, moving from the definition of the objects that play the role of formulas in this setting~--- the equations~--- and of consequence relations between them. An \emph{equation}\index{Equation} over $\lang$ is a pair of $\lang$-formulas $s, t \in \fml$, and we usually denote it by the expression $s \eq t$. The pair $\la \fml^2, L^2 \ra$ is clearly an $\lang$-algebra, that we denote by $\EQ_\lang$ and call the $\lang$-algebra of equations over $\lang$.

In order to define consequence relations over $\Eq_\lang$, we need to use the notions of $\lang$-\emph{homomorphism} and \emph{true equality}. If $\AA$ is an $\lang$-algebra, a homomorphism $h: \AA \lto \mathbf{B}$ from $\la A, L^A \ra$ to $\la B, L^B \ra$ is defined, obviuosly, as a map that preserves the $\lang$-operations. If $h: \Fml \lto \AA$ is a homomorphism and $(s \eq t) \in \Eq_\lang$, then we denote by $h(s \eq t)$ the pair $(h(s), h(t)) \in A^2$, and we refer to it as an \emph{equality}. An equality $(h(s),h(t))$ is said to be \emph{true} if $h(s) = h(t)$.

If $\mathcal{K}$ is a class of $\lang$-algebras, and $E \cup \{\e\}$ is a subset of $\Eq_\lang$, $E \models_\mathcal{K} \e$ means that for all $\AA \in \mathcal{K}$ and all homomorphisms $h: \Fml \lto \AA$, if $h[E]$ is a set of true equalities, then $h(\e)$ is a true equality. It is clear that $\models_\mathcal{K}$ is a substitution invariant consequence relation over $\Eq_\lang$, i.e. it satisfies conditions (\ref{entailin},\ref{entailmp},\ref{entailstructural}). It is well known, see e.g. Corollary 7.2 of~\cite{blokjonsson2}, that $\models_\mathcal{K}$ is finitary iff the class of algebras $\mathcal{K}$ is closed with respect to ultraproducts.

According to~\cite{blokpigozzi}, we set the following definition.
\begin{definition}\label{alglog}
Let $\cat S = \la \lang, \vdash \ra$ be a deductive system. A class $\mathcal{K}$ of $\lang$-algebras is called an \emph{equivalent algebraic semantics}\index{Equivalent!-- algebraic semantics}\index{Algebraic semantics!Equivalent --} for $\cat S$ if there exist a finite set of equations $u_i \eq v_i$, $i \in I$, on a single variable and a finite set of binary definable connectives $\varDelta_j$, $j \in J$, such that for every subset $\Phi \cup \{\psi\}$ of $\fml$ and for every equation $s \eq t$ over $\fml$, the following conditions are verified
\begin{enumerate}
\item[$(i)$] $\Phi \vdash \psi$ iff $\{u_i(\phi)\eq v_i(\phi) \mid \phi \in \Phi \} \models_\mathcal{K} u_i(\psi)\eq v_i(\psi)$, for all $i \in I$,
\item[$(ii)$] $s \eq t \Sse_\mathcal{K} \{u_i(s \mathop{\varDelta_j} t) \eq v_i(s \mathop{\varDelta_j} t)\mid i \in I, j \in J\}$.
\end{enumerate}
A deductive system $\cat S$ is said to be \emph{algebraizable}\index{Algebraizable logic}\index{Logic!Algebraizable --} iff there exists an equivalent algebraic semantics for it.
\end{definition}
It can be shown that conditions $(i)$ and $(ii)$ above can be expressed equivalently as follows: for every set of equations $E \cup \{s \eq t\}$ over $\fml$ and for every $\psi \in \fml$,
\begin{enumerate}
\item[$(iii)$] $E \models_\mathcal{K} s \eq t$ iff $\{u\mathop{\varDelta_j} v \mid u \eq v \in E, j \in J \} \vdash  s \mathop{\varDelta_j} t$, for all $j \in J$.
\item[$(iv)$] $\psi \dashv \vdash \{u_i(\psi) \mathop{\varDelta_j} v_i (\psi) \mid i \in I, j \in J\}$. 
\end{enumerate}
Moreover, if we define the maps $\tau: \fml \lto \wp(\Eq_\lang)$ and $\rho: \Eq_\lang \lto \wp(\fml)$ by $\tau(\psi) = \{u_i(\psi) \eq v_i(\psi) \mid i \in I\}$ and $\rho(s \eq t) = \{s \mathop{\Delta_j} t \mid j \in J\}$, then conditions $(i)$ and $(ii)$ take the more elegant form
\begin{enumerate}
\item[$(i')$] $\Phi \vdash \psi$ iff $\tau[\Phi] \models_\mathcal{K} \tau(\psi)$,
\item[$(ii')$] $\e \Sse_\mathcal{K} \tau\rho(\e)$.
\end{enumerate}

In~\cite{galatostsinakis}, the authors characterize the maps  $\tau: \fml \lto \wp(\Eq_\lang)$ and $\rho: \Eq_\lang \lto \wp(\fml)$ that allow to rewrite conditions $(i),(ii)$ in the form of $(i')$ and $(ii')$. First of all, we call $\tau$ (respectively, $\rho$) \emph{finitary} if for all $\psi \in \fml$ (resp., for all $\e \in \Eq_\lang$), $\tau(\psi)$ (resp., $\rho(\e)$) is a finite set. Analogously, we call $\tau$ (resp., $\rho$) \emph{structural} or \emph{substitution invariant} if it commutes~--- w.r.t. the composition of maps~--- with substitutions.
 
\begin{lemma}\emph{\cite{galatostsinakis}}\label{algandtrans}
For maps $\tau: \fml \lto \wp(\Eq_\lang)$ and $\rho: \Eq_\lang \lto \wp(\fml)$, the following conditions are equivalent.
\begin{enumerate}
\item[$(a)$] The maps $\tau$, $\rho$ are finitary and substitution invariant.
\item[$(b)$] There exists a finite set of equations $u_i \eq v_i$, $i \in I$, on a single variable and a finite set of binary definable connectives $\varDelta_j$, $j \in J$, such that $\tau(\psi) = \{u_i(\psi) \eq v_i(\psi) \mid i \in I\}$ and $\rho(s \eq t) = \{s \mathop{\varDelta_j} t \mid j \in J\}$.
\end{enumerate}
\end{lemma}
\begin{proof}
We just need to prove that $(a)$ implies $(b)$, the converse being already discussed. Let $x, y$ be distinct variables in $V$ and assume that $\tau(x) = \{u_i \eq v_i \mid i \in I\}$ and $\rho(x \eq y) = \{t_j \mid j \in J\}$. Since $\tau$ and $\rho$ are finitary, it follows that $I$ and $J$ are finite. 

If $\psi \in \fml$, let $\kappa_\psi \in \sfm$ be the substitution that sends all variables to $\psi$. Since $\tau$ is substitution invariant we have $\kappa_x(\tau(x)) = \tau(\kappa _x(x)) = \tau(x)$. In other words, if we replace all variables in $\tau(x)$ by $x$, we get $\tau(x)$ back; in other words the equations $u_i \eq v_i$ contain the single variable $x$. Moreover, for all $\psi \in \fml$, we have $\tau(\psi) = \tau(\kappa_\psi(x)) = \kappa_\psi(\tau(x)) = \{\kappa_\psi(u_i(x) \eq v_i(x)) \mid i \in I\} = \{u_i(\psi) \eq v_i(\psi) \mid i \in I\}$.

Let $V_1$ and $V_2$ be two sets that partition the set $V$ of all variables in a way that $x \in V_1$ and $y \in V_2$. For all $s \eq t \in \Eq_\lang$, let $\kappa_{s \eq t} \in \sfm$ be the substitution that sends all variables in $V_1$ to $s$ and all variables in $V_2$ to $t$. Since $\tau$ is substitution invariant, we have $\kappa_{x \eq y}(\rho(x \eq y)) = \rho(\kappa_{x \eq y}(x \eq y)) = \rho(x \eq y)$. In other words, the terms $t_j$ are binary and depend only on the variables $x$ and $y$; we set $t_j = x \mathop{\varDelta_j} y$. Moreover, for all $s \eq t \in \Eq_\lang$, we have $\rho(s \eq t) = \rho(\kappa _{s \eq t}(x \eq y)) = \kappa_{s \eq t}(\rho(x \eq y)) = \{\kappa_{s \eq t}(x \varDelta _j y) \mid j \in J\} = \{s \mathop{\varDelta_j} t \mid i \in I\}$.
\end{proof}
  
\begin{corollary}\emph{\cite{galatostsinakis}}
A deductive system $\la \lang, \vdash \ra$ is algebraizable iff there exist finitary and substitution invariant maps $\tau: \fml \lto \wp(\Eq_\lang)$ and $\rho: \Eq_\lang \lto \wp(\fml)$ and a class of $\lang$-algebras $\mathcal{K}$ such that, for every subset $\Psi \cup \{\psi\}$ of $\fml$ and  $\e \in \Eq_\lang$,
\begin{enumerate}
\item[$(i)$] $\Phi \vdash \psi$ iff $\tau[\Phi] \models_\mathcal{K} \tau(\psi)$,
\item[$(ii)$] $\e \Sse_\mathcal{K} \tau\rho(\e)$.
\end{enumerate}
\end{corollary}

Obviously, the maps $\tau$ and $\rho$ extend to maps $\tau': \wp(\fml) \lto \wp(\Eq_\lang)$ and $\rho': \wp(\Eq_\lang) \lto \wp(\fml)$. Moreover, $\tau'(\Psi)$ and $\rho'(E)$ are finite if $\Psi \in \wp(\fml)$ and $E \in \wp(\Eq_\lang)$ are finite; we will call maps that have this property \emph{finitary}. Also, if  $\Psi \in \wp(\fml)$, $E \in \wp(\Eq_\lang)$ and $\sigma \in \sfm$, then $\sigma[\tau'(\Psi)] = \tau'(\sigma[\Psi])$ and $\sigma[\rho'(E)] = \rho'(\sigma[E])$; we will call such maps \emph{substitution invariant}. It is easy to verify that maps like $\tau'$ and $\rho'$ arise from maps like $\tau$ and $\rho$ iff they preserve unions.

\begin{corollary}\emph{\cite{galatostsinakis}}
A deductive system $\la \lang, \vdash \ra$ is algebraizable iff there exist finitary and substitution invariant maps $\tau: \wp(\fml) \lto \wp(\Eq_\lang)$ and $\rho: \wp(\Eq_\lang) \lto \wp(\fml)$ that preserve unions, and a class of $\lang$-algebras $\mathcal{K}$ such that for every subset $\Psi \cup \{\psi\}$ of $\fml$ and $\e \in \Eq_\lang$,
\begin{enumerate}
\item[$(i)$] $\Phi \vdash \psi$ iff $\tau(\Phi) \models_\mathcal{K} \tau(\psi)$,
\item[$(ii)$] $\e \Sse _\mathcal{K} \tau\rho(\e)$.
\end{enumerate}
\end{corollary}

\begin{example}\label{ba}
Recall that a \emph{Boolean algebra} is an algebra $(B, \vee, \wedge, ', \bot, \top)$ with two binary operations, a unary operation and two constants, satisfying the following equations
$$\begin{array}{llll}
\textrm{(L1$_\vee$)} \ & x \vee y \eq y \vee x & \textrm{(L1$_\wedge$)} \ & x \wedge y \eq y \wedge x \\
\textrm{(L2$_\vee$)} & x \vee (y \vee z) \eq (x \vee y) \vee z & \textrm{(L2$_\wedge$)} & x \wedge (y \wedge z) \eq (x \wedge y) \wedge z \\
\textrm{(L3$_\vee$)} & x \vee x \eq x & \textrm{(L3$_\wedge$)} & x \wedge x \eq x \\
\textrm{(L4$_\vee$)} & x \eq x \vee (x \wedge y) & \textrm{(L4$_\wedge$)} & x \eq x \wedge (x \vee y) \\
\textrm{(D$_\vee$)} & x \wedge (y \vee z) \eq (x \wedge y) \vee (x \wedge z) \quad & \textrm{(D$_\wedge$)} & x \vee (y \wedge z) \eq (x \vee y) \wedge (x \vee z) \\
\textrm{(B$_\vee$)} & x \vee \top \eq \top & \textrm{(B$_\wedge$)} & x \wedge \bot \eq \bot \\
\textrm{(C$_\vee$)} & x \vee x' \eq \top & \textrm{(C$_\wedge$)} & x \wedge x' \eq \bot.
\end{array}$$

Let $\mathcal{BA}$ denote the class of Boolean algebras, and let $\to$ the operation defined on any Boolean algebra by setting $x \to y \eq x' \vee y$. The \emph{strong completeness theorem}\index{Completeness!-- theorem} for $\mathit{CPL}$ states that for every subset $\Phi \cup \{\psi\}$ of $\mathit{Fm_{CPL}}$
\begin{center}
$\Phi \vdash_{CPL} \psi$ iff $\{\phi \eq \top \mid \phi \in \Phi \} \models_\mathcal{BA} \psi \eq \top$.
\end{center}
Conversely, by the \emph{inverse strong completeness theorem}\index{Completeness!Inverse -- theorem} for $\mathit{CPL}$, for every set of equations $E \cup \{s \eq t\}$ over $\mathit{Fm_{CPL}}$,
$$\textrm{$E \models_\mathcal{BA} s \eq t$ iff $\{u \to v, v \to u \mid u \eq v \in E\} \vdash_{\mathit{CPL}} \{s \to t, t \to s \}$.}$$
Furthermore, for every $\psi \in \mathit{Fm_{CPL}}$ and  $s \eq t \in \Eq_{\mathit{CPL}}$,
\begin{eqnarray}
&s \eq t \Sse _\mathcal{BA} \{s \to t \eq \top, t \to s \eq \top\},& \nonumber \\ 
&\psi \dashv \vdash_{\mathit{CPL}} \{\psi \to \top, \top \to \psi\},&  \nonumber
\end{eqnarray}
where $X \dashv \vdash Y$ denotes the conjunction of $X \vdash Y$ and $Y \vdash X$.
  
If for every $\psi \in \mathit{Fm_{CPL}}$ and  $\e = (s \eq t) \in \Eq_{\mathit{CPL}}$ we define $\tau(\psi) = \{\psi \eq \top\}$ and $\rho(\e) = \rho (s \eq t) = \{s \to t, t \to s \}$, then the above take the form
\begin{enumerate}
\item[-] $\Phi \vdash_{\mathit{CPL}} \psi$ iff $\tau(\Phi) \models_\mathcal{BA} \tau(\psi)$,
\item[-] $E \models_\mathcal{BA} \e$ iff $\rho(E) \vdash_{\mathit{CPL}} \rho(\e)$,
\item[-] $\e \Sse_\mathcal{BA} \tau\rho(\e)$,
\item[-] $\psi \dashv \vdash_{\mathit{CPL}} \rho\tau (\psi)$.
\end{enumerate}
Observe that the last two statements mean that the application of $\tau$ and $\rho$ one after the other might not yield the original formula or equation, but it will give a set of formulas or equations that are mutually deducible with the original ones. Then the class $\mathcal{BA}$ of Boolean algebras is an equivalent algebraic semantics for the Classical Propositional Logic presented in Example~\ref{cpl}.
\end{example}

\begin{example}\label{mvalg}
An algebra $\AA = \la A, \oplus, ^*, 0 \ra$ of type (210), is called an \emph{MV-algebra} iff it satisfies the following equations: 
\begin{enumerate}
\item[(MV1)]$x \oplus (y \oplus z) = (x \oplus y) \oplus z$;
\item[(MV2)]$x \oplus y = y \oplus x$;
\item[(MV3)]$x \oplus 0 = x$;
\item[(MV4)]$(x^*)^* = x$;
\item[(MV5)]$x \oplus 0^* = 0^*$;
\item[(MV6)]$(x^* \oplus y)^* \oplus y = (y^* \oplus x)^* \oplus x$.
\end{enumerate}
It has been proved in many different ways~--- see~\cite{chang1,chang2,cignoli,mundici2,panti,roserosser}~--- that the \L ukasiewicz Propositional Logic, presented in Example~\ref{lpl}, is complete with respect to the class $\cat{MV}$ of MV-algebras. Also the inverse completeness theorem holds (see~\cite{mundici} for an extensive study on MV-algebras). As for $\mathit{CPL}$ and Boolean algebras, we can define the operation $\to$ in any MV-algebra by setting $x \to y \eq x^* \oplus y$ and the constant $1 \eq 0^*$. Thus we obtain the maps
\begin{enumerate}
\item[$\tau:$]$\psi \in \mathit{Fm}_{\luk} \lmapsto \{\psi \eq 1\} \in \wp(\Eq_{\luk})$,
\item[$\rho:$]$s \eq t \in \Eq_{\luk} \lmapsto \{s \to t, t \to s\} \in \wp(\mathit{Fm}_{\luk})$,
\end{enumerate}
and we have, for all $\Phi \cup \{\psi\} \in \wp(\mathit{Fm}_{\luk})$ and $E \cup \{\e\} \in \wp(\Eq_{\luk})$,
\begin{enumerate}
\item[-] $\Phi \vdash_{\luk} \psi$ iff $\tau(\Phi) \models_\mathcal{MV} \tau(\psi)$,
\item[-] $E \models_\mathcal{MV} \e$ iff $\rho(E) \vdash_{\luk} \rho(\e)$,
\item[-] $\e \Sse_\mathcal{MV} \tau\rho(\e)$,
\item[-] $\psi \dashv \vdash_{\luk} \rho\tau (\psi)$.
\end{enumerate}
Then the class $\mathcal{MV}$ is an equivalent algebraic semantics for the propositional logic \emph{\L}$\mathit{PL}$ presented in Example~\ref{lpl}.
\end{example}

\begin{example}\label{bckalg}
It is shown in~\cite{blokpigozzi} that $\vdash_{\mathit{BCK}}$ is algebraizable and the $\{\to\}$-subreducts of commutative integral residuated lattices (we will introduce residuated lattices in Section~\ref{reslatquantsec}) form an algebraic semantics for it. The corresponding maps $\tau$ and $\rho$ are given by $\tau(\psi) = \{\psi \eq (\psi \to \psi)\}$ and  $\rho(u \eq v) = \{u \to v, v \to u\}$.
\end{example}
It has been proved by W. J. Blok and D. Pigozzi (see Corollary 5.6 of~\cite{blokpigozzi}) that the propositional logic $S3$ of Example~\ref{s3} is not algebraizable.

As for the case of propositional logics, we define the lattice of theories $\Th_{\models}$ of a consequence relation $\models$ over $\Eq_\lang$ to be the lattice of the sets of equations closed under $\models$, that is the sets $E$ of equations such that $E \models \e$ implies $\e \in E$. The notions of finitarity and substitution invariance have analogues for closure operators and lattices of theories. We discuss the  connections between consequence relations, closure operators and lattices of theories in a more general setting in Section~\ref{abstconsrelsec}. Recalling that the lattices of (equational or propositional) theories are closed under inverse substitutions, i.e. $\sigma^{-1}(t)$ is an element of the lattice for any substitution $\sigma \in \sfm$ and for any element $t$ of the considered lattice of theories, we have the following characterization of algebraizability of a deductive system.
\begin{theorem}\emph{\cite{blokpigozzi}}\label{algeqv}
A deductive system $\la \lang, \vdash \ra$ is algebraizable with equivalent algebraic semantics a quasivariety $\mathcal{K}$ iff there exists an isomorphism $T$ between $\Th_\vdash$ and $\Th_{\models_\mathcal{K}}$ that commutes with inverse substitutions, i.e. such that $T(\sigma^{-1}(t)) = \sigma^{-1}(T(t))$, for any substitution $\sigma$ over $\lang$.
\end{theorem} 

For the sake of being formally coherent and complete, we close this section with a more precise and formal (but equivalent to the one given above) definition of equational deductive system, according to~\cite{blokpigozzi2}. A \emph{quasi-equation}\index{Equation!Quasi --}\index{Quasi-equation} over $\lang$ is the equational concept corresponding to that of inference rule; it is a pair $\la E,\e \ra$ where $E$ is a finite set of equations and $\e$ is a single equation. An equation $\d$ is \emph{directly derivable} from a set $D$ of equations by the quasi-equation $\la E, \e \ra$ if there is a substitution $\sigma$ such that $\sigma \e = \d$ and $\sigma[E] \subseteq D$. A quasi-equation $\la E,\e \ra$ is usually denoted by $\e_1 \& \ldots \&\e_n \Rightarrow \e$, where $\e_1, \ldots, \e_n$ are all the elements of $E$. The \emph{axioms} of an equational deductive systems are simply $\lang$-equations.
\begin{definition}\label{eqdedsys}
An \emph{equational deductive system}\index{Deductive system!Equational --} $\cat S$ over a given language $\lang$, is defined by means of a (possible infinite) set of quasi-equations and axioms. It consists of the pair $\cat S = \la \lang, \vDash \ra$, where $\vDash$ is a subset of $\wp(\Eq_\lang) \times \Eq_\lang$ defined by the following condition: $E \vDash \e$ iff $\e$ is contained in the smallest set of equations that includes $E$ together with all substitution instances of the axioms of $\cat S$, and is closed under direct derivability by the quasi-equations~of~$\cat S$.
\end{definition}

\section{Gentzen-style systems}
\label{seqlogicsec}

As we anticipated in the previous section, we will now generalize the concept of deductive system in such a way that the new definition will include, as special cases, those of propositional and equational deductive systems.

Let $\lang$ be a propositional language and let $m, n \in \N_0$ such that at least one of them is positive. A (\emph{classical, associative}) \emph{sequent over $\lang$ of type $(m,n)$}\index{Sequent}, is a pair $(\Phi, \Psi)$ composed by two sequences of $\lang$-formulas $\Phi = (\phi_1, \phi_2, \dots, \phi_m)$, of length $m$, and $\Psi = (\psi_1, \psi_2, \dots, \psi_n)$, of length $n$. For the sequent $(\Phi, \Psi)$, also the notation $\phi_1, \phi_2, \dots, \phi_m \Rightarrow \psi_1, \psi_2, \dots, \psi_n$ is often used.

Usually, by a \emph{set of sequents} $\mathit{Seq}_\lang$, it is understood a set of sequents that is \emph{closed under type}, namely a set of sequents such that, for all $m, n$, if it contains an $(m,n)$-sequent, then it contains all $(m,n)$-sequents. If $\mathit{Seq}_\lang$ is a set of sequents, then $\Tp(\mathit{Seq}_\lang) \subseteq \N_0 \times \N_0$ denotes the set of all types of the sequents in $\mathit{Seq}_\lang$. If $s = (\Phi, \Psi)$ is a sequent and $\sigma \in \sfm$ is a substitution, the sequent $(\sigma(\Phi), \sigma(\Psi))$ is denoted by $\sigma(s)$.

As first, immediate, examples we observe that the set $\fml$ can be identified with the set of all $(0,1)$-sequents, and the set $\Eq_\lang$ can be identified with the set of all $(1,1)$-sequents.

\begin{definition}\label{seqconsrel}
In analogy with the cases of $\fml$ and $\Eq_\lang$, if $\mathit{Seq}_\lang$ is a set of sequents, we define an \emph{asymmetric consequence relation over $\mathit{Seq}_\lang$} as a subset $\vdash$ of $\wp(\mathit{Seq}_\lang) \times \mathit{Seq}_\lang$ such that conditions (\ref{entailin},\ref{entailmp}) hold. A consequence relation over $\mathit{Seq}_\lang$ is called \emph{finitary} if (\ref{entailfinitary}) hold, and is called \emph{substitution invariant} if (\ref{entailstructural}) is verified.

Again, a \emph{symmetric consequence relation over $\mathit{Seq}_\lang$} is a binary relation over $\wp(\mathit{Seq}_\lang)$ such that (\ref{sentailin}--\ref{sentailjoin}) hold; it is called \emph{finitary} if it satisfies (\ref{sentailfinitary}) and \emph{substitution invariant} if it satisfies (\ref{sentailstructural}).

A deductive system over a set of sequents, that is a set of sequents endowed with a finitary and substitution invariant consequence relation, is also called a \emph{Gentzen-style}\index{Deductive system!Gentzen-style --}\index{Deductive system!-- over a set of sequents} system.
\end{definition}

The notion of algebraizability of a set $\mathit{Seq}$ of sequents closed under type has been defined by Rebagliato and Verd\'u~\cite{rebver}. If $\mathit{Seq}_1$ and $\mathit{Seq}_2$ are sets of sequents over $\lang$, and $\vdash_1$ and $\vdash_2$ are two consequence relations over $\mathit{Seq}_1$ and $\mathit{Seq}_2$, respectively, a \emph{translation}\index{Translation} between $\mathit{Seq}_1$ and $\mathit{Seq}_2$ is a set $\tau = \{\tau_{(m,n)} \mid (m,n) \in \Tp(\mathit{Seq}_1)\}$, where $\tau_{(m,n)}$ is a finite subset of $\mathit{Seq}_2$ in (at most) $m + n$ variables. If $s \in \mathit{Seq}_1$ is an $(m,n)$-sequent, $\tau(s) = \tau_{(m,n)}(s)$ denotes the result of replacing the variables in $\tau_{(m,n)}$ by the $m + n$ formulas of $s$.

\begin{definition}\label{seqeq}
Two consequence relations $\vdash_1$ and $\vdash_2$ over $\mathit{Seq}_1$ and $\mathit{Seq}_2$, respectively, are called \emph{equivalent} (in the sense of Rebagliato and Verd\'u), if there are translations $\tau$ and $\rho$ between $\mathit{Seq}_1$ and $\mathit{Seq}_2$ such that for all subsets $S_1 \cup \{s_1\}$ of $\mathit{Seq}_1$ and all subsets $S_2 \cup \{s_2\}$ of $\mathit{Seq}_2$,
\begin{enumerate}
\item[$(i)$]$S_1 \vdash _1 s_1$ iff $\tau(S_1) \vdash_2 \tau(s_1)$ and
\item[$(ii)$]$s_2 \dashv \vdash_2 \tau \rho (s_2)$.
\end{enumerate}
It follows that 
\begin{enumerate}
\item[$(iii)$]$S_2 \vdash _2 s_2$ iff $\rho(S_2) \vdash_1 \rho(s_2)$ and
\item[$(iv)$]$s_1 \dashv \vdash_1 \rho\tau(s_1)$.
\end{enumerate}
\end{definition}

\begin{lemma}\emph{\cite{galatostsinakis}}
Consider maps $\tau: \mathit{Seq}_1 \lto \wp(\mathit{Seq}_2)$ and $\rho: \mathit{Seq}_2 \lto \wp(\mathit{Seq}_1)$. The following are equivalent.
\begin{enumerate}
\item[$(a)$]The maps $\tau$, $\rho$ are finitary and substitution invariant.
\item[$(b)$]There exist translations $\tau'$ and $\rho'$ between $\mathit{Seq}_1$ and $\mathit{Seq}_2$ such that $\tau'(s_1) = \tau(s_1)$ and $\rho'(s_2) = \rho(s_2)$ for all $s_1 \in \mathit{Seq}_1$ and $s_2 \in \mathit{Seq}_2$.
\end{enumerate}
\end{lemma}
In Chapter~\ref{consrelchap} we will show a unification of all the notions and results defined so far. 

\section{Consequence relations on powersets}
\label{conspwsetsec}

In Sections~\ref{proplogicsec},~\ref{eqlogicsec} and~\ref{seqlogicsec} we have defined, according to the tradition, asymmetric and symmetric consequence relations on sets, respectively, of formulas, equations and sequents over a propositional language. We will see that all these notions can be reformulated in a more general categorical setting but, in order to operate such a generalization, we need notions and results that we will present in Chapters~\ref{orderchap} and~\ref{mqchap}; therefore we postpone it to Chapter~\ref{consrelchap}.

By the way, we can give a glimpse of how this abstraction will work, and make a first step in that direction, by showing a first generalization of consequence relations to powersets. First of all, we define~--- as the reader may expect~--- asymmetric consequence relations on an arbitrary set.

Let $S$ be a set. An \emph{asymmetric consequence relation over $S$} is a subset $\vdash$ of $\wp(S) \times S$ such that, for all subsets $X \cup Y \cup \{x, y, z\}$ of $S$,
\begin{eqnarray}
&&\textrm{if $x \in X$, then $X \vdash x$,} \label{wpentailin} \\
&&\textrm{if $X \vdash y$, for all $y\in Y$, and $Y \vdash z$, then $X \vdash z$.} \label{wpentailmp}
\end{eqnarray}
An asymmetric consequence relation over $S$ is called \emph{finitary}, if for all subsets $X \cup \{x\}$ of $S$,
\begin{eqnarray}
&&\textrm{if $X \vdash x$, then there is a finite subset $X_0$ of $X$ such that $X_0 \vdash x$.} \qquad \label{wpentailfinitary}
\end{eqnarray}
 
Now, in order to generalize the notion of substitution invariance to arbitrary powersets, we need something that somehow generalizes the role of substitutions. Then we can observe that the substitution monoid $\Sfm$ acts on $\fml$, $\Eq_\lang$ and $\mathit{Seq}_\lang$ in the sense that, for all $\sigma_1, \sigma_2 \in \sfm$ and $s$ in either $\fml$, $\Eq_\lang$ or $\mathit{Seq}_\lang$,
\begin{enumerate}
\item[$(i)$]$(\sigma_1 \sigma_2)(s) = \sigma_1(\sigma_2(s))$,
\item[$(ii)$]$\id_{\sfm}(s) = s$.
\end{enumerate}
Conditions $(i)$ and $(ii)$ are precisely the ones that define an \emph{action of a monoid over a set}\index{Action of a monoid}\index{Monoid!-- action}. Indeed, a monoid $\MM = \la M, \cdot, e \ra$ is said to act on a set $S$, if there exists a map $\star: M \times S \lto S$ such that for all $m_1, m_2 \in M$ and $x \in S$,
\begin{enumerate}
\item[$(i)$]$(m_1 \cdot m_2) \star x = m_1 \star (m_2 \star x)$
\item[$(ii)$]$e \star x = x$.
\end{enumerate}

Then we say that an asymmetric consequence relation $\vdash$ on $S$ is \emph{$\MM$-invariant} iff, for all $X \cup \{y\} \subseteq S$ and $m \in M$, $X \vdash y$ implies $\{m \star x \mid x \in X\} \vdash m \star y$.

It is immediate to verify that, for any monoid $\MM \la M, \cdot, e \ra$, the structure $\wp(\MM) = \la \wp(M), \cdot, \{e\} \ra$~--- where $A \cdot B := \{a \cdot b \mid a \in A, b \in B\}$, for all $A, B \in \wp(M)$~--- is a monoid as well. Moreover, if $\MM$ acts on $S$, then $\wp(\MM)$ acts on $\wp(S)$, i.e. there exists a map $\star': \wp(M) \times \wp(S) \lto \wp(S)$ such that for all $A, B \in \wp(M)$ and $X \in \wp(S)$,
\begin{enumerate}
\item[$(i)$]$(A \cdot B) \star' X = A \star' (B \star' X)$
\item[$(ii)$]$\{e\} \star' X = X$,
\end{enumerate}
where $A \star' X = \{a \star x \mid a \in A, x \in X\}$. It is easy to see also that $\star'$ preserves arbitrary unions in both its arguments. In what follows we will denote by $\star$ both the action of a monoid on a set and its extension to the respective powersets.

If $S_1$ and $S_2$ are sets over which a monoid $\MM$ acts (by $\star_1$ and $\star_2$ respectively), a map $\tau: \wp(S_1) \lto \wp(S_2)$ is called \emph{$\wp(\MM)$-invariant}, if for all $A \in \wp(M)$ and $X \in \wp(S_1)$, we have $A \star_2 \tau(X) = \tau (A \star_1 X)$.

With all this setting, and the notations we have just introduced, we can suitably define finitary and structural maps in this case. Assume that $\vdash_1$ and $\vdash_2$ are asymmetric consequence relations on $S_1$ and $S_2$, respectively. If there exist maps $\tau: \wp(S_1) \lto \wp(S_2)$ and $\rho: \wp(S_2) \lto \wp(S_1)$ that preserve unions and such that for every subset $X \cup \{x\}$ of $S_1$ and $y \in S_2$,
\begin{enumerate}
\item[$(i)$]$X \vdash_1 x$ iff $\tau(X) \vdash_2 \tau(x)$,
\item[$(ii)$]$y \dashv \vdash_2 \tau\rho(y)$,
\end{enumerate}
then we will say that $\vdash_1$ and $\vdash_2$ are \emph{similar via $\tau$ and $\rho$}\index{Consequence relation!Similar --s}. We will show in Lemma~\ref{sim} that in this case $\vdash_2$ and $\vdash_1$ are similar via $\rho$ and $\tau$, as well. 
  
Now, let $\MM$ be a monoid that acts on $S_1$ and $S_2$. If $\vdash_1$ and $\vdash_2$ are similar via $\tau$ and $\rho$, and both $\tau$ and $\rho$ are $\wp(\MM)$-invariant, then we say that $\vdash_1$ and $\vdash_2$ are \emph{equivalent via $\tau$ and $\rho$}\index{Consequence relation!Equivalent --s}.

The generalization of symmetric consequence relations to the case of powersets is immediate as well. A \emph{symmetric consequence relation over $S$} is a binary relation $\vdash$ on $\wp(S)$ that satisfies, for all $X, Y, Z \in \wp(S)$,
\begin{eqnarray}
&&\textrm{if $Y \subseteq X$, then $X \vdash Y$,} \label{wpsentailin} \\
&&\textrm{if $X \vdash Y$ and $Y \vdash Z$, then $X \vdash Z$,} \label{wpsentailimp} \\
&&X \vdash \bigcup_{X \vdash Y} Y. \label{wpsentailjoin}
\end{eqnarray}
Note that $\vdash$ satisfies the first two conditions iff it is a pre-order on $\wp(S)$ that contains the relation $\supseteq$.

A symmetric consequence relation over $S$ is called \emph{finitary} provided that, for all $X, Y \in \wp(S)$, if $X \vdash Y$ and $Y$ is finite, then there is a finite subset $X_0$ of $X$ such that $X_0 \vdash Y$. If $\MM$ is a monoid acting on $S$, we will say that $\vdash$ is \emph{$\wp(\MM)$-invariant} if, for all $X, Y \in \wp(S)$ and $A \in \wp(\MM)$,  $X \vdash Y$ implies $A \star X \vdash A \star Y$.

As we anticipated, given an asymmetric consequence relation $\vdash$, we can define its symmetric counterpart~--- that we will denote by $\vdash^s$~--- by $X \vdash^s Y$ iff $X \vdash y$ for all $y \in Y$, for any given $X, Y \in \wp(S)$. Conversely, given a symmetric consequence relation $\vdash$, we define its asymmetric counterpart $\vdash^a$ by $X \vdash x$ iff $X \vdash \{x\}$, for $X \in \wp(S)$ and $x \in S$. 
 
\begin{lemma}\emph{\cite{galatostsinakis}}\label{symasym}
Symmetric consequence relations on $\wp(S)$, where $S$ is a set, are in bijective correspondence with asymmetric consequence relations on $\wp(S)$ via the maps $\vdash \ \lmapsto \ \vdash^a$ and $\vdash \ \lmapsto \ \vdash^s$. Moreover finitarity is preserved under these maps and the same holds for $\wp(\MM)$-invariance (resp., $\MM$-invariance), if a monoid $\MM$ acts on $S$.
\end{lemma}
\begin{proof}
Let $\vdash$ be an asymmetric consequence relation and let $X, Y, Z \in \wp(S)$. If $Y \subseteq X$, then $y \in X$, for all $y \in Y$; so $X \vdash y$, for all $y \in Y$, hence $X \vdash^s Y$. Also, if $X \vdash^s Y$ and $Y \vdash^s Z$, then, for all $z \in Z$, $X \vdash y$ for all $y \in Y$ and $Y \vdash z$. Therefore $X \vdash z$, for all $z \in Z$; hence $X \vdash^s Z$. If $y \in \bigcup_{X \vdash^s Y} Y$, then $y \in Y$, for some $Y$ such that $X \vdash^s Y$; so $X \vdash y$. Consequently, $X \vdash^s  \bigcup_{X \vdash^s  Y} Y$.

Assume, now, that $\vdash$ is a symmetric consequence relation and let $X \cup Y \cup \{x, y, z\} \in \wp(S)$. If $x \in X$, then $\{x\} \subseteq X$, so $X \vdash \{x\}$; i.e. $X \vdash^a x$. If $X \vdash^a y$, for all $y \in Y$, and $Y \vdash ^a z$, then $X \vdash \{y\}$ for all $y \in Y$, and $Y \vdash \{z\}$. Note that $Y \subseteq \bigcup_{X \vdash Z} Z$, so $\bigcup_{X \vdash Z} Z \vdash Y$. Since $X \vdash \bigcup_{X \vdash Z} Z$, we have $X \vdash Y$. Consequently, $X \vdash \{z\}$; i.e. $X \vdash^a z$.
 
Let $\vdash$ be an asymmetric consequence relation and let $X \cup \{x\} \in \wp(S)$. We  have $X \vdash^{sa} x$ iff $X \vdash^s \{x\}$ iff $X \vdash x$. So $\vdash^{sa} = \vdash$. Conversely, let $\vdash$ be a symmetric consequence relation and let $X \cup Y \cup \{x, y, z\} \in \wp(S)$. We have $X \vdash^{as} Y$ iff $X \vdash^a y$, for all $y \in Y$, iff $X \vdash \{y\}$, for all $y \in Y$. Note that $Y = \bigcup_{y \in Y} \{y\} \subseteq \bigcup_{X \vdash Z} Z$, so $\bigcup_{X \vdash Z} Z \vdash Y$. Since $X \vdash \bigcup_{X \vdash Z} Z$, we have that $X \vdash \{y\}$, for all $y\in Y$, implies $X \vdash Y$. Conversely, if $X \vdash Y$, then $X \vdash \{y\}$, for all $y\in Y$, since $Y \vdash \{y\}$, for all $y\in Y$. So $\vdash^{as} = \vdash$.

Let $\vdash$ be a finitary asymmetric consequence relation and let $X \vdash^s Y$, for some finite set $Y$. Then, $X \vdash y$, for all $y \in Y$. By the fact that $\vdash$ is finitary, for all $y \in Y$, there is a finite subset $X_y$ of $X$ such that $X_y \vdash y$. Note that the set $X_0 = \bigcup_{y \in Y} X_y$ is finite and $X_0 \vdash y$, for all $y \in Y$; so $X_0 \vdash^s Y$. Thus, $\vdash^s$ is finitary, as well. Conversely, let $\vdash$ be a finitary symmetric consequence relation and let $X \vdash^a y$. Then $X \vdash \{y\}$ and, by the finitariness of $\vdash$, there exists a finite subset $X_0$ of $X$ such that $X_0 \vdash \{y\}$; i.e. $X_0 \vdash^a y$. Thus, $\vdash^a$ is finitary, as well. 

Now let $\MM$ be a monoid that acts on $S$, and assume that $\vdash$ is $\MM$-invariant. Then, if $X, Y \subseteq S$ are such that $X \vdash y$ for all $y \in Y$, $\{m \star x \mid x \in X\} = \{m\} \star X \vdash m \star y$, for all $m \in M$ and $y \in Y$. By (\ref{wpentailin}), we have $M \star X \vdash m \star y$, for all $m \in M$ and $y \in Y$, so $M \star X \vdash M \star Y$.

Conversely, if $\vdash$ is a $\wp(\MM)$-invariant symmetric consequence relation on $S$, and $X \vdash^a y$, then $X \vdash \{y\}$. But $\vdash$ is $\wp(\MM)$-invariant, so $\{m\} \star X \vdash \{m\} \star \{y\}$; hence $m \star X \vdash^a m \star y$.
\end{proof}
Thanks to Lemma~\ref{symasym}, we can overcome the distinction between asymmetric and symmetric consequence relations. Thus, in what follows, we will refer to \emph{consequence relations} without specifying whether they are symmetric or asymmetric. Last, we remark once more that all the definitions and results we showed for consequence relations over sets of formulas, equations and sequents, are special cases of those we presented in this section.

\section*{References and further readings}
\addcontentsline{toc}{section}{References and further readings}

The reader may have already noticed that the works~\cite{blokjonsson2},~\cite{blokpigozzi} and~\cite{galatostsinakis} have been repeatedly cited in this chapter. Actually, most of the contents of this chapter can be found in many scientific papers and books, but the aforementioned three works have been our main guide in drawing up this preliminaries on propositional logics. Besides those we already cited in the chapter, we point out, as suggested readings, some further works.

Substructural logics are some of the best known examples of deductive systems over sets of sequents. Relevant works on this subject are the book by N. Galatos, P. Jipsen, T. Kowalski and H. Ono,~\cite{nickbook}, and the papers~\cite{galatosono2,galatosono}, by N. Galatos and H. Ono, and~\cite{ono}, by H. Ono. Many examples of algebraizable and non-algebraizable deductive systems can be found in the books~\cite{gerla}, by G. Gerla, and~\cite{hajek}, by P. H\'ajek. Last we cite~\cite{yde}, by P. Blackburn, M. de Rijke and Y. Venema, as a reference book on modal logics.


\chapter{An Outline of Categorical Tools}
\label{catschap}

In the present chapter we will briefly introduce some basic notions and results of Category Theory. Their exposition, far from being exhaustive and detailed, is exploitable to the introduction of the main topics of this thesis, and~--- therefore~--- only the bare necessities will be presented. For these reasons, a reader that is already familiar with the basic notions and constructions of category theory may even skip the whole chapter.

Category Theory was born with the aim of providing a language that was able to describe, with precision, many similar constructions and phenomena that occur in different mathematical fields. For example products of structures can be defined for vector spaces, groups, topological spaces, Banach spaces, automata; free objects are defined for lattices, vector spaces, modules, and so on. Further advantages of Category Theory are, for instance, the symbolism that allows to quickly visualize quite complicated facts by means of diagrams, the concept of ``functor''~--- a vehicle that allows one to transport problems from one area of mathematics to another one~--- and the ``duality principle'', that category theorists also like to call the ``two for the price of one'' principle, since it states, essentially, that every concept is two concepts, and every result is two results.

At the very beginning we need to establish some foundational aspects. The basic concepts that we need are those of \emph{sets} and \emph{classes}\index{Class}. Sets can be thought of as the usual sets of some axiomatic set theory (Zermelo--Fr\ae nkel\index{Set theory!Zermelo--Fr\ae nkel --}, von Neumann--Bernays--G\"odel\index{Set theory!vonNeumann--Bernays--G\"odel --}, Morse--Kelly\index{Set theory!Morse--Kelly --}, etc.) or even of Cantor ``na\"ive''\index{Set theory!Cantor --} set theory. The concept of class has been created to deal with ``large collections of sets''. In particular, we require that:
\begin{enumerate}
\item[(1)] the members of each class are sets,
\item[(2)] for every property $P$ the class of all sets satisfying $P$ can be formed.
\end{enumerate}
Hence there is the largest class: the class of all sets, called the \emph{universe}\index{Class!Universal --}\index{Universe} (or \emph{universal class}) and denoted by $\U$. Classes are precisely the subcollections of $\U$. Thus, given classes $A$ and $B$, we can build such classes as $A \cup B$, $A \cap B$, and $A \times B$. Then it is possible to define functions between classes, equivalence relations on classes, etc. A \emph{family} $(A_i)_{i \in I}$ of sets is a function $A : I \lto \U$ sending $i \in I$ to $A(i) = A_i$, where $I$ is a class. The following conditions are required for convenience
\begin{enumerate}
\item[(3)] if $X_1, X_2, \ldots, X_n$ are classes, then so is the $n$-tuple $(X_1, X_2, \ldots, X_n)$,
\item[(4)] every set is a class (equivalently: every member of a set is a set).
\end{enumerate}
Hence sets are special classes. Classes that are not sets are called \emph{proper classes}\index{Class!Proper --}. They cannot be members of any class. Examples of proper classes are the universe $\U$, the class of all vector spaces, the class of all topological spaces, and the class of all automata.

Notice that, in this setting, condition $(4)$ above gives us the \emph{Axiom of Replacement}\index{Axiom!-- of Replacement}:
\begin{enumerate}
\item[(5)] there is no surjection from a set to a proper class.
\end{enumerate}
Therefore sets are also called \emph{small classes}\index{Class!Small --}, and proper classes are called \emph{large classes}\index{Class!Large --}.

The framework of sets and classes described so far is not strong enough for treating all the theory of categories. Nevertheless it suffices for our purposes and, therefore, we will not go any further.

\section{Categories and functors}
\label{catfuncsec}

\begin{definition}\label{category}
A \emph{category}\index{Category} is a quadruple $\cat C = (O, \hom, \id, \circ)$ consisting of
\begin{enumerate}
\item[$(a)$] a class $O$\index{Object!-- class}\index{Class!Object --}, whose members are called \emph{objects}\index{Object} of $\cat C$,
\item[$(b)$] for each pair $(A,B)$ of $\cat C$-objects, a set $\hom(A,B)$, whose members are called $\cat C$-\emph{morphisms}\index{Morphism} from $A$ to $B$ (the statement ``$f \in \hom(A,B)$'' is expressed more graphically by using arrows; e.g., by statements such as ``$f: A \lto B$ is a
morphism'' or ``$A \stackrel{f}{\lto} B$ is a morphism''),
\item[$(c)$] for each object $A$, a morphism $\id_A: A \lto A$, called the \emph{identity}\index{Morphism!Identity --}\index{Identity} of $A$,
\item[$(d)$] a \emph{composition law} associating with each $\cat C$-morphism $f: A \lto B$ and each $\cat C$-morphism $g: B \lto C$ a $\cat C$-morphism $g \circ f: A \lto C$, called the composite of $f$ and $g$, subject to the following conditions:
\begin{enumerate}
\item[$(i)$] composition is associative; i.e., for morphisms $f: A \lto B$, $g: B \lto C$, and $h: C \lto D$, the equation $h \circ (g \circ f) = (h \circ g) \circ f$ holds,
\item[$(ii)$] identities act as units with respect to composition; i.e., for any morphism $f: A \lto B$, we have $\id_B \circ f = f$ and $f \circ \id_A = f$,
\item[$(iii)$] the sets $\hom(A,B)$ are pairwise disjoint.
\end{enumerate}
\end{enumerate}
\end{definition}

\begin{remark}
Let $\cat C = (O, \hom, \id, \circ)$ be a category.
\begin{itemize}
\item The class $O$ of $\cat C$-objects is usually denoted by $\obj(\cat C)$\index{Object!-- class}\index{Class!Object --}.
\item The class of all $\cat C$-morphisms, denoted by $\Mor(\cat C)$, is defined to be the union of all the sets $\hom(A,B)$ in $\cat C$.
\item If $f: A \lto B$ is a $\cat C$-morphism, we call $A$ the \emph{domain}\index{Morphism!Domain of a --} of $f$ (denoted by $\dom(f)$) and $B$ the \emph{codomain}\index{Morphism!Codomain of a --} of $f$ (denoted by $\cod(f)$). Observe that condition $(iii)$ guarantees that each $\cat C$-morphism has a unique domain and a unique codomain. However, this condition is given for technical convenience only, because whenever all other conditions are satisfied, it is easy to ``force'' condition $(iii)$ by simply replacing each morphism $f \in \hom(A,B)$ by a triple $(A, f, B)$.
\item The composition, $\circ$, is a partial binary operation on the class $\Mor(\cat C)$. For a pair $(f, g)$ of morphisms, $f \circ g$ is defined if and only if the domain of $f$ and the codomain of $g$ coincide.
\item If more than one category is involved, subscripts may be used (as in $\hom_\cat C(A,B)$) for clarification.
\item If $\cat C$ is a category and $A$ is an object of $\cat C$, we should write $A \in \obj(\cat C)$ but, with an abuse of notation, we will often use the simplified expression $A \in \cat C$.
\end{itemize}
\end{remark}

Best known examples of categories are the following ones
\begin{enumerate}
\item[$\mathcal{S}$:] the category whose objects are sets and morphisms are maps between them,
\item[$\mathcal{G}$:] with the class of groups as $\obj(\cat G)$ and group homomorphisms as morphisms,
\item[$\mathcal{OS}$:] the category of ordered sets, where the morphisms are the order preserving maps,
\end{enumerate}
while the following will be of interest for this thesis
\begin{enumerate}
\item[$\Q$:]\index{Category!-- $\Q$, of quantales} the category whose object are \emph{quantales}, i.e. complete residuated lattice-ordered monoids (see Section~\ref{reslatquantsec} for details), and whose morphisms are the maps that preserve arbitrary joins, the monoid operation, the unit and the bottom element,
\item[$\SL$:]\index{Category!-- $\SL$, of sup-lattices} with the class of complete lattices as $\obj(\SL)$ and the maps preserving the bottom element and arbitrary joins as morphisms~--- we will study such category in details in Sections~\ref{suplatsec} and~\ref{sltensorsec},
\item[$\QQ\Mod$:]\index{Category!-- $\QQ\Mod$, of left $\QQ$-modules} given a quantale $\QQ$, the category $\QQ\Mod$ has left $\QQ$-modules as objects and $\QQ$-homomorphisms as morphisms (see Definition~\ref{modules}),
\item[$\Modd\QQ$:]\index{Category!-- $\Modd\QQ$, of right $\QQ$-modules} is analogous to $\QQ\Mod$, $\obj(\Modd\QQ)$ being the class of right $\QQ$-modules.
\end{enumerate}

For any category $\cat C = (O, \hom_\cat C, \id, \circ)$, the \emph{dual}\index{Category!Dual --}\index{Dual!-- category} (or \emph{opposite}) category of $\cat C$ is the category $\cat C^{\op} = (O, \hom_{\cat C^{\op}}, \id, \circ^{\op})$, where $\hom_{\cat C^{\op}}(A,B) = hom_\cat C(B,A)$ and $f \circ^{\op} g = g \circ f$. Thus $\cat C$ and $\cat C^{\op}$ have the same objects and, except for their direction, the same morphisms.

Because of the way dual categories are defined, every statement $s_{\cat C^{\op}}(X)$ concerning an object $X$ in the category $\cat C^{\op}$ can be translated into a logically equivalent statement $s_{\cat C}^{\op}(X)$ concerning the object $X$ in the category $\cat C$. In particular, any property $\cat P$~--- regarding objects or morphisms in a category $\cat C$~--- can be translated in an obvious way into a property $\cat P^{\op}$ in such a way that, if $\cat P$ holds for some objects or morphisms in $\cat C$, then $\cat P^{\op}$ holds for the same objects or morphisms in $\cat C^{\op}$. Let us make an example that will clarify this situation.

Let $X$ be an object of a category $\cat C$ and let $\cat P_{\cat C}(X)$ be the following property: \emph{for all $Y \in \obj(\cat C)$ there exists a unique morphism $f: Y \lto X$}. We can translate $\cat P_{\cat C}(X)$ in the logically equivalent statement $\cat P_{\cat C^{\op}}^{\op}(X)$: \emph{for all $Y \in \obj(\cat C^{\op})$ there exists a unique morphism $f: X \lto Y$}. So, obviously, $\cat P$ holds for a certain object in a category $\cat C$ if and only if $\cat P^{\op}$ holds for the same object in the category $\cat C^{\op}$. Now we can easily prove the

\begin{proposition}[{\em Duality Principle}]\emph{\cite{cats}}\label{dualityprinciple}\index{Duality!-- Principle}
Whenever a property $\cat P$ holds for all categories, then the property $\cat P^{\op}$ holds for all categories.
\end{proposition}
\begin{proof}
Let $\cat C$ be any category. Since the property $\cat P$ holds for all categories, in particular it holds for $\cat C^{\op}$. But $\cat P$ holds for $\cat C^{\op}$ if and only if $\cat P^{\op}$ holds for $(\cat C^{\op})^{\op} = \cat C$. Then the thesis follows from the arbitrary choice of $\cat C$.
\end{proof}

Because of this principle, each result in category theory has two equivalent formulations (which at first glance might seem to be quite different). However, only one of them needs to be proved, since the other one follows by virtue of the Duality Principle.

Often the dual concept\index{Dual!-- concept} $\cat P^{\op}$ of a concept $\cat P$ is denoted by ``co-$\cat P$'' (e.g. products and coproducts). A concept $\cat P$ is called \emph{self-dual}\index{Dual!Self-{--} concept} if $\cat P = \cat P^{\op}$. An example of a self-dual concept is that of identity morphism.

\begin{definition}\label{catiso}
A morphism $f : A \lto B$ in a category $\cat C$ is called an \emph{isomorphism}\index{Isomorphism}\index{Morphism!Iso--} provided there exists a morphism $g : B \lto A$ with $g \circ f = \id_A$ and $f \circ g = \id_B$. Such a morphism $g$ is called the \emph{inverse} of $f$ and is, of course, an isomorphism as well that will be denoted by $f^{-1}$. If there exists an isomorphisms between two objects of a category, such objects are said to be \emph{isomorphic}.

A morphism $f: A \lto B$ is said to be a \emph{monomorphism}\index{Monomorphism}\index{Morphism!Mono--} provided that for all pairs $\xymatrix{C \ar@<+.4ex>[r]^h \ar@<-.4ex>[r]_k & A}$ of morphisms such that $f \circ h = f \circ k$, it follows that $h = k$ (i.e., $f$ is ``left-cancellable'' with respect to composition). Dually, a morphism $f: A \lto B$ is said to be an \emph{epimorphism}\index{Epimorphism}\index{Morphism!Epi--} provided that for all pairs $\xymatrix{B \ar@<+.4ex>[r]^h \ar@<-.4ex>[r]_k & C}$ of morphisms such that $h \circ f = k \circ f$, it follows that $h = k$ (i.e., $f$ is ``right-cancellable'' with respect to composition).
\end{definition}

It is clear from the above definition that the statement ``$f$ is an isomorphism'' is self-dual, i.e., $f$ is an isomorphism in $\cat C$ if and only if $f$ is an isomorphism in $\cat C^{\op}$.

In the definition of isomorphism, we have called the morphism $g$ ``the'' inverse of $f$ even if, a priori, there should be more than one inverse. Indeed, if $h$ is a morphism such that $h \circ f = \id_A$ and $f \circ h = \id_B$, we have $h = \id_A \circ h = (g \circ f) \circ h = g \circ (f \circ h) = g \circ \id_B = g$; the inverse morphism is unique. Moreover, if $f: A \lto B$ and $g: B \lto C$ are two isomorphisms, then it is easy to prove~--- using the definition and the associativity of $\circ$~--- that $g \circ f: A \lto C$ is an isomorphism whose inverse is $f^{-1} \circ g^{-1}$.

Originally it was believed that monomorphisms and epimorphisms would constitute the correct categorical abstractions of the notions ``embeddings of substructures'' and ``projections over quotient structures''~--- respectively~--- that exist in various constructs. However, in many instances the concepts of monomorphism and epimorphism are too weak; e.g., in the category of topological spaces (whose morphisms are continuous maps), monomorphisms are just injective continuous maps and thus need not be embeddings. Later on, stronger notions that more frequently correspond with embeddings and projections have been defined, but satisfactory concepts of ``embeddings'' and ``projections'' seem to be possible only in the setting of constructs (see Definition~\ref{concrete} for the definition of constructs). We will limit our attention to the latter case, when~--- in Section~\ref{concretesec}~--- we will define the notions of injective and projective objects.

Now, we take a more global viewpoint and consider categories themselves as structured entities. The ``morphisms'' between them that preserve their structure are called functors.

\begin{definition}\label{functor}
If $\mathcal{C}$ and $\mathcal{C}'$ are categories, then a \emph{functor}\index{Functor} $F$ from $\mathcal{C}$ to $\mathcal{C}'$ is a function that assigns to each $\mathcal{C}$-object $A$ a $\mathcal{C}'$-object $F(A)$, and to each $\mathcal{C}$-morphism $f: A \lto B$ a $\mathcal{C}'$-morphism $F(f): F(A) \lto F(B)$, in such a way that
\begin{enumerate}
\item[$(i)$] $F$ preserves composition, i.e. $F(f \circ g) = F(f) \circ F(g)$ whenever $f \circ g$ is defined,
\item[$(ii)$] $F$ preserves identity morphisms, i.e. $F(\id_A) = \id_{F(A)}$ for each $\mathcal{C}$-object $A$.

Functors are sometimes called \emph{covariant functors}\index{Functor!Covariant --}. A \emph{contravariant functor}\index{Functor!Contravariant --} from $\cat C$ to $\cat C'$ means a functor from $\cat C^{\op}$ to $\cat C'$.
\end{enumerate}
\end{definition}

A functor $F$ from $\mathcal{C}$ to $\mathcal{C}'$ will be denoted by $F: \mathcal{C} \lto \mathcal{C}'$ or $\mathcal{C} \stackrel{F}{\lto} \mathcal{C}'$. We frequently use the simplified notations $FA$ and $Ff$ rather than $F(A)$ and $F(f)$. Indeed, we sometimes denote the action on both objects and morphisms by $F(A \stackrel{f}{\lto} B) = FA \stackrel{Ff}{\lto} FB$.

A functor $F: \mathcal{C} \lto \mathcal{C}'$ is called an \emph{embedding}\index{Categorical!-- embedding}\index{Embedding!Categorical --} provided it is injective on morphisms, it is called \emph{faithful}\index{Functor!Faithful --} if all the hom-set restrictions $F: \hom_\cat C(A,A') \lto \hom_{\cat C'}(FA, FA')$ are injective and \emph{full}\index{Functor!Full --} if they are surjective; moreover it is called \emph{isomorphism-dense}\index{Functor!Isomorphism-dense --}\index{Isomorphism!{--}-dense functor} provided that for any $\cat C'$-object $B$ there exists some $\cat C$-object $A$ such that $FA$ is isomorphic to $B$. A functor $F: \cat C \lto \cat C'$ is called an \emph{isomorphism}\index{Categorical!-- isomorphism}\index{Isomorphism!Categorical --} if there exists a functor $G: \cat C' \lto \cat C$ such that $G \circ F$ is the identity functor of $\cat C$ and $F \circ G$ is the identity functor of $\cat C'$. In this case, the categories $\cat C$ and $\cat C'$ are called \emph{isomorphic}.

Even if there exists a definition of isomorphism between categories, there is a weaker concept that is much more useful in practice: the one of ``categorical equivalence''. Indeed, in mathematics, two object that are isomorphic can be treated as essentially the same object, and categories does not make exception; nonetheless, most of the categorical properties are preserved under equivalences and equivalences are more than isomorphisms, i.e. every isomorphism is an equivalence but not vice versa. Therefore, this sort of ``weak isomorphism'' deserves its own definition: a functor $F: \cat C \lto \cat C'$ is called an \emph{equivalence}\index{Categorical!-- equivalence} if it is full, faithful and isomorphism-dense. In this case, the categories $\cat C$ and $\cat C'$ are said to be \emph{equivalent}.

The following properties of functors are very easy to prove and, therefore, we limit ourselves to list them omitting their proofs.

\begin{proposition}\label{functprop}
Let $F: \cat C \lto \cat C'$ and $G: \cat C' \lto \cat C''$ be functors.
\begin{enumerate}
\item[$(i)$] All functors preserve isomorphisms, i.e., whenever $f: A \lto A'$ is a $\cat C$-isomorphism, then $F(f)$ is a $\cat C'$-isomorphism.
\item[$(ii)$] The composition $G \circ F: \cat C \lto \cat C''$ defined by $$(G \circ F)(f: A \lto A') = G(Ff): G(FA) \lto G(FA')$$ is a functor.
\item[$(iii)$] A functor is an embedding if and only if it is faithful and injective on objects.
\item[$(iv)$] A functor is an isomorphism if and only if it is full, faithful, and bijective on objects.
\item[$(v)$] If $F$ and $G$ are both isomorphisms (respectively: embeddings, faithful, full), then so is $G \circ F$.
\item[$(vi)$] If $G \circ F$ is an embedding (respectively: faithful), then so is $F$.
\item[$(vii)$] If $F$ is surjective on objects and $G \circ F$ is full, then $G$ is full.
\item[$(viii)$] If $F$ is full and faithful, then for every $\cat C'$-morphism $f : FA \lto FA'$ there exists a unique $\cat C$-morphism $g: A \lto A'$ such that $Fg = f$. Furthermore, $g$ is a $\cat C$-isomorphism if and only if $f$ is a $\cat C'$-isomorphism.
\item[$(ix)$] If $F$ is full and faithful, then it reflects isomorphisms; i.e., whenever $g$ is a $\cat C$-morphism such that $Fg$ is a $\cat C'$-isomorphism, then $g$ is a $\cat C$-isomorphism.
\item[$(x)$] If $F$ is an equivalence, then there exists an equivalence $H: \cat C' \lto \cat C$.
\item[$(xi)$] If $F$ and $G$ are equivalences, then so is $G \circ F$.
\end{enumerate}
\end{proposition}

We have formulated a duality principle related to objects, morphisms, and categories. We now extend this to functors, i.e. we introduce, for any functor, the concept of its dual functor\index{Dual!-- functor} that can be used to formulate the duals of categorical statements involving functors. Given a functor $F: \cat C \lto \cat C'$, the \emph{dual functor} $F^{\op}: \cat C^{\op} \lto \cat C'^{\op}$ is the functor defined by
$$F^{\op}(A \stackrel{f}{\lto} A') = FA \stackrel{Ff}{\lto} FA'.$$

It is immediate to verify that
\begin{proposition}
Each of the following properties of functors is self-dual: ``isomorphism'', ``embedding'', ``faithful'', ``full'', ``isomorphism-dense'' and ``equivalence''.
\end{proposition}

Two categories $\cat C$ and $\cat C'$ are called \emph{dually equivalent} if $\cat C^{\op}$ and $\cat C'$~--- or, that is the same, $\cat C$ and $\cat C'^{\op}$~--- are equivalent.

\section{Concrete categories and constructs}
\label{concretesec}

\begin{definition}\label{concrete}
Let $\cat X$ be a category. A \emph{concrete category}\index{Category!Concrete --} over $\cat X$ is a pair $(\cat C,U)$, where $\cat C$ is a category and $U: \cat C \lto \cat X$ is a faithful functor. Sometimes $U$ is called the \emph{forgetful}\index{Functor!Forgetful --}\index{Functor!Underlying --} (or \emph{underlying}) \emph{functor} of the concrete category and $\cat X$ is called the \emph{base category}\index{Category!Base --} for $(\cat C,U)$. A concrete category over the category $\cat S$ of sets is called a \emph{construct}\index{Construct}.
\end{definition}

For example, any category whose objects are sets with some structure (namely: topological spaces, groups, lattices, etc.), is a construct. The forgetful functor is the one that sends every object into its underlying set and every morphism into itself (the latter regarded as a morphism of sets). Moreover, many constructs can also be seen as concrete categories over another category. For instance, the category of vectorial spaces is both a construct and a concrete category over $\cat G^\textrm{ab}$, the category of Abelian groups; in this case, the forgetful functor does not ``forget'' the whole structure but only the external multiplication.

As we already stated at the beginning of the chapter, our aim is to introduce categorical concepts and results that we need for our purposes. So, since we will only deal with concrete categories and constructs, in what follows we will focus our interest exclusively on them.

Let $(\cat C,U)$ be a concrete category over a category $\cat X$, and let $X \in \obj(\cat X)$. A \emph{structured arrow}\index{Arrow!Structured --} with domain $X$ is a pair $(f,A)$ consisting of a $\cat C$-object $A$ and an $\cat X$-morphism $f: X \lto UA$. A \emph{universal arrow}\index{Arrow!Universal --} over an $\cat X$-object $X$ is a structured arrow $(u,A)$ with domain $X$ that has the following universal property: for any $\cat C$-object $B$ and any structured arrow $(f,B)$ with domain $X$ there exists a unique $\cat C$-morphism $h_f: A \lto B$ such that $U(h_f) \circ u = f$, i.e. such that the diagram
\begin{equation*}
\xymatrix{
X \ar[rr]^u \ar[ddrr]_f && UA \ar[dd]^{U(h_f)} \\
&& \\
&& UB \\
}
\end{equation*}
commutes.

\begin{definition}\label{freeobject}
A \emph{free object}\index{Object!Free --} over an $\cat X$-object $X$ is a $\cat C$-object $A$ such that there exists a universal arrow $(u,A)$ over $X$.
\end{definition}

\begin{remark}\label{universal=unique}
It is immediate to prove that, for any $\cat X$-object $X$, universal arrows over $X$ are essentially unique. In other words, if $(u,A)$ and $(u',A')$ are two universal arrows with domain $X$, then there exists a $\cat C$-isomorphism $i: A \lto A'$ such that $i \circ u = u'$ and $i^{-1} \circ u' = u$. Conversely, if $u: X \lto A$ is a universal arrow and $i: A \lto A'$ is a $\cat C$-isomorphism, then $(i \circ u, A')$ is universal too. This fact means also that a free $\cat C$-object over an $\cat X$-object $X$ is unique up to isomorphisms, and it is usually denoted by $\fr_\cat C(X)$. Moreover, if $\cat C$ is a construct, i.e. if the underlying category $\cat X$ is the one of sets, then the free object over a certain set $X$ depends only on the category $\cat C$ and on the cardinality of $X$. Hence, in this case, free objects are also denoted by $\fr_\cat C(\kappa)$ or $\fr_\cat C^\kappa$, where $\kappa$ is the cardinality of the set $X$.
\end{remark}

Let $\cat C$ be a construct, i.e. a concrete category with $\cat S$~--- the category of sets~--- as base category. A $\cat C$-morphism $f: A \lto B$ is called \emph{initial}\index{Morphism!Initial --}\index{Initial!-- morphism} provided that, for any $\cat C$-object $C$, an $\cat S$-morphism $g: UC \lto UA$ is a $\cat C$-morphism whenever $f \circ g: UC \lto UB$ is a $\cat C$-morphism. An initial morphism $f: A \lto B$ that has a monomorphic underlying $\cat S$-morphism $U(f): UA \lto UB$ is called an \emph{embedding}\index{Embedding}. If $f: A \lto B$ is an embedding, then $(f,B)$ is called an \emph{extension}\index{Extension (of an object)}\index{Object!Extension of an --} of $A$ and $(A, f)$ is called an \emph{initial subobject}\index{Initial!-- subobject} of $B$.

The concepts of final morphism and quotient morphism are dual to the concepts of initial morphism and embedding, respectively. A $\cat C$-morphism $f: A \lto B$ is called \emph{final}\index{Morphism!Final --}\index{Final!-- morphism} provided that, for any $\cat C$-object $C$, an $\cat S$-morphism $g: UB \lto UC$ is a $\cat C$-morphism whenever $g \circ f: UA \lto UC$ is a $\cat C$-morphism. A final morphism $f: A \lto B$ with epimorphic underlying $\cat S$-morphism $U(f): UA \lto UB$ is called a \emph{quotient morphism}\index{Morphism!Quotient --}. If $f: A \lto B$ is a quotient morphism, then $(f,B)$ is called a \emph{final quotient object}\index{Object!Final quotient -}\index{Final!-- quotient object} of A.

\begin{definition}\label{injprojobj}
In a construct $\cat C$, an object $C$ is called \emph{injective} provided that for any embedding $m: B \lto A$ and any $\cat C$-morphism $f: B \lto C$ there exists a $\cat C$-morphism $g: A \lto C$ extending $f$, i.e. such that the diagram
\begin{equation*}
\xymatrix{
A \ar@{-->}[ddrr]^g  && B \ar[dd]_f \ar@{^{(}->}[ll]^m  \\
&&\\
&& C \\
}
\end{equation*}
commutes.

Dually, an object $C$ is called \emph{projective} provided that for any morphism $f: C \lto B$ and any quotient morphism $e: A \lto B$ there exists a morphism $h: C \lto A$ extending $f$, i.e. such that the diagram
\begin{equation*}
\xymatrix{
A \ar@{->>}[rr]^e && B \\
&&\\
&& C \ar[uu]_f \ar@{-->}[lluu]^h\\
}
\end{equation*}
commutes.
\end{definition}
We observe explicitly that the morphisms $g$ and $h$ in the above definition are not required to be uniquely determined by $m$ and $f$ and by $e$ and $f$ respectively.

\section{Products and coproducts}
\label{prodsec}

A basic categorical concept that simultaneously generalizes the concepts of object and morphism is that of source.

\begin{definition}\label{source}
A \emph{source}\index{Source} is a pair $(A, (f_i)_{i \in I})$ consisting of an object $A$ and a family of morphisms $f_i: A \lto A_i$ with domain $A$, indexed by some class $I$. $A$ is called the domain of the source and the family $(A_i)_{i \in I}$ is called the codomain of the source. A source can be denoted also by $(f_i: A \lto A_i)_{i \in I}$ or $(A \stackrel{f_i}{\lto} A_i)_{i \in I}$.
\end{definition}
Given a source $S = (f_i: A \lto A_i)_{i \in I}$ and a morphism $f: B \lto A$, we write $S \circ f$ for denoting the source $(f_i \circ f: B \lto A_i)_{i \in I}$. A source $S = (A, (f_i)_{i \in I})$ is called a \emph{mono-source}\index{Source!Mono-{--}} provided it can be cancelled from the left, i.e., if for any two morphisms $h, k \in \hom(B,A)$ the equation $S \circ h = S \circ k$ implies $h = k$.

Cartesian products of families of sets, direct products of families of vector spaces, topological products of families of topological
spaces, etc., can be regarded as objects together with families of (projection) morphisms emanating from them, i.e., as sources. As such~--- but not as objects alone~--- they can be characterized, up to isomorphisms, by the following categorical definition:
\begin{definition}\label{product}
A source $T = (P \stackrel{p_i}{\lto} A_i)_{i \in I}$ is called a \emph{product}\index{Product!-- of objects} provided that for every source $S = (A \stackrel{f_i}{\lto} A_i)_{i \in I}$ with the same codomain there exists a unique morphism $f: A \lto P$ such that $S = T \circ f$. A product with codomain $(A_i)_{i \in I}$ is called a product of the family $(A_i)_{i \in I}$.
\end{definition}

\begin{lemma}\emph{\cite{cats}}\label{prodmono}
Every product is a mono-source.
\end{lemma}
\begin{proof}
If $T = (P, (p_i)_{i \in I})$ is a product and $h, k \in \hom(A,P)$ verify the equation $T \circ h = T \circ k$, then $S = T \circ h = T \circ k$ is a source with the same codomain as $T$. The uniqueness requirement in the definition of product implies that $h = k$. Hence $T$ is a mono-source.
\end{proof}

\begin{proposition}\emph{\cite{cats}}\label{produnique}
For any family $(A_i)_{i \in I}$ of objects, the product of $(A_i)_{i \in I}$ is essentially unique; i.e., if $T = (P \stackrel{p_i}{\lto} A_i)_{i \in I}$ is a product of $(A_i)_{i \in I}$, then the following hold:
\begin{enumerate}
\item[$(i)$] for each product $T' = (P' \stackrel{p_i'}{\lto} A_i)_{i \in I}$ there exists an isomorphism $h: P' \lto P$ with $T' = T \circ h$,
\item[$(ii)$] for each isomorphism $h: P' \lto P$ the source $T \circ h$ is a product of $(A_i)_{i \in I}$.
\end{enumerate}
\end{proposition}
\begin{proof}
\begin{enumerate}
\item[$(i)$] Since $T$ and $T'$ are products with the same codomain, there exist unique morphisms $h$ and $k$ such that $T' = T \circ h$ and $T = T' \circ k$. Therefore $T' \circ \id_{P'} = T' \circ (k \circ h)$ and $T \circ \id_P = T \circ (h \circ k)$. Since, by Lemma~\ref{prodmono}, $T$ and $T'$ are mono-sources, these equations imply that $\id_{P'} = k \circ h$ and $\id_P = h \circ k$. Hence $h$ is an isomorphism.
\item[$(ii)$] Obvious.
\end{enumerate}
\end{proof}

The above uniqueness result allows us to use the following notations for products (whenever they exist): the product of $(A_i)_{i \in I}$ will be denoted by $\left(\prod_{i \in I} A_i \stackrel{\pi_j}{\lto} A_j\right)_{j \in I}$, and the morphisms $\pi_j$ will be called \emph{projections}\index{Projection}.

We will now define the concepts of sink and coproduct, i.e. the respective dual concepts of source and product.
\begin{definition}\label{sink}
A \emph{sink}\index{Sink} is a pair $((f_i)_{i \in I},A)$ (sometimes denoted by $(f_i: A_i \lto A)_{i \in I}$ or $(A_i \stackrel{f_i}{\lto} A)_{i \in I}$) consisting
of an object $A$~--- called the codomain of the sink~--- and a family of morphisms $f_i: A_i \lto A$ indexed by some class $I$. The family $(A_i)_{i \in I}$ is called the domain of the sink.
\end{definition}
As for sources, if $S = (f_i: A_i \lto A)_{i \in I}$ is a sink and $f: A \lto B$ is a morphism, we will denote by $f \circ S$ the sink $(f \circ f_i: A_i \lto B)_{i \in I}$. Moreover we ``dualize'' in an obvious way the concept of mono-source: a sink $S = ((f_i)_{i \in I},A)$ is called an \emph{epi-sink}\index{Sink!Epi-{--}} provided it can be cancelled from the right, i.e. if, for any two morphisms $h, k \in \hom(A,B)$, the equation $h \circ S = k \circ S$ implies $h = k$.

\begin{definition}\label{coproduct}
A sink $T = (A_i \stackrel{m_i}{\lto} C)_{i \in I}$ is called a \emph{coproduct}\index{Coproduct of objects}\index{Product!Co-{--} of objects} if for every sink $S = (A_i \stackrel{f_i}{\lto} A)_{i \in I}$ with the same domain of $T$ there exists a unique morphism $f: C \lto A$ such that $S = f \circ T$. A coproduct with domain $(A_i)_{i \in I}$ is called a coproduct of the family $(A_i)_{i \in I}$.
\end{definition}

The Duality Principle gives us the following two results for free, since they are the dual statements of, respectively, Lemma~\ref{prodmono} and Proposition~\ref{produnique}.
\begin{proposition}\label{coprodepi}
Every coproduct is an epi-sink.
\end{proposition}
\begin{proposition}\label{coprodunique}
For any family $(A_i)_{i \in I}$ of objects, the coproduct of $(A_i)_{i \in I}$ is essentially unique; i.e., if $T = (A_i \stackrel{m_i}{\lto} C)_{i \in I}$ is a coproduct of $(A_i)_{i \in I}$, then the following hold:
\begin{enumerate}
\item[$(i)$] for each coproduct $T' = (A_i \stackrel{m_i'}{\lto} C')_{i \in I}$ there exists an isomorphism $h: C' \lto C$ with $T' = T \circ h$,
\item[$(ii)$] for each isomorphism $h: C' \lto C$ the sink $T \circ h$ is a coproduct of $(A_i)_{i \in I}$.
\end{enumerate}
\end{proposition}
Then we can, again, introduce a special notation for coproducts (whenever they exist): the coproduct of $(A_i)_{i \in I}$ will be denoted by $\left(A_j \stackrel{\mu_j}{\lto} \coprod_{i \in I} A_i\right)_{j \in I}$, and the morphisms $\mu_j$ will be called \emph{injections}\index{Injection}.

\section*{References and further readings}
\addcontentsline{toc}{section}{References and further readings}

A classical book on Category Theory is the famous ``Categories for the working mathematician'', by Saunders MacLane, published for the first time in 1971 and whose last edition,~\cite{maclane}, is dated 1998. An Italian version of this book,~\cite{maclaneit}, was published in 1977.

By the way, another book that can be recommended is~\cite{cats}, written by Ji\v{r}\'i Ad\'amek, Horst Herrlich and George E. Strecker. This book, besides being a complete and detailed work on Category Theory, has two further important merits. First, it is designed to be used both as a textbook for beginners and as a reference source, and~--- as the authors themselves claim in the preface~--- is organized and written in a ``pedagogical style''. The second (but not less important for students and young researchers) merit is that the book is available for free on the internet, with a GNU Free Documentation License, in its newest edition~\cite{catssito}. Actually the book by Ad\'amek, Herrlich and Strecker has been our reference in drawing up this chapter.

\chapter{Ordered Structures and Residuation Theory}
\label{orderchap}

In Part~\ref{modupart} we will study quantale modules and show their application to deductive systems on propositional languages. We have already anticipated that we can think of quantales and quantale modules as objects that somehow remind rings and ring modules. However, likenesses apart, it is worthwile to point out two important and distinctive features of quantale modules.

First of all, the object classes of the categories of quantales and modules over a fixed quantale are subclasses of the class of complete lattices. This characteristic makes quantale and $\Q$-module categories richer and~--- somehow~--- more flexible than other similarly paired categories (e.g. rings and ring modules, monoids and sets with an action), thanks to the presence of infinitary operations and of an infinitely distributive operation.

In this connection, we point out the second important characteristic: the infinite distributivity of the external product. It is well known that a binary operation in a lattice-ordered structure distributes with respect to any existing join if and only if it is biresiduated with respect to the lattice order of the structure. Such an equivalence allows us to use the distributivity or some Galois connection, depending on our convenience.

From all these considerations it follows that, before going through the study of $\Q$-modules, we need to recall several notions and results of Lattice Theory and Residuation Theory; it will be the subject of this chapter.

In the first section, we will treat residuated maps and adjoint pairs, and we will show the main known properties of these maps and pairs. Sections~\ref{suplatsec} and~\ref{sltensorsec} are directed to providing the reader with a certain knowledge of the category of sup-lattices. For the same reason that obliges us to know enough on Abelian groups, in order to study ring modules, we need to be familiar at least with the basic properties of sup-lattices, before starting the examination of quantale modules. Finally, in Section~\ref{reslatquantsec}, we will give the definition of residuated lattices and quantales. Then we will show several properties of the category of quantales like, e.g., the construction of free objects.

A last note. In constructs, i.e. in concrete categories over the category of sets, monomorphisms, morphisms with injective underlying map and embeddings need not coincide, in general, as well as their respective dual notions. Nonetheless it can be proved that, in all the categories we will encounter in the rest of the thesis, such notions coincide and, therefore, we will use any of them indifferently.

\section{Residuated maps}
\label{resmapsec}

In this section we recall some basic definitions and results on residuated maps.

\begin{definition}\label{residuated map}
Let $\la X, \leq \ra$ and $\la Y, \leq \ra$ be two posets. A map $f: X \lto Y$ is said to be \emph{residuated}\index{Residuated!-- map} iff there exists a map $g: Y \lto X$ such that, for all $x \in X$ and for all $y \in Y$, the following condition holds:
\begin{equation*}
f(x) \leq y \quad \iff \quad x \leq g(y).
\end{equation*}
It is immediate to verify that the map $g$ is uniquely determined; we will call it the \emph{residual map}\index{Residual map} or the \emph{residuum}\index{Residuum} of $f$, and denote it by $f_*$. The pair $(f, f_*)$ is said \emph{adjoint}\index{Adjoint pair}.
\end{definition}

Before discussing the basic properties of adjoint pairs, we recall that, if $\la X, \leq \ra$ is a poset, a map $\g: X \lto X$ is called a \emph{closure operator}\index{Operator!Closure --}\index{Closure!-- operator} iff it is order preserving, extensive and idempotent, i.e. iff for all $x, y \in X$
\begin{enumerate}
\item[$(i)$]$x \leq y$ implies $\g(x) \leq \g(y)$,
\item[$(ii)$]$x \leq \g(x)$,
\item[$(iii)$]$\g \circ \g = \g$.
\end{enumerate}
Dually, a map $\d: X \lto X$ is called a \emph{coclosure operator}\index{Operator!Coclosure --}\index{Operator!Interior --}\index{Closure!Co-- operator}, or an \emph{interior operator}, iff it satisfies
\begin{enumerate}
\item[$(i)$]$\d$ is order preserving,
\item[$(ii)$]$\d(x) \leq x$ for all $x \in X$,
\item[$(iii)$]$\d(\d(x)) = \d(x)$, for all $x \in X$.
\end{enumerate}

The following result is a classical characterization of residuated maps:
\begin{theorem}\label{residchar}
Let $\la X, \leq \ra$ and $\la Y, \leq \ra$ be two posets, and $f: X \lto Y$. The following statements are equivalent:
\begin{enumerate}
\item[$(a)$]$f$ is residuated, with residual $f_*$;
\item[$(b)$]$f$ preserves all existing joins in $X$, i.e., if $Z \subseteq X$ and there exists $\bigvee Z$ in $X$, then there exists, in $Y$ $\bigvee f[Z]$ and $\bigvee f[Z] = f\left[\bigvee Z \right]$;
\item[$(c)$]$f$ is isotone and there exists a unique isotone map $f_*: Y \lto X$ such that
\begin{equation}\label{resid1}
f \circ f_* \leq \id_Y
\end{equation}
and
\begin{equation}\label{resid2}
f_* \circ f \geq \id_X.
\end{equation}
Moreover, the following hold for any adjoint pair $(f,f_*)$:
\item[$(d)$]\begin{equation}\label{resid3}
f \circ f_* \circ f = f
\end{equation}
and
\begin{equation}\label{resid4}
f_* \circ f \circ f_* = f_*.
\end{equation}
\item[$(e)$]$$f_*(y) = \max\{x \in X \mid f(x) \leq y\}.$$
\end{enumerate}
\end{theorem}

\begin{corollary}\label{resclosureinterior}
Let $\la X, \leq \ra$ and $\la Y, \leq \ra$ be two posets and let $(f, f_*)$ be an adjoint pair, with $f: X \lto Y$. Then
\begin{enumerate}
\item[$(i)$]$f_* \circ f$ is a closure operator over $X$;
\item[$(ii)$]$f \circ f_*$ is an interior operator in $Y$.
\end{enumerate}
\end{corollary}
\begin{proof}
It is enough to apply (\ref{resid1}--\ref{resid4}).
\end{proof}

Moreover we have:

\begin{proposition}\label{residprop}
Let $\la X, \leq \ra$ and $\la Y, \leq \ra$ be posets, and let $(f,f_*)$ be an adjoint pair, with $f: X \lto Y$. Then the following hold:
\begin{enumerate}
\item[$(i)$]$f$ preserves all existing joins, i.e. if $\{x_i\}_{i \in I}$ is a family of elements of $X$ such that there exists $\bigvee_{i \in I} x_i$, then also $\bigvee_{i \in I} f(x_i)$ exists and $f\left(\bigvee_{i \in I} x_i\right) = \bigvee_{i \in I} f(x_i)$;
\item[$(ii)$]$f_*$ preserves all existing meets, i.e. if $\{y_j\}_{j \in J}$ is a family of elements of $Y$ such that there exists $\bigwedge_{j \in J} y_j$, then also $\bigwedge_{j \in J} f_*(y_j)$ exists and $f_*\left(\bigwedge_{j \in J} y_j\right) = \bigwedge_{j \in J} f_*(y_j)$;
\item[$(iii)$]$f$ is surjective \quad $\iff$ \quad $f_*$ is injective \quad $\iff$ \quad $f \circ f_* = \id_Y$;
\item[$(iv)$]$f$ is injective \quad $\iff$ \quad $f_*$ is surjective \quad $\iff$ \quad $f_* \circ f = \id_X$.
\end{enumerate}
\end{proposition}

Let $\la X, \leq \ra$, $\la Y, \leq \ra$ and $\la W, \leq \ra$ be posets. A map $f: X \times Y \lto W$ of two variables is said to be \emph{biresiduated}\index{Residuated!Bi-- map} if it is residuated with respect to each variable, i.e. if the following two conditions hold:
\begin{enumerate}
\item[-]for any fixed $\ov y$, there exists a map $g_{\ov y}: W \lto X$ such that
\begin{equation*}
f(x,\ov y) \leq w \qquad \iff \qquad x \leq g_{\ov y}(w);
\end{equation*}
\item[-]for any fixed $\ov x$, there exists a map $h_{\ov x}: W \lto Y$ such that
\begin{equation*}
f(\ov x,y) \leq w \qquad \iff \qquad y \leq h_{\ov x}(w).
\end{equation*}
\end{enumerate}
In an analogous way the notion of residuated map extends to that of $n$-residuated map, for all $n \in \N$.

\section{Sup-lattices}
\label{suplatsec}

A poset $\la L, \leq \ra$ which admits arbitrary joins is called a \emph{sup-lattice}\index{Sup-lattice}. A sup-lattice homomorphism is a map that preserves arbitrary joins. If $S \subseteq L$, the join over $S$ will be denoted indifferently by $\bigvee_{s \in S} s$ or $\bigvee S$. It is easily seen that a sup-lattice admits also arbitrary meets. Indeed, for any $S \subseteq L$, we can consider the set $S' = \{x \in L \mid x \leq s \ \forall s \in S\}$ of the lower bounds of $S$, and we have $\bigvee S' = \bigwedge S$. Then the category $\SL$ of sup-lattices\index{Category!-- $\SL$, of sup-lattices} is the one whose objects are complete lattices and morphisms are maps preserving arbitrary joins. If $\LL_1$ and $\LL_2$ are sup-lattices and $\phi: \LL_1 \lto \LL_2$ is a sup-lattice homomorphism, then clearly $\phi(\bot_1) = \phi\left({}^{\LL_1}\bigvee \varnothing\right) = {}^{\LL_2}\bigvee \varnothing = \bot_2$, where $\bot_i$ is the bottom element of $\LL_i$, $i = 1, 2$. Thus, for a sup-lattice $\LL$, we will use the signature $\LL = \la L, \vee, \bot \ra$.

For any sup-lattice $\LL = \la L, \vee, \bot \ra$, it is possible to define a dual sup-lattice in an obvious way: if we consider the opposite partial order $\geq$ (also denoted by $\leq^{\op}$), then $\LL^{\op} = \la L, \vee^{\op}, \bot^{\op} \ra = \la L, \wedge, \top \ra$, where $\top = \bigvee L$, is a sup-lattice and, clearly, $(\LL^{\op})^{\op} = \LL$.

It is also clear that, given two sup-lattices $\LL_1$ and $\LL_2$, the hom-set $\hom_\SL(\LL_1,\LL_2)$ is a sup-lattice itself, with the order relation~--- and, therefore, the operations~--- defined pointwisely: $f \leq g \iff f(x) \leq g(x) \ \forall x \in L_1$. We will denote it by $\Hom_\SL(\LL_1,\LL_2)$.

If $\LL_1$ and $\LL_2$ are sup-lattices and $f \in \hom_\SL(\LL_1, \LL_2)$, $f$ preserves all suprema in $L_1$. Then, by Theorem~\ref{residchar}, $f$ is a residuated map whose residual, $f_*$, is defined by $f_*(y) = \bigvee\{x \in L_1 \mid f(x) \leq y\}$ for all $y \in L_2$. Moreover, $f_*$ is also a homomorphism between the dual sup-lattices $\LL_2^{\op}$ and $\LL_1^{\op}$; to emphasize this fact, we will denote $f_*$ also by $f^{\op}$. For all $f, g \in \hom_\SL(\LL_1,\LL_2)$, we have:
\begin{enumerate}
\item[-] $(f^{\op})^{\op} = f$,
\item[-] $f \leq g \iff g_* \leq f_* \iff f^{\op} \leq^{\op} g^{\op}$,
\item[-] $(g \circ f)^{\op} = f^{\op} \circ g^{\op}$,
\end{enumerate}
hence
\begin{equation}\label{homsldual}
\Hom_\SL(\LL_1,\LL_2) \cong \Hom_\SL(\LL_2^{\op},\LL_1^{\op}).
\end{equation}
Conversely, if $f: \LL_1 \lto \LL_2$ is a residuated map between sup-lattices, it is a sup-lattice homomorphism between $\LL_1$ and $\LL_2$ by Proposition~\ref{residprop}$(i)$, and its residuum $f_*$ is a sup-lattice homomorphism from $\LL_2^{\op}$ to $\LL_1^{\op}$, by the $(ii)$ of the same proposition.

Then we have that sup-lattice homomorphisms coincide with residuated maps between sup-lattices. Moreover, it follows
\begin{proposition}\label{slselfduality}
The contravariant functor $( \ )^{\op}: \SL \lto \SL$ is a self-duality.
\end{proposition}

We will now describe free objects in the category $\SL$ showing that, for any given set $X$, the free sup-lattice $\frsl(X)$ over the set of generators $X$ is $\wp(\mathbf{X}) = \la \wp(X), \cup, \varnothing \ra$, i.e. the power set of $X$ with the set-union as the join operation and the empty set as the bottom element.

\begin{proposition}\label{freesl}
For any set $X$, the free sup-lattice $\frsl(X)$ generated by $X$ is $\wp(\mathbf{X})$, equipped with the \emph{singleton map} $\sigma: x \in X \lmapsto \{x\} \in \wp(X)$.
\end{proposition}
\begin{proof}
We need to prove that, given an arbitrary sup-lattice $\LL = \la L, \vee, \bot \ra$ and an arbitrary map $f: X \lto L$, there exists a unique morphism $h_f: \wp(\mathbf{X}) \lto \LL$ that extends $f$, i.e. such that $h_f \circ \sigma = f$. Let $h_f$ be the map defined, for all $S \in \wp(X)$, by $h_f(S) = \bigvee_{s \in S} f(s)$. Trivially, for any $\ov x \in X$, $f(\ov x) = \bigvee_{x \in \{\ov x\}}f(x) = h_f(\{\ov x\}) = (h_f \circ \sigma)(x)$, hence $h_f \circ \sigma = f$; moreover, the fact that $h_f$ is a sup-lattice morphism follows easily from how the map is defined and the properties of complete lattices. If $h_f': \wp(\mathbf{X}) \lto \LL$ is another morphism that extends $f$, then we have:
$$h_f'(S) = h_f'\left(\bigcup_{s \in S} \{s\}\right) = h_f'\left(\bigcup_{s \in S} \sigma(s)\right) = \bigvee_{s \in S} (h_f' \circ \sigma)(s) = \bigvee_{s \in S} f(s) = h_f(S).$$
\end{proof}

In particular, by Proposition~\ref{freesl}, the free sup-lattice over one generator $\wp 1$ is isomorphic to $\{\bot, \top\}$ and, clearly, $\wp 1^{\op} \cong \wp 1$. On the other hand, it is also clear that, for an arbitrary sup-lattice $\LL$, a map $f: \wp1 \lto \LL$ is a sup-lattice homomorphism if and only if $f(\bot) = \bot_\LL$. Thus, using also (\ref{homsldual}) and self duality of $\wp1$, we have
$$\LL \cong \Hom_\SL(\wp1,\LL) \cong \Hom_\SL(\LL^{\op},\wp1)$$
and
$$\LL^{\op} \cong \Hom_\SL(\wp1,\LL)^{\op} \cong \Hom_\SL(\LL,\wp1),$$
hence
\begin{equation}\label{homp1l}
\LL \cong \Hom_\SL(\wp1,\LL) \cong \Hom_\SL(\LL^{\op},\wp1) \cong \Hom_\SL(\LL,\wp1)^{\op},
\end{equation}
\begin{equation}\label{homlp1}
\LL^{\op} \cong \Hom_\SL(\wp1,\LL)^{\op} \cong \Hom_\SL(\LL,\wp1) \cong \Hom_\SL(\LL^{\op},\wp1)^{\op}.
\end{equation}

\begin{proposition}\label{coprodsl}
For any family of sup-lattices $\{\LL_i\}_{i \in I}$, the coproduct $\coprod_{i \in I} \LL_i$ is the product $\prod_{i \in I} \LL_i$ equipped with the maps $\mu_i: \LL_i \lto \prod_{i \in I} \LL_i$ that send each $x_i \in L_i$ in the family of $\prod_{i \in I} \LL_i$ in which all the elements are equal to $\bot$ except the $i$-th that is equal to $x_i$. Moreover, if $\pi_i: \prod_{i \in I} \LL_i \lto \LL_i$ is the canonical $i$-th projection, for all $i \in I$, $\pi_i \circ \mu_i = \id_{\LL_i}$.
\end{proposition}
\begin{proof}
First of all let us observe that, for any family $(x_i)_{i \in I}$, $(x_i)_{i \in I} = \bigvee_{i \in I} \mu_i(x_i)$. We want to prove that, given an arbitrary sup-lattice $\LL$ and a family of homomorphisms $f_i: \LL_i \lto \LL$, there exists a unique homomorphism $f: \prod_{i \in I} \LL_i \lto \LL$ such that $f \circ \mu_i = f_i$ for all $i \in I$.

Let $f\left((x_i)_{i \in I}\right) = \bigvee_{i \in I} f_i(x_i)$. For all $i \in I$ and $x_i \in L_i$, we have $(f \circ \mu_i)(x_i) = f\left(\mu_i(x_i)\right) = f_i(x_i) \vee \bigvee_{j \in I \setminus \{i\}} f_j(\bot) = f_i(x_i)$.

Now, let $f': \prod_{i \in I} \LL_i \lto \LL$ be another homomorphism such that $f' \circ \mu_i = f_i$ for all $i \in I$. Then $f'\left((x_i)_{i \in I}\right) = f'\left(\bigvee_{i \in I} \mu_i(x_i)\right) = \bigvee_{i \in I} f'(\mu_i(x_i)) = \bigvee_{i \in I} f_i(x_i) = f\left((x_i)_{i \in I}\right)$. Thus $f' = f$ and also the uniqueness is proved.
\end{proof}

\begin{proposition}\label{projsl}
Every free sup-lattice is projective, and an object $\LL$ in $\SL$ is projective if and only if its dual $\LL^{\op}$ is injective.
\end{proposition}
\begin{proof}
Let $\wp(\mathbf{X})$ be a free sup-lattice, and let $\LL_1$ and $\LL_2$ be two sup-lattice such that there exists a surjective morphism $g: \LL_1 \lto \LL_2$ and a morphism $f: \wp(\mathbf{X}) \lto \LL_2$. Then $f$ the unique homomorphism that extends the map $f_X: x \in X \lmapsto f(\{x\}) \in \LL_2$. If we consider the map $g_* \circ f_X: X \lto \LL_1$, we can extend it to a homomorphism $h: \wp(\mathbf{X}) \lto \LL_1$, and it is immediate to verify that $g \circ h = f$. Then any free sup-lattice is projective.

The second assertion is a trivial application of the self-duality of $\SL$.
\end{proof}

\section{Tensor products of sup-lattices}
\label{sltensorsec}

The present section is dedicated to the tensor product of sup-lattices. We will show its existence and some of its properties. In order to do that, we will start by introducing some notations and results.

Let $\Cl(\LL)$ be the set of all the closure operators over a given sup-lattice $\LL$. It is easily seen that it is a poset with respect to the order relation defined pointwise: for all $\g, \d \in \Cl(\LL)$, $\g \leq \d \iff \g(x) \leq \d(x) \ \forall x \in L$. Furthermore, let us consider the set, $\Quot(\LL)$, of all the quotients of $\LL$ and the one, $\Lambda(\LL)$, of all the subsets of $\LL$ that are closed under arbitrary meets, both partially ordered by the set inclusion.

The following result proves that the posets $\Cl(\LL)^{\op}$, $\Lambda(\LL)$ and $\Quot(\LL)$ are isomorphic and, moreover, that any quotient of a sup-lattice is isomorphic~--- in a precise sense~--- to a subset of the same sup-lattice that is closed under arbitrary meets.

\begin{theorem}\label{quotmeetclosure}
For any sup-lattice $\LL$, $\Cl(\LL)^{\op}$, $\Quot(\LL)$ and $\Lambda(\LL)$ are isomorphic posets. Moreover, every element of $\Lambda(\LL)$ has a sup-lattice structure that is canonically isomorphic to a quotient of $\LL$.
\end{theorem}
\begin{proof}
For all $S \in \Lambda(\LL)$, we can define a map $\g_S: L \lto L$ by setting, for all $x \in L$, $\g_S(x) = \bigwedge\{y \in S \mid x \leq y\}$. It comes straightforwardly from its definition that $\g_S$ is monotone, extensive and idempotent; hence $\g_S \in \Cl(\LL)$ and $\g_S[L] = S$. Then we can define the map
$$\Gamma: \ S \in \Lambda(\LL) \ \lmapsto \ \g_S \in \Cl(\LL).$$
If $S, T \in \Lambda(\LL)$ are two different sets, then there exists $\ov x \in (S \setminus T) \cup (T \setminus S)$. If $\ov x \in S \setminus T$, then $\g_S(\ov x) = \ov x \neq \g_T(\ov x)$; analogously, if $\ov x \in T \setminus S$, then $\g_T(\ov x) = \ov x \neq \g_S(\ov x)$. Hence $S \neq T$ implies $\g_S \neq \g_T$, and $\Gamma$ is injective.

On the other hand, if $\g \in \Cl(\LL)$ and $S$ is an arbitrary subset of $\g[L]$, then $S$ is also a subset of $L$; thus there exists $z = \bigwedge S \in L$. Now, since $z \leq x$ for all $x \in S$, $\g(z) \leq \g(x) = x$ for all $x \in S$. Therefore $\g(z) \leq z$, whence $\g(z) = z \in S$. Then $\g[L]$ is closed under arbitrary meets. Moreover, if we assume that there exists $\ov x \in L$ such that $\g(\ov x) \neq \g_{\g[L]}(\ov x)$, then clearly $\g(\ov x) \gneqq \g_{\g[L]}(\ov x)$, by the definition of $\g_{\g[L]}$. So we have $\ov x \leq \g_{\g[L]}(\ov x)$ and $\g(\g_{\g[L]}(\ov x)) = \g_{\g[L]}(\ov x) < \g(\ov x)$, that contradicts the monotonicity of $\g$. It follows that $\g = \g_{\g[L]}$, for all $\g \in \Cl(\LL)$, and $\Gamma$ is surjective too, hence bijective. Now, if $S \subseteq T \in \Lambda(\LL)$, then clearly $\g_T \leq \g_S$ by the definition of such maps, so $\Gamma$ is order reversing. Analogously it is immediate to verify that also $\Gamma^{-1}$ is order reversing and therefore $\Gamma$ is an isomorphism between $\Lambda(\LL)$ and $\Cl(\LL)^{\op}$.

Now it is easy to see that, if $\g = \g_S$ is a closure operator, then $\mathbf{S} = \la S, \g \circ \vee, \g(\bot)\ra$ is a sup-lattice. Then, given a set $S \in \Lambda(\LL)$, we can define a map $\rho_S: L \lto S$ just by restricting the codomain of $\g_S$ to $S$:
\begin{equation}\label{reflection}
\rho_S: \ x \in L \ \lmapsto \ \g_S(x) = \bigwedge \{y \in S \mid x \leq y\} \in S.
\end{equation}
The map $\rho_S$, also called the \emph{reflection}\index{Sup-lattice!-- reflection} of $L$ over $S$, is the left inverse of the inclusion $\id_{L\restr S}: S \hookrightarrow L$ and is a surjective homomorphism from $\LL$ to $\mathbf{S}$; thus $\mathbf{S}$ is isomorphic to the quotient $\LL/\rho_S$ of $\LL$. In this way we define a map
$$\Pi: \ \mathbf{S} \in \Lambda(\LL) \ \lmapsto \ \LL/\rho_S \in \Quot(\LL),$$
and a map $\Theta = \Pi \circ \Gamma^{-1}$, from $\Cl(\LL)$ to $\Quot(\LL)$.

It is easily seen that $\Pi$ is an isomorphism and, consequently, $\Theta$ is an isomorphism between $\Cl(\LL)^{\op}$ and $\Quot(\LL)$; the theorem is proved.
\end{proof}

\begin{lemma}\emph{\cite{joyaltierney}}\label{quotsl}
Let $\LL$ be a sup-lattice and $R \subseteq L \times L$. Then the set
$$S = \{x \in L \ \mid \ \forall (r_1,r_2) \in R \quad r_1 \leq x \iff r_2 \leq x\}$$
is closed under arbitrary meets, and the quotient of $\LL$ with respect to the congruence $\equiv_R$ generated by $R$ is isomorphic to the sup-lattice $\mathbf{S} = \la S, \g_S \circ \vee, \g_S(\bot)\ra$, i.e. $\Pi[\mathbf{S}] = \LL/\equiv_R$.
\end{lemma}
\begin{proof}
The fact that $S$ is closed under arbitrary meets is an easy consequence of its definition. Now let us consider the reflection $\rho_S: L \lto S$; for any pair $(r_1, r_2) \in R$, we have
$$\rho_S(r_1) = \bigwedge\{x \in S \mid r_1 \leq x\} = \bigwedge\{x \in S \mid r_2 \leq x\} = \rho_S(r_2),$$
thus $\equiv_R$ is contained in the congruence determined by $\rho_S$.

In order to complete the proof, suppose that $\rho_T: L \lto T$ is another reflection such that $\rho_T(r_1) = \rho_T(r_2)$ for all $(r_1,r_2) \in R$. Then, for any $x \in T$, we have
$$r_1 \leq x \iff \rho_T(r_1) \leq x \iff \rho_T(r_2) \leq x \iff r_2 \leq x,$$
whence $T \subseteq S$, i.e. the congruence $\rho_T$ contains the congruence $\rho_S$. The thesis follows.
\end{proof}

Now we are ready to define the tensor product of sup-lattices, show its existence, and then to describe its shape and prove some of its basic properties.
\begin{definition}\label{tensorsl}
Let $\LL_1$, $\LL_2$ and $\LL$ be sup-lattices; a map $f: L_1 \times L_2 \lto L$ is said to be a \emph{bimorphism}\index{Bimorphism!-- of sup-lattices}\index{Sup-lattice!-- bimorphism} if it preserves arbitrary joins in each variable separately:
$$f\left(\bigvee_{i \in I} x_i, y\right) = \bigvee_{i \in I} f(x_i, y) \quad \textrm{and} \quad f\left(x, \bigvee_{j \in J} y_j\right) = \bigvee_{j \in J} f(x,y_j).$$

The \emph{tensor product}\index{Tensor product!-- of sup-lattices}\index{Sup-lattice!-- tensor product} $\LL_1 \tensor \LL_2$, of $\LL_1$ and $\LL_2$, is the codomain of the universal bimorphism $\LL_1 \times \LL_2 \lto \LL_1 \tensor \LL_2$. In other words, we call tensor product of $\LL_1$ and $\LL_2$ a sup-lattice $\LL_1 \tensor \LL_2$, equipped with a bimorphism $\tau: \LL_1 \times \LL_2 \lto \LL_1 \tensor \LL_2$, such that, for any sup-lattice $\LL$ and any bimorphism $f: \LL_1 \times \LL_2 \lto \LL$, there exists a unique homomorphism $h_f: \LL_1 \tensor \LL_2 \lto \LL$ with $h_f \circ \tau = f$.
\end{definition}

\begin{remark}\label{tensorunique}
By definition, if the tensor product of two algebraic structures exists, then it is unique up to isomorphisms.
\end{remark}

\begin{theorem}\label{tensorslexists}
Let $\LL_1$ and $\LL_2$ be sup-lattices. The tensor product $\LL_1 \tensor \LL_2$ exists; it is, up to isomorphisms, the quotient of the free sup-lattice $\wp(\LL_1 \times \LL_2)$ with respect to the congruence generated by the set
\begin{equation}\label{R}
R = \left\{
	\begin{array}{l}
		\left(\left\{\left(\bigvee A, y\right)\right\}, \bigcup_{a \in A}\{(a,y)\}\right) \\
		\left(\left\{\left(x, \bigvee B\right)\right\}, \bigcup_{b \in B}\{(x,b)\}\right) \\
	\end{array} \right\vert
	\left.
	\begin{array}{l}
	A \subseteq L_1, y \in L_2 \\
	B \subseteq L_2, x \in L_1 \\
	\end{array}
	\right\}.
\end{equation}
\end{theorem}
\begin{proof}
Let $\LL$ be any sup-lattice and let $f: \LL_1 \times \LL_2 \lto \LL$ be a bimorphism. Since $f$ is, of course, a map, we can extend it to a homomorphism $h_f: \wp(\LL_1 \times \LL_2) \lto \LL$; thus $h_f \circ \sigma = f$. On the other hand, the fact that $f$ is a bimorphism implies $f\left(\bigvee A, y\right) = \bigvee_{a \in A} f(a,y)$ and $f\left(x,\bigvee B\right) = \bigvee_{b \in B} f(x,b)$, for all $x \in L_1$, $y \in L_2$, $A \subseteq L_1$ and $B \subseteq L_2$. Then, since $h_f$ is a homomorphism, we have
\begin{eqnarray}
&&h_f\left(\left\{\left(\bigvee A,y\right)\right\}\right) = (h_f \circ \sigma)\left(\bigvee A,y\right) \nonumber \\
&=& f\left(\bigvee A,y\right) = \bigvee_{a \in A}f(a,y) \nonumber \\
&=& \bigvee_{a \in A} h_f(\sigma(a,y)) = h_f\left(\bigcup_{a \in A}\sigma(a,y)\right) \nonumber \\
&=& h_f\left(\bigcup_{a \in A}\{(a,y)\}\right) \nonumber
\end{eqnarray}
and, analogously, $h_f\left(\left\{\left(x,\bigvee B\right)\right\}\right) = h_f\left(\bigcup_{b \in B}\{(x,b)\}\right)$.

What above means that the kernel of $h_f$ contains $R$, thus~--- once denoted by $\mathbf{T}$ the quotient sup-lattice $\wp(\LL_1 \times \LL_2)/\equiv_R$ and by $\pi$ the canonical projection of $\wp(\LL_1 \times \LL_2)$ over  it~--- the map
$$h_f' \ : \quad [X]_{\equiv_R} \ \in \ \mathbf{T} \quad \lmapsto \quad h_f(X) \ \in \ \LL$$
is well defined and is a homomorphism. Moreover we have $h_f' \circ \pi \circ \sigma = h_f \circ \sigma = f$, so we have extended the bimorphism $f$ to a homomorphism $h_f'$, and it is easy to verify that the map $\tau = \pi \circ \sigma$ from $\LL_1 \times \LL_2$ to $\mathbf{T}$ is indeed a bimorphism.

The following commutative diagram may clarify the constructions above.
\begin{equation}\label{tensorsldiagram}
\xymatrix{
\LL_1 \times \LL_2 \ar[rr]^\sigma \ar[rddd]_f \ar[rd]_\tau && \wp(\LL_1 \times \LL_2) \ar[lddd]^{h_f} \ar[ld]^{\pi} \\
 & \mathbf{T}  \ar[dd]^{h_f'} & \\
 &&\\
 & \LL & \\
}
\end{equation}

It is useful to remark explicitly that $R$ and $\tau$ do not depend either on the sup-lattice $\LL$ or on the bimorphism $f$. Then we have proved that $\tau$ is the universal bimorphism whose domain is $\LL_1 \times \LL_2$, and that $\mathbf{T}$ is its codomain, i.e. the tensor product $\LL_1 \tensor \LL_2$.
\end{proof}

\begin{definition}\label{tensors}
If $x \in \LL_1$ and $y \in \LL_2$, we will denote by $x \tensor y$ the image of the pair $(x,y)$ under $\tau$, i.e. the congruence class $[\{(x,y)\}]_{\equiv_R}$, and we will call it a \emph{tensor}\index{Tensor!-- in $\SL$}. It can be proved that every element of $\LL_1 \tensor \LL_2$ is a join of tensors, so
$$\LL_1 \tensor \LL_2 = \left\{\bigvee_{i, j} x_i \tensor y_j \ \Big\vert \ \{x_i\}_{i \in I} \subseteq L_1, \{y_j\}_{j \in J} \subseteq L_2\right\}.$$
\end{definition}

\begin{theorem}\label{isohomtenssl}
Let $\LL_1$, $\LL_2$ and $\LL_3$ be sup-lattices. Then
$$\Hom_\SL(\LL_1 \tensor \LL_2, \LL_3) \cong \Hom_\SL(\LL_1,\Hom_\SL(\LL_2,\LL_3)).$$
\end{theorem}
\begin{proof}
If $h$ is a homomorphism from $\LL_1 \tensor \LL_2$ to $\LL_3$, then $h(x \tensor y)$ is clearly an element of $\LL_3$, for every tensor $x \tensor y$. Thus, fixed $x \in \LL_1$, $h$ defines a map
$$h_x: y \in \LL_2 \lmapsto h(x \tensor y) \in \LL_3.$$
Since $h$ is a homomorphism, given a family $\{y_i\}_{i \in I} \subseteq L_2$, we have
$$\begin{array}{l}
h_x\left(\bigvee_{i \in I}y_i\right) \\
= h\left(x \tensor \bigvee_{i \in I} y_i\right) \\
= h\left(\bigvee_{i \in I}(x \tensor y_i)\right) \\
= \bigvee_{i \in I} h(x \tensor y_i) \\
= \bigvee_{i \in I} h_x(y_i),
\end{array}$$
so $h_x \hom_\SL(\LL_2,\LL_3)$, for any fixed $x$. Hence we have a map $h_-: x \in L_1 \lmapsto h_x \in \hom_\SL(\LL_2,\LL_3)$, but~--- again~--- the fact that $h$ is a homomorphism implies $h_-$ to be a homomorphism as well:
$$\begin{array}{l}
h_{\bigvee_{i \in I}x_i}(y) \\
= h\left(\left(\bigvee_{i \in I} x_i\right) \tensor y\right) \\
= h\left(\bigvee_{i \in I} (x_i \tensor y)\right) \\
= \bigvee_{i \in I} h(x_i \tensor y) \\
= \bigvee_{i \in I} h_{x_i}(y)
\end{array},$$
for all $\{x_i\}_{i \in I} \subseteq L_1$ and $y \in L_2$. Besides, we also have
$$\left(\bigvee_{i \in I}h_i\right)(x \tensor y) = \bigvee_{i \in I}h_i(x \tensor y) = \bigvee_{i \in I}{h_i}_x(y) = \left(\bigvee_{i \in I}{h_i}_x\right)(y),$$
for any family $\{h_i\} \subseteq \hom_\SL(\LL_1 \tensor \LL_2, \LL_3)$, and for all $x \in L_1$ and $y \in L_2$.

Therefore we have a homomorphism
\begin{equation}\label{etasl}
\eta: \ \Hom_\SL(\LL_1 \tensor \LL_2, \LL_3) \ \lto \ \Hom_\SL(\LL_1, \Hom_\SL(\LL_2,\LL_3)),
\end{equation}
defined by $\eta(h) = h_-$, i.e. $((\eta(h))(x))(y) = h(x \tensor y)$.

Let us show that $\eta$ has an inverse. If $f \in \hom_\SL(\LL_1,\hom_\SL(\LL_2,\LL_3))$, then the map $f': (x,y) \in \LL_1 \times \LL_2 \lmapsto (f(x))(y) \in \LL_3$ is clearly a bimorphism. Hence there exists a unique homomorphism $h_{f'}: \LL_1 \tensor \LL_2 \lto \LL_3$ such that $h_{f'} \circ \tau = f'$, i.e. such that $h_{f'}(x \tensor y) = f'(x,y) = (f(x))(y)$, for all $x \in L_1$ and $y \in L_2$, and~--- clearly~--- $\eta(h_{f'}) = f$. On the other hand, if $f = \eta(h)$ with $f \in \Hom_\SL(\LL_1,\hom_\SL(\LL_2,\LL_3)$ and $h \in \Hom_\SL(\LL_1 \tensor \LL_2, \LL_3)$, then the uniqueness of the homomorphism that extends the map $f'$, defined above, to $\LL_1 \tensor \LL_2$ ensures us that $h_{f'} = h$. Then we have the inverse homomorphism
$$\eta^{-1}: \ f \in \Hom_\SL(\LL_1,\Hom_\SL(\LL_2,\LL_3)) \ \lmapsto \ h_{f'} \in \Hom_\SL(\LL_1 \tensor \LL_2, \LL_3),$$
and the theorem is proved.
\end{proof}

\begin{theorem}\label{isohomtenssl2}
Let $\LL_1$ and $\LL_2$ be sup-lattices. Then the following identities hold:
\begin{enumerate}
\item[$(i)$]$\Hom_\SL(\LL_1,\LL_2) \cong (\LL_2^{\op} \tensor \LL_1)^{\op}$,
\item[$(ii)$]$\LL_1 \tensor \LL_2 \cong \Hom_\SL\left(\LL_1, \LL_2^{\op}\right)^{\op}$.
\end{enumerate}
\end{theorem}
\begin{proof}
\begin{enumerate}
\item[$(i)$]\begin{equation*}
\begin{array}{lll}
& \quad & \Hom_\SL(\LL_1,\LL_2) \\
\textrm{by (\ref{homsldual})} & & \cong \Hom_\SL(\LL_2^{\op},\LL_1^{\op}) \\
\textrm{by (\ref{homlp1})} & & \cong \Hom_\SL(\LL_2^{\op},\Hom_\SL(\LL_1,\wp1)) \\
\textrm{by Theorem~\ref{isohomtenssl}} & & \cong \Hom_\SL(\LL_2^{\op} \tensor \LL_1,\wp1) \\
\textrm{again by (\ref{homlp1})} & & \cong (\LL_2^{\op} \tensor \LL_1)^{\op}. \\
\end{array}
\end{equation*}
\item[$(ii)$]\begin{equation*}
\begin{array}{lll}
 & \quad & \Hom_\SL(\LL_1, \LL_2^{\op})^{\op} \\
\textrm{by (\ref{homlp1})} & & \cong \Hom_\SL(\LL_1, \Hom_\SL(\LL_2,\wp1))^{\op} \\
\textrm{by Theorem~\ref{isohomtenssl}} & & \cong \Hom_\SL(\LL_1 \tensor \LL_2, \wp1)^{\op} \\
\textrm{again by (\ref{homlp1})} & & \cong ((\LL_1 \tensor \LL_2)^{\op})^{\op} \\
 & & \cong \LL_1 \tensor \LL_2.
\end{array}
\end{equation*}
\end{enumerate}
\end{proof}

\begin{remark}\label{tensorcomm}
Let us observe that if we follow the proof $(i)$ of Theorem~\ref{isohomtenssl2}, replacing the pair $(\LL_1,\LL_2)$ with $(\LL_1,\LL_2^{\op})$ and applying the self-duality $( \ )^{\op}$, we get immediately $\Hom_\SL(\LL_1, \LL_2^{\op})^{\op} \cong \LL_2 \tensor \LL_1$. Therefore it follows from the $(ii)$ that the tensor product of sup-lattices is commutative up to isomorphisms, i.e. $\LL_1 \tensor \LL_2 \cong \LL_2 \tensor \LL_1$ for any pair of sup-lattices $(\LL_1,\LL_2)$.
\end{remark}

\begin{theorem}\label{producttensor}
Let $\LL$ be a sup-lattice and $\{\LL_i\}_{i \in I}$ be a family of sup-lattices. Then
$$\coprod_{i \in I} \LL \tensor \LL_i \cong \LL \tensor \coprod_{i \in I} \LL_i \quad \textrm{and} \quad \prod_{i \in I} \LL \tensor \LL_i \cong \LL \tensor \prod_{i \in I} \LL_i. $$
\end{theorem}
\begin{proof}
Let $\LL'$ be any sup-lattice and, for all $i \in I$, $f_i: \LL \tensor \LL_i \lto \LL'$  and $g_i: \LL' \lto \LL \tensor \LL_i$ be homomorphisms. As remarked in Definition~\ref{tensors}, any element of $\LL \tensor \prod_{i \in I} \LL_i$~--- and then, by Proposition~\ref{coprodsl}, any element of $\LL \tensor \coprod_{i \in I} \LL_i$~--- is a join of tensors, i.e. can be expressed in the form $\bigvee_{j,k} x_j \tensor (y_{ik})_{i \in I}$, with $\{x_j\}_{j \in J} \subseteq L$ and $\{(y_{ik})_{i \in I}\}_{k \in K} \subseteq \prod_{i \in I} L_i$. Now, for any fixed $\ov i \in I$, let us consider the following homomorphisms: 
$$\iota_\LL \tensor \mu_{\ov i}: \ \ \bigvee_{j,k} x_j \tensor y_{\ov i k} \ \in \ \LL \tensor \LL_{\ov i} \ \lmapsto \ \bigvee_{j,k} x_j \tensor \mu_{\ov i}(y_{\ov i k}) \ \in \ \LL \tensor \prod_{i \in I} \LL_i,$$
$$\iota_\LL \tensor \pi_{\ov i}: \ \ \bigvee_{j,k} x_j \tensor \left(y_{ik}\right)_{i \in I} \ \in \ \LL \tensor \prod_{i \in I} \LL_i \ \lmapsto \ \bigvee_{j,k} x_j \tensor \pi_{\ov i}\left((y_{ik})_{i \in I}\right) \ \in \ \LL \tensor \LL_{\ov i};$$
we observe that, for all $\ov i \in I$, $(\iota_\LL \tensor \pi_{\ov i}) \circ (\iota_\LL \tensor \mu_{\ov i}) = \iota_{\LL \tensor \LL_{\ov i}}$ and, if $\ov{i'} \neq \ov i$, $(\iota_\LL \tensor \pi_{\ov{i'}}) \circ (\iota_\LL \tensor \mu_{\ov i})$ is the map $\iota_\LL \tensor \bot_{\ov{i'}}$ that sends the element $\bigvee_{j,k} x_j \tensor y_{\ov i k}$ to $\bigvee_j x_j \tensor \bot$.

What we need to prove is the existence~--- and uniqueness~--- of two homomorphisms $f: \LL \tensor \prod_{i \in I} \LL_i \lto \LL'$ and $g: \LL' \lto \LL \tensor \prod_{i \in I} \LL_i$ such that the following diagrams are commutative.
\begin{equation}\label{prodtensors}
\begin{array}{cc}
\xymatrix{
\LL \tensor \LL_i \ar[rr]^{\iota_\LL \tensor \mu_i} \ar[dd]^{f_i} && \LL \tensor \coprod_{i \in I} \LL_i \ar[ddll]^f\\
&&\\
\LL'&&\\
}
&
\xymatrix{
&& \LL' \ar[ddll]_g \ar[dd]_{g_i} \\
&&\\
\LL \tensor \prod_{i \in I} \LL_i \ar[rr]^{\iota_\LL \tensor \pi_i} && \LL \tensor \LL_i \\
}
\end{array}
\end{equation}

Since it is clear that any element $\bigvee_{j,k} x_j \tensor (y_{ik})_{i \in I}$ of $\LL \tensor \coprod_{i \in I} \LL_i$ can be written as $\bigvee_{\ov i \in I} \bigvee_{j,k} x_j \tensor \mu_{\ov i}\left((y_{ik})_{i \in I}\right)$, the first diagram in (\ref{prodtensors}) can be easily made commutative by setting, for all $\bigvee_{j,k} x_j \tensor (y_{ik})_{i \in I} \in \LL \tensor \coprod_{i \in I} \LL_i$,
$$f\left(\bigvee_{j,k} x_j \tensor (y_{ik})_{i \in I}\right) = \bigvee_{\ov i \in I} \bigvee_{j,k} f_{\ov i}\left(x_j \tensor \mu_{\ov i}\left((y_{ik})_{i \in I}\right)\right).$$
Regarding the second diagram, we define, for all $z \in \LL'$,
$$g(z) = \bigvee_{i \in I} ((\iota_\LL \tensor \mu_i) \circ g_i)(z).$$
The fact that $f$ and $g$ are homomorphisms is easily seen since the tensor product preserves joins in both coordinates and all the maps involved in their definition are homomorphisms. Then, for any fixed $\ov i \in I$ and for any $z \in L'$, we have
$$\begin{array}{l}
((\iota_\LL \tensor \pi_{\ov i}) \circ g)(z) \\
= (\iota_\LL \tensor \pi_{\ov i})\left(\bigvee_{i \in I} ((\iota_\LL \tensor \mu_i) \circ g_i)(z)\right) \\
= \bigvee_{i \in I} (\iota_\LL \tensor \pi_{\ov i})\left(((\iota_\LL \tensor \mu_i) \circ g_i)(z)\right) \\
= \bigvee_{i \in I} \left(((\iota_\LL \tensor \pi_{\ov i}) \circ (\iota_\LL \tensor \mu_i) \circ g_i)(z)\right) \\
= g_{\ov i}(z) \vee \bigvee_{i \in I} \iota_\LL \tensor \bot_i (g_i(z)) \\
= g_{\ov i}(z).
\end{array}$$
The proof of the fact that such $f$ and $g$ are unique is straightforward.

Then $\LL \tensor \coprod_{i \in I} \LL_i$ (respectively: $\LL \tensor \prod_{i \in I} \LL_i$) has the universal property of extending sinks (resp.: sources) whose domain (resp.: codomain) is the family $\{\LL \tensor \LL_i\}_{i \in I}$. The theorem is proved.
\end{proof}

\begin{corollary}
Let $\LL$ be a sup-lattice, $X$ and $Y$ be non-empty sets. Then
\begin{enumerate}
\item[$(i)$]$\wp1 \tensor \LL \cong \LL$;
\item[$(ii)$]$\wp(\mathbf{X}) \tensor \LL \cong \LL^X$;
\item[$(iii)$]$\wp(\mathbf{X}) \tensor \wp(\mathbf{Y}) \cong \wp(\mathbf{X \times Y})$.
\end{enumerate}
\end{corollary}
\begin{proof}
\begin{enumerate}
\item[$(i)$]\begin{equation*}
\begin{array}{lll}
& \quad & \wp1 \tensor \LL \\
\textrm{by Theorem~\ref{tensorcomm}} && \cong \LL \tensor \wp1 \\
\textrm{by Theorem~\ref{isohomtenssl2}} && \cong \Hom_\SL(\LL,\wp1)^{\op} \\
\textrm{by (\ref{homp1l})} && \cong \LL.
\end{array}
\end{equation*}
\item[$(ii)$]Let us denote by $\{\wp1_x\}_{x \in X}$ and $\{\LL_x\}_{x \in X}$ two families of copies, of $\wp1$ and $\LL$ respectively, with set of indices $X$. Then $\wp(\mathbf{X}) \cong \wp1^X \cong \prod_{x \in X} \wp1_x$ and $\LL^X \cong \prod_{x \in X} \LL_x$. We have:
\begin{equation*}
\begin{array}{lll}
&\quad&\wp(\mathbf{X}) \tensor \LL \\
&& \cong \left(\prod_{x \in X} \wp1_x\right) \tensor \LL \\
\textrm{by Theorems~\ref{tensorcomm} and~\ref{producttensor}} && \cong \prod_{x \in X} \wp1_x \tensor \LL \\
\textrm{by $(i)$} && \cong \prod_{x \in X} \LL_x \\
&& \cong \LL^X.
\end{array}
\end{equation*}
\item[$(iii)$]\begin{equation*}
\begin{array}{lll}
&\quad& \wp(\mathbf{X}) \tensor \wp(\mathbf{Y}) \\
\textrm{by $(ii)$} && \cong \wp(\mathbf{Y})^X \\
&& \cong \left(\wp1^Y\right)^X \\
&& \cong \wp1^{X \times Y} \\
&& \cong \wp(\mathbf{X \times Y}).
\end{array}
\end{equation*}
\end{enumerate}
\end{proof}

\section{Residuated lattices and quantales}
\label{reslatquantsec}

A binary operation $\cdot$ on a partially ordered set $\la P, \leq \ra$ is said to be \emph{residuated}\index{Residuated!-- binary operation} iff there exist binary operations $\under$ and $/$ on $P$ such that for all $x, y, z \in P$,
\begin{equation*}
x \cdot y \leq z \quad \textrm{iff} \quad x \leq z/y \quad \textrm{iff} \quad y \leq x \under z.
\end{equation*}
The operations $\under$ and $/$ are referred to as the left and right \emph{residuals}\index{Residual!Left --}\index{Residual!Right --}, or \emph{divisions}\index{Division!Left --}\index{Division!Right --}, of $\cdot$, respectively. In other words, a residuated binary operation over $\la P, \leq \ra$ is a map from $P \times P$ to $P$ that is biresiduated\index{Residuated!Bi-- map}. It follows from the results of Section~\ref{resmapsec} that the operation $\cdot$ is residuated if and only if it is order preserving in each argument and, for all $x, y, z \in P$, the inequality $x \cdot y \leq z$ has a largest solution for $x$ (namely $z/y$) and for $y$ (i.e. $x \under z$). In particular, the residuals are uniquely determined by $\cdot$ and $\leq$. The system $P = \la P, \cdot, \under, /, \leq \ra$ is called a \emph{residuated partially ordered groupoid}\index{Residuated!-- partially ordered groupoid} or \emph{residuated po-groupoid}\index{Residuated!-- po-groupoid}.

In the situations where $\cdot$ is a monoid operation with a unit element $e$ and the partial order is a lattice order, we can add the monoid unit and the lattice operations to the similarity type to get an algebraic structure $\RR = \la R, \vee, \wedge, \cdot, \under, /, e \ra$ called a \emph{residuated lattice-ordered monoid}\index{Residuated!-- lattice-ordered monoid} or \emph{residuated lattice}\index{Residuated!-- lattice} for short. It is not hard to see that $\RL$, the class of all residuated lattices, is a variety and the identities
\begin{equation*}
\begin{array}{lllll}
x \wedge (xy \vee z) / y \approx x, & \quad & x (y \vee z) \approx xy \vee xz, & \quad & (x / y)y \approx x \\
y \wedge x \under (xy \vee z) \approx y, & \quad & (y \vee z) x \approx yx \vee zx, & \quad & y(y \under x) \approx x, \\
\end{array}
\end{equation*}
together with monoid and lattice identities, form an equational basis for it.

In the category $\Q$ of quantales\footnote{In the literature, quantales are often defined as complete residuated po-groupoids, i.e. they are not unital. However, since we deal only with unital quantales, we will use this definition for avoiding notational complications.}, $\obj(\Q)$ is the class of complete residuated lattices and the morphisms are the maps preserving products, the unit, arbitrary joins and the bottom element.

An alternative definition of quantale is the following

\begin{definition}\label{quantale}
A \emph{quantale}\index{Quantale} is an algebraic structure $\QQ = \la Q, \vee, \cdot, \bot, e \ra$ such that
\begin{enumerate}
\item[$(Q1)$]$\la Q, \vee, \bot \ra$ is a sup-lattice,
\item[$(Q2)$]$\la Q, \cdot, e \ra$ is a monoid,
\item[$(Q3)$]$x \cdot \bigvee\limits_{i \in I} y_i = \bigvee\limits_{i \in I} (x \cdot y_i)$ \ and \ $\left(\bigvee\limits_{i \in I} y_i\right) \cdot x = \bigvee\limits_{i \in I} (y_i \cdot x)$ \ for all $x \in Q$, $\{y_i\}_{i \in I} \subseteq Q$.
\end{enumerate}
$\QQ$ is said to be \emph{commutative}\index{Quantale!Commutative --} iff the commutative property holds for the multiplication too.
\end{definition}

The equivalence between the two definitions above is immediate to verify, by means of the completeness of the lattice order in $\QQ$ and the properties of residuated maps. Indeed, since $\la Q, \vee, \bot \ra$ is a sup-lattice, it is possible to define the left and right residuals of $\cdot$:
\begin{equation*}
x \under y = \bigvee\{z \in Q \ \mid \ x \cdot z \leq y\},
\end{equation*}
\begin{equation*}
y / x = \bigvee\{z \in Q \ \mid \ z \cdot x \leq y\}.
\end{equation*}
Obviously, if $\QQ$ is commutative then the left and right divisions coincide:
\begin{equation*}
x \under y = y / x.
\end{equation*}

The following example will be very useful in Chapter~\ref{consrelchap}, and the subsequent results shows also the importance of such special quantales.

\begin{example}\label{powermonoid}
Let $\AA = \la A, \cdot, e \ra$ be a monoid. We define, for all $X, Y \subseteq A$,
\begin{equation}\label{complexmultiplication}
\begin{array}{l}
X \cdot Y = \{x \cdot y \mid x \in X, y \in Y\},\\
X \cdot \varnothing = \varnothing \cdot X = \varnothing.
\end{array}
\end{equation}
It is immediate to verify that the structure $\wp(\AA) = \la \wp(A), \cup, \cdot, \varnothing, \{e\} \ra$ is a quantale, whose product is often called \emph{complex multiplication}.
\end{example}

\begin{proposition}\label{freeqonm}
Let $\Q$ and $\cat M$ be, respectively, the categories of quantales and monoids. Then $\Q$ is a concrete category over $\cat M$. Moreover, given a monoid $\MM$, the free quantale over $\MM$ is $\wp(\MM)$.
\end{proposition}
\begin{proof}
The first part of the thesis is trivial. Indeed it is evident that the functor $U: \Q \lto \cat M$~--- defined as the functor that forgets the join and the bottom element~--- is faithful.

Now let us consider a monoid $\MM = \la M, \cdot, u \ra$, a quantale $\QQ = \la Q, \vee, \cdot, \bot, e \ra$ and a monoid morphism $f: \MM \lto U\QQ$. If we consider the singleton map $\sigma: x \in M \lto \{x\} \in \wp(M)$, it is obviously a monoid morphism and, if we set $h_f: X \in \wp(M) \lto \bigvee_{x \in X} f(x)$, the thesis follows as an easy application of Proposition~\ref{freesl}.
\end{proof}

Let $S$ be a non-empty set. It is well-known that the free monoid on $S$ is $\mathbf{S}^* = \la S^*, \cdot, \varnothing \ra$, where $S^* = \{x_1 \ldots x_n \mid n \in \N, x_i \in S\} \cup \{\varnothing\}$ and the product of two elements in $S^*$ is defined simply as the juxtaposition of them. With these notations, we have immediately the following result.

\begin{proposition}\label{freeqons}
Let $S$ be a non-empty set. The free quantale over $S$ is $\wp(\mathbf{S}^*)$.
\end{proposition}
\begin{proof}
If we consider the category $\cat M$ as a construct, with forgetul functor $U$, and $\Q$ as a concrete category over $\cat M$, with forgetful functor $U'$, it is immediate to verify that $U \circ U'$ is the underlying functor that makes $\Q$ a construct. Then $UU' \QQ = Q$, for any quantale $\QQ$. Now let $\QQ \in \Q$ and $f: S \lto Q$ be a map, and consider the following diagram
$$\xymatrix{
S \ar[rr]^{\id_{S^* \restr S}} \ar[ddrr]_f && S^* \ar[dd]^{Uh_f} \ar[rr]^\sigma && \wp(S^*) \ar[dd]^{UU'h_f'} \\
&&&&\\
	&& Q \ar[rr]_{\id_Q} && Q
},$$
where $\id_{S^* \restr S}$ is the inclusion map, $\sigma$ is the singleton map, $h_f$ is the unique monoid homomorphism that extends $f$ and $h_f'$ is the unique quantale homomorphism that extends $h_f$.

Then the diagram above is easily seen to be commutative, hence $h_f'$ is a homomorphism of quantales that extends $f$; the uniqueness of $h_f'$ comes from that of $h_f$. The assertion is proved.
\end{proof}

Observe that, in $\wp(\mathbf{S}^*)$, $\varnothing$ is the bottom element and $\{\varnothing\}$ is the unit.

\begin{proposition}\label{aritq}
Let $\QQ$ be a quantale. Then, for any $x, y, z \in Q$ and $\{y_i\}_{i \in I} \subseteq Q$,
\begin{enumerate}
\item[$(i)$] $x \bot = \bot x = \bot$,
\item[$(ii)$] if $x \leq y$ then $xz \leq yz$ and $zx \leq zy$,
\item[$(iii)$] if $x \leq y$ then $x/z \leq y/z$, $z \under x \leq z \under y$, $z/y \leq z/x$ and $y \under z \leq x \under z$,
\item[$(iv)$] $(y/x)x \leq y$ and $x (x \under y) \leq y$,
\item[$(v)$] $x/e = e \under x = x$,
\item[$(vi)$] $x/\left(\bigvee_{i \in I} y_i\right) = \bigwedge_{i \in I}(x/y_i)$ and $\left(\bigvee_{i \in I} y_i\right) \under x = \bigwedge_{i \in I} (y_i \under x)$,
\item[$(vii)$] $\left(\bigwedge_{i \in I} y_i\right)/x = \bigwedge_{i \in I}(y_i/x)$ and $x \under \left(\bigwedge_{i \in I} y_i\right) = \bigwedge_{i \in I}(x \under y_i)$,
\item[$(viii)$] $y \under (x \under z) = xy \under z$ and $(z/y)/x = z/xy$.
\end{enumerate}
\end{proposition}
\begin{proof}
See Proposition 2.12 of~\cite{leethesis}.
\end{proof}

\section*{References and further readings}
\addcontentsline{toc}{section}{References and further readings}

In this chapter we have recalled some essential notions of Lattice Theory and Residuation Theory. Thus we must mention, of course, the most important works on these subjects.

The classical book by G. Birkhoff, cited in~\cite{birkhoff}, is still probably the most important reference to what extent Lattice Theory. Meanwhile, the theory of residuated lattices and residuated maps finds its foundations in the pioneer paper of 1939,~\cite{warddilworth}, by M. Ward and R. P. Dilworth and its first systematic treatment in the book~\cite{blythjanowitz} by T. S. Blyth and M. F. Janowitz. The most recent overviews on the state of the art of residuated lattices are the survey paper~\cite{jipsentsinakis} by P. Jipsen and C. Tsinakis, published in 2002, and the book~\cite{nickbook} by N. Galatos, P. Jipsen, T. Kowalski and H. Ono, published in 2007.

Apart from the aforemetioned works, there exists an immense literature on these fields, and it would be impossible to cite all the remarkable papers that deal with lattices and residuated lattices. However, almost every notion and result contained in Sections~\ref{resmapsec}~and~\ref{reslatquantsec} can easily be found in one of the works cited here or in the text itself. 

A deep investigation on quantales is proposed in the book~\cite{rosenthal} by K. I. Rosenthal, and further interesting readings are the works cited in~\cite{abramskyvickers}, by S. Abramsky and S. Vickers,~\cite{pasekarosicky}, by J. Paseka and J. Rosick\'y, and~\cite{resende1,resende2,resende3,resende4}, by P. Resende.

Regarding sup-lattices, the paper~\cite{joyaltierney} by A. Joyal and M. Tierney, cited several times in the text, contains a quick presentation of most of the results presented in Sections~\ref{suplatsec}~and~\ref{sltensorsec}. Unfortunately, except where explicitly marked near the theorem number, there is not any proof or useful reference, in that paper, of the results stated for sup-lattices. Last, some applications of sup-lattices can be found, for instance, in~\cite{resende4}~--- again~--- and in~\cite{resendevickers}, by P. Resende and S. Vickers.


\part{Quantale Modules}
\label{modupart}

\chapter{Algebraic and Categorical Properties of Quantale Modules}
\label{mqchap}

Until this point, we have often driven the reader to think of quantales and quantale modules as object similar to rings and ring modules. Now that quantales have been defined, the reader will agree that the parallel drawn between quantales and rings is not forced; besides, we will see that the axioms defining quantale modules reflect the ones defining ring modules. Nevertheless, there is another kind of structure that is even closer to quantale modules: the semimodule over a semiring.

Semirings are ring-like objects themselves; the only difference between the two structures is that a semiring is not an Abelian group, but only a commutative monoid, with respect to the addition. Then a semimodule over a semiring is a commutative monoid endowed with an action from the semiring that satisfies the same conditions that define ring modules; see, for instance,~\cite{golan} for an introduction to semimodules over semirings. In a certain sense, we can look at a quantale as a special semiring whose addition may have arbitrarily infinite addenda and is idempotent. Consequently, a quantale module will be a semimodule, over this special semiring, such that also its internal operation is idempotent and is defined for an arbitrary number of arguments; moreover, the module action is arbitrarily distributive in both the arguments. Obviously, in the case where the quantale and its module are both finite, they are really~--- respectively~--- a semiring and a semimodule (both additively idempotents).

Another way to see this likeness is to observe that the join-product reduct of a bounded residuated lattice is an idempotent semiring, and we have defined quantales exactly as join-product reducts of a complete residuated lattice. Then, since every logician knows that residuated lattices are almost omnipresent in Logic, this connection, together with the results we will present in Chapter~\ref{consrelchap}, proves once more the importance of studying these kinds of algebraic structures.

This chapter can be considered the algebraic heart of this thesis. Here we will study the categories of quantale modules and, in particular, the three principal classes of operators between modules: homomorphisms, transforms and nuclei. In the first two sections, we will give the basic algebraic and categorical definitions; then, in Sections~\ref{nucleisec} and~\ref{mqtransec}, we will introduce $\Q$-module structural closure operators (nuclei, for short) and $\Q$-module transforms, respectively, also discussing their relations. In Chapter~\ref{consrelchap} we will see that deductive systems and nuclei on a quantale module are in one-one correspondence while in Chapters~\ref{imagechap} and~\ref{ltbchap} applications of $\Q$-module transforms will be shown. In Section~\ref{projinjsec} we will investigate the projective and injective objects in the categories of quantale modules and, again in Chapter~\ref{consrelchap}, we will see that the powersets of the sets of formulas, equations and sequents are all projective modules; this fact has important consequences to what extent the similarity and the equivalence of deductive systems over the same language. In Section~\ref{amalsec} we will prove that the categories of $\Q$-modules have the strong amalgamation property and, last, in Section~\ref{mqtensorsec} we will build the tensor product of $\Q$-modules.

\section{Basic notions}
\label{mqbasnotsec}

\begin{definition}\label{modules}
Let $\QQ$ be a quantale. A (left) \emph{$\QQ$-module}\index{Q-module@$\Q$-module}\index{Quantale!-- module} $\MM$, or a \emph{module over $\QQ$}, is a sup-lattice $\la M, \vee, \bot \ra$ with an external binary operation, called \emph{scalar multiplication}\index{Multiplication!Scalar --},
$$\star: (q,m) \in Q \times M \lmapsto q \star m \in M,$$
such that the following conditions hold:
\begin{enumerate}
\item[$(M1)$]$(q_1 \cdot q_2) \star m = q_1 \star (q_2 \star m)$, for all $q_1, q_2 \in Q$ and $m \in M$;
\item[$(M2)$]the external product is distributive with respect to arbitrary joins in both coordinates, i.e.
\begin{enumerate}
\item[$(i)$]for all $q \in Q$ and $\{m_i\}_{i \in I} \subseteq M$,
$$q \star {}^\MM\bigvee_{i \in I} m_i = {}^\MM\bigvee_{i \in I} q \star m_i,$$
\item[$(ii)$]for all $\{q_i\}_{i \in I} \subseteq Q$ and $m \in M$,
$$\left({}^\QQ\bigvee_{i \in I} q_i\right) \star m = {}^\MM\bigvee_{i \in I} q_i \star m,$$\
\end{enumerate}
\item[$(M3)$]$e \star m = m$.
\end{enumerate}
\end{definition}

The $(M2)$ can be expressed, equivalently, as follows:
\begin{enumerate}
\item[$(M2')$]The scalar multiplication is residuated with respect to the lattice order in $M$, i.e.
\begin{enumerate}
\item[$(i)$]for all $q \in Q$, the map
$$q^\star: m \in M \lmapsto q \star m \in M$$
is residuated,
\item[$(ii)$]for all $m \in M$, the map
$$^{\star}m: q \in Q \lmapsto q \star m \in M$$
is residuated.
\end{enumerate}
\end{enumerate}
Then, from $(M2')$ it follows that, for all $q \in Q$, there exists the residual map $(q^\star)_*$ of $q^\star$, and for all $m \in M$ there exists the residual map $(^\star m)_*$ of $^\star m$. Consequently Theorem~\ref{residchar} implies
\begin{equation}\label{qbulletsharp}
(q^\star)_*: m \in M \lmapsto \bigvee\{n \in M \mid q \star n \leq m\} \in M
\end{equation}
and
\begin{equation}\label{bulletmsharp}
(^\star m)_*: n \in M \lmapsto \bigvee\{q \in Q \mid q \star m \leq n\} \in Q.
\end{equation}
Condition (\ref{qbulletsharp}) defines another external operation over $M$:
$$\ust: (q,m) \in Q \times M \lmapsto q \ust m = (q^\star)_*(m) \in M.$$
Analogously, (\ref{bulletmsharp}) defines a map from $M \times M$ to $Q$:
$$\ost: (m,n) \in M \times M \lmapsto m \ost n = (^\star m)_*(n) \in Q.$$

The proof of the following proposition is straightforward from the definitions of $\star$, $\ust$ and $\ost$ and from the properties of quantales.

\begin{proposition}\emph{\cite{galatostsinakis}}\label{basicmq}
For any quantale $\QQ$ and any $\QQ$-module $\MM$, the following hold.
\begin{enumerate}
\item[$(i)$] The operation $\star$ is order-preserving in both coordinates.
\item[$(ii)$] The operations $\ust$ and $\ost$ preserve meets in the numerator; moreover, they convert joins in the denominator into meets. In particular, they are both order-preserving in the numerator and order reversing in the denominator.
\item[$(iii)$] $(m \ost n) \star n \leq m$.
\item[$(iv)$] $q \star (q \ust m) \leq m$.
\item[$(v)$] $m \leq q \ust (q \star m)$.
\item[$(vi)$] $(q \ust m) \ost n = q \under (m \ost n)$.
\item[$(vii)$] $((m \ost n) \star n) \ost n = m \ost n$.
\item[$(viii)$] $e \leq m \ost m$.
\item[$(ix)$] $(m \ost m) \star m = m$.
\end{enumerate}
\end{proposition}
\begin{proof}
The $(i)$ and the $(ii)$ are immediate consequences of the definitions and of Proposition~\ref{aritq}. Also the $(iii$--$v)$ come straightforwardly from the definitions of $\ost$ and $\ust$.
\begin{enumerate}
\item[$(vi)$]For all $q' \in Q$,
$$\begin{array}{l}
q' \leq q \under (m \ost n) \liff qq' \leq m \ost n \liff (qq') \star n \leq m \liff \\
q \star (q' \star n) \leq m \liff q' \star n \leq q \ust m \liff q' \leq (q \ust m) \ost n,
\end{array}$$
thus the thesis follows from the definitions of $\under$ and $\ust$.
\item[$(vii)$]By the $(i)$, $(m \ost n) \star n \leq m$, thus $((m \ost n) \star n) \ost n \leq m \ost n$. On the other hand, $((m \ost n) \star n) \ost n = \bigvee\{q \in Q \mid q \star n \leq (m \ost n) \star n\}$, but $(m \ost n) \star n \leq (m \ost n) \star n$. Therefore $m \ost n \in \{q \in Q \mid q \star n \leq (m \ost n) \star n\}$, so $m \ost n \leq ((m \ost n) \star n) \ost n$.
\item[$(viii)$]It is trivial.
\item[$(ix)$]The inequality $(m \ost m) \star m \leq m$ is obvious; the inverse inequality follows from the $(viii)$.
\end{enumerate}
\end{proof}

\begin{example}\label{freemodule}
Let $\QQ$ be a quantale and $X$ be an arbitrary non-empty set. We can consider the sup-lattice $\la Q^X, \vee^X, \bot^X \ra$, where $\bot^X$ is the $\bot$-constant function from $X$ to $Q$ and
$$(f \vee^X g)(x) = f(x) \vee g(x) \quad \textrm{ for all } x \in X.$$
Then we can define a scalar multiplication in $Q^X$ as follows:
$$\star : (q,f) \in Q \times Q^X \lmapsto q \star f \in Q^X,$$
with the map $q \star f$ defined as $(q \star f)(x) = q \cdot f(x)$ for all $x \in X$.

It is clear that $Q^X$ is a left $\QQ$-module~--- denoted by $\QQ^X$~--- and, for all $q \in Q$, $f \in Q^X$ and $x \in X$, the following holds:
$$(q \ust f)(x) = q \under f(x).$$
\end{example}

\begin{example}\label{monoidaction}
Let $S$ be a set and $\AA = \la A, \cdot, e \ra$ be a monoid from which an action $\star$ on $S$ is defined (see Section~\ref{conspwsetsec} for details). Then the sup-lattice $\wp(\SS) = \la \wp(S), \cup, \varnothing \ra$ is a module over the free quantale $\wp(\AA) = \la \wp(A), \cup, \cdot, \varnothing, \{e\} \ra$. Indeed we observed, in Section~\ref{conspwsetsec}, that the module action $\star'$ of the monoid $\la \wp(A), \cdot, \{e\} \ra$ on the set $\wp(S)$, preserves arbitrary unions both on the left and on the right.
\end{example}

The definition, and properties, of right $\QQ$-modules are completely analogous. If $\QQ$ is commutative, the concepts of right and left $\QQ$-modules coincide and we will say simply $\QQ$-modules. If a sup-lattice $\MM$ is both a left $\QQ$-module and a right $\QQ'$-module~--- over two given quantales $\QQ$ and $\QQ'$~--- we will say that $\MM$ is a \emph{$(\QQ,\QQ')$-bimodule}\index{Q-bimodule@$\Q$-bimodule} if the following associative law holds:
\begin{equation}\label{bimodule}
(q \star_l m) \star_r q' = q \star_l (m \star_r q'), \quad \textrm{for all} \ m \in M, q \in Q, q' \in Q',
\end{equation}
where, clearly, $\star_l$ and $\star_r$ are~--- respectively~--- the left and right scalar multiplications.

If $\QQ$ is a quantale and $\MM$ is a left (respectively: right) $\QQ$-module, the dual sup-lattice $\MM^{\op}$ is a right (resp.: a left) $\QQ$-module, with the external multiplication $\ust$ defined by (\ref{qbulletsharp}).

\begin{definition}\label{qmhomo}
Let $\QQ$ be a quantale and $\MM_1, \MM_2$ be two left $\QQ$-modules. A map $f: M_1 \lto M_2$ is a $\QQ$-module homomorphism if $f\left({}^{\MM_1}\bigvee_{i \in I} m_i\right) = {}^{\MM_2}\bigvee_{i \in I} f(m_i)$ for any family $\{m_i\}_{i \in I} \subseteq \MM_1$, and $f(q \star_1 m) = q \star_2 f(m)$, for all $q \in Q$ and $m \in M_1$, where $\star_i$ is the external product of $\MM_i$, for $i = 1,2$. The definition of right $\QQ$-module homomorphisms is analogous.
\end{definition}

Thus, given a quantale $\QQ$, the categories $\QQ\Mod$ and $\Modd\QQ$ have, respectively, left and right $\QQ$-modules as objects, and left and right $\QQ$-module homomorphisms as morphisms. If $\QQ$ is commutative, $\QQ\Mod$ and $\Modd\QQ$ coincide, and we denote such category by $\QQ\Mod$.
\begin{remark}\label{notation}
In what follows, in all the definitions and results that can be stated both for left and right modules, we will refer generically to ``modules''~--- without specifying left or right~--- and we will use the notations of left modules.
\end{remark}

\begin{proposition}\label{dualhomomq}
Let $\QQ$ be a quantale, $\MM_1$, $\MM_2$ be two $\QQ$-modules and $f: \MM_1 \lto \MM_2$ be a homomorphism. Then $f$ is a residuated map and the residual map $f_*: M_2 \lto M_1$ is a $\QQ$-module homomorphism between $\MM_2^{\op}$ and $\MM_1^{\op}$.
\end{proposition}
\begin{proof}
Since $f$ is a $\QQ$-module homomorphism, a fortiori it is a sup-lattice homomorphism, hence a residuated map and, as we proved in Section~\ref{suplatsec}, $f_*$ is a sup-lattice homomorphism from $\MM_2^{\op}$ to $\MM_1^{\op}$. What we need to prove, now, is that $f_*\left(q \under_{\star_2} m_2\right) = q \under_{\star_1} f_*(m_2)$. Let $m \in M_1$; we have
$$\begin{array}{llllll}
& m \leq f_*\left(q \under_{\star_2} m_2\right) & \ \liff \ & f(m) \leq q \under_{\star_2} m_2 & \ \liff \ & q \star_2 f(m) \leq m_2 \\
\liff \ & f(q \star_1 m) \leq m_2 & \ \liff \ & q \star_1 m \leq f_*(m_2) & \ \liff \ & m \leq q \under_{\star_1} f_*(m_2), \\
\end{array}$$
and the thesis follows.
\end{proof}

If $\MM_1$ and $\MM_2$ are left modules, then their duals are right modules and vice versa. Then it follows from Proposition~\ref{dualhomomq} that~--- as well as for sup-lattices~--- there is a one-one correspondence between the hom-set of two given modules and the hom-set of their dual modules. We will be back on these properties (and clarify them) in the next section.

Let $\QQ$ be a quantale and let $\MM = \la M, \vee, \bot \ra$ be a $\QQ$-module. A \emph{$\QQ$-submodule}\index{Q-submodule@$\Q$-submodule} $\NN$ of $\MM$ is a sup-sublattice $\la N, \vee, \bot \ra$ of $\la M, \vee, \bot \ra$ that is stable with respect to the external product of $\MM$. It is easy to verify that, for any family $\{\NN_i\}_{i \in I}$ of $\QQ$-submodules of $\MM$, $\left\la \bigcap_{i \in I} N_i, \vee, \bot \right\ra$ is still a $\QQ$-submodule of $\MM$. Thus, given an arbitrary subset $S$ of $M$, we define the $\QQ$-submodule $\left\la \la S \ra, \vee, \bot \right\ra$ \emph{generated by $S$} as the intersection of all the $\QQ$-submodules of $\MM$ containing $S$. Vice versa, given a submodule $\NN$ of $\MM$, we will say that a subset $S$ of $M$ is a \emph{system of generators}\index{System of generators} for $\NN$~--- or that $S$ \emph{generates} $\NN$~--- if $\NN = \la S \ra$.

If $\{\MM_i\}$ is a family of $\QQ$-modules, $\MM$ is a $\QQ$-module and $X$ is a non-empty set, then the \emph{product}\index{Q-module@$\Q$-module!Product of --s} $\prod_{i \in I} \MM_i = \la \prod_{i \in I} M_i, \vee, (\bot_i)_{i \in I}) \ra$ of the family $\{\MM_i\}_{i \in I}$, and $\MM^X = \la M^X, \vee^X, \bot^X \ra$ are clearly $\QQ$-modules with the operations defined pointwise. $\MM^X$ is also called the \emph{power module}\index{Q-module@$\Q$-module!Power --} of $\MM$ by $X$.

\begin{proposition}\label{chsubgen}
Let $\QQ$ be a quantale, $\MM$ a $\QQ$-module and $S \subseteq M$. Then $\la S \ra = \bigvee Q \star S$, where
$$\bigvee Q \star S = \left\{\bigvee_{x \in S} q_x \star x \ \Big| \ \{q_x\}_{x \in S} \in Q^S\right\}.$$
\end{proposition}
\begin{proof}
First of all we observe that $\bot = \bigvee_{x \in S} \bot \star x \in \bigvee Q \star S$; now let $\{y_i\}_{i \in I}$ be an arbitrary family of elements of $\bigvee Q \star S$, with set of indices $I$. Then, for all $i \in I$, $y_i = \bigvee_{x \in S} q_x^i \star x$, for a suitable family $\{q_x^i\}_{x \in S}$ of elements of $Q$. We have
$$\bigvee_{i \in I} y_i = \bigvee_{i \in I} \left(\bigvee_{x \in S} q_x^i \star x\right) = \bigvee_{x \in S} \left(\bigvee_{i \in I} q_x^i \star x\right) = \bigvee_{x \in S} \left(\bigvee_{i \in I} q_x^i\right) \star x \in \bigvee Q \star S.$$
Moreover, for all $y = \bigvee_{x \in S} q_x \star x \in \bigvee Q \star S$ and for all $q \in Q$,
$$q \star y = q \star \bigvee_{x \in S} q_x \star x = \bigvee_{x \in S} q \star (q_x \star x) = \bigvee_{x \in S} (q \cdot q_x) \star x \in \bigvee Q \star S;$$
hence $\bigvee Q \star S$ is a $\QQ$-submodule of $\MM$ and it is clear that $S \subseteq \bigvee Q \star S$. Therefore $\la S \ra \subseteq \bigvee Q \star S$.

On the other hand, if $\NN$ is a $\QQ$-submodule of $\MM$ containing $S$, if we fix an element $q_x \in Q$ for each $x \in S$, the scalar product $q_x \star x$ must be in $N$ for all $x \in S$; thus also $\bigvee_{x \in S} q_x \star x \in N$; hence $\bigvee Q \star S \subseteq \NN$. The arbitrary choice of $\NN$, among the submodules of $\MM$ containing $S$, ensures that $\bigvee Q \star S \subseteq \la S \ra$; the thesis follows.
\end{proof}

Given a quantale $\QQ$, a $\QQ$-module is called \emph{cyclic}\index{Q-module@$\Q$-module!Cyclic --} iff it is generated by a single element $v$. According to the notation introduced in Proposition~\ref{chsubgen}, a cyclic module generated by a certain $v$, will be also denoted by $\QQ \star v$. In Section~\ref{projinjsec} a characterization of cyclic projective modules will be shown. Such modules are very important for the applications to consequence relations and deductive systems since, as proved in Chapter~\ref{consrelchap}, the powersets of the set of formulas and of those of equations over a propositional language are cyclic projective modules over a suitable quantale. Moreover, the powerset of a set of sequents closed under type is the coproduct of cyclic projective modules.

\begin{lemma}\emph{\cite{galatostsinakis}}\label{cyclicchar}
A $\QQ$-module $\MM$ is cyclic with generator $v$ iff $(m \ost v) \star v = m$, for all $m \in M$.
\end{lemma}
\begin{proof}
If $v$ is a generator, then for all $m \in M$, there exists a $q \in Q$ such that $m = q \star v$; so $q \leq (m \ost v)$. We have $m = q \star v \leq (m \ost v) \star v \leq m$, by Proposition~\ref{basicmq}. So, $(m \ost v) \star v = m$. The converse, is obvious.
\end{proof}

\begin{definition}\label{intqmodule}
Let $\QQ$ be a quantale, and let $\MM = \la M, \vee, \bot \ra$ be a $\QQ$-module and $m$ be a fixed element of $M$. If we consider the set $M\upm = \{n \in M \mid m \leq n\}$, we can endow such set with a structure of $\QQ$-module. Indeed it is clear that $M\upm$ is closed both under arbitrary joins and meets; on the other hand, its bottom element is $m$ and we define the external operation $\star\upm$ as
\begin{equation*}
q \star\upm n = m \vee q \star n, \qquad \textrm{for all } n \in M\upm.
\end{equation*}
It is easy to verify that $\MM\upm = \la M\upm, \vee, m \ra$ is a $\QQ$-module, with such external operation. We will call it the (upper) \emph{interval $\QQ$-module}\index{Q-module@$\Q$-module!Interval --} determined by $m$ in $\MM$. For $M\upm$, we will use also the notation $[m, \top]$, indifferently.
\end{definition}

\section{Free modules, hom-sets, products and coproducts}
\label{freemqsec}

Starting from this section, we will investigate several constructions and properties of the categories of quantale modules. According to Remark~\ref{notation}, in all the definitions and statements regarding modules over a non-commutative quantale $\QQ$, whenever we say simply $\QQ$-module or write $\QQ\Mod$, we mean that the definition or the result holds for both left and right modules (suitably reformulated, where necessary).

\begin{proposition}\label{freemq}
For any set $X$, the free $\QQ$-module generated by $X$ is the function module $\QQ^X = \la Q^X, \vee^X, \bot^X \ra$, with join and scalar multiplication defined pointwise, equipped with the map $\chi: x \in X \lmapsto \chi_x \in Q^X$, where $\chi_x$ is defined, for all $x \in X$, by
\begin{equation}\label{chi}
\chi_x(y) = \left\{\begin{array}{ll} \bot & \textrm{if } y \neq x \\ e & \textrm{if } y = x\end{array}\right..
\end{equation}
\end{proposition}
\begin{proof}
Let $\MM = \la M, \vee, \bot \ra$ be any $\QQ$-module and $f: X \lto M$ be an arbitrary map. We shall prove that there exists a unique $\QQ$-module morphism $h_f: \QQ^X \lto \MM$ such that $h_f \circ \chi = f$. First observe that, for any $\alpha \in Q^X$, $\alpha = \bigvee_{x \in X} \alpha(x) \star \chi_x$; then let us set $h_f(\alpha) = \bigvee_{x \in X} \alpha(x) \star f(x)$. For any fixed $\ov x \in X$, $(h_f \circ \chi)(\ov x) = h_f(\chi_{\ov x}) = \bigvee_{x \in X} \chi_{\ov x}(x) \star f(x) = e \cdot f(\ov x) \vee \bigvee_{x \in X \setminus \{\ov x\}} \bot \cdot f(x) = f(\ov x)$, hence $h_f \circ \chi = f$.

Let now $(\alpha_i)_{i \in I}$ be a family of elements of $Q^X$ and observe that $\bigvee_{i \in I} \alpha_i = \bigvee_{i \in I} \left(\bigvee_{x \in X} \alpha_i(x) \star \chi_x\right) = \bigvee_{x \in X} \left(\bigvee_{i \in I} \alpha_i(x) \star \chi_x \right)$. Then
$$\begin{array}{lll}
h_f\left(\bigvee_{i \in I} \alpha_i\right) &=& h_f \left(\bigvee_{i \in I} \left(\bigvee_{x \in X} \alpha_i(x) \star \chi_x \right) \right) \\
 &=& h_f\left(\bigvee_{x \in X} \left(\bigvee_{i \in I} \alpha_i(x) \star \chi_x \right) \right) \\
 &=& h_f\left(\bigvee_{x \in X} \left(\bigvee_{i \in I} \alpha_i(x) \right) \star \chi_x \right) \\
 &=& \bigvee_{x \in X} \left(\bigvee_{i \in I} \alpha_i(x)\right) \star f(x) \\
 &=& \bigvee_{x \in X} \left(\bigvee_{i \in I} \alpha_i(x) \star f(x)\right) \\
 &=& \bigvee_{i \in I} \left(\bigvee_{x \in X} \alpha_i(x) \star f(x)\right) \\
 &=& \bigvee_{i \in I} h_f(\alpha_i), \\
\end{array}$$
thus $h_f$ preserves arbitrary joins; the proof of the fact that it also preserves scalar multiplication is straightforward, and therefore $h_f$ is a $\QQ$-module homomorphism. Further, the uniqueness of $h_f$ can be proved exactly as in the proof of Proposition~\ref{freesl}.
\end{proof}
Obviously, every $\Q$-module is homomorphic image of a free module.

\begin{definition}\label{hommq}
Given $\Q$-modules $\MM$ and $\NN$, we define, on the set of all the homomorphisms from $\MM$ to $\NN$ --- $\hom_\QQ(\MM,\NN)$ --- the following operations and constants:
\begin{enumerate}
\item[-]for all $h_1, h_2 \in \hom_\QQ(\MM,\NN)$, $h_1 \sqcup h_2$ is the map defined by $(h_1 \sqcup h_2)(x) = h_1(x) \vee h_2(x)$, for all $x \in M$,
\item[-]let $\bot^\bot$ and $\top^\top$ be the maps defined, respectively, by $\bot^\bot(x) = \bot_N$ and $\top^\top(x) = \top_N$, for all $x \in M$
\end{enumerate}
and, if $\QQ$ is commutative,
\begin{enumerate}
\item[-]for all $q \in Q$ and $h \in \hom_\QQ(\MM,\NN)$, let $q \diamond h$ be the map defined by $(q \diamond h)(x) = q \star h(x) = h(q \star x)$, for all $x \in M$.
\end{enumerate}
It is easy to see that $\Hom_\QQ(\MM,\NN) = \la \hom_\QQ(\MM,\NN), \sqcup, \bot^\bot \ra$ is a sup-lattice and, if $\QQ$ is a commutative quantale, it is a $\QQ$-module with the external multiplication $\diamond$. If $\NN = \MM$, the sup-lattice (or, in case, the module) of the endomorphisms $\Hom_\QQ(\MM,\MM)$ will be denoted by $\Edm(\MM) = \la \edm(\MM), \sqcup, \bot^\bot \ra$.
\end{definition}

By Proposition~\ref{dualhomomq}, $\Hom_\QQ(\MM,\NN)$ and $\Hom_\QQ(\NN^{\op},\MM^{\op})$ are isomorphic sup-lattices and, if $\QQ$ is a commutative quantale, they are isomorphic as $\QQ$-modules.

\begin{proposition}\label{coprodmq}
Let $\QQ$ be a quantale. For any given family of $\QQ$-modules $\{\MM_i\}_{i \in I}$, the coproduct $\coprod_{i \in I} \MM_i$ is the product $\prod_{i \in I} \MM_i$ equipped with the right inverses $\mu_i: \MM_i \lto \prod_{i \in I} \MM_i$ of the projections $\pi_i: \prod_{i \in I} \MM_i \lto \MM_i$. Moreover, for all $i \in I$, $\pi_i \circ \mu_i = \id_{\MM_i}$.
\end{proposition}
\begin{proof}
First of all, let us observe that, for any fixed index $j \in I$ and for any $x \in M_j$, $\mu_j(x)$ is the family of $\prod_{i \in I} \MM_i$ in which all the elements are equal to $\bot$ except the $j$-th that is equal to $x$. Thus, in particular, for any family $(x_i)_{i \in I}$, $(x_i)_{i \in I} = \bigvee_{i \in I} \mu_i(x_i)$.

Now we want to prove that, given an arbitrary $\QQ$-module $\MM$ and a family of homomorphisms $f_i: \MM_i \lto \MM$, there exists a unique homomorphism $f: \prod_{i \in I} \MM_i \lto \MM$ such that $f \circ \mu_i = f_i$ for all $i \in I$.

Let $f\left((x_i)_{i \in I}\right) = \bigvee_{i \in I} f_i(x_i)$. For all $i \in I$ and $x_i \in M_i$, we have $(f \circ \mu_i)(x_i) = f\left(\mu_i(x_i)\right) = f_i(x_i) \vee \bigvee_{j \in I \setminus \{i\}} f_j(\bot) = f_i(x_i)$.

Now, let $f': \prod_{i \in I} \MM_i \lto \MM$ be another homomorphism such that $f' \circ \mu_i = f_i$ for all $i \in I$. Then $f'\left((x_i)_{i \in I}\right) = f'\left(\bigvee_{i \in I} \mu_i(x_i)\right) = \bigvee_{i \in I} f'(\mu_i(x_i)) = \bigvee_{i \in I} f_i(x_i) = f\left((x_i)_{i \in I}\right)$. Thus $f' = f$ and also the uniqueness is proved.
\end{proof}

\section{Structural closure operators}
\label{nucleisec}

In this section we introduce a class of $\Q$-module operators, the structural closure operators~--- nuclei, for short~--- that we will encounter very often henceforward. Indeed we will see that, apart from the fact that they are strongly related with module morphisms (see Theorem~\ref{homoandnuclei}), they are the key that allows the representation of deductive systems in the framework of $\Q$-modules and they are also involved in the applications to image processing we will present in Chapters~\ref{imagechap} and~\ref{ltbchap}.

\begin{definition}\label{structural}
Let $\QQ$ be a quantale and $\MM$ a $\QQ$-module. A map $\g: M \lto M$ is said to be a \emph{structural operator}\index{Operator!Structural --}\index{Structural!-- operator} on $\MM$ provided it satisfies, for all $m, n \in M$ and $q \in Q$, the following conditions:
\begin{enumerate}
\item[$(i)$]$m \leq \g(m)$;
\item[$(ii)$]$m \leq n$ implies $\g(m) \leq \g(n)$;
\item[$(iii)$]$q \star \g(m) \leq \g(q \star m)$.
\end{enumerate}
We will say that a structural operator $\g$ is a \emph{nucleus}\index{Nucleus}\index{Q-module@$\Q$-module!-- nucleus}, or a \emph{structural closure operator}\index{Structural!-- closure operator}\index{Operator!Structural closure --}\index{Closure!Structural~--~operator}, if it is also idempotent:
\begin{enumerate}
\item[$(iv)$]$\g \circ \g = \g$.
\end{enumerate}
Then a nucleus is a closure operator that satisfies also condition $(iii)$. So it is natural to call \emph{conucleus}\index{Conucleus}\index{Nucleus!Co--}\index{Q-module@$\Q$-module!-- conucleus}\index{Structural!-- coclosure operator}\index{Operator!Structural coclosure --}\index{Closure!Structural co-- operator}\index{Interior!Structural~--~operator} a coclosure operator satisfying the $(iii)$.
\end{definition}
If $\g$ is a nucleus, we will denote by $M_\g$ the $\g$-closed system $\g[M]$ and it is easily seen that $M_\g$ is closed under arbitrary meets. Dually, the image $M^\d$ of a conucleus $\d$ is closed under arbitrary joins. In the following result we give several characterizations of structural closure operators.
\begin{lemma}\emph{\cite{galatostsinakis}}\label{modclop}
Let $\MM$ be a $\QQ$-module and let $\g$ be a closure operator on $\MM$. The following are equivalent:
\begin{enumerate}
\item[$(a)$]$\g$ is structural;
\item[$(b)$]$\g(q \star \g(m)) = \g (q \star m)$, for all $q \in Q$ and $m \in M$;
\item[$(c)$]$\g(m) \ost n = \g(m) \ost \g(n)$, for all $m, n \in M$;
\item[$(d)$]$\g(q \ust m) \leq q \ust \g(m)$, for all $q \in Q$ and $m \in M$;
\item[$(e)$]$q \ust \g(m) \in M_\g$, for all $q \in Q$ and $m \in M$. 
\end{enumerate} 
\end{lemma}
\begin{proof}
It is clear that $(a)$ is equivalent to $(b)$. To show that $(a)$ implies $(c)$, let $m, n \in M$. The inequality $\g(m) \ost \g(n) \leq \g(m) \ost n$ follows from the fact that $n \leq \g(n)$. For the reverse inequality, by the structurality of $\g$, we have
$$(\g(m) \ost n) \star \g(n) \leq \g((\g(m) \ost n) \star n) \leq \g(\g(m)) = \g(m);$$
so $\g(m) \ost n \leq \g(m) \ost \g(n)$.
 
For the converse implication, let $q \in Q$ and $m \in M$. Since $q \star m \leq \g(q \star m)$, we have $q \leq \g(q \star m) \ost m = \g(q \star m) \ost \g(m)$. Thus, $q \star \g(m) \leq \g(q \star m)$.

For the equivalence of $(a)$ and $(d)$, let $q \in Q$ and $m \in M$. We have $q \star \g(q \ust m) \leq \g(q \star (q \ust m)) \leq \g(m)$, by Proposition~\ref{basicmq}$(iv)$. Conversely, $q \star \g(m) \leq q \star \g(q \ust (q \star m)) \leq q \star (q \ust \g(q \star m)) \leq \g(q \star m)$, by Proposition~\ref{basicmq}$(v,iv)$. 
 
To show that $(a)$ implies $(e)$, let $q \in Q$ and $m \in M$. It suffices to show that $\g(q \ust \g(m)) \leq q \ust \g(m)$; i.e. $q \star \g(q \ust \g(m)) \leq \g(m)$. Indeed,
$$q \star \g(q \ust \g(m)) \leq \g(q \star (q \ust \g(m))) \leq \g(\g(m)) \leq \g(m).$$

Conversely implication, let $q \in Q$ and $m \in M$. Since $q \star m \leq \g(q \star m)$, we have $m \leq q \ust \g(q \star m)$. By the hypothesis, it follows that $\g(m) \leq q \ust \g (q \star m)$, hence $q \star \g(m) \leq \g(q \star m)$. 
\end{proof}

Before we prove the next result, recall that~--- according to Proposition~\ref{dualhomomq}~--- any $\Q$-module homomorphism is a residuated map.

\begin{theorem}\label{homoandnuclei}
Let $\QQ$ be a quantale, $\MM$ and $\NN$ be $\QQ$-modules, and $f \in \hom_\QQ(\MM,\NN)$. Then $f_* \circ f$ is a nucleus on $\MM$.

Conversely, if $\g$ is a nucleus on $\MM$, then $M_\g$~--- with the join $\vee_\g = \g \circ \vee$, the external product $\star_\g = \g \circ \star$ and the bottom element $\bot_\g = \g(\bot)$~--- is a $\QQ$-module (denoted by $\MM_\g$) and there exists $f_\g \in \hom_\QQ(\MM, \MM_\g)$ such that ${f_\g}_* \circ f_\g = \g$.
\end{theorem}
\begin{proof}
By Corollary~\ref{resclosureinterior}, if $(f,f_*)$ is an adjoint pair, then $f_* \circ f$ is a closure operator on the domain of $f$; therefore what we need to prove is condition $(iii)$ of Definition~\ref{structural} for $f_* \circ f$.

Since $(f,f_*)$ is an adjoint pair and $f$ is a homomorphism,
\begin{eqnarray}
&& q \star f(m) \leq q \star f(m) \nonumber \\ 
&\liff \quad& q \star (f \circ f_* \circ f)(m) \leq q \star f(m) \nonumber \\
&\liff \quad& f(q \star (f_* \circ f)(m)) \leq q \star f(m) \nonumber \\
&\liff \quad& q \star (f_* \circ f)(m) \leq f_* (q \star f(m)) \nonumber \\
&\liff \quad& q \star (f_* \circ f)(m) \leq (f_* \circ f)(q \star m); \nonumber \\ \nonumber
\end{eqnarray}
then $f_* \circ f$ is a nucleus.

Now let $\g$ be a nucleus on $\MM$. The fact that $\MM_\g$ is a $\QQ$-module comes easily from the fact that it is closed under arbitrary meets, and from how the join and the product have been defined. For the same reasons the map
\begin{equation*}
f_\g: m \in M \lmapsto \g(m) \in M_\g
\end{equation*}
is a $\QQ$-module homomorphism. If we apply Theorem~\ref{residchar} to $f_\g$, we get ${f_\g}_*(\g(m)) = \bigvee \{n \in M \mid f_\g(n) = \g(n) \leq \g(m)\} = \g(m)$, for any element $\g(m) \in M_\g$; thus ${f_\g}_* = \id_{M \restr M_\g}$, whence $({f_\g}_* \circ f_\g)(m) = (\id_{M \restr M_\g} \circ f_\g)(m) = \id_{M \restr M_\g}(\g(m)) = \g(m)$, for all $m \in M$. The theorem is proved.
\end{proof}

In the same hypotheses of the previous theorem, we observe explicitly that, even if $f \circ f_*$ is an interior operator in the sup-lattice $\NN$, it is not~--- in general~--- a conucleus on the $\QQ$-module $\NN$. It is, instead, a nucleus on the $\QQ$-module $\NN^{\op}$. It suffices to notice that, by Proposition~\ref{dualhomomq}, $f_*$ is a $\QQ$-module homomorphism from $\NN^{\op}$ to $\MM^{\op}$ whose residual is $f$, and to apply Theorem~\ref{homoandnuclei}.

We close this section with the following result that, although immediate, will be important in proving the characterization of cyclic projective modules we anticipated.
\begin{proposition}\label{gammacyclic}
Let $\QQ$ be a quantale. If $\g: \QQ \lto \QQ$ is a nucleus on the $\QQ$-module $\QQ$, then $\QQ_\g$ is a cyclic module and it is generated by $\g(e)$.
\end{proposition}

\section{$\Q$-module transforms}
\label{mqtransec}

In this section we introduce the $\Q$-module transforms and we prove some results about them. Then we show that any direct transform is a $\Q$-module homomorphism and a residuated map whose residual is its inverse transform.

On the other hand, $\Q$-module \emph{faithful} transforms~--- i.e. those transforms whose kernel is a coder (see Definition~\ref{coders})~--- have many further interesting properties. We will also see that the \L ukasiewicz transform defined in~\cite{dinolarusso}, and extensively discussed in Chapter~\ref{ltbchap}, is effectively a (orthonormal) $\Q$-module transform.

\begin{definition}\label{qwtransform}
Let $\QQ \in \Q$ and $X, Y$ be non-empty sets and let us consider the free $\QQ$-modules $\QQ^X$ and $\QQ^Y$. We will call a \emph{$\Q$-module transform}\index{Q-module@$\Q$-module!-- transform}\index{Transform!$\Q$-module --} between $\QQ^X$ and $\QQ^Y$, with \emph{kernel}\index{Kernel of a transform}\index{Transform!Kernel of a --} $p$, the operator
$$H_p: Q^X \lto Q^Y$$
defined by
\begin{equation}\label{qwtransformeq}
H_p f(y) = \bigvee_{x \in X} f(x) \cdot p(x,y) \quad \textrm{for all } y \in Y,
\end{equation}
where $p \in Q^{X \times Y}$. Its \emph{inverse transform}\index{Transform!Inverse --} $\Lambda_p: Q^Y \lto Q^X$ will be the map defined by
\begin{equation}\label{qwinverse transformeq}
\Lambda_p g(x) = \bigwedge_{y \in Y} g(y) / p(x,y) \quad \textrm{for all } x \in X.
\end{equation}
\end{definition}

\begin{remark}\label{righttrans}
Recalling that we are using the notations of left modules, we observe that, if we consider $\QQ^X$ and $\QQ^Y$ as right modules, the direct and inverse transforms are defined respectively by
\begin{equation}\label{qwtransformeqright}
H_p f(y) = \bigvee_{x \in X} p(x,y) \cdot f(x) \quad \textrm{for all } y \in Y,
\end{equation}
and
\begin{equation}\label{qwinverse transformeqright}
\Lambda_p g(x) = \bigwedge_{y \in Y} p(x,y) \under g(y) \quad \textrm{for all } x \in X.
\end{equation}
Up to a suitable reformulation, all the results we will present for $\Q$-module transforms hold also for free right modules.
\end{remark}

\begin{theorem}\label{wadjointpair}
Let $\QQ \in \Q$, $X, Y$ be two non-empty sets and $p \in Q^{X \times Y}$. If $H_p$ is the $\Q$-module transform, with kernel $p$, between $\QQ^X$ and $\QQ^Y$, and $\Lambda_p$ is its inverse transform, then the following hold:
\begin{enumerate}
\item[$(i)$]$(H_p,\Lambda_p)$ is an adjoint pair, i.e. $H_p$ is a residuated map and $\Lambda_p = {H_p}_*$;
\item[$(ii)$]$H_p \in \hom_{\QQ\Mod}\left(\QQ^X,\QQ^Y\right)$ and $\Lambda_p \in \hom_{\Modd\QQ}\left(\left(\QQ^Y\right)^{\op},\left(\QQ^X\right)^{\op}\right)$;
\item[$(iii)$]$\Lambda_p \circ H_p$ is a nucleus over $\QQ^X$ and $H_p \circ \Lambda_p$ is a nucleus over $\left(\QQ^Y\right)^{\op}$.
\end{enumerate} 
\end{theorem}
\begin{proof}
\begin{enumerate}
\item[$(i)$] Both $H_p$ and $\Lambda_p$ are clearly isotone; let us prove that $(H_p,\Lambda_p)$ is an adjoint pair by showing that (\ref{resid1}) and (\ref{resid2}) hold. For any $f \in Q^X$ and $g \in Q^Y$, we have:
\begin{eqnarray}
\lefteqn{H_p f \leq g \quad \iff} \nonumber \\
&& H_p f(y) \leq g(y) \quad \forall y \in Y \quad \iff \nonumber \\
&& \bigvee_{x \in X} f(x) \cdot p_Y(x,y) \leq g(y) \quad \forall y \in Y \quad \iff \nonumber \\
&& f(x) \cdot p_Y(x,y) \leq g(y) \quad \forall x \in X, \forall y \in Y \quad \iff \nonumber \\
&& f(x) \leq g(y)/p_Y(x,y) \quad \forall x \in X, \forall y \in Y \quad \iff \nonumber \\
&& f(x) \leq \bigwedge_{y \in Y} g(y)/p_Y(x,y) = \Lambda_p g(x)\quad \forall x \in X. \nonumber
\end{eqnarray}
Hence $H_p f \leq g \iff f \leq \Lambda_p g$; thus, by setting alternatively $g = H_p f$ and $f = \Lambda_p g$, we get respectively (\ref{resid2}) and (\ref{resid1}), and the $(i)$ is proved.
\item[$(ii)$] Since $H_p$ is a residuated map, it is a sup-lattice homomorphism. Moreover it is evident from the definition that $H_p$ preserves the scalar multiplication.
\item[$(iii)$] It follows from the $(ii)$ and Theorem~\ref{homoandnuclei}.
\end{enumerate}
\end{proof}

The following classification of the kernels has a few interesting theoretical implications but it is important for applications to image processing.
\begin{definition}\label{coders}
Let $\QQ \in \Q$, and $X, Y$ be non-empty sets. Let us consider a map $p \in Q^{X \times Y}$; we set the following definitions:
\begin{enumerate}
\item[$(i)$]$p$ is called a \emph{coder}\index{Coder}\index{Transform!Coder of a --} iff there exists an injective map $\e: Y \lto X$ such that $e \leq p(\e(y),y)$ for all $y \in Y$;
\item[$(ii)$]$p$ is said to be \emph{normal}\index{Coder!Normal --} iff there exists an injective map $\e: Y \lto X$ such that $p(\e(y),y) = e$ for all $y \in Y$;
\item[$(iii)$]$p$ is said to be \emph{strong}\index{Coder!Strong --} iff it is normal and
\begin{equation}\label{strong}
p(\e(y_1),y_2) = \bot \quad \textrm{for all } y_1, y_2 \in Y, y_1 \neq y_2;
\end{equation}
\item[$(iv)$]$p$ is said to be \emph{orthogonal}\index{Coder!Orthogonal --} iff $p(x,y_1) \cdot p(x,y_2) = \bot$ for all $y_1, y_2 \in Y$ such that $y_1 \neq y_2$ and for all $x \in X$;
\item[$(v)$]$p$ is said to be \emph{orthonormal}\index{Coder!Orthonormal --} iff it is orthogonal and normal.
\end{enumerate}
If $p$ is a coder, the $\Q$-module transform $H_p$ will be called \emph{faithful}\index{Transform!Faithful --}.
\end{definition}

\begin{remark}\label{orthonormal}
\begin{enumerate}
\item[$(i)$] If $p$ is normal, then it is a coder.
\item[$(ii)$] If $p$ is strong, it is a normal coder by definition.
\item[$(iii)$]If $p$ is an orthonormal map and $\e: Y \lto X$ is an injective map as in Definition~\ref{coders} $(ii)$, for any two arbitrary different elements $y_1, y_2 \in Y$, from $p(\e(y_1),y_1) \cdot p(\e(y_1),y_2) = e \cdot p(\e(y_1),y_2) = \bot$, it follows that $p(\e(y_1),y_2) = \bot$. Then any orthonormal map is a strong coder.
\end{enumerate}
\end{remark}

\begin{definition}\label{qtransform}
Let $\QQ \in \Q$, $X, Y$ be two non-empty sets and $p \in Q^{X \times Y}$ be a coder. Let us consider the faithful $\Q$-module transform $H_p$, with kernel $p$, between $Q^X$ and $Q^Y$. We set the following definitions:
\begin{enumerate}
\item[$(i)$]$H_p$ will be called a \emph{normal transform}\index{Transform!Normal --} iff $p$ is normal; 
\item[$(ii)$]$H_p$ will be called a \emph{strong transform}\index{Transform!Strong --} iff $p$ is strong; 
\item[$(iii)$]$H_p$ will be called an \emph{orthonormal transform}\index{Transform!Orthonormal --} iff $p$ is orthonormal.
\end{enumerate}
\end{definition}

\begin{theorem}\label{sadjointpair}
Let $\QQ \in \Q$ and let $H_p$ be a $\Q$-module strong transform, by the coder $p \in Q^{X \times Y}$, with inverse transform $\Lambda_p$. Then
\begin{equation*}
H_p \circ \Lambda_p = \id_{Q^Y};
\end{equation*}
thus $H_p$ is onto and, by Proposition~\ref{residprop} $(iii)$, $\Lambda_p$ is one-one.
\end{theorem}
\begin{proof}
By Theorem~\ref{wadjointpair} we have $H_p \circ \Lambda_p \leq \id_{Q^Y}$. In order to prove the inverse inequality, let us proceed as follows.

Since $p$ is strong, we can consider an injective map $\e: Y \lto X$ such that (\ref{strong}) holds. Let now $g \in Q^Y$ be an arbitrary function and let us fix an arbitrary $\ov y \in Y$. We have:
\begin{eqnarray}
\lefteqn{(H_p \circ \Lambda_p) g(\ov y) =} \nonumber \\
&& \bigvee_{x \in X}\left(\left(\bigwedge_{y \in Y} g(y) / p(x,y)\right) \cdot p(x,\ov y)\right) \geq \nonumber \\
&& \left(\bigwedge_{y \in Y} g(y) / p(\e(\ov y),y)\right) \cdot p(\e(\ov y),\ov y) = \nonumber \\
&& \left(\bigwedge_{\mathop{}^{y \in Y}_{y \neq \ov y}} g(y) / p(\e(\ov y),y)\right) \wedge \left(g(\ov y) / p(\e(\ov y),\ov y)\right) = \nonumber \\
&& \left(\bigwedge_{\mathop{}^{y \in Y}_{y \neq \ov y}} g(y) / \bot\right) \wedge \left(g(\ov y) / e\right) = \nonumber \\
&& \top \wedge g(\ov y) = \nonumber \\
&& g(\ov y) \nonumber.
\end{eqnarray}
Since the above relations hold for all $g \in Q^Y$ and $\ov{y} \in Y$, the thesis is proved.
\end{proof}

In what follows we will always assume that $Y$ is a subset of $X$ and, if $p \in Q^{X \times Y}$ is a coder, then the map $\e$ is the inclusion map $\id_{X \restr Y}: y \in Y \lmapsto y \in X$.

\begin{lemma}\label{transf=coder}
Let $\QQ \in \Q$, $X$ be a non-empty set, $Y$ be a non-empty subset of $X$ and $p, p' \in Q^{X \times Y}$ be two maps. Then $H_p = H_{p'}$ if and only if $p = p'$.
\end{lemma}
\begin{proof}
Since the other implication is trivial, let us prove that $H_p = H_{p'}$ implies $p = p'$ by showing that, if $p \neq p'$, then $H_p \neq H_{p'}$.

By assumption, there exists a pair $(\ov x,\ov y) \in X \times Y$ such that $p(\ov x, \ov y) \neq p'(\ov x, \ov y)$. Let us consider the map $f \in Q^X$ defined by
$$f(x) = \left\{
			\begin{array}{ll}
			e & \textrm{if } x = \ov x \\
			\bot & \textrm{if } x \in X \setminus \{\ov x\} \\
			\end{array}.
			\right.
$$
It is immediate to verify that $H_p f(\ov y) = p(\ov x, \ov y) \neq p'(\ov x, \ov y) = H_{p'} f(\ov y)$, and the thesis follows.
\end{proof}

The previous result ensures us that a $\Q$-module transform $H_p$ is completely determined by its kernel $p$.

\begin{lemma}\label{projection}
Let $\QQ \in \Q$, $X$ be a non-empty set, $Y$ be a non-empty subset of $X$ and $p \in Q^{X \times Y}$ be a coder. Then, for any fixed $y \in Y$, $H_p f(y) = f(y)$ for all $f \in Q^X$ if and only if
\begin{equation}\label{projcoder}
p(x,y) = \left\{
				\begin{array}{ll}
				e & \textrm{if } x = y \\
				\bot & \textrm{if } x \neq y
				\end{array}.\right.
\end{equation}
\end{lemma}
\begin{proof}
If (\ref{projcoder}) holds, then $H_p f(y) = f(y)$ for all $f \in Q^X$, trivially.

On the other hand, if (\ref{projcoder}) does not hold for $p$, then we distinguish two cases:
\begin{enumerate}
\item[Case 1:]$p(y, y) = q_1 \neq e$, for some $q_1 \in Q$;
\item[Case 2:]there exists $\ov x \in X$ such that $p(\ov x, y) = q_2 \gneqq \bot$.
\end{enumerate}
In the first case, let $f \in Q^X$ be the map defined by
$$f(x) = \left\{
			\begin{array}{ll}
			e & \textrm{if } x = y \\
			\bot & \textrm{if } x \neq y \\
			\end{array}.
			\right.
$$
It is easy to see that $H_p f(y) = q_1 \neq e = f(y)$.

In the second case, let $g \in Q^X$ be the map defined as
$$g(x) = \left\{
			\begin{array}{ll}
			e & \textrm{if } x = \ov x \\
			\bot & \textrm{if } x \neq \ov x \\
			\end{array}.
			\right.
$$
Clearly $H_p g(y) = q_2 \neq \bot = g(y)$, and the lemma is proved.
\end{proof}

Given the sets $X$ and $Y$, we will denote by $\pi_Y$ the coder defined, for all $x \in X$ and for all $y \in Y$, by (\ref{projcoder}) and we will call it a \emph{projective coder}\index{Coder!Projective --}. By Lemma~\ref{projection}, for all $f \in Q^X$, $H_{\pi_Y} f = f_{\restr Y}$, i.e. $H_{\pi_Y}$ is the projection of $Q^X$ on $Q^Y$.

\begin{definition}\label{psupport}
If $p$ is a coder of $Q^{X \times Y}$, let $Y_p' \subseteq Y$ be the set of all the elements $y$ of $Y$ such that $p(x,y)$ is defined by (\ref{projcoder}):
$$Y_p' = \left\{y \in Y \ | \ p(x,y) = \pi_Y(x,y) \right\}.$$
The set $\dot Y_p = Y \setminus Y_p'$ will be called the \emph{support}\index{Coder!Support of a --} of $p$ and the restriction $\ul p = p_{\restr \dot Y_p}$ will be called the \emph{core}\index{Coder!Core of a --} of $p$. If $\dot Y_p = Y$, then $p = \ul p$ and we will say that $p$ is \emph{irreducible}\index{Coder!Irreducible --}; $p$ is \emph{reducible}\index{Coder!Reducible --} if $\dot Y_p \subsetneqq Y$.
\end{definition}

\begin{definition}\label{closurecod}
Given a coder $p \in Q^{X \times Y}$ and a set $Z$ such that $Y \subseteq Z \subseteq X$, let us consider the extension $p^Z$ of the coder $p$ to $X \times Z$, defined as follows:
\begin{equation*}
p^Z(x,z) = \left\{
						\begin{array}{ll}
							p(x,z) & \textrm{if } (x,z) \in X \times Y \\
							\pi_{Z \setminus Y}(x,z) & \textrm{if } (x,z) \in (X \times Z) \setminus (X \times Y) \\
						\end{array}.
						\right.
\end{equation*}
The coder $p^Z$ will be called the \emph{projective extension of $p$ to $Z$}\index{Coder!Projective extension of a --}. In this case, it is clear that $\dot Y_p = \dot Z_{p^Z}$ and $\ul p = \ul p^Z$.

If $Z = X$, we will denote $p^X$ by $\ov p$ and we will call it the \emph{closure}\index{Coder!Closure of a --}\index{Closure!-- of a coder} of $p$. So $\ov p$ is the coder defined by 
\begin{equation*}
\ov p(x,z) = \left\{
						\begin{array}{ll}
							p(x,z) & \textrm{if } (x,z) \in X \times Y \\
							\pi_{X \setminus Y}(x,z) & \textrm{if } (x,z) \in (X \times X) \setminus (X \times Y) \\
						\end{array}.
						\right.
\end{equation*}
Clearly, for any coder $p \in Q^{X \times X}$, $\ov p = p$; therefore, such coders will be called \emph{closed coders}\index{Coder!Closed --}.
\end{definition}

\begin{definition}\label{equiproj}
Let $\QQ \in \Q$, $X$ be a non-empty set, $Y, Z$ be two non-empty subsets of $X$ and $p \in Q^{X \times Y}$, $p' \in Q^{X \times Z}$ be two coders. We will say that $p$ and $p'$ are \emph{equivalent up to projections}~--- and we will write $p \doteq p'$~--- iff $\ul p = \ul p'$, i.e. iff $\dot Y_p = \dot Z_{p'}$ and $p_{\restr \dot Y_p} = p'_{\restr \dot Y_p}$.
\end{definition}

\begin{proposition}\label{eqclosure}
Let $\QQ \in \Q$, $X$ be a non-empty set, $Y, Z$ be two non-empty subsets of $X$ and $p \in Q^{X \times Y}$, $p' \in Q^Z$ be two coders. Then
$$p \doteq p' \quad \iff \quad \ov p = \ov p'.$$
In other words, $p$ and $p'$ are equivalent up to projections if and only if they have the same closure.
\end{proposition}
\begin{proof}
It is trivial.
\end{proof}
The last definitions and Proposition~\ref{eqclosure} are significant, again, for applications. In the next result we invert Theorem~\ref{wadjointpair}, showing that all the homomorphisms between free modules are transforms.

\begin{theorem}\label{repr}
The sup-lattices $\Hom_\QQ(\QQ^X,\QQ^Y)$ and $\QQ^{X \times Y}$ are isomorphic.
\end{theorem}
\begin{proof}
Let
\begin{equation}\label{hbarxy}
\hbar: \ Q^{X \times Y} \ \lto \ \hom_\QQ(\QQ^X,\QQ^Y)
\end{equation}
be the map defined by $\hbar(p) = H_p$, for all $p \in Q^{X \times Y}$; in other words $\hbar$ sends every map $p \in Q^{X \times Y}$ in the transform between $\QQ^X$ and $\QQ^Y$ whose kernel~is~$p$.

The fact that $\hbar$ is injective comes directly from Lemma~\ref{transf=coder}. Moreover it is clear that $\hbar(\bot^{X \times Y}) = \bot^\bot$. Now let $\{k_i\}_{i \in I} \subseteq Q^{X \times Y}$; we must prove that $\hbar\left(\bigvee_{i \in I} k_i \right) = \bigsqcup_{i \in I}\hbar(k_i)$. For any $f \in Q^X$ and for all $y \in Y$, we have
\begin{eqnarray}
\lefteqn{\hbar\left(\bigvee_{i \in I} k_i \right) f(y) = \bigvee_{x \in X} f(x) \cdot \left(\bigvee_{i \in I} k_i \right)(x,y)} \nonumber\\
&&= \bigvee_{x \in X} f(x) \cdot \left(\bigvee_{i \in I} k_i(x,y) \right) = \bigvee_{x \in X} \bigvee_{i \in I} f(x) \cdot k_i(x,y) \nonumber\\
&&= \bigvee_{i \in I} \bigvee_{x \in X} f(x) \cdot k_i(x,y) = \bigvee_{i \in I} \hbar(k_i)f(y) \nonumber\\
&&= \left(\bigvee_{i \in I} \hbar(k_i)f\right)(y) = \left(\bigsqcup_{i \in I} \hbar(k_i)\right)f(y); \nonumber
\end{eqnarray}
whence $\hbar$ is a sup-lattice monomorphism. 

Now we must prove that $\hbar$ is surjective too. Let $h \in \hom_\QQ(\QQ^X,\QQ^Y)$ and, for any $x \in X$, let us consider the map $\chi_x$ defined by (\ref{chi}). Let now $k^h \in Q^{X \times Y}$ be the function defined by
\begin{equation*}
k^h(x,y) = h \chi_x(y), \qquad \textrm{for all } (x,y) \in Q^{X \times Y};
\end{equation*}
then we have
\begin{eqnarray}
\lefteqn{h f(y) = h\left(\bigvee_{x \in X} f(x) \star \chi_x\right)(y)} \nonumber \\
&& = \left(\bigvee_{x \in X} f(x) \star h \chi_x\right)(y) = \bigvee_{x \in X} f(x) \cdot h \chi_x(y) \nonumber \\
&& = \bigvee_{x \in X} f(x) \cdot k^h(x,y) = H_{k^h} f(y), \nonumber \\ \nonumber
\end{eqnarray}
for all $f \in Q^X$ and for all $y \in Y$. It follows that $h = H_{k^h} = \hbar(k^h)$, hence $\hbar$ is a sup-lattice isomorphism whose inverse is~--- obviously~--- defined by 
$$\hbar^{-1} h(x,y) = h \chi_x(y),$$
for all $h \in \hom_\QQ(\QQ^X,\QQ^Y)$ and $(x,y) \in X \times Y$.
\end{proof}

The previous theorem allows us to define a structure of $\QQ$-module on $\hom_\QQ(\QQ^X,\QQ^Y)$, also when $\QQ$ is not commutative, by defining the external multiplication
\begin{equation}\label{asthom}
\begin{array}{cccc}
\star: & Q \times \hom_\QQ(\QQ^X,\QQ^Y) & \lto 		& \hom_\QQ(\QQ^X,\QQ^Y) \\
			& (q,h) 												 & \lmapsto & \hbar(q \star \hbar^{-1}(h))
\end{array},
\end{equation}
in such a way that this $\QQ$-module is isomorphic to $\QQ^{X \times Y}$; we will denote this structure by $\Hom_\QQ^\star(\QQ^X,\QQ^Y)$.

\section{Projective and injective $\Q$-modules}
\label{projinjsec}

In the present section we study projective and injective objects in the categories of quantale modules. A characterization of projective cyclic modules is proposed, together with a result (Theorem~\ref{projlem3}) that allows the identification of a large class of projective modules and, by duality, of a class of injective objects. Moreover, as we mentioned several times, such cyclic projective modules are of great importance in the applications to deductive systems.
\begin{proposition}\label{projmqlem0}
Let $\QQ \in \Q$. Then a left $\QQ$-module $\MM$ is projective iff the right module $\MM^{\op}$ is injective, and vice versa.
\end{proposition}
\begin{proof}
It is a trivial application of the basic properties of $\Q$-module homomorphisms.
\end{proof}

\begin{proposition}\label{freeproj}
Every free $\Q$-module is projective.
\end{proposition}
\begin{proof}
Let $X$ be a non-empty set, and $\MM, \NN$ be $\QQ$-modules such that there exist an epimorphism $g: \MM \lto \NN$ and a homomorphism $f: \QQ^X \lto \NN$. Then $f$ the unique homomorphism that extends the map $f_X: x \in X \lmapsto f(\chi_x) \in \NN$. If we consider the map $g_* \circ f_X: X \lto \MM$, we can extend it to a homomorphism $h: \QQ^X \lto \MM$, and it is immediate to verify that $g \circ h = f$. Then any free module is projective.
\end{proof}

\begin{lemma}\emph{\cite{galatostsinakis}}\label{astarv}
Let $\QQ$ be a quantale, $\MM$ a $\QQ$-module, $v \in M$ and consider the $\QQ$-module $\QQ \star v = \{q \star v \mid q \in Q\}$. Then we have: 
\begin{enumerate}
\item[$(i)$] the map $\g_v: Q \lto Q$, defined by $\g_v(q) = (q \star v) \ost v$ is a structural closure operator;
\item[$(ii)$] $\QQ \star v$ is isomorphic to $\QQ_{\g_v}$.
\end{enumerate}
Consequently, a $\QQ$-module is cyclic iff it is isomorphic to a module $\QQ_\g$, where $\g:\QQ \lto \QQ$ is a structural closure operator.\footnote{Recall that this condition is equivalent to that of being isomorphic to a quotient (or, that is the same, being homomorphic image) of the $\QQ$-module $\QQ$.} 
\end{lemma}
\begin{proof}
$[(i)]$ \quad The map $\g_v$ is clearly extensive; moreover, if $q \leq r$, then $\g_v(q) \leq \g_v(r)$, and $\g_v(\g_v(q)) = \g_v(q)$, by Proposition~\ref{basicmq}$(vii)$. Also, $q \g_v(r) \star v = q((r \star v) \ost v) \star v \leq qr \star v$, so $q \g_v(r) \leq \g_v(qr)$. Thus, $\g_v$ is structural.

$[(ii)]$ \quad For all $x \in Q \star v$, $x = (x \ost v) \star v$, hence $x \ost v = ((x \ost v) \star v) \ost v \in Q_{\g_v}$ by the definition of $\g_v$. Then we set
$$f: q \in Q_{\g_v} \lmapsto q \star v \in Q \star v \qquad \textrm{ and } \qquad g: x \in Q \star v \lmapsto x \ost v \in Q_{\g_v}.$$
For all $x \in Q \star v$, we have $f(g(x)) = (x \ost v) \star v = x$, by Lemma~\ref{cyclicchar}. Also, for all $q \in Q_{\g_v}$, $g(f(q)) = \g_v(q) = q$. So, $f^{-1} = g$. The fact that $f$ and $g$ are homomorphisms is trivial. 
\end{proof}
 
\begin{corollary}\label{cor8.5}
If $\QQ$ is a quantale and $u \in Q$, then the $\QQ$-module $\QQ \cdot u$ is isomorphic to $\QQ_{\g_u}$, where $\g_u: q \in Q \lmapsto qu/u \in Q$.
\end{corollary}

\begin{lemma}\emph{\cite{galatostsinakis}}\label{uidempotent}
Let $\QQ$ be a quantale, $\g: Q \lto Q$ a structural closure operator and $u \in Q$. The following are equivalent.
\begin{enumerate}
\item[$(a)$] $\g(q)u = qu$, for all $q \in Q$, and $\g(u) = \g(e)$;
\item[$(b)$] $\g = \g_u$ and $u = u^2$.
\end{enumerate}
\end{lemma}
\begin{proof}
The fact that $(b)$ implies $(a)$ is trivial. Conversely, from $\g(q)u = qu$, we obtain $\g(q) \leq qu/u = \g_u(q)$, for all $q \in Q$. Also, from $\g(u) = \g(e)$, we obtain, for all $r \in Q$, $\g(ru) = r \cdot_\g \g(u) = r \cdot_\g \g(e) = \g(re) = \g(r)$. Then we have
$$\g_u(q)u \leq qu \lTo \g(\g_u(q)u) \leq \g(qu) \lTo\g(\g_u(q)) \leq \g(q) \lTo\g_u(q) \leq \g(q).$$
Moreover, since $\g = \g_u$, we have $u^2 = \g(u)u = \g(e)u = e \cdot u = u$, which completes the proof.
\end{proof}

\begin{corollary}\emph{\cite{galatostsinakis}}\label{cycliccond}
Let $\MM$ be a cyclic $\QQ$-module with generator $v$. The following conditions are equivalent for an element $u \in Q$.
\begin{enumerate}
\item[$(a)$] $u \star v = v$ and $((q \star v)\ost v)u = qu$, for all $q \in Q$;
\item[$(b)$] $\g_v(q)u = qu$, for all $q \in Q$, and $\g_v(u) = \g_v(e)$;
\item[$(c)$] $\g_v = \g_u$ and $u^2 = u$;
\item[$(d)$] $\MM$ is isomorphic to $\QQ \cdot u$ and $u^2 = u$.
\end{enumerate}
\end{corollary}
\begin{proof}
The equivalence $(a) \liff (b)$ follows from the fact that $\g_v(u) = \g_v(e)$ iff $(u \star v) \ost v = v \ost v$ iff $u \star v = v$, by Proposition~\ref{basicmq}.

The implication $(b) \lTo(c)$ follows from the Lemma~\ref{uidempotent}, while the implication $(c) \lTo(d)$ follows from $\QQ \star v \cong \QQ_{\g_v}$ (Lemma~\ref{astarv}), $\QQ \cdot u \cong \QQ_{\g_u}$ (Corollary~\ref{cor8.5}) and $\g_u = \g_v$.

Finally, $(d) \lTo(a)$ follows from the fact that if $u^2 = u$, then $\QQ \cdot u$ satisfies the $(a)$ with $v = u$.
\end{proof}

\begin{theorem}\emph{\cite{galatostsinakis}}\label{cyclicproj}
A cyclic $\QQ$-module $\MM$, with generator $v$, is projective if and only if there exists an element $u \in Q$ that satisfies the equivalent conditions of Corollary~\ref{cycliccond}.
\end{theorem}
\begin{proof}
Every cyclic module is of the form $\QQ_\g$ for some nucleus $\g: Q \lto Q$, by Lemma~\ref{astarv}. Suppose that $\QQ_\g$ is projective; we will verify condition $(d)$ of Lemma~\ref{cycliccond}. Since $\QQ_\g$ is projective, there exists a module morphism $f$ that completes the diagram below.
$$\xymatrix{
\QQ_\g \ar@{_{(}-->}@<.5ex> [rr]^f \ar@{^{(}->>}@<-.5ex>[ddrr]_{\id_{\QQ_\g}} && \QQ \ar@{->>}[dd]^\g \\
&& \\
&& \QQ_\g
}$$
Now let $u = f(\g(e))$. For all $q \in Q$, we have $\g(q) = \g(q e) = \g(q \g(e)) = q \cdot_\g \g(e)$, so $f(\g(q)) = q \cdot f(\g(e)) = q u$. Consequently, $f[Q_\g] = Q \cdot u$. Moreover, $f$ is injective, by the diagram, so $\QQ_\g \cong \QQ \cdot u$.

We will now show that $u^2 = u$. Indeed, $u^2 = f(\g(e))f(\g(e)) = f(f(\g(e)) \cdot_\g \g(e)) = f(\g(f(\g(e)))) = f(\g(e)) = u$, because $\g \circ f = \id_{\QQ_\g}$. The $(d)$ of Lemma~\ref{cyclicproj} follows.

Then, if $\MM$ is projective, the equivalent conditions of Lemma~\ref{cyclicproj} hold. The opposite implication is trivial, and the theorem is proved.
\end{proof}

\begin{theorem}\label{projlem3}
If $\{\MM_i\}_{i \in I}$ is a family of projective (respectively: injective) modules, then $\prod_{i \in I} \MM_i$ is projective (resp.: injective) as well.
\end{theorem}
\begin{proof}
Before going through the proof, we recall that, by Proposition~\ref{coprodmq}, $\coprod_{i \in I} \MM_i$ and $\prod_{i \in I} \MM_i$ coincide as objects, and they differ from each other just in the direction of the respective families of morphisms: $\{\mu_i\}_{i \in I}$ and $\{\pi_i\}_{i \in I}$. Because of that, we will prove only the projective case, the injective one differing exclusively in the direction of the arrows and in the use of $\{\pi_i\}_{i \in I}$ instead of $\{\mu_i\}_{i \in I}$; a similar proof can be found in~\cite{galatostsinakis}.

Let $\MM = \prod_{i \in I} \MM_i$, $\NN$ and $\PP$ be two $\QQ$-modules, and $f: \MM \lto \PP$ and $g: \NN \lto \PP$ be homomorphisms, with $g$ surjective. For any $\MM_i$, we can consider the homomorphism $f_i = f \circ \mu_i$ and, because of the projectivity of $\MM_i$, there exists a homomorphism $h_i: \MM_i \lto \NN$ such that $g \circ h_i = f_i$. Since $\MM$ is the coproduct of the family $\{\MM_i\}$, there exists a homomorphism $h: \MM \lto \NN$ such that $h \circ \mu_i = h_i$ for all $i \in I$, hence $f \circ \mu_i = f_i = g \circ h_i = g \circ h \circ \mu_i$. Then we have two morphisms $f$ and $g \circ h$ that extends the same family of morphisms $\{\MM_i \stackrel{f_i}{\lto} \NN\}_{i \in I}$ to the coproduct $\prod_{i \in I} \MM_i$ and, by the definition of coproducts, such morphisms must coincide, i.e. $f = g \circ h$. The proof is clarified by the first of the following diagrams, while the second one describes the case of injectivity:
\begin{equation*}
\begin{array}{lr}
\xymatrix{
&&&& \NN \ar @{->>} [dddd]^g \\
&&&&\\
\MM \ar[ddrrrr]_{f} \ar[uurrrr]^h && \MM_i \ar @{^{(}->} [ll]^{\mu_i} \ar[ddrr]_{f_i} \ar[uurr]^{h_i} && \\
&&&&\\
&&&& \PP\\
}
&
\xymatrix{
\NN \ar[ddrr]^{h_i} \ar[ddrrrr]^h &&&& \\
&&&&\\
 && \MM_i \ar @{->>} [rr]^{\pi_i}   && \MM  \\
&&&&\\
\PP \ar @{_{(}->} [uuuu]^g \ar[uurr]_{f_i} \ar[uurrrr]_{f} &&&& \\
}
\end{array}.
\end{equation*}

The thesis follows.
\end{proof}

\section{On the amalgamation property}
\label{amalsec}

Given a category $\cat C$, an \emph{amalgam}\index{Amalgam} in $\cat C$ is a 5-tuple $(A, f, B, g, C)$, where $A, B, C \in \obj(\cat C)$ and $f: A \lto B$, $g: A \lto C$ are injective morphisms.

\begin{definition}\label{ama}
We say that a category $\cat C$ \emph{has the amalgamation property}\index{Amalgamation property} if whenever an amalgam $(A, f, B, g, C)$ is given, with $A \neq \varnothing$, there exists an object $D$ and two injective morphisms $f': B \lto D$ and $g': C \lto D$ such that $f' \circ f = g' \circ g$.

A category $\cat C$ is said to have the \emph{strong amalgamation property}\index{Amalgamation property!Strong --} if it has the amalgamation property with $f'[B] \cap g'[C] = (f' \circ f)[A] = (g' \circ g)[A]$.
\end{definition}

In the next result we prove that any category of $\Q$-modules has the strong amalgamation property. Even if, in the whole thesis, we assumed our quantales to be unital, we observe explicitly that in the proof of Theorem~\ref{amalgam} such an assumption is not used at all. Then also the categories of modules on quantales that are not unital have the strong amalgamation property.

\begin{theorem}\label{amalgam}
$\QQ\Mod$ and $\Modd\QQ$ have the strong amalgamation property, for any quantale $\QQ$.
\end{theorem}
\begin{proof}
As usual, we will prove the assertion for $\QQ\Mod$. Let $\MM$, $\NN$ and $\PP$ be $\QQ$-modules, and $f: \PP \lto \MM$ and $g: \PP \lto \NN$ be two injective morphisms. Then $\PP$, $f[\PP]$ and $g[\PP]$ are isomorphic $\QQ$-modules. Let us consider the coproduct $\MM \times \NN$ with the associated (injective) morphisms
$$\mu_M: m \in \MM \lmapsto (m, \bot_\NN) \in \MM \times \NN, \ \quad \mu_N: n \in \NN \lmapsto (\bot_\MM, n) \in \MM \times \NN,$$
and let $\sim$ be the smallest $\QQ$-module congruence over $\MM \times \NN$ such that $(f(p),\bot_\NN) \sim (\bot_\MM,g(p))$ for all $p \in P$.

It is easy to verify, by direct inspection, that
\begin{equation}\label{simama}
(m,n) \sim (m',n') \quad \liff \quad \begin{array}{l}\exists p, p' \in P: \ \left\{\begin{array}{l} f(p) = m \\ f(p') = m' \\ g(p) = n' \\ g(p') = n\end{array}\right.\\
\text{or } (m,n) = (m',n').\end{array}
\end{equation}
Indeed, (\ref{simama}) trivially defines an equivalence relation. If (\ref{simama}) holds for two families $\{(m_i,n_i)\}_{i \in I}, \{(m_i',n_i')\}_{i \in I} \in M \times N$, with $\{p_i\}_{i \in I}$ and $\{p_i'\}_{i \in I}$ as the corresponding families of elements of $P$, then $\bigvee_{i \in I} m_i = \bigvee_{i \in I} f(p_i) = f\left(\bigvee_{i \in I} p_i\right)$, and similar equalities can be stated for the other families of elements of $M$ and $N$. Then (\ref{simama}) defines a sup-lattice congruence. On the other hand, it is immediate as well to verify that the relation defined by (\ref{simama}) is invariant with respect to the external multiplications, whence it is a module congruence that~--- clearly~--- contains $\sim$.

Conversely, if $(m,n)$ and $(m',n')$ are two pair in $M \times N$ that satisfy the condition of (\ref{simama}), then in particular $(m,\bot_\NN) \sim (\bot_\MM,n')$ and $(m',\bot_\NN) \sim (\bot_\MM,n)$. So
$$(m,n) = (m,\bot_\NN) \vee (\bot_\MM,n) \sim (m',\bot_\NN) \vee (\bot_\MM,n') = (m',n'),$$
whence $\sim$ is really completely described by (\ref{simama}).

Then, since $f$ and $g$ are injective, if $m, m' \in M$ and $n, n' \in N$, $(m,\bot_\NN) \sim (m',\bot_\NN)$ iff $m = m'$ and $(\bot_\MM, n) \sim (\bot_\MM, n')$ iff $n = n'$.

Now, if we set $\RR = (\MM \times \NN)/\sim$, we can define the following $\QQ$-module homomorphisms
$$f': m \in \MM \lmapsto [\mu_M(m)]_\sim \in \RR \quad \textrm{ and } \quad g': n \in \NN \lmapsto [\mu_N(n)]_\sim \in \RR,$$
and it follows from the considerations above that $f'$ and $g'$ are both injective. Moreover, the fact that $f' \circ f = g' \circ g$ is an immediate consequence of how $\sim$ has been defined; hence the categories of quantale modules have the amalgamation property.
\begin{equation}
\xymatrix{
& \RR & \\
&&\\
& \MM \times \NN \ar@{->>} [uu]^{\pi_\sim} & \\
\MM \ar@{_{(}->} [ur]^{\mu_M} \ar@{_{(}->}@<.6ex> [uuur]^{f'} && \NN \ar@{^{(}->} [ul]_{\mu_N} \ar@{^{(}->}@<-.6ex> [uuul]_{g'} \\
& \PP \ar@{_{(}->} [ul]^{f} \ar@{^{(}->} [ur]_{g} &
}
\end{equation}

Now observe that
$$f'[M] = \{[(m,\bot_\NN)]_\sim \mid m \in M\} \quad \textrm{ and } \quad g'[N] = \{[(\bot_\MM,n)]_\sim \mid n \in N\}.$$
Since, obviously, $(f' \circ f)[P] = (g' \circ g)[P] \subseteq f'[M] \cap g'[N]$, in order to prove that $\QQ\Mod$ has the strong amalgamation property, we just need to show that $(f' \circ f)[P] \supseteq f'[M] \cap g'[N]$ or, equivalently, $(g' \circ g)[P] \supseteq f'[M] \cap g'[N]$.

We first observe that $f'[M] = \{[(m,\bot_\NN)]_\sim \mid m \in M\}$ and $g'[N] = \{[(\bot_\MM,n)]_\sim \mid n \in N\}$; then, if $[(m,n)]_\sim \in f'[M] \cap g'[N]$, there exist $m' \in M$ and $n' \in N$ such that $(m,n) \sim (m',\bot_\NN)$ and $(m,n) \sim (\bot_\MM,n')$. Therefore, by (\ref{simama}), there exists $p \in P$ such that $f(p) = m'$ and $g(p) = n'$. Thus
$$[(m,n)]_\sim = (f' \circ f)(p) = (g' \circ g)(p) \in (f' \circ f)[P] = (g' \circ g)[P],$$
i.e. $(f' \circ f)[P] = (g' \circ g)[P] = f'[M] \cap g'[N]$, and the theorem is proved.
\end{proof}

\section{Tensor products}
\label{mqtensorsec}

In this section we will discuss the tensor product of $\Q$-modules. First of all, recall that~--- according to the definition of bimodules~--- a sup-lattice $\MM$ that is a left and a right module over a fixed quantale $\QQ$, is a $\QQ$-bimodule iff $(q_1 \star_l m) \star_r q_2 = q_1 \star_l (m \star_r q_2)$, for all $m \in M$ and $q_1, q_2 \in Q$.

\begin{definition}\label{tensormq}
Let $\QQ$ be a quantale and let $\MM_1 = \la M_1, \vee_1, \bot_1\ra$ be a right $\QQ$-module, $\MM_2 = \la M_2, \vee_2, \bot_2\ra$ a left $\QQ$-module and $\LL = \la L, \vee, \bot\ra$ be a sup-lattice. Then $\MM_1 \times \MM_2$ is a bimodule, where the join is defined componentwisely and the scalar multiplications are defined, for all $(x,y) \in M_1 \times M_2$ and $q \in Q$, as follows:
\begin{enumerate}
\item[]$q \star_l (x,y) = (x,q \star_2 y)$ (left multiplication),
\item[]$(x,y) \star_r q = (x \star_1 q,y)$ (right multiplication),
\end{enumerate}
$\star_i$ being the scalar multiplication in $\MM_i$, for $i = 1,2$. In what follows we shall omit subscripts if there will not be any danger of confusion.

A map $f: M_1 \times M_2 \lto L$ is said to be a \emph{$\QQ$-bimorphism}\index{Q-module@$\Q$-module!-- bimorphism}\index{Bimorphism!-- of $\Q$-modules} if it preserves arbitrary joins in each variable separately
\begin{eqnarray}
&& f\left(\bigvee_{i \in I} x_i, y\right) = \bigvee_{i \in I} f(x_i,y), \nonumber \\
&& f\left(x,\bigvee_{j \in J} y_j\right) = \bigvee_{j \in J} f(x,y_j), \nonumber \\
\end{eqnarray}
and also the following condition holds
\begin{equation}\label{bimproduct}
f(x, q \star y) = f(x \star q,y).
\end{equation}

The \emph{tensor product}\index{Q-module@$\Q$-module!-- tensor product}\index{Tensor product!-- of $\Q$-modules} $\MM_1 \tensor_\QQ \MM_2$, of the $\QQ$-modules $\MM_1$ and $\MM_2$, is the codomain of the universal $\QQ$-bimorphism $\MM_1 \times \MM_2 \lto \MM_1 \tensor_\QQ \MM_2$. In other words, we call tensor product of $\MM_1$ and $\MM_2$ a sup-lattice $\MM_1 \tensor_\QQ \MM_2$, equipped with a $\QQ$-bimorphism $\tau: \MM_1 \times \MM_2 \lto \MM_1 \tensor_\QQ \MM_2$, in such a way that, for any sup-lattice $\LL$ and any $\QQ$-bimorphism $f: \MM_1 \times \MM_2 \lto \LL$, there exists a unique sup-lattice homomorphism $k_f: \MM_1 \tensor_\QQ \MM_2 \lto \LL$ such that $k_f \circ \tau = f$.
\end{definition}

The proof of the following theorem is completely analogous to the one of Theorem~\ref{tensorslexists}. Nonetheless we believe it is useful to present it anyway.

\begin{theorem}\label{tensormqexists}
Let $\MM_1$ be a right $\QQ$-module and $\MM_2$ be a left $\QQ$-module. Then the tensor product $\MM_1 \tensor_\QQ \MM_2$ of the $\QQ$-modules $\MM_1$ and $\MM_2$ exists; it is, up to isomorphisms, the quotient $\wp(\MM_1 \times \MM_2)/\equiv_R$ of the free sup-lattice generated by $M_1 \times M_2$ with respect to the (sup-lattice) congruence relation generated by the set $R'$:
\begin{equation}\label{R'}
R' = \left\{
	\begin{array}{l}
		\left(\left\{\left(\bigvee A, y\right)\right\}, \bigcup_{a \in A}\{(a,y)\}\right) \\
		\left(\left\{\left(x, \bigvee B\right)\right\}, \bigcup_{b \in B}\{(x,b)\}\right) \\
		\left(\{(x \star q, y)\}, \{(x,q \star y)\}\right) \\
	\end{array} \right\vert
	\left.
	\begin{array}{l}
	A \subseteq M_1, y \in M_2 \\
	B \subseteq M_2, x \in M_1 \\
	q \in Q \\
	\end{array}
	\right\}.
\end{equation}
\end{theorem}
\begin{proof}
Let $\LL$ be any sup-lattice and let $f: \MM_1 \times \MM_2 \lto \LL$ be a $\QQ$-bimorphism. Since $f$ is, of course, a map, we can extend it to a sup-lattice homomorphism $h_f: \wp(\MM_1 \times \MM_2) \lto \LL$; thus $h_f \circ \sigma = f$, where $\sigma: M_1 \times M_2 \lto \wp(M_1 \times M_2)$ is the singleton map. On the other hand, the fact that $f$ is a $\QQ$-bimorphism implies $f\left(\bigvee A, y\right) = \bigvee_{a \in A} f(a,y)$, $f\left(x,\bigvee B\right) = \bigvee_{b \in B} f(x,b)$ and $f(x \star q, y) = f(x, q \star y)$, for all $x \in M_1$, $y \in M_2$, $A \subseteq M_1$, $B \subseteq M_2$ and $q \in Q$. Then, since $h_f$ is a sup-lattice homomorphism, as in Theorem~\ref{tensorslexists}, we have $h_f\left(\left\{\left(\bigvee A,y\right)\right\}\right) = h_f\left(\bigcup_{a \in A}\{(a,y)\}\right)$ and $h_f\left(\left\{\left(x,\bigvee B\right)\right\}\right) = h_f\left(\bigcup_{b \in B}\{(x,b)\}\right)$. Moreover, we have
\begin{eqnarray}
&&h_f(\{(x \star q,y)\}) = (h_f \circ \sigma)(x \star q,y) = f(x \star q,y) \nonumber \\
&=& f(x,q \star y) = (h_f \circ \sigma)(x,q \star y) = h_f(\{(x,q \star y)\}). \nonumber \\ \nonumber
\end{eqnarray}

What above means that the kernel of $h_f$ contains $R'$, thus~--- once denoted by $\mathbf{T'}$ the quotient sup-lattice $\wp(\MM_1 \times \MM_2)/\equiv_{R'}$ and by $\pi'$ the canonical projection of $\wp(\MM_1 \times \MM_2)$ over  it~--- the map
$$k_f \ : \quad [X]_{\equiv_{R'}} \ \in \ \mathbf{T'} \quad \lmapsto \quad h_f(X) \ \in \ \LL$$
is well defined and is a sup-lattice homomorphism. Moreover we have $k_f \circ \pi' \circ \sigma = h_f \circ \sigma = f$, so we have extended the $\QQ$-bimorphism $f$ to a sup-lattice homomorphism $k_f$, and it is easy to verify that the map $\tau' = \pi' \circ \sigma$ from $\MM_1 \times \MM_2$ to $\mathbf{T'}$ is indeed a bimorphism.

The following commutative diagram may clarify the constructions above.
\begin{equation}\label{tensormqdiagram}
\xymatrix{
\MM_1 \times \MM_2 \ar[rr]^\sigma \ar[rddd]_f \ar[rd]_{\tau'} && \wp(\MM_1 \times \MM_2) \ar[lddd]^{h_f} \ar[ld]^{\pi'} \\
 & \mathbf{T'}  \ar[dd]^{k_f} & \\
 &&\\
 & \LL & \\
}
\end{equation}

It is useful to remark explicitly that $R'$ and $\tau'$ do not depend either on the sup-lattice $\LL$ or on the $\QQ$-bimorphism $f$. Then we have proved that $\tau'$ is the universal bimorphism whose domain is $\MM_1 \times \MM_2$, and that $\mathbf{T'}$ is its codomain, i.e. the tensor product $\MM_1 \tensor_\QQ \MM_2$ of the $\QQ$-modules $\MM_1$ and $\MM_2$.
\end{proof}

As for the tensor product of sup-lattices, if $x \in M_1$ and $y \in M_2$, we will denote by $x \tensor y$ the image of the pair $(x,y)$ under $\tau'$, i.e. the congruence class $[\{(x,y)\}]_{\equiv_{R'}}$, and we will call it a \emph{$\QQ$-tensor}\index{Tensor!-- in $\QQ\Mod$ and $\Modd\QQ$}\index{Q-tensor@$\QQ$-tensor} or, if there will not be danger of confusion, simply a \emph{tensor}. It is clear, then, that every element of $\MM_1 \tensor_\QQ \MM_2$ is a join of tensors, so
$$\MM_1 \tensor_\QQ \MM_2 = \left\{\bigvee_{i \in I} x_i \tensor y_i \ \Big\vert \ x_i \in M_1, y_i \in M_2\right\}.$$

Let now $\QQ$ and $\RR$ be two quantales, if $\MM_1$ is an $\RR$-$\QQ$-bimodule and $\MM_2$ is a left $\QQ$-module, then the tensor product $\MM_1 \tensor_\QQ \MM_2$ naturally inherits a structure of left $\RR$-module from the one defined on $\MM_1$:
\begin{equation*}
\circledast_l: \ \left(r, \bigvee_{i \in I} x_i \tensor y_i\right) \in \RR \times \left(\MM_1 \tensor \MM_2\right) \ \lto \ \bigvee_{i \in I} (r \star x_i) \tensor y_i \in \MM_1 \tensor \MM_2.
\end{equation*}
Indeed it is trivial that $\circledast$ distributes over arbitrary joins in both coordinates; on the other hand, the external associative law comes straightforwardly from the fact that $\MM_1$ is a left $\RR$-module. Analogously, if $\MM_1$ is a right $\QQ$-module and $\MM_2$ is a $\QQ$-$\RR$-bimodule, then the tensor product $\MM_1 \tensor_\QQ \MM_2$ is a right $\RR$-module with the scalar multiplication defined, obviously, as
\begin{equation*}
\circledast_r: \ \left(\bigvee_{i \in I} x_i \tensor y_i, r\right) \in \left(\MM_1 \tensor \MM_2\right) \times \RR \ \lto \ \bigvee_{i \in I} x_i \tensor (y_i \star r) \in \MM_1 \tensor \MM_2.
\end{equation*}
Therefore, it also follows that, if $\SS$ is another quantale such that $\MM_1$ is an $\RR$-$\QQ$-bimodule and $\MM_2$ is a $\QQ$-$\SS$-bimodule, then $\MM_1 \tensor_\QQ \MM_2$ is an $\RR$-$\SS$-bimodule. In particular, if $\QQ$ is a commutative quantale, any tensor product of $\QQ$-modules is a $\QQ$-module itself.

In the case of modules over a commutative quantale, as the tensor product of sup-lattices, also the tensor product of modules is commutative, i.e. $\MM_1 \tensor_\QQ \MM_2 \cong \MM_2 \tensor_\QQ \MM_1$.

\begin{proposition}\label{tensormqcomm}
Let $\QQ$ be a commutative quantale and $\MM_1$, $\MM_2$ be $\QQ$-modules. Then the tensor products $\MM_1 \tensor_\QQ \MM_2$ and $\MM_2 \tensor \MM_1$ are isomorphic $\QQ$-modules.
\end{proposition}
\begin{proof}
It is immediate to verify that the maps
$$\alpha: (x,y) \in \MM_1 \times \MM_2 \lmapsto y \tensor x \in \MM_2 \tensor_\QQ \MM_1$$
and
$$\beta: (y,x) \in \MM_2 \times \MM_1 \lmapsto x \tensor y \in \MM_1 \tensor_\QQ \MM_2$$
are $\QQ$-bimorphism. Then there exist two homomorphisms,
$$\phi: \MM_1 \tensor_\QQ \MM_2 \lto \MM_2 \tensor_\QQ \MM_1 \quad \textrm{and} \quad \psi: \MM_2 \tensor_\QQ \MM_1 \lto \MM_1 \tensor_\QQ \MM_2,$$
extending $\alpha$ and $\beta$ respectively. Trivially, each of these homomorphisms is the inverse of the other one.
\end{proof}

All these properties of the tensor product of modules will allow us to show some relations between tensor products and hom-sets. First we need the following lemma.

\begin{lemma}\label{diam}
Let $\QQ$ and $\RR$ be quantales. If $\MM_1$ is a $\QQ$-$\RR$-bimodule and $\MM_2$ is a left $\QQ$-module, then $\Hom_\QQ(\MM_1,\MM_2)$ is a left $\RR$-module with the external product $\diamond_r$ defined, for $r \in R$, $h \in \hom_\QQ(\MM_1,\MM_2)$ and $x \in M_1$, by
\begin{equation}\label{homleftq}
(r \diamond_l h)(x) = h(x \star_\RR r),
\end{equation}
$\star_\RR$ denoting the right external product of $\MM_1$.

Analogously, if $\MM_1$ is an $\RR$-$\QQ$-bimodule and $\MM_2$ is a right $\QQ$-module, then $\Hom_\QQ(\MM_1,\MM_2)$ is a right $\RR$-module with the external product $\diamond_r$ defined, for $r \in R$, $h \in \hom_\QQ(\MM_1,\MM_2)$ and $x \in M_1$, by
\begin{equation}\label{homrightq}
(h \diamond_r r)(x) = h(r \star_\RR x),
\end{equation}
$\star_\RR$ denoting the left external product of $\MM_1$.
\end{lemma}
\begin{proof}
We will consider only the first case, the latter being completely analogous.

Given a scalar $r \in R$ and a $\QQ$-homomorphism $h$, it is immediate to verify that the map $r \diamond_l h$ sends the bottom element in the bottom element and preserves any join. The fact that $r \diamond_l h$ preserves also the right multiplication comes from the fact that $\MM_1$ is bimodule; indeed, for any $q \in Q$ and $x \in M_1$, we have
$$\begin{array}{l}
(r \diamond_l h)(q \star_\QQ x) \\
= h((q \star_\QQ x) \star_\RR r) = h(q \star_\QQ (x \star_\RR r)) = q \star_\QQ h(x \star_\RR r) \\
= q \star ((r \diamond_l h)(x)).
\end{array}$$
\end{proof}

\begin{theorem}\label{isohomtensmq}
Let $\QQ$ and $\RR$ be quantales and let $\MM_1$ be a right $\QQ$-module, $\MM_2$ a $\QQ$-$\RR$-bimodule and $\MM_3$ a right $\RR$-module. Then, if we consider the right $\RR$-module $\MM_1 \tensor_\QQ \MM_2$ and the right $\QQ$-module $\Hom_\RR(\MM_2,\MM_3)$, we have
$$\Hom_\RR(\MM_1 \tensor_\QQ \MM_2, \MM_3) \cong_\SL \Hom_\QQ(\MM_1,\Hom_\RR(\MM_2,\MM_3)),$$
where $\cong_\SL$ means they are isomorphic as sup-lattices.
\end{theorem}
\begin{proof}
If $h$ is an $\RR$-homomorphism from $\MM_1 \tensor_\QQ \MM_2$ to $\MM_3$, then $h(x \tensor y)$ is clearly an element of $\MM_3$, for every tensor $x \tensor y$. Thus, fixed $x \in \MM_1$, $h$ defines a map
$$h_x: y \in M_2 \lmapsto h(x \tensor y) \in M_3.$$
Since $h$ is an $\RR$-homomorphism, given a family $\{y_i\}_{i \in I} \subseteq M_2$, an element $y \in M_2$ and a scalar $r \in R$, we have
$$\begin{array}{l}
h_x\left(\bigvee_{i \in I}y_i\right) \\
= h\left(x \tensor \bigvee_{i \in I} y_i\right) \\
= h\left(\bigvee_{i \in I}(x \tensor y_i)\right) \\
= \bigvee_{i \in I} h(x \tensor y_i) \\
= \bigvee_{i \in I} h_x(y_i)
\end{array}$$
and
$$h_x(y \star r) = h(x \tensor (y \star r)) = h((x \tensor y) \star r) = h(x \tensor y) \star r = h_x(y) \star r,$$
so $h_x \in \hom_\RR(\MM_2,\MM_3)$, for any fixed $x$. Hence we have a map $h_-: x \in M_1 \lmapsto h_x \in \hom_\RR(\MM_2,\MM_3)$, but the fact that $h$ is also a sup-lattice homomorphism implies $h_-$ to be a sup-lattice homomorphism as well:
$$\begin{array}{l}
h_{\bigvee_{i \in I}x_i}(y) \\
= h\left(\left(\bigvee_{i \in I} x_i\right) \tensor y\right) \\
= h\left(\bigvee_{i \in I} (x_i \tensor y)\right) \\
= \bigvee_{i \in I} h(x_i \tensor y) \\
= \bigsqcup_{i \in I} h_{x_i}(y)
\end{array},$$
for all $\{x_i\}_{i \in I} \subseteq M_1$ and $y \in M_2$. Moreover, if $q \in Q$, by (\ref{R'}) and (\ref{homrightq}),
$$h_{x \star q}(y) = h((x \star q) \tensor y) = h(x \tensor (q \star y)) = h_x(q \star y) = (h_x \diamond_r q)(y),$$
for all $x \in M_1$, $y \in M_2$, $q \in Q$, so $h_-$ is a $\QQ$-homomorphism. Besides, we also have
$$\left(\bigsqcup_{i \in I}h_i\right)(x \tensor y) = \bigvee_{i \in I}h_i(x \tensor y) = \bigvee_{i \in I}{h_i}_x(y) = \left(\bigsqcup_{i \in I}{h_i}_x\right)(y),$$
for any family $\{h_i\} \subseteq \hom_\RR(\MM_1 \tensor_\QQ \MM_2, \MM_3)$, for all $x \in M_1$, for all $y \in M_2$.

Therefore we have a sup-lattice homomorphism
\begin{equation}\label{zetamq}
\zeta: \ \Hom_\RR(\MM_1 \tensor_\QQ \MM_2, \MM_3) \ \lto \ \Hom_\QQ(\MM_1, \Hom_\RR(\MM_2,\MM_3)),
\end{equation}
defined by $\zeta(h) = h_-$, i.e. $((\zeta(h))(x))(y) = h(x \tensor y)$.

Let us show that $\zeta$ has an inverse. If $f \in \hom_\QQ(\MM_1,\hom_\RR(\MM_2,\MM_3))$, then the map $f': (x,y) \in \MM_1 \times \MM_2 \lmapsto (f(x))(y) \in \MM_3$ is clearly a $\QQ$-bimorphism. Hence there exists a unique homomorphism $h_{f'}: \MM_1 \tensor_\QQ \MM_2 \lto \MM_3$ such that $h_{f'} \circ \tau' = f'$, i.e. such that $h_{f'}(x \tensor y) = f'(x,y) = (f(x))(y)$, for all $x \in M_1$ and $y \in M_2$, and~--- clearly~--- $\zeta(h_{f'}) = f$. On the other hand, if $f = \zeta(h)$ with $f \in \Hom_\QQ(\MM_1,\hom_\RR(\MM_2,\MM_3)$ and $h \in \Hom_\RR(\MM_1 \tensor_\QQ \MM_2, \MM_3)$, then the uniqueness of the homomorphism that extends the map $f'$, defined above, to $\MM_1 \tensor_\QQ \MM_2$ ensures us that $h_{f'} = h$. Then we have the inverse sup-lattice homomorphism
$$\zeta^{-1}: \ f \in \Hom_\QQ(\MM_1,\Hom_\RR(\MM_2,\MM_3)) \ \lmapsto \ h_{f'} \in \Hom_\RR(\MM_1 \tensor_\QQ \MM_2, \MM_3),$$
and the theorem is proved.
\end{proof}

With an analogous proof, we have
\begin{theorem}\label{isohomtensmq'}
Let $\QQ$ and $\RR$ be quantales and let $\MM_1$ be a $\RR$-$\QQ$-bimodule, $\MM_2$ a left $\QQ$-module and $\MM_3$ a left $\RR$-module. Then, if we consider the left $\RR$-module $\MM_1 \tensor_\QQ \MM_2$ and the left $\QQ$-module $\Hom_\RR(\MM_1,\MM_3)$, we have
$$\Hom_\RR(\MM_1 \tensor_\QQ \MM_2, \MM_3) \cong_\SL \Hom_\QQ(\MM_2,\Hom_\RR(\MM_1,\MM_3)),$$
where $\cong_\SL$ means they are isomorphic as sup-lattices.
\end{theorem}
\begin{proof}
If $h$ is an $\RR$-homomorphism from $\MM_1 \tensor_\QQ \MM_2$ to $\MM_3$, then $h(x \tensor y)$ is clearly an element of $\MM_3$, for every tensor $x \tensor y$. Thus, fixed $y \in \MM_2$, $h$ defines a map
$$h_y: x \in M_1 \lmapsto h(x \tensor y) \in M_3.$$
Since $h$ is an $\RR$-homomorphism, given a family $\{x_i\}_{i \in I} \subseteq M_1$, an element $x \in M_1$ and a scalar $r \in R$, we have
$$\begin{array}{l}
h_y\left(\bigvee_{i \in I}x_i\right) \\
= h\left(\left(\bigvee_{i \in I} x_i\right) \tensor y \right) \\
= h\left(\bigvee_{i \in I}(x_i \tensor y)\right) \\
= \bigvee_{i \in I} h(x_i \tensor y) \\
= \bigvee_{i \in I} h_y(x_i)
\end{array}$$
and
$$h_y(r \star x) = h((r \star x) \tensor y) = h(r \star (x \tensor y)) = r \star h(x \tensor y) = r \star h_y(x),$$
so $h_y \in \hom_\RR(\MM_1,\MM_3)$, for any fixed $y$. Hence we have a map $h_-: y \in M_2 \lmapsto h_y \in \hom_\RR(\MM_1,\MM_3)$, but the fact that $h$ is also a sup-lattice homomorphism implies $h_-$ to be a sup-lattice homomorphism as well:
$$\begin{array}{l}
h_{\bigvee_{i \in I}y_i}(x) \\
= h\left(x \tensor \left(\bigvee_{i \in I} y_i\right)\right) \\
= h\left(\bigvee_{i \in I} (x \tensor y_i)\right) \\
= \bigvee_{i \in I} h(x \tensor y_i) \\
= \bigsqcup_{i \in I} h_{y_i}(x)
\end{array},$$
for all $\{y_i\}_{i \in I} \subseteq M_2$ and $x \in M_1$. Moreover, if $q \in Q$, by (\ref{R'}) and (\ref{homleftq}),
$$h_{q \star y}(x) = h(x \tensor (q \star y)) = h((x \star q) \tensor y) = h_y(x \star q) = (q \diamond_l h_y)(x),$$
for all $x \in M_1$, $y \in M_2$, $q \in Q$, so $h_-$ is a $\QQ$-homomorphism. Besides, we also have
$$\left(\bigsqcup_{i \in I}h_i\right)(x \tensor y) = \bigvee_{i \in I}h_i(x \tensor y) = \bigvee_{i \in I}{h_i}_y(x) = \left(\bigsqcup_{i \in I}{h_i}_y\right)(x),$$
for any family $\{h_i\} \subseteq \hom_\RR(\MM_1 \tensor_\QQ \MM_2, \MM_3)$, for all $x \in M_1$, for all $y \in M_2$.

Therefore we have a sup-lattice homomorphism
\begin{equation}\label{zeta'mq}
\zeta': \ \Hom_\RR(\MM_1 \tensor_\QQ \MM_2, \MM_3) \ \lto \ \Hom_\QQ(\MM_2, \Hom_\RR(\MM_1,\MM_3)),
\end{equation}
defined by $\zeta'(h) = h_-$, i.e. $((\zeta'(h))(x))(y) = h(x \tensor y)$.

Let us show that $\zeta'$ has an inverse. If $f \in \hom_\QQ(\MM_2,\hom_\RR(\MM_1,\MM_3))$, then the map $f': (x,y) \in \MM_1 \times \MM_2 \lmapsto (f(y))(x) \in \MM_3$ is clearly a $\QQ$-bimorphism. Hence there exists a unique homomorphism $h_{f'}: \MM_1 \tensor_\QQ \MM_2 \lto \MM_3$ such that $h_{f'} \circ \tau' = f'$, i.e. such that $h_{f'}(x \tensor y) = f'(x,y) = (f(y))(x)$, for all $x \in M_1$ and $y \in M_2$, and~--- clearly~--- $\zeta'(h_{f'}) = f$. On the other hand, if $f = \zeta'(h)$ with $f \in \Hom_\QQ(\MM_2,\hom_\RR(\MM_1,\MM_3)$ and $h \in \Hom_\RR(\MM_1 \tensor_\QQ \MM_2, \MM_3)$, then the uniqueness of the homomorphism that extends the map $f'$, defined above, to $\MM_1 \tensor_\QQ \MM_2$ ensures us that $h_{f'} = h$. Then we have the inverse sup-lattice homomorphism
$$\zeta'^{-1}: \ f \in \Hom_\QQ(\MM_2,\Hom_\RR(\MM_1,\MM_3)) \ \lmapsto \ h_{f'} \in \Hom_\RR(\MM_1 \tensor_\QQ \MM_2, \MM_3),$$
and the theorem is proved.
\end{proof}

\begin{corollary}\label{isohomtensmqcor}
Let $\QQ$ be a commutative quantale and let $\MM_1$, $\MM_2$ and $\MM_3$ be $\QQ$-modules. Then
$$\begin{array}{l}
\Hom_\QQ(\MM_1 \tensor_\QQ \MM_2, \MM_3) \\
\cong \Hom_\QQ(\MM_1,\Hom_\QQ(\MM_2,\MM_3)) \\
\cong \Hom_\QQ(\MM_2,\Hom_\QQ(\MM_1,\MM_3)).
\end{array}$$
\end{corollary}
\begin{proof}
It is an easy application of Theorems~\ref{isohomtensmq} and~\ref{isohomtensmq'} or, equivalently, of Theorem~\ref{isohomtensmq} and Proposition~\ref{tensormqcomm}.
\end{proof}

\begin{lemma}\label{homqm}
Let $\QQ$ be a quantale and $\MM$ be a $\QQ$-module. Then, considered $\QQ = \la Q, \vee, \bot \ra$ as a module over itself, we have
\begin{equation*}
\Hom_\QQ(\QQ,\MM) \cong_\SL \MM.
\end{equation*}
\end{lemma}
\begin{proof}
First of all we observe that, for any fixed $m \in M$, the map $f_m: q \in Q \lto q \star m \in M$ is trivially a $\QQ$-module homomorphism. Then we can consider the map $\alpha: m \in M \lto f_m \in \hom_\QQ(\QQ,\MM)$, which is clearly a sup-lattice homomorphism.

Let us consider also the map $\beta: f \in \hom_\QQ(\QQ,\MM) \lto f(e) \in M$. Again, it is immediate to verify that $\beta$ is a sup-lattice homomorphism. But we also have:
$$((\alpha \circ \beta)(f))(q) = (\alpha(f(e)))(q) = f_{f(e)}(q) = q \star f(e) = f(q),$$
for all $f \in \hom_\QQ(\QQ,\MM)$ and $q \in Q$, and
$$(\beta \circ \alpha)(m) = \beta(f_m) = f_m(e) = e \star m = m,$$
for all $m \in M$.

Thus $\alpha \circ \beta = \id_{\hom_\QQ(\QQ,\MM)}$ and $\beta \circ \alpha = \id_M$, i.e. $\alpha$ is an isomorphism whose inverse is $\beta$, and the thesis follows.
\end{proof}

As a consequence of the previous result, the $\QQ$-module structure defined on $\hom_\QQ(\QQ,\MM)$ by Lemma \ref{diam} is isomorphic to $\MM$.

\begin{lemma}\label{hommqdual}
Let $\QQ$ be a commutative quantale and $\MM$ be a $\QQ$-module. Then, once considered the $\QQ$-module $\QQ^{\op} = \la Q, \wedge, \top \ra$, we have
\begin{equation*}
\Hom_\QQ(\MM,\QQ^{\op}) \cong \MM^{\op}.
\end{equation*}
\end{lemma}
\begin{proof}
By Lemma~\ref{homqm}, $\Hom_\QQ(\QQ,\MM^{\op}) \cong \MM^{\op}$. But, by Proposition~\ref{dualhomomq}, if $\QQ$ is a commutative quantale, $\Hom_\QQ(\MM,\NN)$ and $\Hom_\QQ(\NN^{\op},\MM^{\op})$ are isomorphic $\QQ$-modules, for any pair of modules $\MM$ and $\NN$. Then
$$\Hom_\QQ(\MM,\QQ^{\op}) \cong \Hom_\QQ(\QQ,\MM^{\op}) \cong \MM^{\op}.$$
\end{proof}

\begin{theorem}\label{isohomtensmq2}
Let $\QQ$ be a commutative quantale and $\MM_1$ and $\MM_2$ be $\QQ$-modules. Then the following identities hold:
\begin{enumerate}
\item[$(i)$]$\Hom_\QQ(\MM_1,\MM_2) \cong (\MM_1 \tensor_\QQ \MM_2^{\op})^{\op}$,
\item[$(ii)$]$\MM_1 \tensor_\QQ \MM_2 \cong \Hom_\QQ\left(\MM_1, \MM_2^{\op}\right)^{\op}$.
\end{enumerate}
\end{theorem}
\begin{proof}
\begin{enumerate}
\item[$(i)$]\begin{equation*}
\begin{array}{lll}
&\quad & \Hom_\QQ(\MM_1,\MM_2) \\
\textrm{by Lemma~\ref{hommqdual}} && \cong \Hom_\QQ(\MM_1,\Hom_\QQ(\MM_2^{\op},\QQ^{\op})) \\
\textrm{by Corollary~\ref{isohomtensmqcor}} && \cong \Hom_\QQ(\MM_1 \tensor_\QQ \MM_2^{\op},\QQ^{\op}) \\
\textrm{again by Lemma~\ref{hommqdual}} && \cong (\MM_1 \tensor_\QQ \MM_2^{\op})^{\op}. \\
\end{array}
\end{equation*}
\item[$(ii)$]\begin{equation*}
\begin{array}{lll}
&\quad & \Hom_\QQ(\MM_1, \MM_2^{\op})^{\op} \\
\textrm{by Lemma~\ref{hommqdual}} && \cong \Hom_\QQ(\MM_1, \Hom_\QQ(\MM_2,\QQ^{\op}))^{\op} \\
\textrm{by Corollary~\ref{isohomtensmqcor}} && \cong \Hom_\QQ(\MM_1 \tensor_\QQ \MM_2, \QQ^{\op})^{\op} \\
\textrm{again by Lemma~\ref{hommqdual}} && \cong ((\MM_1 \tensor_\QQ \MM_2)^{\op})^{\op} \\
 && \cong \MM_1 \tensor_\QQ \MM_2.
\end{array}
\end{equation*}
\end{enumerate}
\end{proof}

\begin{theorem}\label{producttensormq}
Let $\QQ$ be a quantale, $\MM$ be a right $\QQ$-module and $\{\MM_i\}_{i \in I}$ be a family of left $\QQ$-modules. Then
$$\coprod_{i \in I} \MM \tensor_\QQ \MM_i \cong_\SL \MM \tensor_\QQ \coprod_{i \in I} \MM_i \quad \textrm{and} \quad \prod_{i \in I} \MM \tensor_\QQ \MM_i \cong_\SL \MM \tensor_\QQ \prod_{i \in I} \MM_i.$$
If $\MM$ is an $\RR$-$\QQ$-bimodule (respectively, if $\QQ$ is commutative), the isomorphisms are $\RR$-module (resp., $\QQ$-module) isomorphisms.
\end{theorem}
\begin{proof}
Let $\MM'$ be any $\QQ$-module and, for all $i \in I$, $f_i: \MM \tensor_\QQ \MM_i \lto \MM'$  and $g_i: \MM' \lto \MM \tensor_\QQ \MM_i$ be homomorphisms. As we remarked when we defined the $\QQ$-tensors, any element of $\MM \tensor_\QQ \prod_{i \in I} \MM_i$~--- and then, by Proposition~\ref{coprodmq}, any element of $\MM \tensor_\QQ \coprod_{i \in I} \MM_i$~--- is a join of tensors, i.e. can be expressed in the form $\bigvee_{j,k} x_j \tensor (y_{ik})_{i \in I}$, with $\{x_j\}_{j \in J} \subseteq L$ and $\{(y_{ik})_{i \in I}\}_{k \in K} \subseteq \prod_{i \in I} L_i$. Now, for any fixed $\ov i \in I$, let us consider the following homomorphisms: 
$$\iota_\MM \tensor \mu_{\ov i}: \bigvee_{j,k} x_j \tensor y_{\ov i k} \in \MM \tensor_\QQ \MM_{\ov i} \ \lmapsto \ \bigvee_{j,k} x_j \tensor \mu_{\ov i}(y_{\ov i k}) \in \MM \tensor_\QQ \prod_{i \in I} \MM_i,$$
$$\iota_\MM \tensor \pi_{\ov i}: \bigvee_{j,k} x_j \tensor \left(y_{ik}\right)_{i \in I} \in \MM \tensor_\QQ \prod_{i \in I} \MM_i \ \lmapsto \ \bigvee_{j,k} x_j \tensor \pi_{\ov i}\left((y_{ik})_{i \in I}\right) \in \MM \tensor_\QQ \MM_{\ov i};$$
we observe that, for all $\ov i \in I$, $(\iota_\MM \tensor \pi_{\ov i}) \circ (\iota_\MM \tensor \mu_{\ov i}) = \iota_{\MM \tensor_\QQ \MM_{\ov i}}$ and, if $\ov{i'} \neq \ov i$, $(\iota_\MM \tensor \pi_{\ov{i'}}) \circ (\iota_\MM \tensor \mu_{\ov i})$ is the map $\iota_\MM \tensor \bot_{\ov{i'}}$ that sends the element $\bigvee_{j,k} x_j \tensor y_{\ov i k}$ to $\bigvee_j x_j \tensor \bot$.

What we need to prove is the existence~--- and uniqueness~--- of two homomorphisms $f: \MM \tensor_\QQ \prod_{i \in I} \MM_i \lto \MM'$ and $g: \MM' \lto \MM \tensor_\QQ \prod_{i \in I} \MM_i$ such that the following diagrams are commutative.
\begin{eqnarray}\label{prodtensorsmq1}
\xymatrix{
\MM \tensor_\QQ \MM_i \ar[rr]^{\iota_\MM \tensor \mu_i} \ar[dd]^{f_i} && \MM \tensor_\QQ \coprod_{i \in I} \MM_i \ar[ddll]^f\\
&&\\
\MM'&&\\
}
\\
\label{prodtensorsmq2}
\xymatrix{
&& \MM' \ar[ddll]_g \ar[dd]_{g_i} \\
&&\\
\MM \tensor_\QQ \prod_{i \in I} \MM_i \ar[rr]^{\iota_\MM \tensor \pi_i} && \MM \tensor_\QQ \MM_i \\
}
\end{eqnarray}

Since it is clear that any element $\bigvee_{j,k} x_j \tensor (y_{ik})_{i \in I}$ of $\MM \tensor_\QQ \coprod_{i \in I} \MM_i$ can be written as $\bigvee_{\ov i \in I} \bigvee_{j,k} x_j \tensor \mu_{\ov i}\left((y_{ik})_{i \in I}\right)$, the diagram in (\ref{prodtensorsmq1}) can be easily made commutative by setting, for all $\bigvee_{j,k} x_j \tensor (y_{ik})_{i \in I} \in \MM \tensor_\QQ \coprod_{i \in I} \MM_i$,
$$f\left(\bigvee_{j,k} x_j \tensor (y_{ik})_{i \in I}\right) = \bigvee_{\ov i \in I} \bigvee_{j,k} f_{\ov i}\left(x_j \tensor \mu_{\ov i}\left((y_{ik})_{i \in I}\right)\right).$$
Regarding the diagram in (\ref{prodtensorsmq2}), we define, for all $z \in \MM'$,
$$g(z) = \bigvee_{i \in I} ((\iota_\MM \tensor \mu_i) \circ g_i)(z).$$
The fact that $f$ and $g$ are homomorphisms is easily seen since the tensor product preserves joins in both coordinates and all the maps involved in their definition are homomorphisms. Then, for any fixed $\ov i \in I$ and for any $z \in L'$, we have
$$\begin{array}{l}
((\iota_\MM \tensor \pi_{\ov i}) \circ g)(z) \\
= (\iota_\MM \tensor \pi_{\ov i})\left(\bigvee_{i \in I} ((\iota_\MM \tensor \mu_i) \circ g_i)(z)\right) \\
= \bigvee_{i \in I} (\iota_\MM \tensor \pi_{\ov i})\left(((\iota_\MM \tensor \mu_i) \circ g_i)(z)\right) \\
= \bigvee_{i \in I} \left(((\iota_\MM \tensor \pi_{\ov i}) \circ (\iota_\MM \tensor \mu_i) \circ g_i)(z)\right) \\
= g_{\ov i}(z) \vee \bigvee_{i \in I} \iota_\MM \tensor \bot_i (g_i(z)) \\
= g_{\ov i}(z).
\end{array}$$
The proof of the fact that such $f$ and $g$ are unique is straightforward.

Then $\MM \tensor_\QQ \coprod_{i \in I} \MM_i$ (respectively: $\MM \tensor_\QQ \prod_{i \in I} \MM_i$) has the universal property of extending sinks (resp.: sources) whose domain (resp.: codomain) is the family $\{\MM \tensor_\QQ \MM_i\}_{i \in I}$.

Now assume that $\MM$ is an $\RR$-$\QQ$-bimodule, $\MM'$ is an $\RR$-module and that the $f_i$'s and the $g_i$'s are $\RR$-module homomorphisms. Then the maps $\iota_\MM \tensor \mu_{\ov i}$'s and $\iota_\MM \tensor \pi_{\ov i}$'s are easily seen to be $\RR$-module homomorphisms and, consequently, the same holds for $f$ and $g$, for how such maps have been defined. Therefore the isomorphisms proved are isomorphisms of $\RR$-modules and the commutative case follows as a consequence. The theorem is proved.
\end{proof}

\begin{corollary}
Let $\QQ$ be a commutative quantale, $\MM$ be a $\QQ$-module, and $X$ and $Y$ be non-empty sets. Then
\begin{enumerate}
\item[$(i)$]$\QQ \tensor_\QQ \MM \cong \MM$;
\item[$(ii)$]$\QQ^X \tensor_\QQ \MM \cong \MM^X$;
\item[$(iii)$]$\QQ^X \tensor_\QQ \QQ^Y \cong \QQ^{X \times Y}$.
\end{enumerate}
\end{corollary}
\begin{proof}
\begin{enumerate}
\item[$(i)$]\begin{equation*}
\begin{array}{lll}
& \quad & \QQ \tensor_\QQ \MM \\
\textrm{by Theorem~\ref{tensormqcomm}} && \cong \MM \tensor_\QQ \QQ \\
\textrm{by Theorem~\ref{isohomtensmq2}} && \cong \Hom_\QQ(\MM,\QQ^{\op})^{\op} \\
\textrm{by Lemma~\ref{homqm}} && \cong \MM.
\end{array}
\end{equation*}
\item[$(ii)$]Let us denote by $\{\QQ_x\}_{x \in X}$ and $\{\MM_x\}_{x \in X}$ two families of copies, of $\QQ$ and $\MM$ respectively, with set of indices $X$. Then $\QQ^X \cong \prod_{x \in X} \QQ_x$ and $\MM^X \cong \prod_{x \in X} \MM_x$. We have:
\begin{equation*}
\begin{array}{lll}
& \quad & \QQ^X \tensor_\QQ \MM \\
&& \cong \left(\prod_{x \in X} \QQ_x\right) \tensor_\QQ \MM \\
\textrm{by Theorems~\ref{tensormqcomm} and~\ref{producttensormq}} && \cong \prod_{x \in X} \QQ_x \tensor_\QQ \MM \\
\textrm{by $(i)$} && \cong \prod_{x \in X} \MM_x \\
&& \cong \MM^X.
\end{array}
\end{equation*}
\item[$(iii)$]\begin{equation*}
\begin{array}{lll}
& \quad & \QQ^X \tensor_\QQ \QQ^Y \\
\textrm{by $(ii)$} && \cong \left(\QQ^Y\right)^X \\
&& \cong \QQ^{X \times Y}.
\end{array}
\end{equation*}
\end{enumerate}
\end{proof}

Let now $\RR$ be a quantale and $\QQ$ be a subquantale of $\RR$. If $\MM$ is a left $\QQ$-module, we can use the tensor product in order to extend, in a precise sense, the $\QQ$-module $\MM$ to an $\RR$-module. Indeed, if we consider $\RR$ as an $\RR$-$\QQ$-bimodule, the $\QQ$-tensor product $\RR \tensor_\QQ \MM$ is a left $\RR$-module (hence, also a left $\QQ$-module).

Let $m \in M$; for any $q \in Q$,
$$q \circledast (e \tensor m) = q \tensor m = (e \cdot q) \tensor m = e \tensor (q \star m).$$
So the set $e \tensor_\QQ \MM = \{e \tensor m \mid m \in M\}$~--- that clearly generates $\RR \tensor_\QQ \MM$ as $\RR$-module~--- is a $\QQ$-submodule of $\RR \tensor_\QQ \MM$, homomorphic image of $\MM$. Indeed the map
\begin{equation*}
\begin{array}{llll}
e \tensor \iota_\MM : \quad & \MM \quad & \lto \quad 	& \RR \tensor_\QQ \MM \\
										 				&	m  				& \lmapsto		& e \tensor m \\
\end{array}
\end{equation*}
is a $\QQ$-module homomorphism.

If $\MM = \QQ^X$ is a free module, the tensor product is isomorphic to the free $\RR$-module over the same basis: $\RR \tensor_\QQ \QQ^X \cong \RR^X$. Indeed the map $\phi: (r,f) \in \RR \times \QQ^X \lmapsto (r \cdot f(x))_{x \in X} \in \RR^X$ is clearly a $\QQ$-bimorphism, and the homomorphism that extends $\phi$ to $\RR \tensor_\QQ \QQ^X$ is
$$k_\phi: \bigvee_{i \in I} r_i \tensor f_i \in \RR \tensor_\QQ \QQ^X \lmapsto \bigvee_{i \in I} \left(r_i \cdot f_i(x)\right)_{x \in X} \in \RR^X.$$
Then it is easy to verify that $k': g \in \RR^X \lmapsto \bigvee_{x \in X} g(x) \tensor \chi_x \in \RR \tensor_\QQ \QQ^X$ is a homomorphism and it is the inverse of $k_\phi$. Every element of $\RR \tensor_\QQ \QQ^X$ can be written in a unique way as $\bigvee_{x \in X} r_x (e \tensor \chi_x)$, i.e. $\RR \tensor_\QQ \QQ^X$ is the free $\RR$-module generated by the set $\{e \tensor \chi_x \mid x \in X\}$, equipotent to $X$.

In general, if $\MM$ is a left $\QQ$-module, $X$ is a set of generators for $\MM$ and $\RR$ is a quantale containing $\QQ$, then the left $\RR$-module $\RR \tensor_\QQ \MM$ is generated by $e \tensor_\QQ X = \{e \tensor x \mid x \in X\}$. As a consequence we have that, if $\MM$ is cyclic, generated by $v$, then $\RR \tensor_\QQ \MM$ is cyclic as well, and it is generated by $e \tensor v$. But we have more:
\begin{theorem}\label{projtensor}
If $\MM$ is a projective cyclic $\QQ$-module and $\QQ$ is a subquantale of a quantale $\RR$, then $\RR \tensor_\QQ \MM$ is a projective cyclic $\RR$-module. Consequently, if $\MM$ is the coproduct of projective cyclic modules, then $\RR \tensor_\QQ \MM$ is the coproduct of projective cyclic modules as well, thus projective itself.
\end{theorem}
\begin{proof}
Let $\MM = \la v \ra$ be cyclic and projective. We have already observed that $\RR \tensor_\QQ \MM$ is cyclic and it is generated by $e \tensor v$; by Corollary~\ref{cycliccond} and Theorem~\ref{cyclicproj}, there exists an idempotent element $u \in \QQ$ such that $\g_v = \g_u$. On the other hand, $u$ is also an idempotent element of $\RR$, and we have
$$\RR \tensor_\QQ \MM \cong \RR \tensor_\QQ (\QQ \star v) \cong \RR \tensor_\QQ (\QQ \cdot u) \cong \RR \cdot u \cong \RR_{\g_u},$$
where $\g_u: r \in R \lmapsto ru \in R$, then $\g_{e \tensor v} = \g_u$.

So it is satisfied condition $(c)$ of Corollary~\ref{cor8.5}, hence $\RR \tensor_\QQ \MM$ is projective by Theorem~\ref{cyclicproj}.

The second assertion is a consequence of the first one, together with Theorems~\ref{projlem3}~and~\ref{producttensormq}. The theorem is proved.
\end{proof}

\begin{corollary}\label{cyclictensoriso}
Let $\RR$ be a quantale, $\QQ$ and $\QQ'$ be two subquantales of $\RR$, and $u \in \RR$ be an idempotent element contained in both $\QQ$ and $\QQ'$. Then the (cyclic projective) $\RR$-modules $\RR \tensor_\QQ (\QQ \cdot u)$ and $\RR \tensor_{\QQ'} (\QQ' \cdot u)$ are isomorphic.
\end{corollary}
\begin{proof}
By Theorem~\ref{projtensor}, $\RR \tensor_\QQ (\QQ \cdot u)$ and $\RR \tensor_{\QQ'} (\QQ' \cdot u)$ are both isomorphic to $\RR \cdot u$.
\end{proof}

\section{Restriction of scalars}
\label{changesec}

\begin{lemma}\label{indmod}
Let $\QQ$ and $\RR$ be quantales and $h: \QQ \lto \RR$ a quantale homomorphism. Then $h$ induces a structure of $\QQ$-module on any $\RR$-module.

In particular, $h$ induces structures of $\QQ$-bimodule, $\RR$-$\QQ$-bimodule and $\QQ$-$\RR$-bimodule on $\RR$ itself.
\end{lemma}
\begin{proof}
Let $\NN = \la N, \vee, \bot \ra$ be an $\RR$-module with scalar multiplication $\star$. It is easy to verify that
\begin{equation}\label{starh}
\star_h: (q, x) \in Q \times N \lmapsto h(q) \star x \in N
\end{equation}
makes $\NN$ into a $\QQ$-module, henceforth denoted by $\NN_h$. Since $\RR$ is a bimodule over itself, the second part of the assertion follows immediately.
\end{proof}

The operation performed in (\ref{starh}) is well-known in the theory of ring modules as \emph{restricting the scalars along $h$}. In fact it defines a functor
\begin{equation}\label{subh}
\begin{array}{cccc}
H: & \RR\Mod & \lto & \QQ\Mod \\
& \NN & \lmapsto & \NN_h
\end{array}
\end{equation}
having both a right and a left adjoint. This property was already pointed out by Joyal and Tierney \cite{joyaltierney} in the commutative case. In the general case, however, the situation is precisely the same, as shown by the following result.

\begin{theorem}\label{adjfunct}
The functor $H$ defined in (\ref{subh}) has both a left adjoint $H_l$ and a right adjoint $H_r$.
\end{theorem}
\begin{proof}
For any $\MM \in \QQ\Mod$, viewing $\RR$ as a $\RR$-$\QQ$-bimodule, we can construct the tensor product $\RR \tensor_{\QQ} \MM$ which is a left $\RR$-module. Hence
\begin{equation}\label{hl}
\begin{array}{cccc}
H_l: & \QQ\Mod & \lto & \RR\Mod \\
			& \MM			& \lmapsto & \RR \tensor_{\QQ} \MM
\end{array}
\end{equation}
is the left adjoint of $H$. In order to prove that, we need to show that, for any $\QQ$-module $\MM$ and any $\RR$-module $\NN$, there exists a natural bijection between $\hom_\RR(\RR \tensor_{\QQ} \MM, \NN)$ and $\hom_\QQ(\MM,\NN_h)$. The first hom-set is isomorphic, as a sup lattice, to $\hom_\QQ(\MM, \hom_\RR(\RR,\NN))$, by Theorem~\ref{isohomtensmq'}; on the other hand, by Lemma~\ref{homqm}, $\hom_\RR(\RR,\NN) \cong_{\SL} \NN$ and such an isomorphism is a $\QQ$-module isomorphism (with $\NN_h$ instead of $\NN$) for how the $\QQ$-module structure is induced on $\hom_\RR(\RR,\NN)$. Hence the two hom-sets are isomorphic sup-lattices, and $H_l$ is the left adjoint of $H$. 

The right adjoint is defined by
\begin{equation}\label{hr}
\begin{array}{cccc}
H_r: & \QQ\Mod & \lto & \RR\Mod \\
			 & \MM & \lmapsto & \Hom_\QQ(\RR_h,\MM),
\end{array}
\end{equation}
where the left $\RR$-module structure on $\Hom_\QQ(\RR_h,\MM)$ is the one introduced in Lemma~\ref{diam}. This part of the proof is analogous to the case of $H_l$. Indeed, for any $\QQ$-module $\MM$ and any $\RR$-module $\NN$, by Theorem~\ref{isohomtensmq'}, $\hom_\RR(\NN,H_r(\MM))$ --- namely $\hom_\RR(\NN, \Hom_\QQ(\RR_h,\MM))$ --- is isomorphic, as a sup-lattice, to $\hom_\QQ((\RR \tensor_\RR \NN)_h, \MM)$; on the other hand, since every tensor $r \tensor n \in \RR \tensor_\RR \NN$ can be rewritten in the form $e \tensor r \star n$, such a tensor product is easily seen to be isomorphic to $\NN_h$. Therefore $\hom_\RR(\NN,H_r(\MM))$ is a sup-lattice isomorphic to $\hom_\QQ(\NN_h, \MM)$ and the theorem is proved.
\end{proof}

\section*{References and further readings}
\addcontentsline{toc}{section}{References and further readings}

The notion of quantale module is relatively recent. Even if they appear~--- though ``in disguise''--- in the paper~\cite{joyaltierney} by A. Joyal and M. Tierney, a real consideration of these structures has begun only a few years ago, and the number of mathemticians that work on this topic is still limited.

As far as we know, this thesis represents a first systematic theoretic treatment of the categories of $\Q$-modules while, on the contrary, the contents of Chapters~\ref{consrelchap},~\ref{imagechap} and~\ref{ltbchap} are just a part of the applications of $\Q$-modules that are already present in literature. To what extent such applications, they appeared in~\cite{dinolarusso,dinolarusso2}, written with A. Di Nola, besides the work~\cite{galatostsinakis} already cited many times.

Further readings on $\Q$-modules, or on subjects somehow connected with them, that can be suggested are the papers~\cite{baltag} by A. Baltag, B. Coecke and M. Sadrzadeh,~\cite{kruml} by D. Kruml,~\cite{paseka} by J. Paseka,~\cite{resende5} by P. Resende, and~\cite{stubbe} by I. Stubbe.

\chapter{Deductive Systems on Quantale Modules}
\label{consrelchap}

In this chapter we will see that most of the algebraic and categorical results of the previous chapters have reflections on Mathematical Logic. First of all, we will recall several results obtained by N. Galatos and C. Tsinakis, and presented in~\cite{galatostsinakis}, suitably rewritten in the language of quantale modules, that establish once and for all that deductive systems can be represented as special quantale modules, and consequence relations are essentially $\Q$-module nuclei.

After that, in Section~\ref{transsec}, we will return to the concrete situations of deductive systems on propositional languages. Here we will clarify (at least for the case of propositional languages) some concepts regarding the comparison of deductive systems that, although very common in the literature of Mathematical Logic, are still rather vague, and whose meaning usually relies on logicians' intuition. Once these notions have been clarified, we will generalize them in such a way that the role of the lattices of theories will become leading, to what extent this issue, with respect to that of languages.

Such generalizations are amenable to an algebraic characterization that, under additional hypotheses, yield results~--- that are mainly of algebraic nature as well--- on the corresponding classical notions.

\section{Consequence relations on sup-lattices}
\label{abstconsrelsec}

We have seen, in Section~\ref{conspwsetsec}, that symmetric consequence relations are binary relations on the powerset of a set. In this section~--- again, according to~\cite{galatostsinakis}~--- we introduce the notion of a consequence relation on an arbitrary sup-lattice, and show that consequence relations on a given sup-lattice are in bijective correspondence with closure operators on it. Also in this case, a suitable notion of \emph{finitarity} can be defined, for closure operators, in such a way that the correspondence, between consequence relations and closure operators, associates each finitary relation to a finitary operator and vice versa.

Moreover, in Section~\ref{conspwsetsec} we have seen that the notion of substitution can be formalized, in the case of powersets, as a module action from the powerset of a monoid. But we know, from Example~\ref{monoidaction}, that this is a paradigmatic example of quantale module. Then, if a sup-lattice is a module over a quantale, the substitution invariance for consequence relations will be formalized, in this case, as invariance under the action of the quantale. Finally, the closure operators associated to these structural consequence relations will be precisely the nuclei over the given $\Q$-module.

\begin{definition}\label{consrelsl}
Let $\MM = \la M, \vee, \bot \ra$ be a sup-lattice. A (\emph{symmetric}) \emph{consequence relation} on $\MM$ is a binary relation $\vdash$ on $\MM$ that satisfies, for all $x, y, z \in M$,
\begin{enumerate}
\item[$(i)$] if $y \leq x$, then $x \vdash y$;
\item[$(ii)$] if $x \vdash y$ and $y \vdash z$, then $x \vdash z$;
\item[$(iii)$] $x \vdash \bigvee_{x \vdash y} y$.
\end{enumerate}
Note that $\vdash$ satisfies the first two conditions iff it is a pre-order on $\MM$ that contains the relation $\geq$.
\end{definition}

A subset $X$ of $M$ is called \emph{directed} in $\MM$ if it upward directed in the poset $\la M, \leq \ra$, namely if for all $x, y \in X$, there exists a $z \in X$ such that $x, y \leq z$. An element $x$ of $\MM$ is called \emph{compact} if, for all directed $Y \subseteq M$, $x \leq \bigvee Y$ implies $x \leq y$, for some $y \in Y$. Equivalently, $x$ is compact if for all $Z \subseteq M$, if $x \leq \bigvee Z$, then there is a finite subset $Z_0$ of $Z$ such that $x \leq \bigvee Z_0$. For every subset $X$ of $M$, we denote by $K_\MM(X)$ the set of compact elements of $\MM$ that are contained in $X$. We write $K_\MM$ for $K_\MM(M)$.

\begin{definition}\label{slfinitary}
A consequence relation on $\MM$ is called \emph{finitary}, if for all $x, y \in M$, if $x \vdash y$ and $y$ is compact, then there is a compact element $x_0 \in M$ such that $x_0 \leq x$ and $x_0 \vdash y$.
\end{definition}

A closure operator $\g$ on a sup-lattice $\MM$ is called \emph{finitary}, if it preserves directed joins; i.e., for all directed $X$, $\g(\bigvee X) = \bigvee \g[X]$. A \emph{finitary sup-lattice}\index{Finitary!-- sup-lattice}\index{Sup-lattice!Finitary --} is a sup-lattice in which every element is a join of compact elements; in particular, $x = \bigvee K_\MM([\bot,x])$, for all $x \in M$.\footnote{Sup-lattices and closure operators that we call ``finitary'', are usually called ``algebraic'' in literature. This different nomenclature is made necessary by other uses of the term ``algebraic'' in this area.}

Notoriously, the compact elements of a powerset $\wp(S)$ are precisely the finite subsets of $S$. So, in the case where $M = \wp(S)$ for some set $S$, both notions of consequence relation and finitary consequence relation give back the ones defined for powersets.

As we anticipated, to define substitution invariance of a consequence relation on a sup-lattice, we need to assume that it is endowed with a module action with certain features we recognize to be exactly those of a module action of a quantale on the given sup-lattice. Then we have the following definition.
\begin{definition}\label{mqstructural}
A consequence relation $\vdash$ on a $\QQ$-module $\MM$ is called \emph{structural}, if $x \vdash y$ implies $q \star x \vdash q \star y$, for all $x, y \in M$ and $q \in Q$.
\end{definition}
Again, if $M = \wp(S)$, for some set $S$, and $Q = \wp(\varSigma)$, where $\varSigma$ is a monoid that acts on $S$, the notions of structurality and of substitution invariance of a consequence relation coincide.

Let $S_1$, $S_2$ be arbitrary sets, and let $\vdash_1$, $\vdash_2$ be consequence relations on $\wp(S_1)$ and $\wp(S_2)$, respectively. We have seen that the maps $\tau: \wp(S_1) \lto \wp(S_2)$ and $\rho: \wp(S_2) \lto \wp(S_1)$ involved in the definition of similarity of $\vdash_1$ and $\vdash_2$ were assumed to preserve unions (see Section~\ref{conspwsetsec}). We have also noted that this is a necessary and sufficient condition for these maps to extend maps from the sets $S_1$ and $S_2$ to the powersets $\wp(S_1)$ and $\wp(S_2)$ respectively. The generalization of this notion in the setting of sup-lattices is that of a map that preserves arbitrary joins, hence that of a sup-lattice homomorphism or~--- equivalently~--- of a residuated map between sup-lattices.

Given a consequence relation $\vdash$ on a sup-lattice $\MM$, we define the map
\begin{equation}\label{gammavdash}
\begin{array}{lccc}
\g_\vdash: & M & \lto & M \\
& x & \lmapsto & \bigvee_{x \vdash y} y.
\end{array}
\end{equation}
Conversely, given a closure operator $\g: \MM \lto \MM$, we define a binary relation $\vdash_\g$ on $M$, by
\begin{equation}\label{vdashgamma}
x \vdash_\g y \quad \textrm{ iff } \quad y \leq \g(x), \qquad \textrm{for all } \ x, y \in M.
\end{equation}

\begin{lemma}\emph{\cite{galatostsinakis}}\label{conrel2clop}
Consequence relations on a sup-lattice $\MM$ are in bijective correspondence with closure operators on $\MM$ via the maps $\vdash \lmapsto \g_\vdash$ and $\g \lmapsto \vdash_\g$. If $\MM$ is finitary, then finitarity is preserved under this correspondence.

If $\MM$ is a $\QQ$-module, for some quantale $\QQ$, then structurality is preserved under the correspondence, i.e. the consequence relation $\vdash$ is structural iff $\g_\vdash$ is structural.
\end{lemma}
\begin{proof}
In the whole proof, $\vdash$ and $\g$ will denote, respectively, a consequence relation and a closure operator, with $\g_\vdash$ and $\vdash_\g$ defined respectively by (\ref{gammavdash}) and (\ref{vdashgamma}) above.

First we prove that $\g_\vdash$ is a closure operator. For all $m \in M$, we have $m \leq m$, so $m \vdash m$, hence $m \leq \bigvee_{m \vdash n} n = \g_\vdash(m)$. If $m \leq n$, then $n \vdash m$; so, for all $p \in M$, $m \vdash p$ implies $n \vdash p$. Consequently, $\{p \in M \mid m \vdash p\} \subseteq \{p \in M \mid n \vdash p\}$; thus $\g_\vdash(m) = \bigvee\{p \in M \mid m \vdash p\} \leq \bigvee \{p \in M \mid n \vdash p\} = \g_\vdash(n)$. If $n \leq \g_\vdash \g_\vdash(m)$, then $\g_\vdash \g_\vdash(m) \vdash n$.

By Definition~\ref{consrelsl}$(iii)$, $m \vdash  \g_\vdash(m)$ and $\g_\vdash(m) \vdash \g_\vdash \g_\vdash(m)$, so $m \vdash n$; hence, $n \leq \bigvee _{m \vdash p} p = \g_\vdash(m)$. Then $\g_\vdash \g_\vdash(m) \leq \g_\vdash(m)$, and $\g_\vdash$ is a closure operator.

As a second step, we show that $\vdash_\g$ is a consequence relation. If $n \leq m$, then $n \leq \g(m)$, so $m \vdash _\g n$. If $m \vdash_\g n$ and $n \vdash_\g p$, then $n \leq \g(m)$ and $p \leq \g(n)$, so $p \leq \g(m)$; i.e. $m \vdash_\g p$. Note that for all $m \in M$, $\bigvee_{m \vdash_\g n} n = \bigvee\{n \in M \mid n \leq \g(m)\}$ is equal to $\g(m)$. Moreover, $m \vdash_\g \bigvee_{m \vdash_\g n} n$ since $\g(m) \leq \g(m)$.

Now we are ready to prove that $\vdash_{\g_\vdash} = \vdash$ and $\g_{\vdash_\g} = \g$. For all $m \in M$, we have $\g_{\vdash_\g}(m) = \bigvee_{m \vdash_\g n} n = \bigvee_{n \leq \g(m)} n = \g(m)$. So, $\g_{\vdash_\g} = \g$. Conversely, for all $m, n \in M$, we have $m \vdash_{\g_\vdash} n$ iff $n \leq \g_\vdash (m)$ iff $n \leq \bigvee _{m \vdash p} p$. Note that $m \vdash \bigvee_{m \vdash p} p$ and that $n \leq \bigvee_{m \vdash p} p$ implies $\bigvee_{m \vdash p} p \vdash n$. Therefore, $n \leq \bigvee_{m \vdash p} p$ implies $m \vdash n$. On the other hand, $m \vdash n$ obviously implies $n \leq \bigvee_{m \vdash p} p$; hence $m \vdash_{\g_\vdash} n$ iff $m \vdash n$, for all $m, n \in M$, i.e. $\vdash_{\g_\vdash} = \vdash$.

Let $\g$ and $\MM$ be both finitary. If $m \vdash_\g n$ and $n$ is compact, then $n \leq \g(m) = \bigvee\g\left[K_{\MM}([\bot,m])\right]$. So, $n \leq \g(m_0)$, for some $m_0 \in K_{\MM}([\bot,m])$. In other words, there is a compact element $m_0 \leq m$ such that $m_0 \vdash_\g n$. It follows the finitarity of $\vdash_\g$.

Conversely, let $\vdash$ and $\MM$ be both finitary; we will show that $\g_\vdash$ is finitary too. Assume that $X$ is directed; we want to show that $\g_\vdash(\bigvee X) = \bigvee \g_\vdash[X]$. First note that $\g_\vdash(\bigvee X)$ is an upper bound of $\g_\vdash[X]$. To show that it is the least upper bound, let $p \in M$ be such that $\g_\vdash(m) \leq p$, for all $m \in X$. If $n$ is a compact element such that $n \leq \g_\vdash(\bigvee X)$, then $\bigvee X \vdash n$. Since $\vdash$ is finitary, there exists a compact element $m_0 \leq \bigvee X$, such that $m_0 \vdash n$; i.e. $n \leq \g_\vdash (m_0)$. Now, since $m_0 \leq \bigvee X$, $X$ is directed and $m_0$ is compact, there exists $m \in X$ such that $m_0 \leq m$. So, $n \leq \g_\vdash (m_0) \leq \g_\vdash(m) \leq p$. Then, from the fact that $n \leq p$, for all compact elements $n \leq \g_\vdash(\bigvee X)$, and the finitarity of $\MM$, we have $\g_\vdash(\bigvee X) \leq p$.

Last, we prove that structurality is preserved by (\ref{gammavdash},\ref{vdashgamma}). Let $\MM$ be a module over a certain quantale $\QQ$, and assume $\g$ to be structural. If $m \vdash_\g n$, then $n \leq \g(m)$; so $q \star n \leq \g (q \star m)$, by the structurality of $\g$. Consequently, $q \star m \vdash_\g q \star n$. Conversely, assume that $\vdash_\g$ is structural. We have $m \vdash_\g \g(m)$, so $q \star m \vdash_\g q \star \g(m)$, by the structurality of $\vdash_\g$. Hence, $q \star \g(m) \leq \g (q \star m)$.
\end{proof}

Now we can abstract also the notion of theory of a consequence relation and show that, also in the case of sup-lattices, the set of theories of a consequence relation can be structured as a lattice.

Let $\vdash$ be a consequence relation on a sup-lattice $\MM$. An element $t$ of $\MM$ is called a \emph{theory} of $\vdash$ if $t \vdash x$ implies $x \leq t$. Note that if $t$ is a theory, then $x \leq t$ and $x \vdash y$ imply $y \leq t$. We denote the set of theories of $\vdash$ by $\th_\vdash$.

\begin{lemma}\emph{\cite{galatostsinakis}}
If $\vdash$ is a consequence relation on the sup-lattice $\MM$, then $\th_\vdash = M_{\g_\vdash}$.
\end{lemma}
\begin{proof}
Let $t \in \th_\vdash$ and set $\g = \g_\vdash$. We will show that $t \in M_{\g_\vdash}$, i.e. that $\g(t) = t$. We have $\g(t) \leq \g(t)$, so $t \vdash \g(t)$; since $t$ is a theory, $\g(t) \leq t$. The other inequality holds because $\g$ is extensive.

Conversely, assume that $\g(t) = t$, and let $m \in M$ such that $t \vdash m$. Then $m \leq \g(t) = t$.
\end{proof}

Then, according to the last two results and Theorem~\ref{homoandnuclei}, for any sup-lattice $\MM$, and any consequence relation $\vdash$ on it, $\Th_\vdash = \la \th_\vdash, \g_\vdash \circ \vee, \g_\vdash(\bot)\ra$ is a sup-lattice and, if $\MM$ is a $\QQ$-module and $\vdash$ is structural, then $\Th_\vdash$ is a $\QQ$-module (a quotient of $\MM$). Moreover, every homomorphic image of a sup-lattice $\MM$ is the sup-lattice of the theories of a consequence relation on $\MM$, and every homomorphic image of a $\QQ$-module $\MM$ is the $\QQ$-module of the theories of a structural consequence relation on $\MM$.

\section{Similarity and equivalences of two consequence relations}
\label{eqtransec}

In this section we define the notions of representation, similarity and equivalence between two closure operators or two consequence relations. The approach is, again, in the wake of~\cite{galatostsinakis}, and generalizes the corresponding notions in~\cite{blokjonsson2}.

Let $\g$ and $\d$ be closure operators on the sup-lattices $\MM$ and $\NN$, respectively. A \emph{(non-structural) representation}\index{Consequence relation!Representation of --s} of $\g$ in $\d$ is a sup-lattice homomorphism $f: \MM_\g \lto \NN_\d$. A representation $f: \MM_\g \lto \NN_\d$ of $\g$ in $\d$ is said to be \emph{induced} by the homomorphism $h: \MM \lto \NN$, if $f \circ \g = \d \circ h$.
$$\xymatrix{
\MM \ar[r]^h \ar @{->>}[d]_\g & \NN  \ar @{->>}[d]^\d \\
\MM_\g \ar[r]^f               		& \NN_\d	
}$$
A non-structural representation is called \emph{conservative} if $f$ is injective, \emph{non-conservative} otherwise.

In view of the correspondence between consequence relations and closure operators, we will denote an arbitrary consequence relation on a sup-lattice $\MM$ by $\vdash_\g$ with the understanding that $\g$ is the associated closure operator. 
 
We say that a consequence relation $\vdash_\g$ is \emph{represented} in the consequence relation $\vdash_\d$ if the associated closure operator $\g$ is represented in $\d$; the representation of $\vdash_\g$ in $\vdash_\d$ is \emph{induced} by a homomorphism $h: \MM \lto \NN$, if the representation of the corresponding closure operators is induced by $h$. Corollary~\ref{crrep} shows that $\vdash_\g$ is represented in $\vdash_\d$ via $h$ if and only if for all  $x, y \in M$,
$$x \vdash_\g y \quad \textrm{iff} \quad h(x) \vdash_\d h(y).$$

Recall that, by Lemma~\ref{conrel2clop}, a closure operator $\g$ on an finitary sup-lattice $\MM$ is finitary iff $\vdash_\g$ is finitary, that is iff, for all $x, y \in M$, if $y \leq \g(x)$ and $y$ is compact, then there exists a compact element $x_0 \leq x$ such that $y \leq \g(x_0)$.

\begin{lemma}\emph{\cite{galatostsinakis}}
Let $\MM$ and $\NN$ be sup-lattices and let $h \in \hom_\SL(\MM, \NN)$. If $\d$ is a closure operator on $\NN$, then the map $\d^h = h_* \circ \d \circ h: M \lto M$ is a closure operator on $\MM$. If $\MM$, $h$ and $\d$ are finitary, then $\d^h$ is finitary as well. If $\MM$ and $\NN$ are $\QQ$-modules, $h \in \hom_\QQ(\MM, \NN)$ and $\d$ is a structural closure operator on $\NN$, then $\d^h$ is structural. 
\end{lemma}
\begin{proof}
Note that $\d: \NN \lto \NN_\d$ is a homomorphism, and the inclusion map $i_\d = \id_{N \restr N_\d}$ is its dual homomorphism (or, that is the same, its residual map), so $\d \circ h: \MM \lto \NN_\d$ is a homomorphism as well, with residual $h_* \circ i_\d = {h_*}_{\restr \NN_\d}$, by the basic properties of sup-lattice morphisms (see Section~\ref{suplatsec}). 
$$\xymatrix{
\MM \ar[rr]<.7ex>^h && \NN  \ar@{->>}[dd]<.7ex>^\d \ar[ll]<.7ex>^{h_*} \\
&&\\
                		&& \NN_\d	\ar@<.5ex>[uull]^{{h_*}_{\restr \NN_\d}}	\ar@{^{(}->}[uu]<.7ex>
}$$
Therefore, $\d^h = h_* \circ \d \circ h = {h_*}_{\restr \NN_\d} \circ \d \circ h$ is a closure operator on $\MM$.

Assume that $\MM$, $h$ and $\d$ are finitary. If $y \leq \d^h(x)$, for some compact element $y$, then $y \leq h_* \circ \d \circ h(x)$, so $h(y) \leq \d \circ h(x)$. Since $h$ is finitary and $y$ is compact, $h(y)$ is compact. Furthermore, since $\d$ is finitary, there is a compact element $x' \leq h(x)$ such that $h(y) \leq \d(x')$. Since $\MM$ is finitary, $x = \bigvee K_{\MM}([\bot,x])$, so $h(x) = \bigvee h\left[K_{\MM}([\bot,x])\right]$. Since $x' \leq h(x)$, there exists a compact element $x_0 \leq x$ such that $x' \leq h(x_0)$. Consequently, $h(y) \leq \d \circ h (x_0)$, hence $y \leq h_* \circ \d \circ h (x_0) = \d^h (x_0)$, for some compact element $x_0 \leq x$. Thus, $\d^h$ is finitary.

For all $q \in Q$ and $x \in M$, we have
$$\begin{array}{l}
h(q \star \d^h (x)) = q \star h(\d^h(x)) = q \star ((h \circ h_* \circ \d \circ h)(x)) \leq \\
q \star (\d(h(x)) \leq \d(q\star h(x)) = \d(h(q \star x)),
\end{array}$$
so $q \star \d^h(x) \leq (h_* \circ \d \circ h)(q \star x) = \d^h(q \star x)$.
\end{proof}

We will call $\d^h$ the \emph{$h$-transform} of $\d$. Similarly, we can define the $h$-transform of a consequence relation $\vdash$ on $\NN$ to be the relation $\vdash^h$ on $\MM$ defined by $x \vdash^h y$ iff $h(x) \vdash h(y)$, for all $x,y \in P$. The following lemma shows that the $h$-transform of a consequence relation is a consequence relation and the associated closure operator is the $h$-transform of the original relation.

\begin{lemma}\emph{\cite{galatostsinakis}}\label{delta2tau}
Let $\MM$ and $\NN$ be sup-lattices, $h: \MM \lto \NN$ a homomorphism and $\d$ a closure operator on $\NN$. The following statements are equivalent:
\begin{enumerate}
\item[$(a)$]$\g = \d^h$,
\item[$(b)$]for all $x, y \in M$, $x \vdash_\g y$ iff $h(x) \vdash_\d h(y)$,
\item[$(c)$]$M_\g = h^{-1}[M_\d]$
\end{enumerate}
\end{lemma}
\begin{proof}
Assume $(a)$ holds; then for all $x, y \in M$, we have $x \vdash_{\d^h} y$ iff $y \leq h_* \d h(x)$ iff $h(y) \leq \d h(x)$ iff  $h(y) \vdash_\d h(x)$. Conversely, for all $x, y \in M$, we have $y \leq \g(x)$ iff $x \vdash_\g y$ iff $h(x) \vdash_\d h(y)$ iff $h(y) \leq \d h(x)$ iff $y \leq h_* \d h(x)$. Consequently, $\g = \d^h$.

For the equivalence between $(a)$ and $(c)$, note that $M_\g = h^{-1}[M_\d]$ means that, for all $x \in M$, $x = \g(x)$ iff $\d h(x) = h(x)$. Moreover, for all $x \in M$, we have $\d h(x) = h(x)$ iff $\d h(x) \leq h(x)$ iff $h_* \d h(x) \leq x$ iff $\d^h(x) = x$, since $\d$ and $\d^h$ are closure operators. Consequently, $(c)$ holds iff $\g$ and $\d^h$ have the same fixed elements; i.e. $\g = \d^h$.
\end{proof}

\begin{lemma}\emph{\cite{galatostsinakis}}\label{jnrep}
Let $\MM$ and $\NN$ be sup-lattices, $h: \MM \lto \NN$ a homomorphism and $\d$ a closure operator on $\NN$. 
\begin{enumerate}
\item[$(i)$]The map $f = \d \circ h_{\restr \MM_{\d^h}}: \MM_{\d^h} \lto \NN_\d$ is a sup-lattice homomorphism whose residuum is $f_* = {h_*}_{\restr \NN_\d} = \d^h \circ {h_*}_{\restr \NN_\d}: \left(\NN_\d\right)^{\op} \lto \left(\MM_{\d^h}\right)^{\op}$. 
\item[$(ii)$]$f$ is a representation of $\d^h$ in $\d$ induced by $h$.
\item[$(iii)$]$\d^h$ is the only closure operator on $\MM$ that is represented in $\d$ under a representation induced by $h$.
\item[$(iv)$]If $\MM$, $h$ and $\d$ are finitary, and $\d$ reflects existing joins, then $f$ is finitary, as well.
\item[$(v)$]If $\MM$ and $\NN$ are $\QQ$-modules, $h: \MM \lto \NN$ is a $\QQ$-module morphism and $\d$ is a structural closure operator on $\NN$, then $f$ is structural. 
\end{enumerate}
\end{lemma}

\begin{corollary}\emph{\cite{galatostsinakis}}\label{crrep}
Let $\MM$ and $\NN$ be posets and let $\vdash_\g$ and $\vdash_\d$ be consequence relations on $\MM$ and $\NN$, respectively. Then, $\vdash_\g$ is represented in $\vdash_\d$ via a homomorphism $h: \MM \lto \NN$, iff for all $x, y \in M$, we have  $x \vdash_\g y$ iff $h(x) \vdash_\d h(y)$.
\end{corollary}
\begin{proof}
It is a direct consequence of Lemma~\ref{delta2tau} and Lemma~\ref{jnrep}$(iii)$.
\end{proof}

It is easy to see that $\vdash_\g$ is represented in $\vdash_\d$ by $f: \Th_{\vdash_\g} \lto \Th_{\vdash_\d}$ means that $f$ is residuated and for all $x, y \in M$, 
$$x \vdash_\g y \quad \textrm{ iff } \quad (f \circ \g)(x) \vdash_\d (f \circ \g)(y).$$

Indeed, if $\vdash_\g$ is represented in $\vdash_\d$ by $f$, then $x \vdash_\g y$ iff $y \leq \g(x)$ iff $f(y) \leq f(\g(x))$ (since $f$ preserves and reflects order) iff $f(y) \leq \d(f(\g(x)))$ iff $f(\g(x)) \vdash_\d f(\g(y))$. Conversely, to show that $f$ reflects order, let $f(\g(y)) \leq f(\g(x))$. Then $f(\g(y)) \leq \d(f(\g(x)))$ that is $f(\g(x)) \vdash_\d f(\g(y))$; so $x \vdash_\g y$ that is $\g(y) \leq \g(x)$.

Let $\g$ and $\d$ be closure operators on the sup-lattices $\MM$ and $\NN$, respectively. A \emph{similarity}\index{Consequence relation!Similarity of --s} between $\g$ and $\d$ is an isomorphism $f: \MM_\g \lto \NN_\d$. If there exists a similarity between $\g$ and $\d$, then $\g$ and $\d$ are called \emph{similar}. A similarity $f$ between $\g$ and $\d$ is said to be \emph{induced} by the homomorphisms $h: \MM \lto \NN$ and $k: \NN \lto \MM$, if $f \circ \g = \d \circ h$ and $f^{-1} \circ \d = \g \circ k$. In this case we will say that $\g$ and $\d$ are similar \emph{via} $h$ and $k$.
$$\xymatrix{
\MM \ar[rr]<.7ex>^h \ar @{->>}[dd]_\g && \NN  \ar @{->>}[dd]^\d \ar[ll]<.7ex>^k  \\
&&\\
\MM_\g \ar@{^{(}->>}[rr]<.7ex>^f               		&& \NN_\d	\ar@{^{(}->>}[ll]<.7ex>^{f^{-1}}
}$$
It is clear that $f$ is a similarity 	between $\g$ and $\d$ iff $f$ is a representation	of $\g$ in $\d$, $f$ is a bijection and
$f^{-1}$ is a representation of $\d$ in $\g$.

A consequence relation $\vdash_\g$ is called \emph{similar} to the consequence relation $\vdash_\d$ (via a homomorphism $h$) if $\g$ is similar to $\d$ (via $h$).

\begin{lemma}\emph{\cite{galatostsinakis}}\label{sim}
Let $\g$ and $\d$ be closure operators on the sup-lattices $\MM$ and $\NN$, respectively. The following statements are equivalent.
\begin{enumerate}
\item[$(a)$]$\g$ and $\d$ are similar via (a similarity induced by) $h$ and $k$. 
\item[$(b)$]$\g = \d^h$  and $\d \circ h \circ k = \d$.
\item[$(c)$]$\d = \g^k$ and $\g \circ k \circ h = \g$.
\end{enumerate}
\end{lemma}
\begin{proof}
We will show the equivalence of the first two statements; the equivalence of the first to the third will follow by symmetry. The forward direction follows from Lemma~\ref{jnrep}$(iii)$ and the definition of similarity ($\d \circ h \circ k = f \circ \g \circ k = f \circ f^{-1} \circ \d = \d$). For the converse, assume that $\g = \d^h$ and $\d \circ h \circ k = \d$. Let $f$ be the representation of $\g = \d^h$ in $\d$ given in Lemma~\ref{jnrep}$(i)$. We have $f \circ \g = \d \circ h$, by Lemma~\ref{jnrep}$(ii)$.
$$\xymatrix{
\MM \ar[rr]<.7ex>^h \ar @{->>}[dd]_\g && \NN  \ar @{->>}[dd]^\d \ar[ll]<.7ex>^k  \\
&&\\
\MM_\g \ar[rr]<.7ex>^f               		&& \NN_\d	\ar@{-->}[ll]<.7ex>^{f^{-1}}
}$$
To show that $f$ is onto, let $y \in N_\d$ and set $x = \g k(y) \in M_\g$. We have $f(x) = (f \circ \g \circ k)(y) = (\d \circ h \circ k)(y) = \d(y) = y$. Consequently, $f$ is an order-isomorphism and $\g$ and $\d$ are similar. To show that the similarity $f$ is induced by $h$ and $k$, we need only prove that $f^{-1} \circ \d = \g \circ k$, or equivalently that $\d = f \circ \g \circ k$. This is true, because $\d = \d \circ h \circ k$ and $f \circ \g = \d \circ h$.
\end{proof}

\begin{corollary}\emph{\cite{galatostsinakis}}\label{crsim}
Let $\MM$ and $\NN$ be sup-lattices and let $\vdash_\g$ and $\vdash_\d$ be consequence relations on $\MM$ and $\NN$, respectively. Then, $\vdash_\g$ is similar to $\vdash_\d$ via the homomorphisms $h: \MM \lto \NN$ and $k: \NN \lto \MM$, iff the following two conditions hold:
\begin{enumerate}
\item[$(i)$]for all $x, y \in M$, we have  $x \vdash_\g y$ iff $h(x) \vdash_\d h(y)$,
\item[$(ii)$]for all $z \in N$, $z \dashv \vdash_\d (h \circ k)(z)$.
\end{enumerate}
\end{corollary}
\begin{proof}
It is easy to see that $\d \circ h \circ k = \d$ iff for all $z \in N$, $z \dashv \vdash_\d (h \circ k)(z)$. Now, the corollary follows from Lemma~\ref{sim}$(b)$ and Corollary~\ref{crrep}.
\end{proof}

Let $\MM$ and $\NN$ be $\QQ$-modules and let $\g$ and $\d$ be structural closure operators on $\MM$ and $\NN$, respectively. An \emph{equivalence} between $\g$ and $\d$ is a module isomorphism $f: \MM_\g \lto \NN_\d$. Note that an equivalence is just a structural similarity. Moreover, $f^{-1}$ is also structural. If such an isomorphism exists then $\g$ and $\d$ are called \emph{equivalent}. If the equivalence is induced by module morphisms $h: \MM \lto \NN$ and $k: \NN \lto \MM$, then $\g$ and $\d$ are called \emph{equivalent via $h$ and $k$}.

\begin{theorem}\emph{\cite{galatostsinakis}}
Let $\MM$ and $\NN$ be $\QQ$-modules, and let $\g$ and $\d$ be structural closure operators on $\MM$ and $\NN$, respectively. If $\g$ and $\d$ are similar via the structural translators $h$ and $k$ then they are equivalent via $h$ and $k$.
\end{theorem}
\begin{proof}
It suffices to show that the similarity $f$ of $\g$ in $\d$ is structural. Indeed, for all $q \in Q$ and $x \in M_\g$, we have
$$\begin{array}{ll}
f(q \star_\g x) & = (f \circ \g)(q \star x) = (\d \circ h)(q \star x) = \d(q \star h(x)) \\
& = \d(q \star (\d \circ h)(x)) = q \star_\d (\d \circ h)(x) = q \star_\d (f \circ \g)(x) \\
& = q \star_\d f(x),
\end{array}$$
since $\g(x) = x$. The thesis follows.
\end{proof}

\section{Equivalences induced by translators}
\label{eqtranssec}

The next step, facing the problem of comparing two deductive systems on different languages, consists of finding an answer to this question: in which cases the equivalences of consequence relations, defined in the previous section, are induced by translators? It is easy to see that we cannot always build the homomorphisms $h$ and $k$ that close the following diagram
\begin{equation}\label{S}
\xymatrix{
\MM \ar@{-->}@<.7ex>[rr]^h \ar @{->>}[dd]_\g && \NN  \ar @{->>}[dd]^\d \ar@{-->}@<.7ex>[ll]^k  \\
&&\\
\MM_\g \ar@{^{(}->>}@<.7ex>[rr]^f               		&& \NN_\d	\ar@{^{(}->>}@<.7ex>[ll]^{f^{-1}}.
}
\end{equation}
Nevertheless, we will show that such a construction is possible in all standard situations including the powersets of formulas, equations and sequents.

In this section we present conditions on modules, due again to N. Galatos and C. Tsinakis, under which every equivalence is induced by translators. More precisely, we will see that, if $\MM$ is a $\QQ$-module that satisfies these conditions, then the map $h$, for the diagram above can be built. Analogously, the map $k$ will exist if $\NN$ satisfies the same conditions.

\begin{lemma}\label{8.1}
The $\QQ$-modules $\MM$ for which all squares of type (\ref{S}) can be completed are the projective objects of $\QQ\Mod$.
\end{lemma}
\begin{proof}
See Lemma 8.1 of~\cite{galatostsinakis}.
\end{proof}

\begin{lemma}\emph{\cite{galatostsinakis}}\label{formeqproj}
$\wp(\Fml)$ and $\wp(\EQ_\lang)$ are projective cyclic $\wp(\Sfm)$-modules.
\end{lemma}
\begin{proof}
Let $v = \{x_1\}$, and $u = \{\kappa_{x_1}\}$, where $\kappa_{x_1}$ is the substitution that maps all variables to $x_1$. Then the thesis, for $\wp(\Fml)$, is a trivial consequence of Corollary~\ref{cycliccond} and Theorem~\ref{cyclicproj}. Of course, the choice of $x_1$ is completely arbitrary, any variable would have served the scope as well.

For the module $\wp(\EQ_\lang)$, we consider two arbitrary distinct variables, say $x_1$ and $x_2$, and we partition the set $\V$ in two disjoint sets $V_1$ and $V_2$. Then we can consider the substitution $\kappa_{x_1 \eq x_2}$ that sends every element of $V_1$ to $x_1$ and every element of $V_2$ to $x_2$. Then, if we take $v = \{x_1 \eq x_2\}$ and $u = \{\kappa_{x_1 \eq x_2}\}$, we obtain the thesis applying, again, Corollary~\ref{cycliccond} and Theorem~\ref{cyclicproj}.
\end{proof}
 
Moreover, Lemma~\ref{formeqproj} extends to the following more general result 
\begin{theorem}\label{seqproj}
The $\wp(\Sfm)$-module $\wp(\mathbf{Seq}_\lang)$ is a coproduct of cyclic projective modules. Consequently it is projective.
\end{theorem}
\begin{proof}
The proof is similar to that of Lemma~\ref{formeqproj}. Indeed, it is easy to see that for all $(m,n) \in \Tp(\mathit{Seq}_\lang)$, the set of sequents of type $(m,n)$ is cyclic, and generated by the sequent $\{x_1, \ldots, x_m \Rightarrow x_{m+1}, \ldots, x_{m+n}\}$. Once the set $\V$ has been partitioned in $m+n$ sets $\{V_j\}_{j=1}^{m+n}$, the idempotent scalar that generates a $\wp(\Sfm)$-submodule of $\wp(\Sfm)$ is $\kappa_{m,n}$, the substitution that sends each variable in $V_j$ to $x_j$, for all $j \leq m + n$. Then $\wp(\mathbf{Seq}_\lang)$ is easily seen to be the coproduct of the cyclic projective modules generated by $\{x_1, \ldots, x_m \Rightarrow x_{m+1}, \ldots, x_{m+n}\}$, with $(m,n) \in \Tp(\mathit{Seq}_\lang)$. See Theorem 8.12 of~\cite{galatostsinakis} for details.
\end{proof}


\section{Interpretations between abstract deductive systems}
\label{interabs}

\begin{definition}\label{absinter}
Let $\QQ$ and $\RR$ be quantales, and $\la \MM, \g \ra$ and $\la \NN, \d \ra$ be two deductive systems, over $\QQ$ and $\RR$ respectively. A \emph{non-structural representation} is a sup-lattice homomorphism $f: \MM_\g \lto \NN_\d$. If $f$ is an isomorphism, it is called a \emph{similarity}.

A \emph{structural representation} $\rho: \la \QQ, \MM_\g \ra \lto \la \RR, \NN_\d \ra$ of $\la \QQ, \MM_\g \ra$ into $\la \RR, \NN_\d \ra$ is a pair $(h, f)$ constituted by a quantale homomorphism $h: \QQ \lto \RR$ and a $\QQ$-module homomorphism $f: \MM_\g \lto (\NN_\d)_h$.

A representation is called
\begin{itemize}
\item \emph{conservative} if $f$ is injective, \emph{non-conservative} otherwise;
\item \emph{faithful} if it is structural and conservative;
\item an \emph{equivalence} if $f$ is bijective and there exists a quantale homomorphism $k: \RR \lto \QQ$ such that $f$ is an isomorphisms of both $\QQ$-modules and $\RR$-modules.
\end{itemize}
If $\rho = (h, f)$ is a structural (respectively: $f$ is a non-structural) representation and there exists a $\QQ$-module (resp.: a sup-lattice) homomorphism $\iota$ such that the diagram
\begin{equation}\label{absinterdiag}
\xymatrix{
\MM \ar@{-->}[rr]^{\iota} \ar@{->>} [dd]_{\g} & & \NN_h \ar@{->>} [dd]^\d\\
&&\\
\MM_\g \ar[rr]_{f} & & (\NN_\d)_h.
\\
}
\end{equation}
commutes, $\iota$ is called an \emph{interpretation}. In this case, we will often use the term ``interpretation'' also for the representation $\rho$ (resp.: for $f$).
\end{definition}

Recalling that, in a lattice $\LL$, a \emph{completely join prime element} is an element $x$ such that, whenever $x \leq \bigvee Y$, for some $Y \subseteq L$, then $x \leq y$ for some $y \in Y$, we set the following further definition.

\begin{definition}\label{i-trans}
Let $\QQ$ and $\RR$ be quantales, $\SS$ a subquantale of $\QQ$, $\TT$ a subquantale of $\RR$ and $i: \SS \lto \TT$ an isomorphism. A quantale homomorphism $h: \QQ \lto \RR$ is called a \emph{translation relatively to $i$}, or an \emph{$i$-translation}, if it satisfies the following conditions:
\begin{enumerate}
\item[$(i)$] $h$ sends all the completely join prime elements of $\QQ$ to completely join prime elements of $\RR$;
\item[$(ii)$] the counterimage, via $h$, of any multiplicative idempotent element of $\RR$, is either empty or composed only of idempotent elements of $\QQ$;
\item[$(iii)$] $h^{-1}(i(x)) = \{x\}$ for all $x \in S$.
\end{enumerate}
\end{definition}

According to Theorem~\ref{adjfunct}, the existence of a representation $\rho$ of $\la \QQ, \MM_\g \ra$ into $\la \RR, \NN_\d \ra$ canonically defines an adjoint and co-adjoint functor $H := ( \ )_h$ from $\RR\Mod$ to $\QQ\Mod$ and, then, two $\RR$-module morphisms
$$H_lf: \RR \tensor_\QQ \MM_\g \lto \NN_\d \quad \textrm{ and } \quad H_rf: \NN_\d \lto \Hom_\QQ(\RR_h, \MM_\g);$$
moreover, if $\rho$ is an interpretation, $\iota$ defines two further $\RR$-module morphisms
$$H_l\iota: \RR \tensor_\QQ \MM \lto \NN \quad \textrm{ and } \quad H_r\iota: \NN \lto \Hom_\QQ(\RR_h, \MM).$$
As the authors observed in \cite{galatostsinakis}, the $\QQ$-modules (respectively: the sup-lattices) $\MM$ for which any representation is an interpretation, i.e. for which any diagram of type (\ref{absinterdiag}) can be completed, are precisely the projective objects of $\QQ\Mod$ (resp.: of $\SL$); hence, in the concrete cases of deductive systems of formulas, equations or sequents, any representation is actually an interpretation. This will be made more explicit in the next sections.

\section{Translations and quantale morphisms}
\label{transsec}

In this section we start facing the problem of comparing two propositional logics having different underlying languages.

Let $\lang = \la L, \nu \ra$ and $\lang' = \la L', \nu' \ra$ be two propositional languages. We recall from Section~\ref{proplogicsec} that, given a propositional language $\lang$ and a denumerable set of variables $\V = \{x_n \mid n \in \N\}$, the $\lang$-formulas are defined recursively by means of the following conditions:
\begin{enumerate}
\item[(F1)] every propositional variable is an $\lang$-formula,
\item[(F2)] every constant symbol is a formula,
\item[(F3)] if $f$ is a connective of arity $\nu(f) > 0$ and $\phi_1, \ldots, \phi_{\nu(f)}$ are $\lang$-formulas, then $f(\phi_1, \ldots, \phi_{\nu(f)})$ is an $\lang$-formula,
\item[(F4)] all $\lang$-formulas are built by iterative applications of (F1), (F2) and (F3).
\end{enumerate}

Let $\lang = \la L, \nu \ra$ be a propositional language. If $n \in \N_0$ and $f: \fml^n \lto \fml$ is a map, $f$ is called a \emph{derived connective in $\lang$} if there exists a formula $\phi_f = \phi_f(x_1, \ldots, x_n) \in \fml$ in the $n$ variables $x_1, \ldots, x_n$ such that $f(\psi_1, \ldots, \psi_n) = \phi_f(x_1/\psi_1, \ldots, x_n/\psi_n)$, for all $\psi_1, \ldots, \psi_n \in \fml$. In this case, we also say that $f$ is a connective \emph{derived from the ones in $L$}. In particular, if $n = 0$, $f$ is a \emph{derived constant}, i.e. a formula in $\fml$ containing only constants and no variables.

Starting from derivable connectives, we want to define a concept of language translation.

\begin{definition}\label{ltranslation}
Let $\lang = \la L, \nu \ra$ and $\lang' = \la L', \nu' \ra$ be two propositional languages and assume that, for each connective $f \in L$, there exists a derived connective $f'$ in $\lang'$ of arity $\nu(f)$. If we denote by $L''$ the set of such derived connectives, the structure $\la \fmll, L'' \ra$ turns out to be an $\lang$-algebra, i.e. an algebra of the same type of $\la \fml, L\ra$. In this case, a map $\tau: \fml \lto \fmll$ is called a \emph{language translation} of $\lang$ into $\lang'$ if
\begin{enumerate}
\item[$(i)$] $\tau^{-1}(x) = \{x\}$ for any variable $x$,
\item[$(ii)$] $\tau$ is an $\lang$-homomorphism, that is
$$\tau(f(\phi_1, \ldots, \phi_{\nu(f)})) = f'(\tau(\phi_1), \ldots, \tau(\phi_{\nu(f)})),$$
for all $f \in L$ and $\phi_1, \ldots, \phi_{\nu(f)} \in \fml$.
\end{enumerate}
\end{definition}

Now, assuming that there exists a language translation of $\lang$ into $\lang'$, let us denote by $\Fmll''$ the $\lang$-algebra $\la \fmll, L'' \ra$. We have the following results.

\begin{lemma}\label{sigmainend}
The monoid of substitutions $\Sfmm$ (over the language $\lang'$) of $\fmll$ is a submonoid of the $\lang$-endomorphism monoid $\mathbf{End}_\lang(\Fmll'')$ of the $\lang$-algebra $\la \fmll, L'' \ra$.
\end{lemma}
\begin{proof}
The inclusion $\sfmm \subseteq \operatorname{End}_\lang(\Fmll'')$ comes easily from the fact that the connectives in $L''$ are derived by those in $L'$, so they are preserved by any substitution of $\fmll$. Then it is clear that $\Sfmm$ is actually a submonoid of $\mathbf{End}_\lang(\Fmll'')$ for it contains the identity map and is closed under composition.
\end{proof}

It may be important to underline, regarding the previous lemma, that the converse inclusion does not hold in general. For example, if $g_1$ and $g_2$ are two connectives of $\lang'$ of the same arity~---~say $n$~---~that are not involved in any of the formulas that define the connectives in $L''$, then a map that sends each variable to itself and $g_1(x_1, \ldots, x_n)$ to $g_2(x_1, \ldots, x_n)$ can be extended to an $\lang$-endomorphism of $\Fmll''$ that is not a substitution in $\lang'$.

\begin{lemma}\label{qinq'}
Let $\lang = \la L, \nu \ra$, $\lang' = \la L', \nu' \ra$ and $\tau$ be a language translation of $\lang$ into $\lang'$. Then the following hold:
\begin{enumerate}
\item[$(i)$]$\tau$ defines a monoid homomorphism $\ov \tau$ from $\Sfm$ to $\Sfmm$ (and, thus, a quantale homomorphism, still denoted by $\ov \tau$, from $\wp(\Sfm)$ to $\wp(\Sfmm)$);
\item[$(ii)$]$\ov \tau$ is injective if and only if so is $\tau$;
\item[$(iii)$]the sup-lattice homomorphism between $\wp(\Fml)$ and $\wp(\Fmll)$ that extends $\tau$ is a $\wp(\Sfm)$-module homomorphism, where the $\wp(\Sfm)$-module structure of $\wp(\Fmll)$ is the one induced by $\ov \tau$.
\end{enumerate}
\end{lemma}
\begin{proof}
\begin{enumerate}
\item[$(i)$]For all $\sigma \in \sfm$, let $\sigma'$ be the substitution uniquely determined by the map $\tau \circ \sigma_{\restr \V} \in \fmll^\V$ and $\ov \tau: \sigma \in \sfm \lmapsto \sigma' \in \sfmm$.

Obviously $\ov \tau(\id_{\fml}) = \id_{\fmll}$. Now let $\sigma_1, \sigma_2 \in \sfm$; we want to show that $\ov \tau(\sigma_2 \circ \sigma_1) = \ov \tau(\sigma_2) \circ \ov \tau(\sigma_1)$ for all $\sigma_1, \sigma_2 \in \sfm$. First of all recall that, for any formula $\phi(x_1, \ldots, x_n) \in \fml$ in the variables $x_1, \ldots, x_n$, by the definition of language  translation, $\tau(\phi)$ is a formula in the same variables, and we will denote it as $\phi'(x_1, \ldots, x_n)$.

We consider an arbitrary variable $x$ and set $\sigma_1(x) = \phi(x_1, \ldots, x_n)$ and $\sigma_2(x_i) = \psi_i(x_{i1}, \ldots, x_{ik_i})$, $i = 1, \ldots, n$, with $\phi, \psi_1, \ldots, \psi_n \in \fml$. Then we have
$$\begin{array}{l}
  \ov \tau(\sigma_2 \circ \sigma_1)(x) \\
= \tau((\sigma_2 \circ \sigma_1)(x)) \\
= \tau(\sigma_2(\phi(x_1, \ldots, x_n))) \\
= \tau(\phi(\psi_1(x_{11}, \ldots, x_{1k_1}), \ldots, \psi_n(x_{n1}, \ldots, x_{nk_n}))) \\
= \phi'(\tau(\psi_1(x_{11}, \ldots, x_{1k_1})), \ldots, \tau(\psi_n(x_{n1}, \ldots, x_{nk_n}))) \\
= \phi'(\tau(\sigma_2(x_1)), \ldots, \tau(\sigma_2(x_n))) \\
= \phi'((\ov \tau(\sigma_2))(x_1), \ldots, (\ov \tau(\sigma_2))(x_n)) \\
= \ov \tau(\sigma_2)(\phi'(x_1, \ldots, x_n)) \\
= \ov \tau(\sigma_2)(\tau(\phi(x_1, \ldots, x_n))) \\
= \ov \tau(\sigma_2)(\tau(\sigma_1(x))) \\
= \ov \tau(\sigma_2)(\ov \tau(\sigma_1)(x)) \\
= (\ov \tau(\sigma_2) \circ \ov \tau(\sigma_1))(x).
\end{array}$$
Then the $(i)$ follows from the arbitrary choice of $\sigma_1, \sigma_2 \in \sfm$ and $x \in \V$.

\item[$(ii)$]If $\tau$ is injective, then $\ov \tau$ is obviously injective too. On the other hand, if $\tau$ is not injective, then there exist two different formulas $\phi, \psi \in \fml$ such that $\tau(\phi) = \tau(\psi)$. Therefore, if we consider the two substitutions $\sigma_\phi$ and $\sigma_\psi$ that send a variable $x$ respectively to $\phi$ and $\psi$, both acting as the identity on $\V \setminus \{x\}$, we have two different substitutions whose images under $\ov \tau$ coincide. Hence $\ov \tau$ is injective if and only if $\tau$ is injective.
\item[$(iii)$]The last property follows immediately from the definition of $\ov \tau$; we will denote by $\tau$ both the map between $\fml$ and $\fmll$ and its extension to the corresponding free sup-lattices.
\end{enumerate}
\end{proof}

Now that we know that any language translation induces a quantale homomorphism, as a next step, we want to find a converse property, namely we want a characterization of quantale morphisms induced by language translations among all the morphisms of quantales of type $\wp(\Sfm)$.

In order to do that, we first observe that both $\Sfm$ and $\Sfmm$ have a submonoid isomorphic to $\la \V^{\V}, \circ, \id_{\V} \ra$. Hence $\wp(\Sfm)$ and $\wp(\Sfmm)$ have two isomorphic subquantales that we will denote respectively by $\VV$ and $\VV'$; we shall denote by
\begin{equation}\label{mu}
\mu: \VV \lto \VV'
\end{equation}
the canonical isomorphism between them. Henceforth, by $\mu$, we will always mean the isomorphism in (\ref{mu}).

Another important remark concerns the completely join-prime elements. indeed, it is well-known that, if $X$ is a non-empty set, the completely join prime elements of the lattice $\la \wp(X), \cap, \cup \ra$ are precisely the atoms, i.e. the singletons. Therefore, a homomorphism between $\wp(\Sfm)$ and $\wp(\Sfmm)$ that preserves completely join-prime elements is necessarily the extension of a monoid homomorphism from $\Sfm$ to $\Sfmm$.

As we are going to show, quantale morphisms induced by language translations are precisely the translations, in the sense of Definition~\ref{i-trans}, relative to the isomorphism $\mu$. Before proving this result, let us introduce another notion and a lemma: if $B$ is a subset of a monoid $\AA$, an element $a \in A$ is called \emph{right $B$-absorbing} if $ba = a$ for all $b \in B$. The following result is trivial.

\begin{lemma}\label{abs}
Let $\AA$ and $\AA'$ be monoids, $B \subseteq A$ and $g: \AA \lto \AA'$ a monoid homomorphism. Then, for any right $B$-absorbing element $a \in A$, $g(a)$ is a right $g[B]$-absorbing element of $\AA'$.
\end{lemma}

\begin{theorem}\label{transchar}
Let $h: \wp(\Sfm) \lto \wp(\Sfmm)$ be a homomorphism. Then $h$ is induced by a language translation of $\lang$ into $\lang'$ if and only if it is a translation relatively to $\mu$.
\end{theorem}
\begin{proof}
One implication is trivial: if $h = \ov\tau$ for a translation $\tau: \fml \lto \fmll$, then it is a translation relatively to $\mu$ by Definition~\ref{ltranslation} and Lemma~\ref{qinq'}.

On the other hand, let us assume that $h$ is a $\mu$-translation. First of all, by Definition~\ref{i-trans} and the above remarks, $h$ is the natural extension to powersets of a monoid homomorphism (still denoted by $h$) from $\Sfm$ to $\Sfmm$.

Now recall that $\wp(\Fml)$ and $\wp(\Fmll)$ can be viewed equivalently as $\wp(\Sfm) \kappa_x$ and $\wp(\Sfmm) \kappa_x$ respectively, where $x$ is any fixed variable and $\kappa_x$ is the substitution that sends all variables to $x$. We set
$$\tau: \sigma \kappa_x \in \fml \lmapsto h(\sigma) \kappa_x \in \fmll,$$
and we shall prove that $\tau$ is a language translation and $h = \ov\tau$. Before we continue, we wish to underline that $\tau$ is independent from the fixed $x$, in the sense that, if $x \neq y$, $\tau': \sigma \kappa_y \in \fml \lmapsto h(\sigma) \kappa_y \in \fmll$, and $\sigma \kappa_x$ and $\sigma' \kappa_y$ yield the same formula of $\lang$, then $\tau(\sigma \kappa_x)$ and $\tau(\sigma' \kappa_y)$ correspond to the same formula of $\lang'$.

As a first step, assume that $c$ is a constant of $\lang$, and let $\sigma_c$ be a substitution that sends $x$ to $c$. For any substitution $\sigma \in \sfm$, $\sigma \sigma_c \kappa_x = \sigma_c \kappa_x$; we assume, by contradiction, that the formula $h(\sigma_c) \kappa_x$ contains a variable $y$. Let $\sigma$ and $\sigma'$ be the substitutions, in $\sfm$ and $\sfmm$ respectively, that send $y$ to another variable $z \neq y$ and fix all the other variables. Then $\sigma \in \VV$, $\sigma' \in \VV'$ and $h(\sigma) = \sigma'$; on the other hand $\tau(\sigma \sigma_c \kappa_x) = \sigma' h(\sigma_c) \kappa_x \neq h(\sigma_c) \kappa_x = \tau(\sigma_c \kappa_x)$, and this is absurd since $\sigma \sigma_c \kappa_x = \sigma_c \kappa_x$. Therefore $h(\sigma_c) \kappa_x$ cannot contain variables and $\lang'$ must have a definable constant.\footnote{We observe explicitly that, in order to define constants, a language must have at least a primitive constant.}

The case of unary connectives could be treated within the general case; nonetheless we prove it separately in order to give the reader a better clue of the argument.

Let $f$ be a unary connective of $\lang$ and $\sigma_f$ a substitution that sends $x$ to $f(x)$. Now let us consider the following subset of $\VV$:
$$B = \{\sigma \in V \mid \sigma(x) = x \textrm{ and } \sigma(y) \neq y, \forall y \in \V \setminus \{x\}\}.$$
Such a set is easily seen to be non-empty and it is clear that $\sigma_f \kappa_x$ (i.e. $f(x)$) is right $B$-absorbing, hence $\sigma \sigma_f \kappa_x = \sigma_f \kappa_x$ for all $\sigma \in B$. Then, if $B' = h[B]$, by Lemma~\ref{abs}, $\tau(\sigma_f \kappa_x) = h(\sigma_f) \kappa_x$ is right $B'$-absorbing, which means essentially that $\tau(\sigma_f \kappa_x)$ contains at most the unique variable $x$. On the other hand, if $\tau(\sigma_f \kappa_x)$ is a constant, then it is a multiplicative idempotent element of $\wp(\sfmm)$, while $\sigma_f \kappa_x \sigma_f \kappa_x \neq \sigma_f \kappa_x$, i.e. $\sigma_f \kappa_x$ is not a multiplicative idempotent. But this is impossible by condition $(ii)$ of Definition~\ref{i-trans}; so $\tau(\sigma_f \kappa_x)$ cannot be a constant, hence it is a formula in the single variable $x$.

Now let $f$ be a connective of arity $n > 1$ and $\sigma_f$ be the substitution that sends $x$ to $f(x_1, \ldots, x_n)$, with $x, x_1, \ldots, x_n$ distinct variables, and acts like the identity on $\V \setminus \{x\}$; let also $X = \V \setminus \{x_1, \ldots, x_n\}$ and consider the subset of $\VV$ defined as follows:
$$B = \{\sigma \in V \mid \sigma(x_i) = x_i, \forall i = 1, \ldots, n \textrm{ and } \sigma(y) \neq y, \forall y \in X\}.$$
$B$ is clearly non-empty and $\sigma_f \kappa_x$ (i.e. $f(x_1, \ldots, x_n)$) is right $B$-absorbing, hence $\sigma \sigma_f \kappa_x = \sigma_f \kappa_x$ for all $\sigma \in B$. As in the case of unary connectives, if $B' = h[B]$, by Lemma~\ref{abs}, $\tau(\sigma_f \kappa_x) = h(\sigma_f) \kappa_x$ is right $B'$-absorbing, which means essentially that $\tau(\sigma_f \kappa_x)$ contains at most the variables $x_1, \ldots, x_n$. Assuming that there exists $i \leq n$ such that $x_i$ is not in $\tau(\sigma_f \kappa_x)$, we can consider the  substitution $\alpha$ that sends $x_i$ to $x$ and acts as the identity on $\V \setminus \{x_i\}$. Then $\tau(\sigma_f \alpha \sigma_f \kappa_x) = h(\sigma_f \alpha \sigma_f) \kappa_x$ is multiplicative idempotent while $\sigma_f \alpha \sigma_f \kappa_x$ is not, since it is in fact the formula
$$f(x_1, \ldots, x_{i-1},f(x_1, \ldots, x_n), x_{i+1}, \ldots, x_n).$$
Again, this is impossible because $h$ is a $\mu$-translation, therefore $\tau(\sigma_f \kappa_x)$ contains precisely the variables $x_1, \ldots, x_n$.

Now we must prove that $\tau$ is a language translation and $h = \ov\tau$. Condition $(i)$ of Definition~\ref{ltranslation} is an obvious consequence of the fact that $h$ is a $\mu$-translation: any variable $y$ corresponds only to substitutions in $\VV$ (e.g. $\kappa_y$ or $\kappa_y \kappa_x$) that, under this hypothesis, are invariant w.r.t. $h$. Regarding Definition~\ref{ltranslation}$(ii)$ observe that, for any connective $f \in \lang$, $\tau(f(x_1, \ldots, x_n))$ is a formula $f'(x_1, \ldots, x_n)$ in the variables $x_1, \ldots, x_n$ and, therefore, for some $\sigma \in \sfm$,
$$\begin{array}{l}
\tau(f(\phi_1, \ldots, \phi_n)) = \tau(\sigma \sigma_f \kappa_x) = h(\sigma) h(\sigma_f) \kappa_x =\\
h(\sigma)(f'(x_1, \ldots, x_n)) = f'(h(\sigma)(x_1), \ldots, h(\sigma)(x_n)) = \\
f'(h(\sigma)\kappa_{x_1}, \ldots, h(\sigma)\kappa_{x_n}) = f'(\tau(\phi_1), \ldots, \tau(\phi_n)).
\end{array}$$

Last, in order to show that $h = \ov\tau$, we need to use again the fact that $h$ is a $\mu$-translation. Indeed, since substitutions are completely and univocally determined by their restriction to $\V$, we can represent any $\sigma \in \sfm$ by the family $\{\sigma \kappa_{x_i} \kappa_x\}_{i \in \N} = \{\sigma(x_i)\}_{i \in \N}$. Then $h(\sigma)$ is completely determined by
$$\{h(\sigma)(x_i)\}_{i \in \N} = \{h(\sigma \kappa_{x_i} \kappa_x)\}_{i \in \N} = \{h(\sigma) \kappa_{x_i} \kappa_x\}_{i \in \N} = \tau \circ \sigma_{\restr \V},$$
that is $h = \ov\tau$. The theorem is proved. 
\end{proof}

In the light of Theorem~\ref{transchar} we will call simply ``translations'' the $\mu$-translations in the case of concrete deductive systems.

In what follows we will often denote the \emph{domains} of deductive systems, i.e. sets of formulas, equations and sequents on which the consequence relation is defined, with the same letter, in italic character, of the respective system; so, for example, the domain of the system $\cat S = \la \lang, \vdash_\g \ra$ shall be denoted by $S$. Moreover, we will consider all of them as sets of sequents by identifying $\fml$ with the set of sequents $S$ with $\Tp(S) = \{(0,1)\}$ and $\Eq_\lang$ with the one such that $\Tp(S) = \{(1,1)\}$. We shall also borrow part of the terminology used in \cite{blokpigozzi2} and \cite{raftery}: for any given natural number $k$, by a \emph{$k$-formula} $\vec\phi$ we mean a sequence of $k$ formulas $(\phi_1, \ldots, \phi_k)$ and by a $k$-variable $\vec x$ we understand a sequence $(x_1, \ldots, x_k)$ of distinct propositional variables.

Hence a sequent of type $(m,n)$ will be denoted by $\vec\phi \To \vec\psi$, where $\vec\phi$ is an $m$-formula and $\vec\psi$ an $n$-formula. A single formula $\phi \in \fml$ shall be identified with the sequent $\varnothing \To \phi$.

\begin{definition}\label{faithful}
Let $\lang = \la L, \nu \ra$ and $\lang' = \la L', \nu' \ra$ be two propositional languages and $\cat S = \la \lang, \vdash_\g \ra$ and $\cat T = \la \lang', \vdash_\d \ra$ be two deductive systems. We say that $\cat S$ is \emph{(faithfully) interpretable} in $\cat T$ if there exist a language translation $\tau$ of $\lang$ in $\lang'$ and a residuated map (or, equivalently, a sup-lattice homomorphism) $\iota: \wp(S) \lto \wp (T)$ such that, for all $\Phi, \Psi \subseteq S$ and $\Sigma \subseteq \Sfm$,
\begin{enumerate}
\item[$(i)$]$\Phi \vdash_\g \Psi$ if and only if $\iota(\Phi) \vdash_\d \iota(\Psi)$,
\item[$(ii)$]$\iota(\Sigma \star \Phi) = \ov\tau(\Sigma) \star \iota(\Phi)$.
\end{enumerate}
The function $\iota$ is called the \emph{interpretation}.

If there exist two faithful interpretations $\iota: \cat S \lto \cat T$ and $\iota': \cat T \lto \cat S$ such that
\begin{equation}\label{equiv}
\Phi \dashv\vdash_\g(\iota' \circ \iota)(\Phi) \quad \textrm{ for all } \Phi \in \wp(S),
\end{equation}
then $\cat S$ and $\cat T$ are said to be \emph{equivalent}.  
\end{definition}

\begin{theorem}\label{faiththm}
Let $\cat S = \la \lang, \vdash\ra$ and $\cat T = \la \lang', \vdash_\d\ra$ be two propositional deductive systems. Then $\cat S$ is interpretable in $\cat T$ if and only if there exists a $\mu$-translation $\ov \tau: \wp(\Sfm) \lto \wp(\Sfmm)$ and an injective homomorphism of $\wp(\Sfm)$-modules $f: \wp(\SS)_\g \lto (\wp(\mathbf{T})_{\d})_{\ov\tau}$.
\end{theorem}
\begin{proof}
Assume that $\cat S$ is interpretable in $\cat T$. The existence of $\ov \tau$ has been proved in Lemma~\ref{qinq'}. So let
$$f: \wp(S)_\gamma \lto \wp(T)_{\d}$$
be defined as follows
\begin{equation}\label{f}
f(\g(\Phi)) = \d(\iota(\Phi)), \quad \textrm{for all $\Phi \in \wp(S)$}.
\end{equation}
If $\Phi, \Psi \in \wp(S)$ are such that $\g(\Phi) = \g(\Psi)$, then $\Phi \vdash_\g\Xi$ for all $\Xi \subseteq \g(\Psi)$, hence $\iota(\Phi) \vdash_\d \iota(\Xi)$ which means that $\d(\iota(\Psi)) \subseteq \d(\iota(\Phi))$. The converse inclusion can be proved analogously, so $\d(\iota(\Psi)) = \d(\iota(\Phi))$ and $f$ is a well defined function.

Now, in order to prove that $f$ is injective, let us consider $\Phi, \Psi \in \wp(S)$ such that $\g(\Phi) \neq \g(\Psi)$; we can assume, without losing generality, that there exists $\phi \in \g(\Phi) \setminus \g(\Psi)$. Then $\Psi \nvdash_\g \phi$ and this implies that $\iota(\Psi) \nvdash_\d \iota(\phi)$. It follows $f(\g(\Phi)) \neq f(\g(\Psi))$ and then $f$ is injective.

Let $\{\Phi_i\}_{i \in I} \subseteq \wp(S)$. We have
$$\begin{array}{l}
f\left({}^\g\bigvee_{i \in I}\g(\Phi_i)\right) = f\left(\g\left(\bigcup_{i \in I} \Phi_i\right)\right) \\
= \d\left(\iota\left(\bigcup_{i \in I} \Phi_i\right)\right) = \d\left(\bigcup_{i \in I}\iota(\Phi_i)\right) \\
= {}^{\d}\bigvee_{i \in I} \d\left(\iota(\Phi_i)\right) = {}^{\d}\bigvee_{i \in I} f(\g(\Phi_i)),
\end{array}$$
whence $f$ is a sup-lattice homomorphism.

Last we need to prove that $f$ is a $\wp(\Sfm)$-module homomorphism, so let $\Sigma \in \wp(\sfm)$ and $\Phi \in \wp(S)$. We have
$$\begin{array}{l}
f(\Sigma \star_\g \g(\Phi)) = f(\g(\Sigma \star \Phi)) = \d(\iota(\Sigma \star \Phi)) \\
= \d(\ov \tau(\Sigma) \star \iota(\Phi)) = \ov \tau(\Sigma) \star_{\d} \d(\iota(\Phi)) \\
= \ov \tau(\Sigma) \star_{\d} f(\Phi).
\end{array}$$

Conversely, let us assume the existence of $\ov\tau$ and $f$. By Theorem~\ref{transchar}, $\ov\tau$ is the natural extension of a language translation $\tau: \fml \lto \fmll$; on the other hand we have the following diagram of $\wp(\Sfm)$-module morphisms
\begin{equation}\label{concinterdiag}
\xymatrix{
\wp(\SS) \ar@{-->}[rr]^{\iota} \ar@{->>} [dd]_{\g} & & \wp(\TT)_{\ov\tau} \ar@{->>} [dd]^\d\\
&&\\
\wp(\SS)_\g \ar[rr]_{f} & & (\wp(\TT)_{\ov\tau})_\d,
\\
}
\end{equation}
that can be completed with a morphism $\iota$ because $\wp(SS)$ is a projective module. Moreover $\iota$ obviously satisfies Definition \ref{faithful}, and the assertion is proved.
\end{proof}

\begin{theorem}\label{equivthm}
Let $\cat S = \la \lang, \vdash_\g \ra$ and $\cat T = \la \lang', \vdash_\d\ra$ be two propositional deductive systems. If $\cat S$ and $\cat T$ are equivalent, then there exist two quantale homomorphisms $\ov \tau: \wp(\Sfm) \lto \wp(\Sfmm)$, $\ov \tau': \wp(\Sfmm) \lto \wp(\Sfm)$, and $\wp(\SS)_\g$ and $\wp(\mathbf{T})_{\d}$ are isomorphic both as $\wp(\Sfm)$-modules and as $\wp(\Sfmm)$-modules.
\end{theorem}
\begin{proof}
The existence of $\ov\tau$ and $\ov\tau'$ is an immediate consequence of Definition~\ref{faithful} and Lemma~\ref{qinq'}; moreover, by Theorem~\ref{faiththm}, we have an injective homomorphism of $\wp(\Sfm)$-modules $f: \wp(\SS)_\g \lto \wp(\TT)_{\d}$ and an injective homomorphism of $\wp(\Sfmm)$-modules $f': \wp(\TT)_{\d} \lto \wp(\SS)_\g$. By (\ref{equiv}), $\g(\Phi) = \g(\iota'(\iota(\Phi)))$ for all $\Phi \in \wp(S)$. Therefore, using the definition of $f$ and $f'$ given in (\ref{f}), it follows
$$\g(\Phi) = \g(\iota'(\iota(\Phi))) = f'(\d(\iota(\Phi))) = f'(f(\g(\Phi))),$$
for all $\Phi \in \wp(S)$, whence $f' \circ f = \id_{\wp(S)_\g}$. For a classical set-theoretic result, $f'$ is surjective and therefore it is a $\wp(\Sfmm)$-isomorphism; on the other hand, by the uniqueness of inverse fuctions, $f$ must be its inverse $\wp(\Sfmm)$-isomorphism. Since $f$ is also a $\wp(\Sfm)$-isomorphism, for the same reasons, $f'$ is its inverse $\wp(\Sfm)$-isomorphism, and the theorem is proved.
\end{proof}

\begin{definition}\label{noncons}
Let $\lang = \la L, \nu \ra$ and $\lang' = \la L', \nu' \ra$ be two propositional languages and $\cat S = \la \lang, \vdash_\g\ra$ and $\cat T = \la \lang', \vdash_\d \ra$ be two deductive systems and assume that there exist a translation $\tau$ of $\lang$ in $\lang'$ and a residuated map $\iota: \wp(S) \lto \wp(T)$ such that
\begin{enumerate}
\item[$(i)$]$\Phi \vdash_\g \Psi$ implies $\iota(\Phi) \vdash_\d \iota(\Psi)$,
\item[$(ii)$]$\iota(\Sigma \star \Phi) = \ov\tau(\Sigma) \star \iota(\Phi)$.
\end{enumerate}

for all $\Phi, \Psi \in \wp(S)$. Then $\iota$ is called a \emph{non-conservative interpretation} of $\cat S$ in $\cat T$.
\end{definition}

\begin{theorem}\label{nonconsthm}
Let $\cat S = \la \lang, \vdash\ra$ and $\cat T = \la \lang', \vdash_\d\ra$ be two propositional deductive systems. If there exists a non-conservative interpretation of $\cat S$ in $\cat T$ then there exist a quantale homomorphism $\ov \tau: \wp(\Sfm) \lto \wp(\Sfmm)$ and an homomorphism of $\wp(\Sfm)$-modules $f: \wp(\SS)_\g \lto \wp(\mathbf{\TT})_{\d}$.
\end{theorem}
\begin{proof}
The proof is analogous to the one of Theorem~\ref{faiththm} except from the fact that, in this case, $f$ is not necessarily injective since the ``if and only if'' in Definition~\ref{faithful}$(i)$ is replaced by a single implication in Definition \ref{noncons}(ii). Hence it is still possibile only to prove the existence of the morphism $f$.
\end{proof}

\begin{definition}\label{weak}
Let $\lang = \la L, \nu \ra$ and $\lang' = \la L', \nu' \ra$ be two propositional languages and $\cat S = \la \lang, \vdash_\g\ra$ and $\cat T = \la \lang', \vdash_\d \ra$ be two deductive systems. If there exists a translation $\tau$ of $\lang$ in $\lang'$ and a non-contradictory extension $\vdash_\e$ of $\vdash_\d$ (i.e. consequence relation that is stronger that $\vdash_\d$) such that $\cat S$ is interpretable via in $\cat T' = \la \lang', \vdash_\e\ra$, then $\cat S$ is said to be \emph{weakly interpretable} in $\cat T$.
\end{definition}

\begin{theorem}\label{weakthm}
Let $\cat S = \la \lang, \vdash\ra$ and $\cat \lang' = \la E, \vdash_\d\ra$ be two propositional deductive systems. If there exists a weak interpretation of $\cat S$ in $\cat T$ then there exist a quantale homomorphism $\ov \tau: \wp(\Sfm) \lto \wp(\Sfmm)$, a $\wp(\Sfmm)$-module structural closure operator $\e \geq \d$ on $\wp(\TT)$ and an injective homomorphism of $\wp(\Sfm)$-modules $f: \wp(\SS)_\g \lto \wp(\TT)_\e$.
\end{theorem}
\begin{proof}
It follows immediately from Definition~\ref{weak} and Theorem~\ref{faiththm}.
\end{proof}

\section{Non-structural interpretations}
\label{nonstrsec}

Recalling that $\fml$ and $\Eq_\lang$ are sets of sequents on $\lang$ such that $\Tp(\fml) = \{(0,1)\}$ and $\Tp(\Eq_\lang) = \{(1,1)\}$, in what follows, for a deductive system on a propositional language $\lang$, we will keep on using the notations introduced in the previous section, so we shall denote by $D$ its domain, i.e. the set of formulas, equations or sequents on which the consequence relation is defined.

\begin{definition}\label{nonstrinter}
Let $\lang$ and $\lang'$ be propositional languages, and $\cat S = \la \lang, \vdash_\g\ra$ and $\cat T = \la \lang', \vdash_\d \ra$ deductive systems.

We say that $\cat S$ is \emph{non-structurally interpretable} in $\cat T$ if there exists a residuated map $\iota$ of $\wp(S)$ into $\wp(T)$ such that, for all $\Phi, \Psi \in \wp(S)$,
\begin{equation}
\Phi \vdash_\g\Psi \quad \textrm{if and only if} \quad \iota(\Phi) \vdash_\d \iota(\Psi).
\end{equation}
\end{definition}

\begin{remark}\label{faithfulcase}
Obviously a faithful interpretation of $\cat S$ into $\cat T$ is also a non-structural one. In the rest of this section, in order to avoid repetitive and not interesting specifications, we will assume that the consequence relations involved are always non-trivial, i.e. their associated nuclei are different from the identity map.
\end{remark}

\begin{theorem}\label{nonstrinterchar}
Let $\cat S = \la \lang, \vdash_\g\ra$ and $\cat T = \la \lang', \vdash_\d \ra$ be two deductive systems, and $\g = \g_\vdash$, $\d = \g_{\vdash_\d}$ be the structural closure operators associated to $\vdash$ and $\vdash_\d$ respectively. Then $\cat S$ is non-structurally interpretable in $\cat T$ if and only if there exists an injective homomorphism of sup-lattices
\begin{equation}\label{nonstreq}
f: \wp(\SS)_\g \lto \wp(\TT)_{\d}.
\end{equation}
\end{theorem}
\begin{proof}
Let us denote still by $\iota: \wp(S) \lto \wp(T)$ the sup-lattice homomorphism univocally determined by $\iota: \phi \in S \lmapsto \{\iota(\phi)\} \in \wp(T)$. By hypothesis, for all $\Phi, \Psi \in \wp(S)$, $\Psi \subseteq \g(\Phi)$ iff $\iota(\Psi) \subseteq \iota(\g(\Phi))$; on the other hand, we have
$$\begin{array}{lll}
\iota(\g(\Phi)) & = & \iota(\{\Psi \in \wp(S) \mid \Phi \vdash_\g\Psi\}) \\
& = & \{\iota(\Psi) \in \iota(\wp(S)) \mid \iota(\Phi) \vdash_\d \iota(\Psi)\} \\
& \subseteq & \{\Psi' \in \wp(T) \mid \iota(\Phi) \vdash_\g\Psi'\} \\
& = & \d(\iota(\Phi)),
\end{array}$$
hence $\Psi \subseteq \g(\Phi)$ iff $\iota(\Psi) \subseteq \d(\iota(\Phi))$. Then
$$f: \ \Phi \in \wp(S)_\g \ \lmapsto \ \d(\iota(\Phi)) \in \wp(T)_{\d}$$
is a well-defined map and $f \circ \g = \d \circ \iota$. Indeed, we have the following commutative diagram
\begin{equation}\label{interdiag}
\xymatrix{
\wp(S) \ar[rr]^{\iota} \ar@{->>} [dd]_{\g} & & \wp(T) \ar@{->>} [dd]^{\d}\\
&&\\
\wp(S)_\g \ar@{-->} [rr]_{f} & & \wp(T)_{\d}
\\
}.
\end{equation}

If $\Phi, \Psi \in \wp(S)_\g$ and $\Phi \neq \Psi$, there exists $\Gamma \subseteq \Phi \setminus \Psi$ or $\Delta \subseteq \Psi \setminus \Phi$. We can assume the first case, without losing generality. From $\Phi \vdash_\g\Gamma$, it follows $\iota(\Gamma) \subseteq \d(\iota(\Phi)) = f(\Phi)$; analogously, from $\Psi \nvdash \Gamma$ it follows $\iota(\Gamma) \nsubseteq \d(\iota(\Psi)) = f(\Psi)$. Thus $\Phi \neq \Psi$ implies $f(\Phi) \neq f(\Psi)$, whence $f$ is injective. Moreover, if $\{\Phi_i\}_{i \in I} \subseteq \wp(S)_\g$, we have
\begin{eqnarray}
\lefteqn{f\left({}^\g\bigvee_{i \in I} \Phi_i\right) = f\left(\g\left(\bigcup_{i \in I} \Phi_i\right)\right)} \nonumber \\
&=& \d\left(\iota\left(\bigcup_{i \in I} \Phi_i \right)\right) = \d\left(\bigcup_{i \in I} \iota(\Phi_i)\right) \nonumber \\
&=& {}^{\d}\bigvee_{i \in I} \d(\iota(\Phi_i)) = {}^{\d}\bigvee_{i \in I} f(\g(\Phi_i)) \nonumber \\
&=& {}^{\d}\bigvee_{i \in I} f(\Phi_i)
\end{eqnarray}
hence $f$ is a sup-lattice injective morphism.

Conversely, assume that such a morphism exists. Then we have the following diagram of sup-lattice homomorphisms
\begin{equation}\label{interdiag2}
\xymatrix{
\wp(\SS) \ar@{->>} [dd]_\g \ar@{-->}[rr]^\iota & & \wp(\TT) \ar@{->>} [dd]_{\d}\\
&&\\
\wp(\SS)_\g \ar[rr]_{f} & & \wp(\TT)_{\d} \ ,\\
}
\end{equation}
that can be closed by a homomorphism of sup-lattices $\iota$ because $\wp(\SS)$ is free, thus projective by Proposition~\ref{projsl}. If $\Phi, \Psi \in \wp(S)$, we have
$$\begin{array}{ll}
& \Phi \vdash_\g\Psi \\
\liff & \g(\Psi) \subseteq \g(\Phi) \\
\liff & (f \circ \g)(\Psi) \subseteq (f \circ \g)(\Phi) \\
\liff & (\d \circ \iota)(\Psi) \subseteq (\d \circ \iota)(\Phi) \\
\liff & \iota(\Phi) \vdash_\d \iota(\Psi). \\
\end{array}$$
Note that the second equivalence above comes from the fact that $f$ is injective. 
\end{proof}

As in the case of structural interpretations, we can define the following variations of the concept of non-structural interpretation.
\begin{definition}\label{nswncinter}
Let $\lang$ and $\lang'$ be propositional languages, and $\cat S = \la \lang, \vdash_\g\ra$ and $\cat T = \la \lang', \vdash_\d \ra$ deductive systems.

A map $\iota$ of $S$ to $T$ is called a \emph{non-structural non-conservative interpretation} of $\cat S$ in $\cat T$ if, for all $\Phi, \Psi \in \wp(S)$,
\begin{equation}
\Phi \vdash_\g \Psi \quad \textrm{implies} \quad \iota[\Phi] \vdash_\d \iota[\Psi].
\end{equation}
It is called a \emph{non-structural weak interpretation} of $\cat S$ in $\cat T$ if there exists a non-contradictory extension $\vdash_\e$ of $\vdash_\d$ such that $\iota$ is a non-structural interpretation of $\cat S$ into $\cat T' = \la \lang', \vdash_\e\ra$.

$\cat S$ and $\cat T$ are said to be \emph{non-structurally equivalent}, or \emph{similar}, if there exist a non-structural interpretation $\iota: S \lto T$ and a non-structural interpretation $\iota': T \to S$ such that $\Phi \dashv\vdash_\g(\iota' \circ \iota)[\Phi]$ for all $\Phi \subseteq S$.
\end{definition}

\begin{corollary}\label{nsncwinterchar}
Let $\cat S = \la \lang, \vdash_\g \ra$ and $\cat T = \la \lang', \vdash_\d \ra$ be two deductive systems, and $\g$, $\d$ be the structural closure operators associated to $\vdash_\g$ and $\vdash_\d$ respectively. Then the following hold:
\begin{enumerate}
\item[$(i)$] $\cat S$ is non-structurally non-conservatively interpretable in $\cat T$ if and only if there exists a sup-lattice homomorphism $f: \wp(\SS)_\g \lto \wp(\TT)_{\d}$;
\item[$(ii)$] $\cat S$ is non-structurally weakly interpretable in $\cat T$ if and only if there exists a $\wp(\Sfmm)$-module structural closure operator $\e \geq \d$ on $\wp(\TT)$ and an injective homomorphism of sup-lattices $f: \wp(\SS)_\g \lto \wp(\TT)_{\e}$;
\item[$(iii)$] $\cat S$ and $\cat T$ are non-structurally equivalent if and only if there exists an isomorphism of sup-lattices $f: \wp(\SS)_\g \lto \wp(\TT)_{\d}$.
\end{enumerate}
\end{corollary}
\begin{proof}
The proof is analogous~---~mutatis mutandis~---~to the one of Theorem~\ref{nonstrinterchar}.
\end{proof}

\section*{References and further readings}
\addcontentsline{toc}{section}{References and further readings}

In this chapter we have essentially added some contributions to the newborn algebraic and categorical theory of consequence relations, presented in~\cite{galatostsinakis}. Since, as we underlined several times, it is the first work where such techniques have been applied to consequence relations and deductive systems, it is impossible to suggest works with the same approach to deductive systems, since they (probably) do not exist.

Works that may be useful for studying in depth (or extend further) the results of this chapter are those~--- some of which have been already cited~--- by W. J. Blok and B. J\'onsson~\cite{blokjonsson1,blokjonsson2}, W. J. Blok and D. Pigozzi~\cite{blokpigozzi,blokpigozzi2}, D'Ottaviano and Feitosa \cite{dott,feitosa}, N. Galatos and H. Ono~\cite{galatosono,galatosono2}, N. Galatos, P. Jipsen, T. Kowalski and H. Ono~\cite{nickbook}, A. Pynko~\cite{pynko}, J. Raftery~\cite{raftery}, J. Rebagliato and V. Verd\'u~\cite{rebver}.


\part{$\Q$-module Transforms in Image Processing}
\label{applpart}

\chapter{Fuzzy Image Compression and Mathematical Morphology}
\label{imagechap}

In this chapter we will show how certain techniques of image processing, even having different scopes, can be grouped together under the common ``algebraic roof'' of $\Q$-module transforms.

The theory of \emph{fuzzy relation equations}\index{Fuzzy!-- relation equation},~\cite{dinolasessa}, is widely used in many applications and particularly in the field of image processing (see, for example,~\cite{hirotanobuhara,nobuharapedrycz1,hirotapedrycz}). As a matter of fact, fuzzy relations fit the problem of processing the representation of an image as a matrix with the range of its elements previously normalized to $[0,1]$. In such techniques, however, the approach is mainly experimental and the algebraic context is seldom clearly defined.

A first unification of fuzzy image processing has been proposed by I. Perfilieva in~\cite{perfilieva}, with an approach that is analytical rather than algebraic. Moreover, the field of applications of the operators (called \emph{Fuzzy transforms}\index{Fuzzy!-- transform}\index{Transform!Fuzzy --}) defined in~\cite{perfilieva} is limited to the real unit interval, $[0,1]$, endowed with the usual order relation and a triangular norm.

Actually, most of the fuzzy algorithms of image processing, make use of join-product operators, and they can be seen as approximate discrete solutions of fuzzy relation equations of the form $A(x,z) = \bigvee_y B(x,y) \cdot C(y,z)$. So it is natural to think of them as examples of $\Q$-module transforms. Indeed we will see in Section~\ref{ordapprsec} that the class of $\Q$-module transforms contains all these operators and much more. 

Further classes of operators that fall within $\Q$-module transforms are those of \emph{mathematical morphological operators}. Mathematical morphology is a technique for image processing and analysis whose birth can be traced back to the book~\cite{matheron}, of 1975, by G. Matheron, and whose establishment is due mainly to the works by J. Serra and H. J. A. M. Heijmans.

Essentially, mathematical morphological operators analyse the objects in an image by ``probing'' them with a small geometric ``model-shape'' (e.g., line segment, disc, square) called the \emph{structuring element}. These operators are defined on spaces having both a complete lattice order (set inclusion, in concrete applications) and an external action from another ordered structure (the set of translations); they are also usually coupled in adjoint pairs. A description of such operators in terms of $\Q$-module transforms can easily be anticipated.

We will begin the chapter by recalling, in Section~\ref{tnormsec}, some definitions and basic facts on triangular norms. Afterwards, rather than dwelling upon technical details, we will try to give the basic ideas of how fuzzy transforms and mathematical morphological operators work, respectively in Section~\ref{fuzalgsec} and Section~\ref{mathmorphsec}. Last, in Section~\ref{ordapprsec}, we will see how $\Q$-module transforms suffice to describe all those techniques.

\section{Left-continuous t-norms and their residua}
\label{tnormsec}

A binary operation $\ast: [0,1]^2 \lto [0,1]$ is called a \emph{triangular norm}\index{Triangular norm}\index{T-norm}, \emph{t-norm} for short, provided it verifies the following conditions
\begin{enumerate}
\item[-] commutativity: $x \ast y = y \ast x$;
\item[-] monotonicity: $x \ast y \leq z \ast y$ if $x \leq z$ and $x \ast y \leq x \ast z$ if $y \leq z$;
\item[-] associativity: $x \ast (y \ast z) = (x \ast y) \ast z$;
\item[-] $1$ is the neutral element: $1 \ast x = x = x \ast 1$.
\end{enumerate}
A t-norm $\ast$ is called \emph{left-continuous}\index{Triangular norm!Left-continuous --}\index{T-norm!Left-continuous --} if, for all $\{x_n\}_{n \in \N}, \{y_n\}_{n \in \N} \in [0,1]^\N$,
$$\left(\bigvee_{n \in \N} x_n\right) \ast \left(\bigvee_{n \in \N} y_n\right) = \bigvee_{n \in \N} (x_n \ast y_n).$$
In this case, clearly, $\ast$ is biresiduated and its residuum (unique, since $\ast$ is commutative) is given by
$$x \to y = \bigvee\{z \in [0,1] \mid z \ast x \leq y\}.$$

Although t-norms are the fuzzy logical analogous of the conjunction of classical logic, here we are mainly interested to them as algebraic operations. The defining conditions of t-norms are exactly the same that define a partially ordered Abelian monoid on the real unit interval $[0,1]$. Therefore some authors call t-norm also the monoidal operation of any partially ordered Abelian monoid; then, in this case, the concept of left-continuity can be substituted by the requirement that the Abelian po-monoid is actually a commutative residuated lattice.

\begin{example}\label{tnorms}
Here we list the best known examples of t-norm
\begin{enumerate}
\item[]\emph{Minimum}, or \emph{G\"odel} t-norm\index{Triangular norm!G\"odel --}\index{T-norm!G\"odel --}: \quad $x \wedge y = \min\{x,y\}$. It is the standard semantics for conjunction in \emph{G\"odel fuzzy logic}\index{Fuzzy!G\"odel -- logic}. Besides that, it occurs in most t-norm based fuzzy logics as the standard semantics for the so-called \emph{weak conjunction}. Last, it is the largest t-norm in the sense that, for all $x, y \in [0,1]$ and for any t-norm $\ast$, $x \ast y \leq x \wedge y$. Its residuum is $x \to_\wedge y = \left\{\begin{array}{ll} y & \textrm{if $x > y$} \\ 1 & \textrm{if $x \leq y$} \end{array}\right.$.
\item[]
\item[]\emph{Product} t-norm\index{Triangular norm!Product --}\index{T-norm!Product --}: \quad $x \cdot y$ (the ordinary product of real numbers). It is the standard semantics for strong conjunction in \emph{product fuzzy logic}\index{Fuzzy!Product -- logic}. Its residuum is defined as $x \to_\cdot y = \left\{\begin{array}{ll} y / x & \textrm{if $x > y$} \\ 1 & \textrm{if $x \leq y$} \end{array}\right.$.
\item[]
\item[]\emph{\L ukasiewicz} t-norm\index{Triangular norm!\L ukasiewicz --}\index{T-norm!\L ukasiewicz --}: \quad $x \odot y = \max\{0, x + y - 1\}$. The name comes from the fact that the t-norm is the standard semantics for strong conjunction in \L ukasiewicz fuzzy logic; it is smaller than the product t-norm. Its residuum is $x \to_{\luk} y = \min\{1, 1 - x + y\}$.
\item[]
\item[]\emph{Generalized \L ukasiewicz} t-norm\index{Triangular norm!Generalized \L ukasiewicz --}\index{T-norm!Generalized \L ukasiewicz --}: \quad $x \odot_p y = \sqrt[p]{\max\{0, x^p + y^p - 1\}}$, where $p$ is a fixed natural number. Its residuum is given by $x \to_{\luk_p} y = \min\{1, \sqrt[p]{1 - x^p + y^p}\}$.
\item[]
\item[]\emph{Nilpotent minimum} t-norm\index{Triangular norm!Nilpotent minimum --}\index{T-norm!Nilpotent minimum --}: \quad $x \bullet y = \left\{\begin{array}{ll} \min(x,y) & \textrm{if $x + y > 1$} \\ 0 & \textrm{otherwise} \end{array}\right.$. It is a standard example of a t-norm which is left-continuous, but not continuous, and its residuum is $x \to_\bullet y = \left\{\begin{array}{ll} \max\{1 - x, y\} & \textrm{if $x > y$} \\ 1 & \textrm{if $x \leq y$}\end{array}\right.$.
\end{enumerate}
\end{example}

\section{Fuzzy algorithms for image compression and reconstruction}
\label{fuzalgsec}

In the literature of image compression, the fuzzy approach is based essentially on the theory of fuzzy relation equations, deeply investigated by A. Di Nola, S. Sessa, W. Pedrycz and E. Sanchez in~\cite{dinolasessa}. The underlying idea is very easy: a grey-scale image is basically a matrix in which every element represents a pixel and its value, included in the set $\{0, \ldots, 255\}$ in the case of a 256-bit encoding, is the ``grey-level'', where $0$ corresponds to black, $255$ to white and the other levels are, obviously, as lighter as they are closer to $255$. Then, if we normalize the set $\{0, \ldots, 255\}$ by dividing each element by $255$, grey-scale images can be modeled equivalently as fuzzy relations, fuzzy functions (i.e. $[0,1]$-valued maps) or fuzzy subsets of a given set.

As we anticipated, we will neither cover the wide literature on this subject, nor show how such techniques have been developed in the last years (also because it would be a thankless task). Here we rather want to point out the connection with our work, and the best way to show it is to present the first attempt of unifying all (or most of) these techniques in a common algebraic framework, namely the \emph{fuzzy transforms expressed by residuated lattice operations}, introduced by I. Perfilieva in~\cite{perfilieva}.

By a \emph{fuzzy partition}\index{Fuzzy!-- partition} of the real unit interval $[0,1]$, we mean a $n$-tuple of fuzzy subsets $A_1, \ldots, A_n$, with $n \geq 2$, identified with their membership functions $A_i: [0,1] \lto [0,1]$ satisfying the following \emph{covering property}
\begin{equation}\label{covering}
\textrm{for all $x \in [0,1]$ there exists $i \leq n$ such that} \quad A_i(x) > 0.
\end{equation}
The membership functions $A_1, \ldots, A_n$ are called the \emph{basic functions} of the partition. There is assumed to exist a a finite subset $P \subset [0,1]$, consisting of \emph{nodes} $p_1, \ldots, p_l$ where $l$ is a sufficiently large natural number. Moreover, we assume that $P$ is \emph{sufficiently dense} with respect to the fixed partition, i.e.
\begin{equation}\label{suffdense}
\textrm{for all $i \leq n$ there exists $j \leq l$ such that $A_i(p_j) > 0$}.
\end{equation}
\begin{definition}\label{ftransf}
Let $f \in [0,1]^P$ and $A_1, \ldots, A_n$, $n < l$, be basic functions of a fuzzy partition of $[0,1]$. We say that the $n$-tuple $(F_1^\ua, \ldots, F_n^\ua)$ is the $F^\ua$-transform of $f$ with respect to $A_1, \ldots, A_n$ if, for all $k \leq n$,
\begin{equation}\label{fuptransf}
F^\ua_k = \bigvee_{j=1}^l (A_k(p_j) \ast f(p_j)).
\end{equation}
We say that the $n$-tuple $(F_1^\da, \ldots, F_n^\da)$ is the $F^\da$-transform of $f$ with respect to $A_1, \ldots, A_n$ if, for all $k \leq n$,
\begin{equation}\label{fdowntransf}
F^\da_k = \bigvee_{j=1}^l (A_k(p_j) \to_* f(p_j)).
\end{equation}
\end{definition}
\begin{definition}\label{finvtransf}
Let $f \in [0,1]^P$, $A_1, \ldots, A_n$, with $n < l$, be basic functions of a fuzzy partition of $[0,1]$, and $(F_1^\ua, \ldots, F_n^\ua)$ be the $F^\ua$-transform of $f$ with respect to $A_1, \ldots, A_n$ if, for all $k \leq n$. The map defined, for all $j \leq l$, by
\begin{equation}\label{fupinvtransf}
f^\ua(p_j) = \bigwedge_{k=1}^n (A_k(p_j) \to_* F^\ua_k)
\end{equation}
is called the \emph{inverse $F^\ua$-transform} of $f$.

Let $(F_1^\da, \ldots, F_n^\da)$ be the $F^\da$-transform of $f$ with respect to $A_1, \ldots, A_n$ if, for all $k \leq n$. The map defined, for all $j \leq l$, by
\begin{equation}\label{fdainvtransf}
f^\da(p_j) = \bigvee_{k=1}^n (A_k(p_j) \ast F^\da_k)
\end{equation}
is called the \emph{inverse $F^\da$-transform} of $f$.
\end{definition}

Apart from the definitions above, several results on such tranforms are presented in the cited paper. We will not list them here because they are special cases of more general results that we presented in Sections~\ref{resmapsec}~and~\ref{mqtransec}.

\section{Dilation and erosion in mathematical morphology}
\label{mathmorphsec}

In~\cite{goutsiasheijmans}, the authors state (quoted verbatim):
\begin{quote}
\small{The basic problem in mathematical morphology is to design nonlinear operators that extract relevant topological or geometric information from images. This requires development of a mathematical model for images and a rigorous theory that describes fundamental properties of the desirable image operators.}
\end{quote}

Then, not surprisingly, images are modeled, in the wake of tradition and intuition, as subspaces or subsets of a suitable space $E$, which is assumed to possess some additional structure (topological space, metric space, graph, etc.), usually depending on the kind of task at hand. We have seen that, in the case of digital image compression, the image space is often modeled as the set of all the functions from a set~--- the set of all the pixels~--- to the real unit interval $[0,1]$. Then, depending on several ``experimental'' factors, the properties of $[0,1]$ involved may be the usual operations, the order relation, t-norms and so on.

In mathematical morphology, the family of binary images is given by $\wp(E)$, where $E$ is, in general, $\R^n$ or $\Z^n$, for some $n \in \N$. In the first case we have continuous binary images, otherwise we are dealing with discrete binary images. The basic relations and operations between images of this type are essentially those between sets, namely set inclusions, unions, or intersections. As a first example, we can consider an image $X$ that is \emph{hidden} by another image $Y$. Then we can formalize this fact by means of set inclusion: $X \subseteq Y$. Analogously, if we simultaneously consider two images $X$ and $Y$, what we see is their union $X \cup Y$; the \emph{background} of an image $X$ is its complement $X^c$ in the whole space, and the part of an image $Y$ that is not covered by another image $X$ is the set difference $Y \setminus X = Y \cap X^c$.

It is easily anticipated, then, that the lattices are the algebraic structures required for abstracting the ideas introduced so far. Nonetheless, keeping in mind the models $\R^n$ and $\Z^n$, it is possible to introduce the concepts of \emph{translation}\index{Translation!-- of an image} of an image and \emph{translation invariance} of an operator, by means of the algebraic operation of sum.

The reader may recognize the following definitions as those of residuated map and its residual, and of adjoint pair.
\begin{definition}\label{dilero}
Let $\LL$, $\MM$ be complete lattices.  A map $\d: L \lto M$ is called a \emph{dilation}\index{Dilation} if it distributes over arbitrary joins, i.e., if $\d\left({}^\LL\bigvee_{i \in I} x_i\right) = {}^\MM\bigvee_{i \in I} \d(x_i)$, for every family $\{x_i\}_{i \in I} \subseteq L$. A map $\e: M \lto L$ is called an \emph{erosion}\index{Erosion} if it distributes over arbitrary meets, i.e., if $\e\left({}^\MM\bigwedge_{i \in I} y_i\right) = {}^\LL\bigwedge_{i \in I} \e(y_i)$, for every family $\{y_i\}_{i \in I}$ of elements of $M$.

Two maps $\d: L \lto M$ and $\e: M \lto L$ are said to form an \emph{adjunction}\index{Adjunction}, $(\d,\e)$, between $\LL$ and $\MM$ if $\d(x) \leq y \liff x \leq \e(y)$, for all $x \in L$ and $y \in M$.
\end{definition}
Notice that the notation used in mathematical morphology is slightly different. Indeed, an adjoint pair is presented with the residuated map in the second coordinate and its residual in the first. Here, in order to avoid confusion, we keep on using the notations introduced in Chapter~\ref{orderchap}. So we may reformulate the definition above by considering the sup-lattice reducts of $\LL$ and $\MM$, and saying that $\d: L \lto M$ is a dilation if it is a sup-lattice homomorphism between $\LL$ and $\MM$. Dually, an erosion $\e: M \lto L$ is a sup-lattice homomorphism between $\MM^{\op}$ and $\LL^{\op}$. Then a dilation $\d$ and an erosion $\e$ form an adjunction if $\e = \d_*$.

Assume that $\d: \LL \lto \MM$ is a dilation. For $x \in L$, we can write
\begin{equation}\label{dilation}
\d(x) = \bigvee_{y \leq x} \d(y),
\end{equation}
where we have used the fact that $\d$ distributes over join. Every dilation defined on $\LL$ is of the form (\ref{dilation}), and the adjoint erosion is
given by
\begin{equation}\label{erosion}
\e(y) = \bigvee_{\d(x) \leq y} x.
\end{equation}
In the case of powersets, if $\d$ is a dilation between $\wp(E)$ and $\wp(F)$, where $E$ and $F$ are nonempty sets. For $X \subseteq E$, we can write
\begin{equation}\label{wpdilation}
\d(X) = \bigcup_{x \in X} \d(\{x\}),
\end{equation}
and the adjoint erosion is, for all $Y \subseteq F$,
\begin{equation}\label{wperosion}
\e(Y) = \{x \in E \mid \d(\{x\}) \subseteq Y\} = \bigcup_{\d(X) \subseteq Y} X.
\end{equation}

Next, we consider the special case when the operators are translation invariant. In this case, the sets $\d(\{x\})$ are translates of a fixed set, called the \emph{structuring element}\index{Structuring element}, by $\{x\}$. Let $E$ be $\R^n$ or $\Z^n$, and consider the complete lattice $\wp(E)$; given an element $h \in E$, we define the $h$-\emph{translation} $\tau_h$ on $\wp(E)$ by setting, for all $X \in \wp(E)$,
\begin{equation}\label{translation}
\tau_h(X) = X + h = \{x + h \mid x \in X\},
\end{equation}
where the sum is intended to be defined coordinatewise.

An operator $f: \wp(E) \lto \wp(E)$ is called \emph{translation invariant}\index{Translation!-- invariant operator}\index{Operator!Translation invariant --}, \emph{T-invariant} for short, if $\tau_h \circ f = f \circ \tau_h$ for all $h \in E$. It can be proved that every T-invariant dilation on $\wp(E)$ is given by
\begin{equation}\label{tdilation}
\d_A(X) = \bigcup_{x \in X} A + x,
\end{equation}
and every T-invariant erosion is given by
\begin{equation}\label{terosion}
\e_A(X) = \{y \in E \mid A + y \subseteq X\} = \{y \in E \mid y \in X + \breve A\},
\end{equation}
where $A$ is an element of $\wp(E)$, called the \emph{structuring element}, and $\breve A = \{- a \mid a \in A\}$ is the reflection of $A$ around the origin.

Now we observe that the expressions for erosion and dilation in (\ref{tdilation}) and (\ref{terosion}) can also be written, respectively, as
\begin{equation}\label{tdilation2}
\d_A(X)(y) = \bigvee_{x \in E} A(y - x) \wedge X(x)
\end{equation}
and
\begin{equation}\label{terosion2}
\e_A(Y)(x) = \bigwedge_{y \in E} A(y - x) \to Y(y),
\end{equation}
where each subset $X$ of $E$ is identified with its membership function
$$X: x \in E \lmapsto \left\{
\begin{array}{ll}
1 & \textrm{if } x \in X \\
0 & \textrm{if } x \in X^c,
\end{array}\right. \in \{0, 1\}$$
and $X \to Y =: X^c \vee Y$.
Moving from these expressions, and recalling that $\wedge$ is a biresiduated commutative operation (that is, a continuous t-norm) whose residuum is $\to$, it is possible to extend these operations from the complete lattice of sets $\wp(E) = \{0,1\}^E$ to the complete lattice of fuzzy sets $[0,1]^E$, by means of continuous t-norms and their residua. What we do, concretely, is extending the morphological image operators of dilation and erosion, from the case of binary images, to the case of grey-scale images.

So let $\ast$ be a continuous t-norm and $\to$ be its residuum; a grey-scale image $X$ is a fuzzy subset of $E$, namely a map $X: E \lto [0,1]$. Given a fuzzy subset $A \in [0,1]^E$, called a \emph{fuzzy structuring element}\index{Structuring element!Fuzzy --}\index{Fuzzy!-- structuring element}, the operator
\begin{equation}\label{fdilation}
\d_A(X)(y) = \bigvee_{x \in E} A(y - x) \ast X(x)
\end{equation}
is a translation invariant dilation on $[0,1]^E$, and the operator
\begin{equation}\label{ferosion}
\e_A(X)(x) = \bigwedge_{y \in E} A(y - x) \to X(y)
\end{equation}
is a translation invariant erosion on $[0,1]^E$.

\begin{figure}[htp]
\centering
	\includegraphics[width=1.5in]{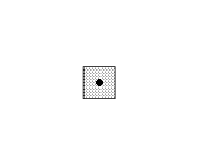}\caption{Structuring element: $A$} \label{strfig}
\end{figure}

\begin{figure}[htp]
\centering
	\includegraphics[width=1.5in]{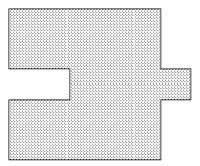} \caption{Original shape: $X$} \label{originalfig}
\end{figure}

\begin{figure}[htp]
\centering
	\includegraphics[width=1.5in]{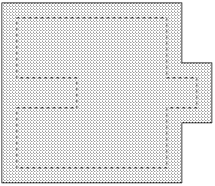} \caption{Dilation: $\d_A(X)$} \label{dilfig}
\end{figure}

\begin{figure}[htp]
\centering
	\includegraphics[width=1.5in]{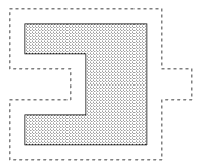} \caption{Erosion: $\e_A(X)$} \label{erofig}
\end{figure}

\begin{figure}[htp]
\centering
	\includegraphics[width=1.5in]{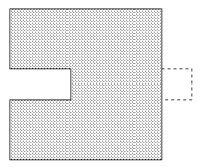} \caption{Opening: $(\d_A \circ \e_A)(X)$} \label{opefig}
\end{figure}

\begin{figure}[htp]
\centering
	\includegraphics[width=1.5in]{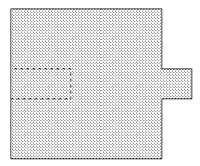} \caption{Closing: $(\e_A \circ \d_A)(X)$} \label{clofig}
\end{figure}

\begin{figure}[htp]
\centering
	\includegraphics[width=1.5in]{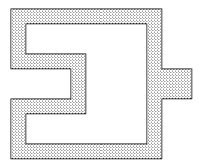} \caption{Outlining: $X - \e_A(X)$} \label{outfig}
\end{figure}

Combining the operators of dilation and erosion by means of the usual algebraic operations in $[0,1]$ it is possible to define new operators, e.g. \emph{outlining} and \emph{top-hat transform}. Their treatment is beyond the scope of this thesis, hence we will not present them in details; however we show one of them among the examples\footnote{Figs.~\ref{strfig}--\ref{outfig} have been taken and edited from~\cite{sito3}.} in Figs.~\ref{strfig}--\ref{outfig}.

Some major references in the area of mathematical morphology are~\cite{goutsiasheijmans,heijmans1,heijmans,heijmansronse,serra2}, besides the aforementioned~\cite{serra1}.

\section{A unified approach by $\Q$-module transforms}
\label{ordapprsec}

The operators defined in the previous sections have a familiar form. Indeed they are all special cases of $\Q$-module transforms between free modules over the quantale reducts of residuated lattice structures defined on the real unit interval $[0,1]$. We will now analyse them in detail.

Let us consider the $F^\ua$ transforms of Definition~\ref{ftransf}. Its domain is $[0,1]^l$ and its codomain is $[0,1]^n$ with $n \leq l$. We get immediately that a $\Q$-module transform
$$H_k: f \in [0,1]^l \lmapsto \bigvee_{j = 1}^l f(j) \ast k(j,{}_-) \in [0,1]^n$$
is an $F^\ua$ transform iff the kernel $k$ satisfies condition (\ref{suffdense}) rewritten as
\begin{equation}\label{qsuffdense}
\textrm{for all $i \leq n$ there exists $j \leq l$ such that $k(j,i) > 0$}.
\end{equation}
Obviously, the inverse $F^\ua$ transform of $H_k$ is right
$$\Lambda_k: g \in [0,1]^n \lmapsto \bigwedge_{i = 1}^n k({}_-,i) \to_* g(i) \in [0,1]^l,$$
i.e. the inverse $\Q$-module transform of $H_k$. The case of $F^\da$ transforms is dual to that of $F^\ua$, in the sense that the direct $F^\da$ transform is an inverse $\Q$-module transform, thus a homomorphism between the duals of free modules, and the inverse transform has the shape of a $\Q$-module transform. In other words, for $F^\da$ transforms we assume $l \leq n$ and condition 
\begin{equation}\label{qsuffdense2}
\textrm{for all $j \leq l$ there exists $i \leq n$ such that $k(j,i) > 0$};
\end{equation}
then $\Lambda_k$ above is the direct $F^\da$ transform and $H_k$ is its inverse.

We have already observed in Section~\ref{mathmorphsec} that the dilations are precisely the sup-lattice homomorphisms, while the erosions are their residua. In order to faithfully represent dilations and erosions that are translation invariant as $\Q$-module transforms from a free $[0,1]$-module to itself, we make the further assumption that the set over which the free module is defined has the additional structure of Abelian group.

So, let $\mathbf{X} = \la X, +, -, 0 \ra$ be an Abelian group, $\ast$ a t-norm on $[0,1]$, and consider the free $[0,1]$-module $[0,1]^X$. For any element $k \in [0,1]^X$, we define the two variable map $\ov k: (x,y) \in X \times X \lmapsto k(y - x) \in [0,1]$. Then, for all $k \in [0,1]^X$, the translation invariant dilation, on $[0,1]^X$, whose structuring element is $k$, is precisely the $\Q$-module transform $H_{\ov k}$, with the kernel $\ov k$ defined above. Obviously, the translation invariant erosion whose structuring element is $k$ is $\Lambda_{\ov k}$.

Then the representation of both fuzzy transforms and pairs dilation--erosion as quantale module transforms is trivial. Actually, what we want to point out here is that, if we drop the assumption that our quantale is defined on $[0,1]$, the classes of transforms defined in this section become much wider. The purpose of this consideration is not to suggest a purely speculative abstractions but, rather, to underline that suitable generalizations of these operators exist already and they may be useful provided their underlying ideas are extended to other kind of tasks. Indeed the aim of fuzzy transforms is to approximate maps that take values in $[0,1]$; hence the area of application of the whole class of $\Q$-module transform, as approximating operators, can be easily enlarged. On the contrary, the idea of dilating and eroding a shape, in order to analyse it, has not yet found an appropriate concrete extension to situations where $[0,1]$ must be replaced by a non-integral quantale. Nonetheless, we strongly believe (and we are working in this direction) that $\Q$-module dilations and erosions will soon find concrete tasks for being fruitfully applied.

\section*{References and further readings}
\addcontentsline{toc}{section}{References and further readings}

The literature on fuzzy image processing is extremely wide. Here we list, for example,~\cite{nobuharapedrycz1,hirotapedrycz,kwakpedrycz,nobuharapedrycz2,nobuharapedrycz3} by W. Pedrycz and others,~\cite{loiasessa} by V. Loia and S. Sessa,~\cite{nachtegael} by M. Nachtegael, D. Van der Weken, D. Van De Ville, E. E. Kerre, W. Philips and I. Lemahieu,~\cite{nobuharatakama} by H. Nobuhara, Y. Takama, K. Hirota, and~\cite{perfilieva} by I. Perfilieva.

An algebraic approach to fuzzy image processing, by means of semimodules over the semiring reducts of the MV-algebra $[0,1]$, together with the algorithm presented in Chapter~\ref{ltbchap}, can be found in~\cite{dinolarusso}, written with A. Di Nola.

The subject of fuzzy relation equations and their applications is deeply investigated in the book~\cite{dinolasessa}, by A. Di Nola A., S. Sessa, W. Pedrycz and E. Sanchez.

Regarding mathematical morphology, we point out the wide production of H. J. A. M. Heijmans on this topic. For example, useful introductory works are the book~\cite{heijmans} and the papers~\cite{goutsiasheijmans}, with J. Goutsias,~\cite{heijmansronse}, with C. Ronse, and~\cite{heijmans}; they have been our main references in drawing up Section~\ref{mathmorphsec}. Then, besides the books by G. Matheron,~\cite{matheron}, and J. Serra,~\cite{serra1}, already cited, we suggest the collection~\cite{serra2}.

\chapter[The \L TB Algorithm for Image Compression and Reconstruction]{An Example: the \L TB Algorithm for Image Compression and Reconstruction}
\label{ltbchap}

So far we have established a great amount of theoretic tools, but a reader interested more in the applications to image processing would feel a sense of dissatisfaction if we would not present a concrete application of such tools. Well, it is the content of this chapter.

In Section~\ref{luktransec} we present an example of $\Q$-module orthonormal transform, defined between modules over the quantale reduct of the MV-algebra $[0,1]$. The orthonormal coder traces its origin back to the construction of normal forms for \L ukasiewicz logic, proposed by A. Di Nola and A. Lettieri in~\cite{dinolalettieri}.

The description of the \L TB~--- \L ukasiewicz Transform Based~--- algorithm for image compression and reconstruction, is the content of Section~\ref{ltbsec}.

What we show here could seem, at a first glance, yet another fuzzy algorithm for image compression and reconstruction. Nonetheless, although also the concrete results are rather promising, what we want to underline, with this example, is the connection between some algebraic results showed in Chapter~\ref{mqchap} and certain properties of the \L TB algorithm. At the end of Section~\ref{applltbsec}, we will show several consequences of Theorem~\ref{wadjointpair} and of Theorem~\ref{sadjointpair}, i.e. of the fact that the \L ukasiewicz transform is orthonormal.

Last, in Section~\ref{ltbjpgsec}, we show and comment a comparison between \L TB and JPEG. The comparison is based both on the respective computational costs and on several efficiency indices, generally used for this kind of confrontations.

\section{The \L ukasiewicz Transform}
\label{luktransec}

In this section we will introduce the \L ukasiewicz transform as a $\Q$-module orthonormal transform between free $\Q$-modules over the quantale reduct $[0,1]_{\luk} = \la [0,1], \vee, \odot, 0, 1 \ra$ of the MV-algebra $\la [0,1], \oplus, \lnot, 0 \ra$. Even though the \L ukasiewicz transform has been introduced, in~\cite{dinolarusso}, mainly with the aim of building an algorithm for compression and reconstruction of digital images, the way it is defined has a strong logical motivation. In~\cite{dinolalettieri}, the authors propose a normal form for formulas of \L ukasiewicz propositional logic, making use of formulas having the property of being canonically associated to the so-called ``simple McNaughton functions''. The coder that determines the \L ukasiewicz transform is a map in two variables defined in such a way that, if we fix the second one, we obtain a map in a single variable that is either a simple McNaughton function or the pointwise meet of one positive and one negative simple McNaughton function. Then, in order to make the genesis of \L ukasiewicz transform clear, we need to recall some definitions and results on MV-algebras and the aforementioned normal forms in \L ukasiewicz logic. 

In the category~--- $\ell\cat G^\textrm{ab}_u$\index{Category!-- $\ell\cat G^\textrm{ab}_u$, of $\ell u$-groups}~--- of \emph{lattice ordered Abelian groups with a strong unit} (\emph{$\ell u$-groups}, for short)\index{lu-groups@$\ell u$-groups} the objects are Abelian groups (we will use the additive signature) endowed with a lattice order that is compatible with the group structure, and with a positive Archimedean element~--- i.e. an element $u$ such that $0 < u$ and, for any other element $x$ of the group, there exists a natural number $n$ with $x \leq nu$~--- called a \emph{strong unit}\index{Unit!Strong --}. The morphisms in $\ell\cat G^\textrm{ab}_u$ are $\ell$-group morphisms~--- i.e. maps that are simultaneously group and lattice homomorphisms~--- that preserve the strong unit.

In~\cite{mundici2} the authors define a functor $\varGamma$\index{Functor!G --@$\varGamma$ --} between the category $\ell\cat G^\textrm{ab}_u$ and the one of MV-algebras; they also prove that $\varGamma$ is a categorical equivalence. Without going into details, we just recall that the image under $\varGamma$ of an $\ell u$-group $\mathbf{G} = \la G, +, -, \vee, \wedge, 0, u \ra$ is the MV-algebra $\la [0,u], \oplus, \lnot, 0 \ra$ where $[0,u] = \{x \in G \mid 0 \leq x \leq u\}$, $x \oplus y := (x + y) \wedge u$ and $\lnot x := u - x$.


Let $\AA = \varGamma(\mathbf{G})$ be an MV-algebra. A finite sequence of elements of $\AA$, $(a_0, \ldots, a_{n-1})$ is called a \emph{partition of the unit}\index{Unit!Partition of the --} if $a_0 + \cdots + a_{n-1} = u$.

In~\cite{dinolalettieri} the authors define a critically separating class of formulas
$$\mathbb{S} = \{\pi^a_b \mid a \in \N, b \in \Z\},$$
of one variable $v$, in \L ukasiewicz logic. In order to simplify the notations, let us assume the following stipulations: for formulas $\phi$ and $\psi$ we set, as usual, $\phi \oplus \psi := \lnot \phi \rightarrow \psi$ and $\phi \odot \psi := \lnot (\phi \rightarrow \lnot \psi)$. Moreover from the above notations, for every positive integer $a$, we set $a.\phi := \underbrace {\phi\oplus \ldots \oplus \phi}_{a \textrm{times}}$.

Let $a \in \N$ and $b \in \Z$ and set:
\begin{enumerate}
\item[]if $b < 0$, then $\pi^{a}_{b} (v) = v \oplus \neg v$;
\item[]if $b \geq a$, then $\pi^{a}_{b} (v) = v \odot \neg v$;
\item[]if $0 \leq b \leq a - 1$, then:
\begin{enumerate}
\item[]$\pi^a_0(v) = a.v$,
\item[]$\pi^a_1(v) = \bigoplus^{a-1}_{i=1}F_{0i}(v)$,
\item[]\hspace{1.5cm}$\ldots$
\item[]$\pi^a_b(v) = \bigoplus^{a-1}_{i=b}F_{0,1,\ldots,b-1,i}(v)$,
\item[]\hspace{1.5cm}$\ldots$
\item[]$\pi^a_{a-1}(v) = F_{0,1,\ldots,a-2,a-1}(v)$,
\end{enumerate}
\end{enumerate}
with $F_{0,1,\dots,b-1,i}(v)$ defined by
\begin{enumerate}
\item[] for every integer $i > 0$, $F_{0,i}(v) = v \odot i.v$,
\item[] for every integer $i > 1$, $F_{0,1,i}(v) = (F_{0,1}(v) \oplus \cdots \oplus F_{0,i-1}(v)) \odot F_{0,i}(v)$,
\end{enumerate}
and, by induction,
\begin{enumerate}
\item[] for every integer $i$ such that $i > b$,
\item[] $F_{0,1,\ldots,b,i}(v)=(F_{0,1,\ldots, b-1,b}(v)\oplus \cdots \oplus F_{0,1,\ldots,b-1,i-1}(v)) \odot F_{0,1,\dots,b-1,i}(v)$.
\end{enumerate}

Now let us denote by $f_{\pi^{a}_{b}(v)}$ the McNaughton functions corresponding to the formulas in $\mathbb{S}$; then, if we fix $n \in \N$ ($n >1$) and set
$$p_k(x) = f_{\pi^{n-1}_{k-1}(v)}(x) \wedge \left(f_{\pi^{n-1}_k(v)}(x)\right)^*, \qquad k = 0, \ldots, n-1,$$
we get a sequence of functions $(p_0, \ldots, p_{n-1})$ in $[0,1]^{[0,1]}$.

In analytical form we have
\begin{equation}\label{partition1}
p_0(x)=
\begin{cases}
-(n-1)x+1 & \text{if $0\leq x\leq \frac{1}{n-1}$ }\\
0 & \text{otherwise}
\end{cases},
\end{equation}

\begin{equation}\label{partition2}
p_{n-1}(x)=
\begin{cases}
(n - 1)x - (n - 2) & \text{if $\frac{n-2}{n-1}\leq x\leq 1$ }\\
0 & \text{otherwise}
\end{cases}
\end{equation}

and, for $k = 1, \ldots, n-2$,
\begin{equation}\label{partition3}
p_k(x)=
\begin{cases}
(n - 1)x - (k - 1) & \text{if $\frac{k-1}{n-1} \leq x \leq \frac{k}{n-1}$}\\
-(n - 1)x + k + 1 & \text{if $\frac{k}{n-1} \leq x \leq \frac{k+1}{n-1}$}\\
0 & \text{otherwise}
\end{cases}.
\end{equation}

Let us recall that the MV-algebra $[0,1]^{[0,1]}$ is the image, by the functor $\varGamma$, of the $\ell u$-group $\langle \R^{[0,1]}, +, \vee, \wedge, \mathbf{0}, \mathbf{1}\rangle$, where $+$, $\vee$ and $\wedge$ are defined pointwise, as usual, and $\mathbf{0}$, $\mathbf{1}$ are the maps constantly equal to $0$ and $1$ respectively.

In what follows, we will denote by $\I_m$ the set $\{0, \ldots, m - 1\}$, for all $m \in \N$; moreover, if $X$ and $Y$ are sets, we will often use expressions like $X^m$ and $X^{Y \times m}$ instead of $X^{\I_m}$ and $X^{Y \times \I_m}$ respectively.

\begin{proposition}\emph{\cite{dinolalettieri}}\label{partition}
The sequence $(p_0, \ldots, p_{n-1})$ is a partition of the unit in the MV-algebra $[0,1]^{[0,1]}$, having the property $p_k \odot p_h = 0$ for $k \neq h$.
\end{proposition}
\begin{proof}
Let us fix an index $k < n$. The thesis is an easy consequence of the following considerations:
\begin{enumerate}
\item[$(i)$]if $x = x_k$, then $p_k(x) = 1$ and $p_h(x) = 0$, for $h \neq k$;
\item[$(ii)$]if $k > 0$ and $x_{k-1} < x < x_k$, then $p_{k-1}^*(x) = p_k(x) \neq 0, 1$, and $p_h(x) = 0$ for $h \neq k - 1, k$;
\item[$(iii)$]for any $x \in [0,1]$, $\sum\limits_{k=0}^{n-1} p_k(x) = 1$, then $p_0 + \ldots + p_{n-1} = \mathbf{1}$.
\end{enumerate}
\end{proof}

\begin{corollary}\label{porthonorm}
Let $n$ be a fixed natural number and $X$ be an element of $\{[0,1]\} \cup \{\I_m\}_{m > n}$. Then the map $p \in [0,1]^{X \times n}$ defined, for all $(x,k) \in X \times \I_n$, by
\begin{equation}\label{pluk}
p(x,k) = \left\{
	\begin{array}{ll}
		p_k(x) & \textrm{if } X = [0,1] \\
		p_k\left(\frac{x}{m-1}\right) & \textrm{if } X = \I_m\\
	\end{array}\right.
\end{equation}
is an orthonormal coder.
\end{corollary}

\begin{definition}\label{lukasiewicz transform}
Let $n$ be a fixed natural number and $X$ be an element of $\{[0,1]\} \cup \{\I_m\}_{m > n}$. We denote by $[0,1]_{\luk}^X$ and $[0,1]_{\luk}^n$ the free $[0,1]_{\luk}$-modules generated, respectively, by $X$ and $\{0, \ldots, n-1\}$. We call \emph{\L ukasiewicz transform}\index{Lukasiewicz transform@\L ukasiewicz transform}\index{Transform!Lukasiewicz --@\L ukasiewicz --} of order $n$ the $\Q$-module orthonormal transform determined by the coder $p$ defined in (\ref{pluk}):
\begin{equation}\label{hn}
H_n(f)(k) = \bigvee_{x \in [0,1]} f(x) \odot p(x,k),
\end{equation}
for all $k \in \{0, \ldots, n-1\}$ and $f \in [0,1]^X$. Obviously, the \emph{\L ukasiewicz inverse transform} is the map defined, for all $x \in X$ and $g \in [0,1]^n$, by
\begin{equation}\label{lambdan}
\Lambda_n(g)(x) = \bigwedge_{k=0}^{n-1} p(x,k) \to_{\luk} g(k).
\end{equation}
\end{definition}

\section{The \L TB algorithm}
\label{ltbsec}

The \L ukasiewicz transform has been defined for maps $f: X \longrightarrow [0,1]$, where $X = [0,1]$ or $X = \I_m$ for some $m \in \N$; then the first step of its application to image processing consists of ``adapting'' the image to the domain of our operator. In other words, each image must be seen as a $[0,1]$-valued map defined on $X$.

We will consider 8-bit greyscale and RGB colour images. A greyscale image of sizes $m \times n$ is an $m \times n$ matrix with values on the set $\{0, \ldots, 255\}$, while an RGB colour image with the same sizes is encoded as a set of three $m \times n$ matrices (one for each colour channel: Red Green Blue) with values on the same set. Therefore the first step consists of normalizing the values of the matrix, or the three matrices, into $[0,1]$, i.e. we simply divide each value by 255. Since the application of the process on RGB images consists just in three parallel applications of the same process used for greyscale images, we will describe it only for the case of a single matrix.

First of all, we choose the sizes $m'$ and $n'$ of the compressed image. Then let $d_m$ be a common divisor of $m$ and $m'$, and $d_n$ be a common divisor of $n$ and $n'$ ($1 < d_m < m'$, $1 < d_n < n'$); we set also the following notations: $a = m/d_m$, $b = n/d_n$, $c = m'/d_m$ and $d = n'/d_n$.\footnote{Here we are implicitly assuming the existence of such $d_m$ and $d_n$. Indeed it is necessary here only to ensure that the dimensions of the compressed image are proportional to those of the original one but, in the unusual case where such assumption fails (i.e. if one of the original dimensions is a prime number), then it is possible to use some technical tricks.}

Let $f \in [0,1]^{m \times n}$ be a matrix (hence an image); we divide $f$ in block~--- i.e. submatrices~--- of type $a \times b$ denoted by $f^{i,j}$, $i \in \I_{d_m}$ and $j \in \I_{d_n}$. Then each $f^{ij}$ is an element of $[0,1]^{a \times b}$ and
\begin{eqnarray}
\lefteqn{f = \left(f^{ij}\right)_{i \in \I_{d_m}}^{j \in \I_{d_n}} } \nonumber \\
&=&\left(\begin{array}{cccc}
				f_{0 \ 0}   & \ldots & f_{0 \ n-1} \\
				f_{1 \ 0}   & \ldots & f_{1 \ n-1} \\
								  & \ddots & 					 \\
				f_{m-1 \ 0} & \ldots & f_{m-1 \ n-1} \\
				\end{array}\right)  \label{matrix} \\
&=&\Big(f^{ij}_{kh}\Big)_{(k,h) \in \I_a \times \I_b}^{(i,j) \in \I_{d_m}\times \I_{d_n}}. \nonumber \\ \nonumber
\end{eqnarray} 
It is easy to see that we can rewrite each block $\left(f^{i j}_{kh}\right)_{(k,h) \in \I_a \times \I_b}$ as an $a \cdot b$ vector $\mathbf{g}^{i j} = (g^{i j}_0, \ldots, g^{i j}_{ab-1})$ by setting, for all $k = 0, \ldots, ab-1$, $g^{i j}_k = f_{q(k,b)r(k,b)}$, where $q(k,b)$ and $r(k,b)$ are, respectively, the quotient and the remainder of the euclidean division $k/b$.

Now we can apply the \L ukasiewicz transform of order $cd$ to each vector $\mathbf{g}^{ij}$, thus obtaining $d_m d_n$ vectors $H_{cd}\left(\mathbf{g}^{ij}\right) \in [0,1]^{cd}$ that can be first turned back into matrices, and then recomposed by giving an $m' \times n'$ matrix: the compressed image.

Eventually, the compressed image can be treated as the original one, and each of its resulting vector can be processed by means of $\Lambda_{cd}$, thus giving the reconstructed image of dimensions $m \times n$.

\section{Applying \L TB algorithm to grey and RGB colour images}
\label{applltbsec}

In order to test the method above, we have extracted and processed several images from~\cite{sito}: the grey images Bridge and Testpat.1k, and the RGB colour ones Mandrill, Lena, Peppers and Redhead.

We have tested three processes of compression and reconstruction; in these processes we have divided the fuzzy matrix (or matrices, in the case of RGB images) associated to the images in square blocks of sizes $a \times b = 2 \times 2$, $4 \times 4$ and $8 \times 8$, respectively compressed to blocks of sizes $c \times d = 2 \times 1$, $2 \times 2$ and $5 \times 5$ by means of the formulas (\ref{matrix}) and (\ref{hn}).

The respective compression rates are obviously $\rho = (2 \cdot 2) / (2 \cdot 1) = 0.5$, $\rho = (4 \cdot 4) / (2 \cdot 2) = 0.25$ and $\rho = (8 \cdot 8) / (5 \cdot 5) \approxeq 0.39$. The blocks we obtained have been afterward decompressed to blocks of the respective original sizes, using the formula (\ref{lambdan}), hence recomposed.

\begin{figure}[ht]
\centering
	\includegraphics[width=3.24in]{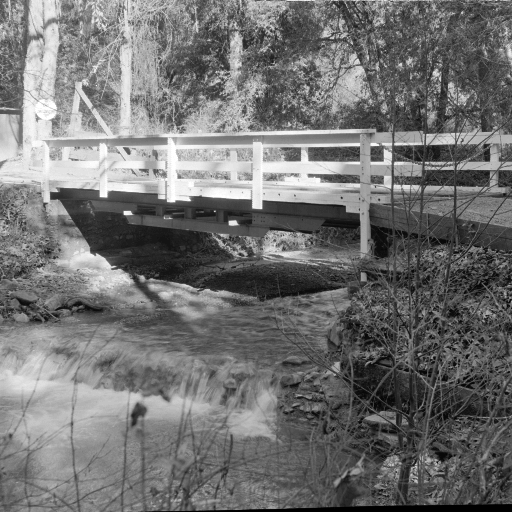}
		\caption{Bridge: original image}
			\label{bridge}
\end{figure}

\begin{figure}[bh]
\centering
	\includegraphics[width=3.24in]{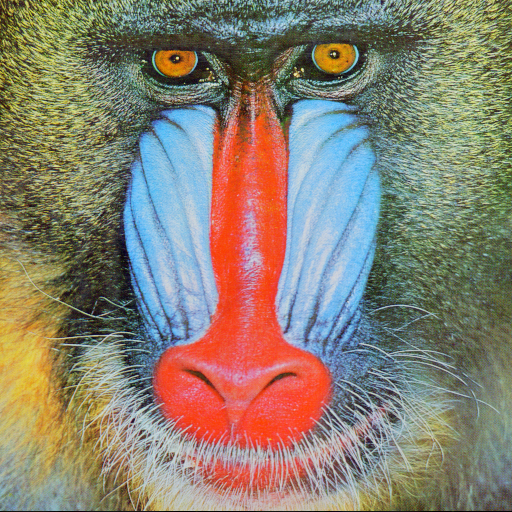}
\caption{Mandrill: original image}
\label{mandrill}
\end{figure}

\begin{figure}[hbpt]
\centering
	\includegraphics[width=3.24in]{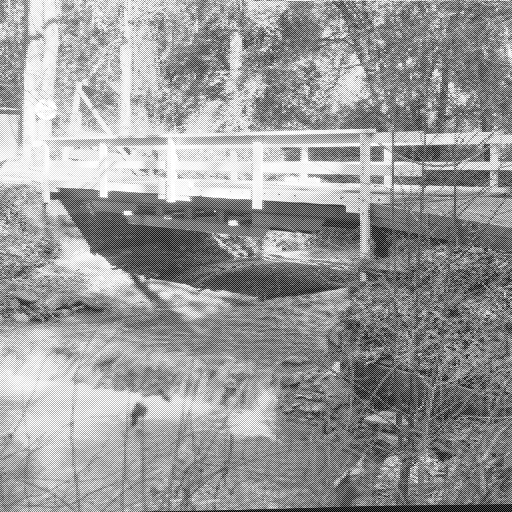}
	\caption{Bridge processed by \L TB, $\rho = 0.5$}
			\label{22-21bridge}
\end{figure}

\begin{figure}[hbp]
\centering
	\includegraphics[width=3.24in]{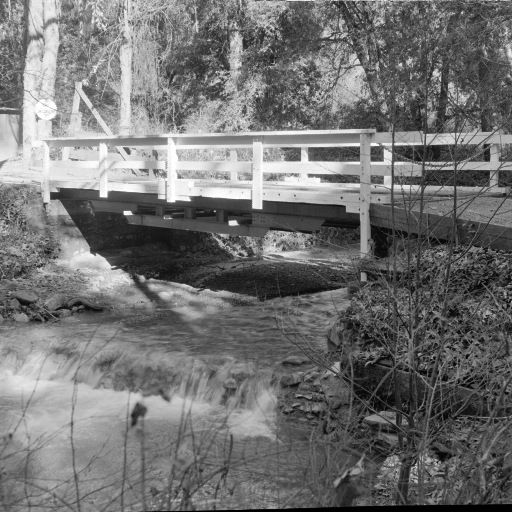}
		\caption{Bridge processed by JPEG, $\rho = 0.5$}
			\label{bridgejpg05}
\end{figure}

\begin{figure}[hbp]
\centering
	\includegraphics[width=3.24in]{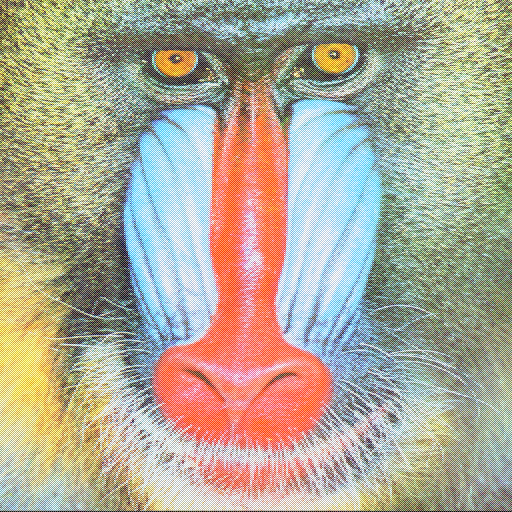}
\caption{Mandrill processed by \L TB, $\rho = 0.5$}
\label{22-21mandrill}
\end{figure}

\begin{figure}[hbp]
\centering
	\includegraphics[width=3.24in]{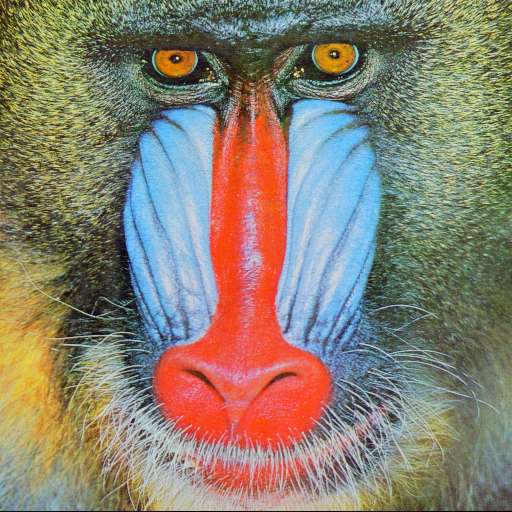}
\caption{Mandrill processed by JPEG, $\rho = 0.5$}
\label{mandrilljpg05}
\end{figure}

\begin{figure}[hbp]
\centering
	\includegraphics[width=3.24in]{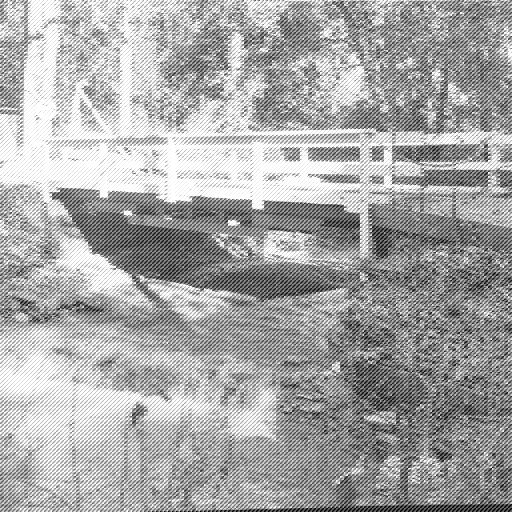}
		\caption{Bridge processed by \L TB, $\rho = 0.25$}
			\label{images/44-22bridge}
\end{figure}

\begin{figure}[hbp]
\centering
	\includegraphics[width=3.24in]{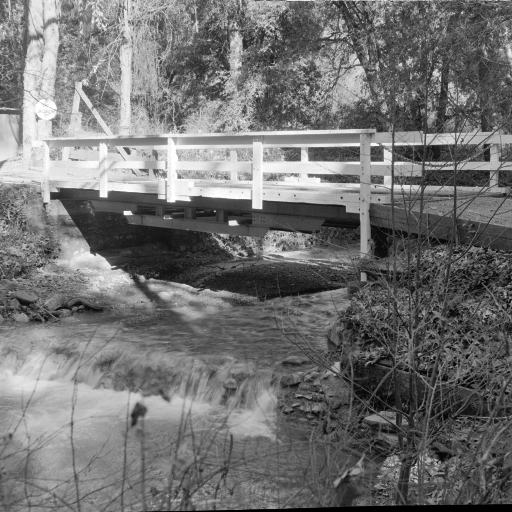}
		\caption{Bridge processed by JPEG, $\rho = 0.25$}
			\label{bridgejpg025}
\end{figure}

\begin{figure}[hbp]
\centering
	\includegraphics[width=3.24in]{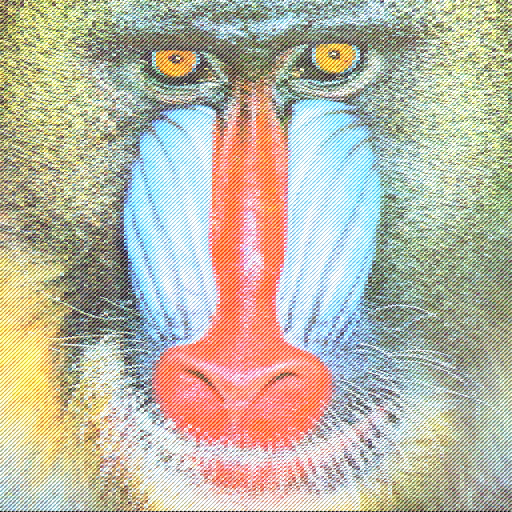}
\caption{Mandrill processed by \L TB, $\rho = 0.25$}
\label{44-22mandrill}
\end{figure}

\begin{figure}[hbp]
\centering
	\includegraphics[width=3.24in]{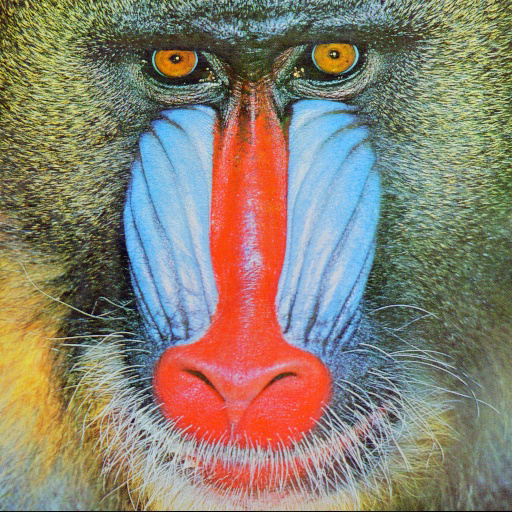}
\caption{Mandrill processed by JPEG, $\rho = 0.25$}
\label{mandrilljpg025}
\end{figure}

In Appendix~\ref{tabfig} we show the images Bridge and Mandrill in their original shape (Figures~\ref{bridge} and~\ref{mandrill}) and after the compression/reconstruction processes with ratios 0.5 and 0.25, compared with the JPEG images with the same compression ratios: Figures~\ref{22-21bridge}--\ref{mandrilljpg025}. Moreover, in Tables~\ref{values05}--\ref{values025}, we list some numerical test indices~--- namely Peak Signal to Noise Ratio (PSNR) and Root Mean Square Error (RMSE)~--- for all of the images we processed. Last, Table~\ref{runningtimes} shows a comparison between the execution times of \L TB and JPEG.

Now we will list some properties of the \L TB algorithm that are direct consequences of the results presented in Chapter~\ref{mqchap}. First, Theorem~\ref{wadjointpair} implies the following:
\begin{enumerate}
\item[-] If the image $i_1$ is pixelwise brighter than or equal to the image $i_2$, then the processed image $(\Lambda \circ H)(i_1)$ is pixelwise brighter than the processed image $(\Lambda \circ H)(i_2)$.
\item[-] For any image $i$, the processed image $(\Lambda \circ H)(i)$ is pixelwise brighter than or equal to $i$.
\item[-] Further applications of the pair compression/reconstruction to any already processed image are lossless.  
\item[-] The \L TB algorithm is invariant under the action of a homogeneous darkening filter applied by the \L ukasiewicz t-norm\footnote{Here, by a \emph{homogeneous darkening filter}, we mean a ``flat'' image, i.e. a constant map. Its application by the \L ukasiewicz t-norm is the pixelwise product of the filter by the image.}.
\end{enumerate}
By Theorem~\ref{sadjointpair}, we have
\begin{enumerate}
\item[-] If we apply the inverse transform $\Lambda$ first, and then the direct one, $H$, we only zoom in and out the image with no errors introduced.
\end{enumerate}

\section{Comparing \L TB with JPEG}
\label{ltbjpgsec}

\subsection*{Computability}

The coding/decoding algorithms are usually compared by means of their execution times and the values of some parameters (PSNR, RMSE, MSE).

The comparison between the \L TB algorithm and JPEG is heavily conditioned by their underlying implementation. Indeed in the process of compression and reconstruction the computational time, for each block of sizes $8 \times 8$, is characterized by the execution of the inverse DCT/DCT. The standard implementation of DCT determines an asymptotic computational time, of the DCT on one block, that is $O(a \cdot b)$, where $a$ and $b$ are, respectively, the numbers of rows and columns of the block.

If we process $8 \times 8$ blocks, the asymptotic time is not relevant anymore and it is more convenient to look at the number of operations executed. The standard implementation requires in general $1024$ products and $896$ sums for computing the DCT on a block of these sizes, but there exist several optimized DCT implementations (FastDCT et al., see for instance~\cite{hung}) that reduce significantly these numbers. For example, the FastDCT proposed in~\cite{feig} requires only $54$ products, $464$ sums and $6$ arithmetical shifts, giving the same result.

Furthermore we should add the time and operations required for other components of the application: quantization, downsampling and entropic encoding.

If we set $x = a \cdot b \cdot c \cdot d$ and $y$ as the number of colour channels of the image (one for grey images, three for RGB images), the \L TB algorithm computes, for the compression of one block, $x \cdot y$ products, $x \cdot y$ sums, $x \cdot y$ comparisons and, at most, $x \cdot y$ assignments. If, for instance, we set $a = b = 3$, $c = d = 2$ and $y = 3$, then the whole compression algorithm requires~--- for each block~--- $108$ products, $108$ sums, $108$ comparisons and at most $108$ assignments. All these values can be still reduced by means of a suitable advanced implementation.

It follows from these considerations that the \L TB algorithm requires an execution time much shorter with respect to JPEG. Nevertheless the values in Tab.~\ref{runningtimes} show that the JPEG application used for our tests (FreeImage Library 3.8, in~\cite{sito2}) is faster for some images. This fact depends on the sampling scheme in blocks of sizes $8 \times 8$, that supports an optimization of JPEG's implementations for several CPU architectures. This is in particular the case of all CPUs supporting MMX, SSE, SSE2, SSE3 and the ones with SIMD (Single Instructions Multiple Data) architecture, using a 64-bit sampling.

On the other hand, the source code of our CoDec has been realized with simple C-like optimizations, since our purpose was just showing the feasibility of this approach and the possible results. So this comparison should be read by also considering the possibility of improving the application overworking SIMD architectures' optimizations.

Some improving techniques could be provided, for example, adopting a sampling scheme similar to the one adopted by the JPEG algorithm, i.e. a scheme enabling the algorithm to work on vectors whose size is a multiple of the length of the machine word. For instance, assuming a CPU with SSE support and 32-bit architecture, it could be $a = b = 8$, $c = d = 4$. This sampling scheme would reduce the required $x \cdot y$ products and comparisons to $\frac{x \cdot y}{4}$ where~--- as we already stated~--- $x = a \cdot b \cdot c \cdot d$ and $y$ is the number of colour channels of the image.

\subsection*{Numerical indices}

With regards to the numerical comparison between JPEG and \L TB, even though there is still a consistent gap, Theorem~\ref{wadjointpair} proves that the \L TB algorithm possesses an interesting property: unlike JPEG, an iterative application of the algorithm on the same image is lossy just for the first process and lossless for the subsequent ones. In other words, once we have compressed and reconstructed an image, we can apply the same process again, on the reconstructed image, obtaining exactly the same compressed and the same reconstructed images (with fixed sizes of the blocks).

Moreover we must also underline that JPEG is composed of two parts, a lossy compression method and a lossless one, while the \L TB is only lossy. This fact draw a possible direction for further studies.

\section*{References and further readings}
\addcontentsline{toc}{section}{References and further readings}

The contents of this chapter can be found in~\cite{dinolarusso} and~\cite{dinolarusso2}. Further examples of fuzzy algorithms for image processing have been already cited at the end of the previous chapter.

\conclusion

In this dissertation we have proposed an investigation of the basic categorical and algebraic properties of quantale modules, and we have shown that certain operators between objects in these categories find important applications in Mathematical Logic and Image Processing. Here we list our main contributions.

In Chapter~\ref{mqchap}, besides recalling basic and known notions and properties on quantale modules, we have shown --- in Sections~\ref{nucleisec} and~\ref{mqtransec} --- the properties of $\Q$-module structural closure operators and $\Q$-module transforms, and their connection with $\Q$-module morphisms. In Section~\ref{amalsec} we proved that the categories of $\Q$-modules have the strong amalgamation property and, last, in Section~\ref{mqtensorsec} we showed the existence of tensor products of $\Q$-modules, also proving that the property, of a module, of being the coproduct of cyclic projectives is inherited by its tensor product with an extension of its quantale of scalars.

These results have been applied to deductive systems in Chapter~\ref{consrelchap}. Once we have recalled the results obtained by N. Galatos and C. Tsinakis, in Section~\ref{transsec} we presented a definition of interpretation between deductive systems, where the concept of translation is separated and independent from the existence of any relationship between the deductive apparatuses. Then we showed how a translation and an interpretation can be algebraically represented in terms of quantale and module morphisms. In Section~\ref{nonstrsec} we introduced the notions of non-structural interpretation of deductive systems, showing an algebraic characterization of this concept, by means of the results of Chapter~\ref{mqchap}.

In Chapter~\ref{imagechap} we have shown that certain operators used for digital image compression and analysis are special cases of $\Q$-module transforms and, last, we have presented in Chapter~\ref{ltbchap} a concrete realization of a $\Q$-module transform for compressing and reconstructing digital images.

Although the results seem to be promising, especially for how easily they can be applied, we cannot pretend~--- of course~--- the applications presented to be not open to further significant developments and improvements.

To what extent Logic, the representation of deductive systems by means of quantale modules is only at its first step~--- the propositional level~--- but all the evidences indicate that an extension to higher order languages should be possible and fruitful as well. Indeed, if the results established in Section~\ref{transsec} will be suitably extended to deductive systems of any type, it would result in the possibility of defining categories whose objects are deductive systems, and then of applying the powerful tool of Category Theory to many logical tasks.

On the other hand, as we already observed in Chapter~\ref{imagechap}, the approach via quantale modules allowed us to group together, in a unique formal context, algorithms that act on digital images in completely different ways and have been proposed for dealing with problems different in nature. Apart from the obvious (and eternal) issue of improving the results of applications, the main open problem is the following: currently, the $\Q$-modules we really encounter in these situations are exclusively $[0,1]$-modules, a very special class of modules, hence such a formal context will be redundant from this point of view, until its applications will be extended to a wider class of tasks in data management. This is probably the most important challenge in this connection. Last, we also need to take into account that winning this challenge would naturally give rise to a further issue, namely the necessity of numerically (or, anyhow, objectively) estimate results of the applications by introducing a sort of measure on quantale modules.

\appendix
\addcontentsline{toc}{part}{Appendix}

\chapter{Tables and Figures}
\label{tabfig}

\twocolumn


\onecolumn

\onecolumn

\onecolumn

\onecolumn

\onecolumn

\onecolumn

\begin{table}[htp]
\caption{Numerical values for grey and RGB colour images ($\rho = 0.5$)}
\centering
\begin{tabular}{|c|c|c|c|c|}
\hline
\multicolumn{1}{|c|}{\raisebox{-1.50ex}[0cm][0cm]{\!Image\!}}
& \multicolumn{1}{|c|}{JPEG}
& \multicolumn{1}{|c|}{\L TB}
& \multicolumn{1}{|c|}{JPEG}
& \multicolumn{1}{|c|}{\L TB}
\cr & RMSE  &  RMSE  & PSNR  & PSNR \\ \hline
 Bridge   &  $2.4985$ & $58.4469$ &	$40.1773$	& $12.7956$\\ \hline
 Testpat.1k  &  $0.0833$ & $66.2213$ & $69.7131$ & $11.7108$\\ \hline
 Lena   &  $2.2606$ & $53.1803$ &	$41.0464$	& $13.6158$\\ \hline
 Mandrill  &  $6.4739$ & $55.8967$ & $31.9075$ & $13.1831$\\ \hline
 Peppers   &  $2.0813$ & $56.7046$ &	$41.7641$	& $13.0584$\\ \hline
 Redhead   &  $0.4454$ & $57.7974$ &	$55.1554$	& $12.8926$\\ \hline 
\end{tabular}
\label{values05}
\end{table}

\begin{table}[htp]
\caption{Numerical values for grey and RGB colour images ($\rho = 0.25$)}
\centering
\begin{tabular}{|c|c|c|c|c|}
\hline
\multicolumn{1}{|c|}{\raisebox{-1.50ex}[0cm][0cm]{\!Image\!}}
& \multicolumn{1}{|c|}{JPEG}
& \multicolumn{1}{|c|}{\L TB}
& \multicolumn{1}{|c|}{JPEG}
& \multicolumn{1}{|c|}{\L TB}
\cr & RMSE  &  RMSE  & PSNR  & PSNR \\ \hline
 Bridge   &  $6.1120$ & $68.3310$ &	$32.4071$	& $11.4384$\\ \hline
 Testpat.1k  &  $0.0833$ & $74.2191$ & $69.7131$ & $10.7205$\\ \hline
 Lena   &  $2.4819$ & $61.6912$ &	$40.2351$	& $12.3263$\\ \hline
 Mandrill  &  $7.1230$ & $65.7072$ & $31.0775$ & $11.7785$\\ \hline
 Peppers   &  $2.1147$ & $65.4880$ &	$41.6257$	& $11.8076$\\ \hline
 Redhead   &  $0.4454$ & $67.0408$ &	$55.1554$	& $11.6040$\\ \hline 
\end{tabular}
\label{values025}
\end{table}

\begin{table}[htp]
\caption{Running times (in $ms$)}
\centering
\begin{tabular}{|c|c|c|c|c|c|c|}
\hline
\multicolumn{1}{|c|}{\raisebox{-1.50ex}[0cm][0cm]{\!Image\!}}
& \multicolumn{1}{|c|}{JPEG}
& \multicolumn{1}{|c|}{\L TB}
& \multicolumn{1}{|c|}{JPEG}
& \multicolumn{1}{|c|}{\L TB}
& \multicolumn{1}{|c|}{JPEG}
& \multicolumn{1}{|c|}{\L TB}
\cr & $\rho = 0.5$ & $\rho = 0.5$ & $\rho = 0.39$ & $\rho = 0.39$ & $\rho = 0.25$ & $\rho = 0.25$\\ \hline
 Bridge   &  $80.98$ & $41.85$ &  $51.68$ & $290.81$ &  $43.64$ & $59.19$ \\ \hline
 Testpat.1k  &  $150.47$ & $123.30$ &  $82.13$ & $4888.75$ &  $4211.14$ & $181.17$\\ \hline
 Lena  &  $184.40$ & $98.68$ &  $173.47$ & $1835.67$ &  $141.04$ & $178.40$\\ \hline
 Mandrill  &  $187.29$ & $105.62$ &  $98.65$ & $6287.89$ &  $113.55$ & $4250.62$\\ \hline
 Peppers  &  $150.61$ & $85.96$ &  $92.95$ & $1551.73$ &  $98.28$ & $157.85$\\ \hline
 Redhead & $98.33$ & $132.97$ & $93.53$ & $2491.04$ & $97.05$ & $194.70$\\ \hline
\end{tabular}
\label{runningtimes}
\end{table}


\backmatter

\printindex

\end{document}